\documentclass[12pt]{article}
\usepackage{amsmath,amsthm,amsfonts}
\usepackage{graphicx}
\usepackage{enumerate}
\usepackage{natbib}
\usepackage{url} % not crucial - just used below for the URL

\usepackage[normalem]{ulem}
\usepackage{bbm}

\usepackage[english]{babel}
\usepackage{apptools}
\usepackage{paralist}
\usepackage{enumitem,kantlipsum}
\usepackage{algorithm,algcompatible}
\usepackage{hyperref}

\newcommand{\blind}{1}

\newtheorem{theorem}{Theorem}[section]
\newtheorem{corollary}{Corollary}[section]
\newtheorem{lemma}{Lemma}[section]
\newtheorem{proposition}{Proposition}[section]
\newtheorem{definition}{Definition}%[section]
\newtheorem{assumption}{Assumption}%[section]
\newtheorem{remark}{Remark}%[section]

\newcommand{\N}{\mathbb{N}}
\newcommand{\R}{\mathbb{R}}

\newcommand{\E}{\mathrm{E}}

\newcommand{\mf}{\mathbf}
\newcommand{\vf}{\mathbf}

\DeclareMathOperator*{\argmin}{arg\,min}

\begin{document}

\def\spacingset#1{\renewcommand{\baselinestretch}%
{#1}\small\normalsize} \spacingset{1}

%%%%%%%%%%%%%%%%%%%%%%%%%%%%%%%%%%%%%%%%%%%%%%%%%%%%%%%%%%%%%%%%%%%%%%%%%%%%%%

\if1\blind
{
  \title{\bf Large covariance matrix estimation by penalized log-det heuristics}
  \author{Enrico Bernardi\hspace{.2cm}\\
    %The authors gratefully acknowledge \textit{please remember to list all relevant funding sources in the unblinded version}}\hspace{.2cm}\\
    Department of Statistical Sciences, University of Bologna\\
    and \\
    Matteo Farn\`e \thanks{We thank the participants to the Conference "Mathematical Methods of Modern Statistics 3", held in Luminy (France) in June 2022, where a preliminary version of this work was presented, for their encouragement and constructive comments, and the Supervisory Statistics Division
    of the European Central Bank for allowing us to make research on confidential data. We have no competing interests to declare.}\\
    Department of Statistical Sciences, University of Bologna}
  \maketitle
} \fi

\if0\blind
{
  \bigskip
  \bigskip
  \bigskip
  \begin{center}
    {\LARGE\bf Large covariance matrix estimation by penalized log-det heuristics}
\end{center}
  \medskip
} \fi

\bigskip
\begin{abstract}
This paper provides a comprehensive estimation framework for large covariance matrices via a log-det heuristics augmented by a nuclear norm plus $\ell_{1}$-norm penalty.
We develop the model framework, which includes high-dimensional approximate factor models with a sparse residual covariance.
We prove that the aforementioned log-det heuristics is locally convex with a Lipschitz-continuous gradient, so that a proximal gradient algorithm may be stated to numerically solve the problem while controlling the threshold parameters.
The proposed optimization strategy recovers in a single step both the covariance matrix components and the latent rank and the residual sparsity pattern with high probability,
and performs systematically not worse than the corresponding estimators employing Frobenius loss in place of the log-det heuristics.
The error bounds for the ensuing low rank and sparse covariance matrix estimators are established, and the identifiability conditions for the latent geometric manifolds are provided, improving existing literature.
The validity of outlined results is highlighted by an exhaustive simulation study and a financial data example involving Euro Area banks.
\end{abstract}

\noindent%
\vfill

\newpage
\spacingset{1.9} % DON'T change the spacing!

%%%% Main text entry area:
\section{Introduction}
\label{Intro}

\subsection{Motivation}\label{motiv}
Estimating high-dimensional covariance or precision matrices has become a
crucial task nowadays, due to the increasing availability of datasets composed of
a large number of variables $p$ compared to the sample size $n$ in many fields, like
economics, finance, biology, genetics, health, climatology, and social sciences.
The consistency of estimated covariance matrices is a prerequisite
to perform several statistical procedures in high dimensions like principal component analysis (PCA),
cluster analysis, graphical model inference, among others.
Recent books on this relevant topic are \cite{pourahmadi2013high,zagidullina2021high},
while recent comprehensive reviews include \cite{fan2016overview,wainwright2019high,lam2020high,ledoit2021shrinkage}.
%
%portfolio allocation,
%
The amplitude of techniques developed to overcome the estimation issues
in high dimensions now provides several state-of-the art solutions,
but also leaves some room to further improve the estimation process in many directions.

Although the theory of covariance matrix estimation for Gaussian data with low and fixed dimension
was developed in the fifties by pioneeristic contributions \citep{anderson1958introduction},
it became soon apparent that the sample covariance matrix $\mf{\Sigma}_n$
is not a reliable estimator of the true covariance matrix $\mf{\Sigma}^{*}$
when the dimension $p$ is large. %{even if it is smaller than $n$}.
%when $p/n$ does not converge to $0$,  where ``$\rightarrow 1^{-}$'' means ``tends to $1$ from the left'',
%and is not a consistent estimator of the true covariance matrix $\mf{\Sigma}^{*}$ when $p/n$ does not converge to $0$
%
As explained in \cite{ledoit2004well}, in fact, sample eigenvalues may be
severely over-dispersed in that case,
%when $p/n\rightarrow 1^{-}$,
as they follow the Mar$\check{\mathrm{c}}$enko-Pastur law
\citep{marchenko1967distribution}. This leads $\mf{\Sigma}_n$ to become numerically unstable as $p\to\infty$.
A large $p$ also easily leads to identifiability issues for $\mf{\Sigma}^{*}$,
as the number of parameters to be recovered grows quadratically in $p$ in absence of further assumptions.
In addition, if $p/n \geq 1$, $p-n$ sample eigenvalues are null,
thus irremediably affecting the invertibility of $\mf{\Sigma}_n$.

A relevant approach to large covariance matrix estimation passes by the %second
assumption of a specific structure for $\mf{\Sigma}^{*}$, in order to drastically reduce the number of parameters.
One option is to assume a factor model structure.
The spiked covariance model, introduced in \cite{johnstone2001distribution}
and recovered in \cite{johnstone2009consistency},
is a successful attempt of this kind.
Linear eigenvalue shrinkage under the same model has been proposed in \cite{donoho2018optimal}.
Another option is to assume some type of sparsity.
Under that assumption, with respect to the underlying structure,
it has been proposed to recover $\mf{\Sigma}^{*}$ by
applying hard-thresholding \citep{bickel2008regularized}, soft-thresholding \citep{bickel2008covariance},
generalized thresholding \citep{rothman2009generalized}, or adaptive thresholding \citep{cai2011adaptive}.
Penalized maximum likelihood (ML) has also been used to estimate a sparse covariance \citep{bien2011sparse}
or precision \citep{friedman2008sparse} matrix.

%tapering \cite{furrer2007estimation}
These two types of algebraic structure enforced in $\mf{\Sigma}^{*}$ may both be too restrictive.
In fact, as pointed out by \cite{giannone2021economic} for economic data, it is likely that
the sparsity assumption is too strong in high dimensions, as the interrelation structure among the variables
is actually more dense than sparse.
At the same time, a strict factor model does not allow for any idiosyncratic covariance structure
that may catch specific pairs of extra-correlated variables beyond the factors.
As a consequence, it became clear that conditional sparsity with respect to an underlying approximate factor model \citep{chamberlain1983arbitrage}
can effectively merge factor model and sparsity assumptions,
thus being simultaneously a parsimonious and flexible approach.
A covariance matrix estimator assuming conditional sparsity is POET \citep{fan2013large},
that proposes to threshold the principal orthogonal complement to obtain a consistent solution.
%The behaviour of PCA and ML estimators of exact or approximate factor models in high dimensions has been studied in \cite{bai2003inferential,bai2012statistical,bai2016maximum}.

\subsection{Approach}
Conditional sparsity can be imposed by assuming for $\mf{\Sigma}^{*}$ a low rank plus sparse decomposition, that is \begin{equation}\mf{\Sigma}^{*}=\mf{L}^{*}+\mf{S}^{*}=\mf{B}\mf{B}'+\mf{S}^{*},\label{mod_l+s}\end{equation}
%where $B$ is any $p\times r$ matrix such that
where $\mf{L}^{*}=\mf{B}\mf{B}'=\mf{U}_L \mf{\Lambda}_L \mf{U}_L'$,
with $\mf{U}_L$ $p\times r$ semi-orthogonal matrix %containing
%as columns the eigenvectors of $\mf{L}^{*}$
and $\mf{\Lambda}_L$ $r\times r$ diagonal positive definite matrix, and
$\mf{S}^{*}$ is a positive definite and element-wise sparse matrix,
containing only $s \ll {p(p-1)}/{2}$ off-diagonal non-zero elements.
%
%of the corresponding eigenvalues.
%and $\mf{S}^{*}$ is a $p \times p$ matrix with full rank. % semi-orthogonal %diagonal
%%
Structure \eqref{mod_l+s} has become the reference model
under which several high-dimensional covariance matrix estimators work.
It can be recovered by nuclear norm plus $l_1$ penalization, that is by solving
\begin{equation}
\bigl(\widehat{\mf{L}}_n,\widehat{\mf{S}}_n\bigr)=
%\arg\min_{{\mf{L},\mf{S},\\ \mf{L} \succeq 0, \mf{S} \succ 0,  \mf{L}+\mf{S} \succ 0}}
%\mathcal{L}(\mf{\Sigma},\mf{\Sigma}_n)+\mathcal{P}(\mf{L},\mf{S}),\label{obj_all}
\argmin_{\mf{L},\mf{S}\in\R^{p \times p},\mf{\Sigma}=\mf{L}+\mf{S}}
\mathcal{L}(\mf{\Sigma},\mf{\Sigma}_n)+\mathcal{P}(\mf{L},\mf{S}),\label{obj_all}
\end{equation}
where $\mathcal{L}(\mf{\Sigma},\mf{\Sigma}_n)$ is a smooth loss function,
measuring the distance between $\mf{\Sigma}\in\R^{p \times p}$, which is an admissible covariance matrix estimate,
and the sample covariance matrix $\mf{\Sigma}_n$, and
$\mathcal{P}(\mf{L},\mf{S})=\psi \Vert\mf{L}\Vert_{*}+\rho \Vert\mf{S}\Vert_{1}$ is a non-smooth composite penalty,
where $\Vert\mf{L}\Vert_{*}=\sum_{i=1}^p \lambda_i(\mf{L})$ is the nuclear norm of $\mf{L}$,
i.e. the sum of the eigenvalues of $\mf{L^{*}}$, and $\Vert\mf{S}\Vert_{1}=\sum_{i=1}^{p} \sum_{j=1}^p \vert\mf{S}_{ij}\vert$
is the $\ell_1$-norm of $\mf{S}$, while $\psi$ and $\rho$ are non-negative threshold parameters.

%The nuclear norm was first proposed in \cite{fazel2001rank} as an alternative to PCA.
%In \cite{fazel2002matrix}, it is proved that $\psi \Vert\mf{L}\Vert_{*}+\rho \Vert\mf{S}\Vert_{1}$ is the tightest convex relaxation of the original non-convex penalty $\psi\mathrm{rk}(\mf{L})+\rho\Vert \mf{S}\Vert _{0}$, which involves the rank of $\mf{L}$ and the number of non-zeros in $\mf{S}$,
%but yields to a NP-hard problem.
%%
%\cite{donoho2006most} proves that $\ell_1$-norm minimization gives the sparsest solution to most large underdetermined linear systems, while \cite{recht2010guaranteed} proves that the nuclear norm minimization provides guaranteed rank minimization under a set of linear equality constraints.
%%
%\cite{candes2009near} shows that $\ell_1$-norm minimization selects the best linear model in a wide range of situations.
%%
%The nuclear norm has instead been used to solve large matrix completion problems, like in
%\cite{srebro2004maximum,candes2010power,mazumder2010spectral,hastie2015matrix}.
%%
%Nuclear norm plus $\ell_1$-norm minimization has been exploited in \cite{candes2011robust} to provide %($\mathcal{P}(\mf{L},\mf{S})$)
%a robust version of PCA under grossly corrupted or missing data, and in \cite{chandrasekaran2011rank} to provide exact recovery for structure \eqref{mod_l+s} in absence of noise.

Heuristics \eqref{obj_all} has generated a stream of literature where new covariance matrix estimators in high dimensions are derived.
\cite{agarwal2012noisy} ensures optimal rates for the solutions of \eqref{obj_all} under model \eqref{mod_l+s} via a purely analytical approach, providing the approximate recovery of $\mf{\Sigma}^{*}$ and a bounded non-identifiability radius for $\mf{L}^{*}$.
In \cite{chandrasekaran2012latent}, a latent graphical model structure, based on a sparse minus low rank decomposition for $\mf{\Sigma}^{*}$,
is learnt by solving a problem of type \eqref{obj_all}, providing both parametric and algebraic consistency
(see Section \ref{Cons_both} for the definitions).
\cite{luo2011high} derives LOREC estimator under model \eqref{mod_l+s} for $\mf{\Sigma}^{*}$
via the objective in \eqref{obj_all} with $\mathcal{L}(\mf{\Sigma},\mf{\Sigma}_n)=\mathcal{L}^{(F)}(\mf{\Sigma},\mf{\Sigma}_n)$ and
$\mathcal{L}^{(F)}(\mf{\Sigma},\mf{\Sigma}_n)=0.5\Vert \mf{\Sigma}_n-(\mf{L}+\mf{S})\Vert_{F}^2$. A similar heuristics
is employed in \cite{barigozzi2024algebraic} to recover the spectral density matrix in high dimensions.
%%

%\subsection{Contribution}
In \cite{farne2020large,farne2024large}, UNALCE estimator is proposed. UNALCE is based on problem \eqref{obj_all}
with $\mathcal{L}(\mf{\Sigma},\mf{\Sigma}_n)=\mathcal{L}^{(F)}(\mf{\Sigma},\mf{\Sigma}_n)$
%and $\mathcal{L}^{(F)}(\mf{\Sigma},\mf{\Sigma}_n)=0.5\Vert \mf{\Sigma}_n-(\mf{L}+\mf{S})\Vert_{F}^2$,
like LOREC, but overcomes LOREC deficiencies thanks to a random matrix theory result that
holds under a wide range of approximate factor models and accounts for the high-dimensional case $p \geq n$.
UNALCE is both algebraically consistent in the sense of \cite{chandrasekaran2012latent} and
parametrically consistent in Frobenius and spectral norm.
%
%In EJS(2022), the theory of \cite{farne2020large} is further refined
%to effectively estimate large-dimensional factor models.
Although UNALCE is the optimal estimator in finite sample
when $\mathcal{L}(\mf{\Sigma},\mf{\Sigma}_n)=\mathcal{L}^{(F)}(\mf{\Sigma},\mf{\Sigma}_n)$
in terms of Frobenius loss, there is room to further improve it by replacing
$\mathcal{L}^{(F)}(\mf{\Sigma},\mf{\Sigma}_n)$ with a different loss.

The Frobenius loss optimizes in fact the
entry by entry performance of the overall covariance matrix estimate
$\widehat{\mf{\Sigma}}_n=\widehat{\mf{L}}_n+\widehat{\mf{S}}_n$. A loss able to explicitly control
the spectrum estimation quality might be desirable, in order to ensure
%unlike eigenvalue shrinkage,
the recovery of latent rank and residual sparsity pattern,
%algebraic consistency of $\bigl(\widehat{\mf{L}}_n,\widehat{\mf{S}}_n\bigr)$,
while simultaneously optimizing the estimated spectrum in terms of distance from the true one.
The loss $\mathcal{L}^{(ld)}(\mf{\Sigma},\mf{\Sigma}_n)=0.5\ln \det (\mf{I}_p+\mf{\Delta}_{n}\mf{\Delta}_{n}')$,
where $\mf{\Delta}_{n}=\mf{\Sigma}-\mf{\Sigma}_n$ and $\mf{\Sigma}=\mf{L}+\mf{S}$, is a possible one satisfying
these needs, because it is controlled by the individual eigenvalues of $\mf{\Delta}_n$, since
\begin{equation}\ln \det (\mf{I}_p+\mf{\Delta}_{n}\mf{\Delta}_{n}')={\ln} \prod_{i=1}^{p}\lambda_i(\mf{I}_p+\mf{\Delta}_{n}\mf{\Delta}_{n}')
\leq\sum_{i=1}^{p}\ln(1+\lambda^2_i(\mf{\Delta}_n))
\leq \sum_{i=1}^{p}\lambda^2_i(\mf{\Delta}_n),\label{eigen}
\end{equation}
thus providing an intrinsic eigenvalue correction.
Inequality \eqref{eigen} holds because $\ln(1+\lambda^2_i(\mf{\Delta}_n))\leq\lambda^2_i(\mf{\Delta}_n),$ since $\ln(x) \leq x-1$ for any $x>0$.
It trivially holds when $\mf{\Delta}_n \mf{\Delta}_n'$ is diagonal, and it can be proved in the general case by
diagonalizing the positive semi-definite symmetric matrix $\mf{\Delta}_n\mf{\Delta}_n'$.
Our estimator pair therefore becomes
\begin{equation}\label{obj}
%\bigl(\widehat{\mf{L}}_n,\widehat{\mf{S}}_n\bigr)=
%\arg\min_{{\mf{L},\mf{S},\\ \mf{L} \succeq 0, \mf{S} \succ 0,  \mf{L}+\mf{S} \succ 0}}
%%\phi(\mf{L},\mf{S})=
%\frac{1}{2}\ln \det (\mf{I}_p+\mf{\Delta}_{n}\mf{\Delta}_{n}')
%+\psi \Vert\mf{L}\Vert_{*}+\rho \Vert\mf{S}\Vert_{1}.
\bigl(\widehat{\mf{L}}_n,\widehat{\mf{S}}_n\bigr)=
\argmin_{\mf{L},\mf{S}\in\R^{p \times p},\mf{\Sigma}=\mf{L}+\mf{S}}
%\phi(\mf{L},\mf{S})=
\mathcal{L}^{(ld)}(\mf{\Sigma},\mf{\Sigma}_n)
+\mathcal{P}(\mf{L},\mf{S}).
\end{equation}
Problem \eqref{obj} minimizes a loss controlled by the Euclidean norm of the eigenvalues of the matrix $\mf{\Delta}_{n}$, while simultaneously minimizing the latent rank and the residual support size, whose the nuclear norm and the $\ell_1$-norm are the respective tightest convex relaxations \citep{fazel2002matrix}. A relevant challenge is that $\mathcal{L}^{(ld)}(\mf{\Sigma},\mf{\Sigma}_n)$ is locally convex, i.e. it is convex into a specific range for $\mf{\Delta}_{n}$. In this paper, we exploit this fact to propose a proximal gradient algorithm to compute \eqref{obj}
%with $\mathcal{L}(\mf{\Sigma},\mf{\Sigma}_n)=\mathcal{L}^{(ld)}(\mf{\Sigma},\mf{\Sigma}_n)$,
and to prove algebraic and parametric consistency for the estimates of $\mf{L}^{*}$, $\mf{S}^{*}$, and $\mf{\Sigma}^{*}$
obtained in this way.

\subsection{Contribution}

%\far{START HERE}

Our model setting imposes for the data vector an approximate factor model with sparse residual covariance. Unlike the original approximate factor model of \cite{chamberlain1983arbitrage},
the residual covariance matrix $\mf{S}^{*}$ is assumed to be positive definite, and its estimate $\widehat{\mf{S}}_n$ obtained through \eqref{obj} also is with high probability.
We then assume $\mf{S}^{*}$ to be element-wise sparse.
This means that, conditionally on the factors, some variables are assumed to be still correlated and their correlation is captured by the non-zero elements of $\mf{S}^{*}$.
This residual sparsity structure is very suitable in high dimensions and quite more general than the restrictions proposed in \cite{bai2012statistical} and \cite{bai2016maximum}.
Nevertheless, the sparsistency of $\widehat{\mf{S}}_n$  (i.e. the exact recovery of the sparsity pattern of $\mf{S}^{*}$) is ensured with high probability.

We follow the algebraic approach of \cite{chandrasekaran2012latent}, consisting in solving a variety-constrained version of \eqref{obj},
which is proven to be equivalent to the unconstrained version of \eqref{obj}. However, instead of learning a latent graphical model from the sample precision matrix as in \cite{chandrasekaran2012latent},
here we learn a factor model from the sample covariance matrix. Problem \eqref{obj} satisfies the Restricted Strong Convexity of \cite{10.1214/12-STS400}
for all $p\in\N$ as $n \to \infty$, so that the error bounds in Frobenius norm of \cite{agarwal2012noisy} also hold for $\widehat{\mf{L}}_n$, $\widehat{\mf{S}}_n$, $\widehat{\mf{\Sigma}}_n$ (see the Supplement Section D).

An important feature of our approach is that we allow for weak factors, i.e. factors whose pervasiveness diverges slower than $O(p)$ as $p \to \infty$,
unlike factor models in \cite{bai2003inferential,bai2012statistical,bai2016maximum}, for instance. A related crucial aspect regards the latent rank selection.
Some well-known criteria to select the latent rank under pervasive approximate factor models are proposed in \cite{bai2002determining}.
Those criteria are widely used in the literature, for instance in \cite{fan2013large} to construct POET estimator
and in \cite{bai2019rank}, which explores pervasive approximate factor model estimation via nuclear norm minimization.
Recently, weak approximate factor model estimation by PCA has been explored in \cite{bai2023approximate},
but no criterion to recover the latent rank is proposed therein. In this paper, the latent rank is instead consistently estimated
by solving the single optimization problem \eqref{obj}, thanks to the nuclear norm term. %%due cinesi paper Stat Sinica review

Such feature was already present in \cite{chandrasekaran2012latent} to learn the number of latent variables driving the marginal graphical model.
It has later been exploited in \cite{farne2020large,farne2024large} to learn the number of latent factors under a similar factor model setting.
However, compared to \cite{farne2020large,farne2024large}, in this paper we replace the Frobenius loss $\mathcal{L}^{(F)}(\mf{\Sigma},\mf{\Sigma}_n)$
with the log-det loss $\mathcal{L}^{(ld)}(\mf{\Sigma},\mf{\Sigma}_n)$, which is proven to work not worse, due to its ability to approximate underlying eigenvalues.
This paper also generalizes the distributional assumptions on factors and residuals, only imposing the existence of second moments.

Summing up, our approach is able to retrieve both the latent rank and the residual sparsity pattern with probability tending to $1$ for all $p \in \N$ as $n\to\infty$, starting from the sample covariance matrix generated from an approximate factor model, under generalized pervasiveness of latent factors and intermediate residual covariance sparsity. This is achieved by solving the single problem \eqref{obj}, which allows the simultaneous recovery of low rank and sparse components, as well of latent rank and residual sparsity pattern, %without the need to apply any further criterion. Furthermore, the residual component and sparsity pattern are recovered
without the need to apply a second-stage algorithm after a first stage recovering low rank component and latent rank.
%Problem \eqref{obj} is identified with increasing probability as $p \to \infty$.

%\far{END HERE}

\subsection{Paper structure}
The remainder of the paper is structured as follows.
In Section \ref{model} we define our model, detailing the necessary assumptions and summarizing the main results.
Section \ref{Cons_both} defines algebraic and parametric consistency.
Section \ref{Cons_alg} establishes the algebraic consistency of $\bigl(\widehat{\mf{L}}_n,\widehat{\mf{S}}_n\bigr)$.
Section \ref{Cons_par} establishes the parametric consistency of $\bigl(\widehat{\mf{L}}_n,\widehat{\mf{S}}_n\bigr)$.
Section \ref{practice} provides some computational details.
Section \ref{real} describes a real data example.
%explains how the adoption of $\mathcal{L}^{(ld)}(\mf{\Sigma},\mf{\Sigma}_n)$ as smooth term in \eqref{obj}
%impacts on algebraic consistency.
%
%and provides a probabilistic framework for it.
%
Concluding remarks follow in Section \ref{concl}.
The Supplement contains the statements and proofs of the technical results on the sample covariance matrix (Section A),
the mathematical and probabilistic analysis of $\mathcal{L}^{(ld)}(\mf{\Sigma},\mf{\Sigma}_n)$ (Section B),
the remaining proofs of mathematical results (Section C),
a comparison with \cite{agarwal2012noisy} (Section D),
the computational aspects of problem \eqref{obj} (Section E),
and a simulation study (Section F).
%\cite{suppl}.
%Appendix \ref{proofs}.

%inverse Frobenius loss.
%maybe generalized to Moore-Penrose?
%algorithm performance on the original problem
%recoverable eigenvalues

%{\textbf{A spectrum-based large covariance matrix estimator}}

%\far{\section{The log-det heuristics based on the estimation error matrix}}

%Before stating them, we need to introduce some matrix norm notation.

%% TO DO %%

%%\partial in lieu of $\delta$
%%^{*} via, dentro '
%nuclear: check index i!
%
%Move algorithm further (Section 6)
%MC criterion
%Threshold selection
%Re-scale or not?
%Un-shrinkage!
%
%Metric: $l_2$ error
%Eigenvalue error
%condition numbers, traces, norms?

%%

\textbf{Notation}.
Given a $p\times p$ symmetric positive semi-definite matrix $\mf{M}$,
we denote by $\lambda_i(\mf{M})$, $i\in\{1,\ldots, p\}$,
the eigenvalues of $\mf{M}$ in descending order.
To indicate that $\mf{M}$ is positive definite or semi-definite we use the notations $\mf{M} \succ 0$ or $\mf{M} \succeq 0$, respectively.
%Given a Hermitian matrix $\mf{M}$, its complex conjugate is denoted by $\mf{M}^{\dag}$.
Then, we recall the following norm definitions:
\begin{inparaenum}
\item Element-wise:
\begin{inparaenum}
\item[(i)]
$\ell_{0}$-norm: $\Vert \mf{M}\Vert_{0}=\sum_{i=1}^p \sum_{j=1}^p \mathbbm{1}( \mf{M}_{ij}\ne 0)$, which is the total number of non-zeros;
\item[(ii)] $\ell_{1}$-norm: $\Vert \mf{M}\Vert_{1}=\sum_{i=1}^p \sum_{j=1}^p \vert \mf{M}_{ij}\vert$;
\item[(iii)] Frobenius norm: $\Vert \mf{M}\Vert_{F}=\sqrt{\sum_{i=1}^p \sum_{j=1}^p  \mf{M}_{ij}^2}$;
\item[(iv)] Maximum norm: $\Vert \mf{M}\Vert_{\infty}=\max_{i\leq p, j \leq p} \vert \mf{M}_{ij}\vert$.
\end{inparaenum}
\item Induced by vector:
\begin{inparaenum}
\item[(i)] $\Vert \mf{M}\Vert_{0,v}=\max_{i \leq p} \sum_{j \leq p} \mathbbm{1}( \mf{M}_{ij} \ne 0)$,
which is the maximum number of non-zeros per row--column,
defined as the maximum `degree' of $ \mf{M}$; %`degree'
\item[(ii)] $\Vert \mf{M}\Vert_{1,v}=\max_{j \leq p} \sum_{i \leq p} \vert \mf{M}_{ij}\vert$;
%\item $\Vert \mf{M} \Vert_{\infty,v}= \max_{i \leq p} \sum_{j \leq p} \vert \mf{M}_{ij}\vert $;
\item[(iii)] Spectral norm: $\Vert \mf{M} \Vert_{2}=\lambda_{1}(\mf{M})$.
%\item $\ell_{1}$-norm: $\Vert \mf{M}\Vert_{1}=\sum_{i=1}^{p} \sum_{j=1}^{p} |\mf{M}_{ij}|$.
\end{inparaenum}
\item Schatten:
\begin{inparaenum}
\item[(i)] Nuclear norm of $ \mf{M}$, here defined as the sum of the eigenvalues of $ \mf{M}$:
$\Vert \mf{M}\Vert_{*}=\sum_{i=1}^p \lambda_i(\mf{M})$.%, which is.
\end{inparaenum}
\end{inparaenum}

The minimum nonzero off-diagonal element of $\mf{M}$ in absolute value is denoted as
$$\Vert \mf{M} \Vert_{\mathrm{min,off}}=\min_{{1\le i,j \leq p, \;
\mbox{s.t.} \, i \ne j \; \mbox{and} \; \mf{M}_{ij} \ne 0}}{\vert \mf{M}_{ij} \vert}.$$
Given a $p$-dimensional vector $\vf{v}$, we denote by:
$\Vert \vf{v} \Vert=\sqrt{\sum_{i=1}^p \vf{v}_i^2}$ the Euclidean norm of $\vf{v}$;
$\Vert \vf{v} \Vert_{\infty}=\max_{i=1,\ldots,p}{|\vf{v}_i|}$ the maximum norm of $\vf{v}$;
${v}^{i}_{k}$ the $i$-th component of the indexed vector $\vf{v}_k$.
%`modello all'inizio!
% radice di \Lambda_L!
Given two matrices $\mf{M}_{1}$ and $\mf{M}_{2}$ of equal size, we call $\mathcal{A}$ the addition operator, such that $\mathcal{A}(\mf{M}_{1},\mf{M}_{2})=\mf{M}_{1}+\mf{M}_{2}$,
and $\mathcal{A^{\dag}}$ the adjoint operator, such that $\mathcal{A^{\dag}}(\mf{M}_{1})=(\mf{M}_{1},\mf{M}_{1})$.
In light of these definitions, given two manifolds $\mathcal{T}'$ and $\Omega$, their orthogonal complements $\mathcal{T}'^{\perp}$ and $\Omega^{\perp}$,
their Cartesian sum $\mathcal{Y}=\mathcal{T}'\oplus\Omega$ and its orthogonal complement $\mathcal{Y}^{\perp}$,
the following identities hold:
\begin{inparaenum}
\item[] $\mathcal{A}^{\dag}\mathcal{A}(\mf{L},\mf{S})=(\mf{L}+\mf{S},\mf{L}+\mf{S})$;
\item[] $\mathbb{P}_{\mathcal{Y}}\mathcal{A}^{\dag}\mathcal{A}\mathbb{P}_{\mathcal{Y}}(\mf{L},\mf{S})=(\mathbb{P}_{\Omega}\mf{L}+\mf{S},\mf{L}+\mathbb{P}_{\mathcal{T}'}\mf{S})$;
\item[] $\mathbb{P}_{\mathcal{Y}^\perp}\mathcal{A}^{\dag}\mathcal{A}\mathbb{P}_{\mathcal{Y}}(\mf{L},\mf{S})=(\mathbb{P}_{\Omega^{\perp}}\mf{L}+\mf{S},\mf{L}+\mathbb{P}_{\mathcal{T}'^{\perp}}\mf{S})$,
\end{inparaenum}
where $\mathbb{P}$ is the projection operator, so that $\mathbb{P}_{\Omega}\mf{L}$ and $\mathbb{P}_{\Omega^{\perp}}\mf{L}$ are the orthogonal projections of $\mf{L}$ onto ${\Omega}$ and ${\Omega^{\perp}}$, while $\mathbb{P}_{\mathcal{T}'}\mf{S}$ and $\mathbb{P}_{\mathcal{T}'^{\perp}}\mf{S}$ are the orthogonal projections of $\mf{S}$ onto ${\mathcal{T}'}$ and ${\mathcal{T}'^{\perp}}$.

Given two sequences $A_\ell$ and $B_\ell$, $\ell \rightarrow \infty$,
we write $A_\ell=O(B_\ell)$ ($A_\ell \preceq B_\ell$), or $B_\ell=O(A_\ell)$ ($A_\ell \succeq B_\ell$), or $A_\ell \simeq B_\ell$,
if there exists a positive real $C$ independent of
$\ell$ such that $A_\ell/B_\ell\leq C$, or $B_\ell/A_\ell\leq C$, or
%More, we write $A_\ell \succeq B_\ell$,
%if $A_\ell/B_\ell \to \infty$ as $\ell \to \infty$,
%we write $A_\ell \preceq B_\ell$,
%if $B_\ell/A_\ell \to \infty$ as $\ell \to \infty$,
$A_\ell/B_\ell\leq C$ and $B_\ell/A_\ell\leq C$, respectively.
%%for all $p$ as $n\to\infty$
%%logica
%%linguaggi
Similarly,
we write $A_\ell=o(B_\ell)$ ($A_\ell \prec B_\ell$), or $B_\ell=o(A_\ell)$ ($A_\ell \succ B_\ell$),
if there exists a positive real $C$ independent of $\ell$ such that $A_\ell/B_\ell < C$ or $B_\ell/A_\ell < C$, respectively.
%More, we write $A_\ell \succeq B_\ell$,
%if $A_\ell/B_\ell \to \infty$ as $\ell \to \infty$,
%we write $A_\ell \preceq B_\ell$,
%if $B_\ell/A_\ell \to \infty$ as $\ell \to \infty$,
%and we write $A_\ell \simeq B_\ell$ if there exists a positive real $C$ independent of
%$\ell$ such that $A_\ell/B_\ell\leq C$ and $B_\ell/A_\ell\leq C$.
%%for all $p$ as $n\to\infty$

We indicate by $\mathrm{P}(.)$ the probability of the event in parentheses.
We denote by $O_P$ the big-O and by $o_P$ the small-O in probability.

\section{Theoretical foundations}\label{model}

\subsection{Model definition}\label{def}
We assume for the data the following factor model structure:
\noindent
\begin{equation}
\mathbf{x}_{k}=\mf{B}\mathbf{f}_{k}+\mathbf{\epsilon}_{k}\label{mod},
\end{equation}
where $\mf{B}$ is a $p \times r$ semi-orthogonal loading matrix
such that $\mf{B}'\mf{B}=\mf{\Lambda}_r$, with $\mf{\Lambda}_r$ $r \times r$ diagonal matrix and $\mf{B}'=[\vf{b}_{1} \; \ldots \; \vf{b}_p]$,
$\vf{f}_{k} \sim (\vf{0}_r,\mf{I}_r)$ is a $r\times 1$ random vector,
$\vf{\epsilon}_{k} \sim (\vf{0}_p, \mf{S^{*}})$ is a $p\times 1$ random vector, and $k=1,\ldots,n$ iterates over the samples.
Assuming that $\mathrm{E}(\vf{f}\vf{\epsilon}')=\mf{0}_{r \times p}$,
we obtain that \begin{equation}\mf{\Sigma}^{*}=\mathrm{E}(\vf{x}\vf{x}')=\mf{B}\mf{B}'+\mf{S^{*}}=\mf{L^{*}}+\mf{S^{*}},\label{cov}\end{equation}
where $\mf{L^{*}}=\mf{B}\mf{B}'$ is positive semi-definite with rank $r<p$, and
$\mf{S^{*}}$ is positive definite with $s<{p(p-1)}/{2}$ non-zero off-diagonal elements.
Equation \eqref{cov} is the low rank plus sparse decomposition of the covariance matrix of $\vf{x}$, $\mf{\Sigma}^{*}$.
We define the unbiased sample covariance matrix as $\mf{\Sigma}_n=n^{-1}\sum_{k=1}^n \vf{x}_k \vf{x}_k'$.

\subsection{Factorial structure}\label{assumption}

%Incorporate assumptions of the paper EJS(2022).
%Roughly speaking, latent eigenvalues spiked with $p^{\alpha}$,
%number of residual non-zeros spiked with $p^{\delta}$, $\delta<\alpha$.

%Previous assumptions: copy and paste in a condensed form.

\begin{assumption}\label{tails}
%Denoting by ${f}^{i}_{k}$ the $i-$th component of the vector $\vf{f}_{k}$,
%$\vf{f}_k \sim iid MNV(\vf{0}_r,\mf{I}_r)$, $\vf{\epsilon}_k \sim iid MNV(\vf{0}_p,\mf{S}^{*})$, $k=1,\ldots,n$.
In model \eqref{mod}, $\mathrm{E}(\vf{f})=\vf{0}_r$, $\mathrm{V}(\vf{f})=\mf{I}_r$,
%$ X, {\epsilon}$
$\mathrm{E}(\vf{\epsilon})=\vf{0}_p$, $\mathrm{V}(\vf{\epsilon})=\mf{S}^{*}$,
$\lambda_p(\mf{S}^{*})>0$, $\mathrm{E}(\vf{f}\vf{\epsilon}')=\vf{0}_{r \times p}$,
and there exist $\delta_{f},M_{f},\delta_{\epsilon},M_{\epsilon}>0$ independent of $p$
such that, for any $k\leq n$, $i \leq r$, $j\leq p$:
$$
\E(\vert{f}^{i}_{k}\vert^{2(1+\delta_f)}) \leq M_{f}, \qquad \E(\vert{\epsilon}^{j}_{k}\vert^{2(1+\delta_\epsilon)}) \leq M_{\epsilon}.
$$
%\Pr(\vert{f}_{k,i}\vert>l) \leq \exp\{-(l/{b_{1}})^{c_1}\}, \qquad \Pr(\vert{\epsilon}_{k,j}\vert>l) \leq \exp\{-(l/b_{2})^{c_2}\}. \nonumber
\end{assumption}
Assumption \ref{tails} renders model \eqref{mod} an approximate factor model, since
the residual covariance matrix $\mf{S}^{*}$ is positive definite, while factors and residuals are uncorrelated random vectors.
The moment conditions, borrowed from \cite{bickel2008covariance}, allow to control the sample size requirement to prove consistency for the sample covariance matrix $\mf{\Sigma}_n$.

\begin{assumption}\label{eigenvalues}
\begin{inparaenum}
\item [(i)] The eigenvalues of the $r\times r$ matrix $p^{-\alpha_{1}}\mf{B}' \mf{B}$
are %\sout{bounded away from $0$ for all $p\in\N$}
such that $\lambda_i(\mf{B}' \mf{B}) \simeq p^{\alpha_i}$, $i=1,\ldots,r$,
for some $0.5 < {\alpha_r} \leq \ldots \leq {\alpha_{1}} \leq 1$;
\item [(ii)] $\Vert \vf{b}_j\Vert _{\infty}=O(1)$ and $r$ is finite and independent of $p$ for all $j=1,\ldots,p$ and $p\in\N$.
\end{inparaenum}
\end{assumption}

Assumption \ref{eigenvalues} prescribes different speeds of divergence for latent eigenvalues,
and imposes a finite latent rank $r$. The left limit $O(p^{1/2})$ in part (i) is imposed on latent eigenvalues to preserve the
factor model structure as $p \to \infty$, because $p^{1/2}/\lambda_r(\mf{B}' \mf{B})=o(1)$ as $p \to \infty$. Part (ii) could actually be relaxed to cope with
$r=O(\ln(p))$, but we avoid it for the sake of simplicity. The maximum loading magnitude is imposed to be bounded
as $p \to \infty$, in order to control the sample covariance matrix error in maximum norm.

\begin{assumption}\label{sparsity}
\begin{inparaenum}
For all $p \in \N$, there exist $\delta_{1} \in (0,0.5]$ and $\delta_{2}>0$ independent of $p$ such that
\item [(i)] $\Vert \mf{S}^{*}\Vert_{0,v}\leq \delta_{2} p^{\delta_{1}}$;
%=\max_{1\le i \leq p}\sum_{j =1}^p \mathbbm{1}(\mf{S}^{*}_{ij} \ne 0)
\item [(ii)] $\Vert \mf{S}^{*} \Vert_{\infty} = O(1)$;
\item[(iii)] $p^{1-\delta_1}\Vert\mf{S}^{*}\Vert_{\mathrm{min,off}} = o(1)$;
\item[(iv)] $\sum_{j=1}^p{\mf{S}^{*}_{jj}}=o(p^{\alpha_1})$.
%\item [(ii)] $\delta_{1}'>0$, with $\delta_{1}'\leq \delta_{1}+\kappa_{0}$, and $\delta_{2}'>0$, such that $\Vert \mf{S}^{*} \Vert_{1,v} \leq \delta_{2}' p^{\delta_{1}'}$.
\end{inparaenum}
\end{assumption}
Assumption \ref{sparsity} explicitly controls the sparsity pattern of $\mf{S}^{*}$,
in that part (i) imposes a bound on the maximum number of non-zeros per row, part (ii) imposes the maximum element magnitude to be $O(1)$,
part (iii) requires the minimum magnitude of nonzero entries to disappear as $n \to \infty$, part (iv) ensures the factor model to be meaningful as $p \to \infty$,
i.e. $\mathrm{tr}(\mf{S}^{*})=o(\mathrm{tr}(\mf{L}^{*}))$, because
$\mathrm{tr}(\mf{S}^{*})=\sum_{j=1}^p{\mf{S}^{*}_{jj}}=o(p^{\alpha})=o(rp^{\alpha})$, and $\mathrm{tr}(\mf{L}^{*})=O(rp^{\alpha})$.
Moreover, parts (i) and (ii) together establish the traditional eigengap between $\lambda_r(\mf{B}' \mf{B})$ and $\lambda_1(\mf{S}^{*})$ as $p\to\infty$ (see \cite{chamberlain1983arbitrage}),
because $\Vert \mf{S}^{*}\Vert_{2} \leq \Vert\mf{S}^{*}\Vert_{1,v} \leq \Vert \mf{S}^{*}\Vert_{0,v} \Vert \mf{S}^{*} \Vert_{\infty} \leq \delta_{2} p^{\delta_{1}}$,
and $\delta_{1} \leq {1}/{2} < {\alpha_r}$ by Assumption \ref{eigenvalues}(i), while
parts (iii) and (iv) are able to ensure that the impact of the factorial component in \eqref{mod} always dominates over the residual component in the sample covariance matrix $\mf{\Sigma}_n$.
%Part (i) is also needed to ensure the identifiability of the sparsity pattern in $\mf{S}^{*}$,
%together with Assumption \ref{alg} (see later). \far{The sparsity framework of Assumption \ref{sparsity} explicitly controls the number of non-zeros in the residual component and their position.}

%sub-Gaussian tails,
%which also ensures that all the moments of $\vf{f}$ and $\vf{\epsilon}$ exist, thus allowing to apply to $\vf{f}$ and $\vf{\epsilon}$ large deviation probabilistic results.
%Note that the magnitude of residual covariances (i.e., the off-diagonal entries of $\mf{S}^{*}$) is controlled by Assumption \ref{sparsity}.

\subsection{Geometry}\label{geometry}

In order to ensure the effectiveness of the composite penalty
$\psi \Vert\mf{L}\Vert_{*}+\rho \Vert\mf{S}\Vert_{1}$
in recovering the latent rank $\mathrm{rk}(\mf{L}^{*})=r$ and
the residual number of nonzeros $\vert\mathrm{supp}(\mf{S}^{*})\vert=s$
(where $\mathrm{supp}(\mf{S}^{*})$ is the support of $\mf{S}^{*}$, denoting the location of nonzeros,
%\sout{orthogonal complement}
%\sout{$ker(\mf{S}^{*})$}
and $\vert\mathrm{supp}(\mf{S}^{*})\vert$ is its dimension),
we need to control the geometric manifolds containing
$\mf{L}^{*}$ and $\mf{S}^{*}$.
As in \cite{chandrasekaran2011rank}, we assume $\mf{L}^{*}\in\mathcal{L}(r)$
and $\mf{S}^{*}\in\mathcal{S}(s)$, where
\begin{eqnarray}
&&\mathcal{L}(r) =  \{\mf{L} \mid \mf{L} \succeq 0, {\mf{L}}={\mf{U}\mf{D}\mf{U}'},
\mf{U} \in \R^{p \times r}, \mf{U}'\mf{U}=\mf{I}_r, \mf{D} \in \R^{r \times r} \mathrm{diagonal}\},\label{var:L}\\
&&\mathcal{S}(s) =  \{\mf{S}\in \R^{p\times p} \mid \mf{S} \succ 0, \vert \mathrm{supp}(\mf{S})\vert \leq s\}.\label{var:S}
\end{eqnarray}
$\mathcal{L}(r)$ is the algebraic variety of matrices with at most rank $r$,
$\mathcal{S}(s)$ is the algebraic variety of (element-wise) sparse matrices with
at most $s$ non-zero elements, and the two varieties
$\mathcal{L}(r)$ and $\mathcal{S}(s)$
can be disentangled if $\mf{L}^{*}$ is far from being sparse,
and $\mf{S}^{*}$ is far from being low rank.
For this reason, \cite{chandrasekaran2011rank} defines
the tangent spaces $\mathcal{T}(\mf{L}^{*})$ and $\Omega(\mf{S}^{*})$ to $\mathcal{L}(r)$ and $\mathcal{S}(s)$ as follows:
\begin{eqnarray}
%\mathcal{T}(\mf{L}^{*})=\{\mf{M}\in\R^{p \times p}|\mf{M}=\mf{U} \mf{Y}_1^\top+\mf{Y}_2 \mf{U}^\top \mid \mf{Y}_1,\mf{Y}_2
%\in \R^{p \times r}\}, \mf{U}'\mf{U}=\mf{I}_r, \mf{U}'\mf{L}^{*}\mf{U}\in\R^{r \times r} \mbox{diagonal}\nonumber\\
\mathcal{T}(\mf{L}^{*})=\{\mf{M}\in \R^{p \times p} \mid \mf{M}=\mf{U} \mf{Y}_1'+\mf{Y}_2 \mf{U}',
\mf{Y}_1, \mf{Y}_2 \in \R^{p \times r}, \mf{U}\in \R^{p\times r}, \mf{U}' \mf{U}=\mf{I}_r;\nonumber\\
\mf{U}' \mf{L}^{*} \mf{U} \in \R^{r \times r} \mbox{diagonal}\},\qquad
\Omega(\mf{S}^{*})=\{\mf{N} \in \R^{p \times p} \mid \mathrm{supp}(\mf{N})\subseteq \mathrm{supp}(\mf{S}^{*})\},\nonumber
%, \far{\mathrm{diag}(\mf{N})=\vf{0}}
\end{eqnarray}
and proposes the following rank-sparsity measures:
\begin{eqnarray}
\xi(\mathcal{T}(\mf{L}^{*})) = \max_{\mf{M} \in \mathcal{T}(\mf{L}^{*}), \Vert\mf{M}\Vert_{2} \leq 1} {\Vert\mf{M}\Vert_\infty},
\qquad%\label{xi}\\
\mu(\Omega(\mf{S}^{*})) = \max_{\mf{N} \in \Omega(\mf{S}^{*}),\Vert\mf{N}\Vert_\infty \leq 1}\
{\Vert\mf{N}\Vert_{2}}.\nonumber%\label{mu}.%\pause \vspace{0.2cm}
\end{eqnarray}
Note that $\xi(\mathcal{T}(\mf{L}^{*}))$ is normalized to attain maximum $1$, irrespectively of the underlying assumptions
on latent eigenvalues.

%\far{Spiegare perchè la definizione di $\xi(\mathcal{T}(\mf{L}^{*}))$ ha senso anche se gli autovalori di $\mf{L}^{*}$ divergono con $p^{\alpha_1}$.}
%\far{Spiegare perchè la condizione $\far{\mathrm{diag}(\mf{N})=\vf{0}}$, dopo aver rivisto la funzione \eqref{obj} rimuovendo dalla norma $\ell_1$ di $\mf{S}$ la diagonale di $\mf{S}$,
%e rivedendo "Overview of results" di conseguenza.}

Then, according to \cite{chandrasekaran2011rank}, the identifiability condition to be satisfied
requires a bound on $\xi(\mathcal{T}(\mf{L}^{*}))\mu(\Omega(\mf{S}^{*}))$.
%(see Section \ref{Cons_both} for more details).
For this reason, recalling from \cite{chandrasekaran2011rank} that
%\begin{eqnarray}
%\mathrm{inc}(\mf{L}^{*}) \leq \xi(\mathcal{T}(\mf{L}^{*})) \leq  2 \mathrm{inc}(\mf{L}^{*}),\nonumber\\%\label{dualL}\\
%\mathrm{deg}_{min}(\mf{S}^{*}) \leq \mu(\Omega(\mf{S}^{*})) \leq \mathrm{deg}_{max}(\mf{S}^{*}),\nonumber%\label{dualS}
%\end{eqnarray}
%where
\begin{compactenum}
\item[$\bullet$] $\mathrm{inc}(\mf{L}^{*})=\max_{j=1,\ldots,p}\Vert \mathbb{P}_{{L}^{*}} \vf{e}_j \Vert$, where $\vf{e}_j$ is the $j$-th canonical basis vector, and $\mathbb{P}_{{L}^{*}}$ is the projection operator onto the row--column space of $\mf{L}^{*}$, with $\sqrt{{r}/{p}}
    \leq \mathrm{inc}(\mf{L}^{*}) \leq 1$;
\item[$\bullet$] $\mathrm{deg}_{min}(\mf{S}^{*})=\min_{1\le i \leq p} \sum_{j =1}^p \mathbbm{1}(\mf{S}^{*}_{ij} \ne 0)$ and
$\mathrm{deg}_{max}(\mf{S}^{*})=\Vert \mf{S}^{*}\Vert_{0,v}$,%=\max_{1\le i \leq p} \sum_{j =1}^p \mathbbm{1}(\mf{S}^{*}_{ij} \ne 0)
\end{compactenum}
with
%\begin{eqnarray}
$\rm{inc}(\mf{L}^{*}) \leq \xi(\mathcal{T}(\mf{L}^{*})) \leq  2 \rm{inc}(\mf{L}^{*})$ and %\;\mbox{and}\;\nonumber\\
$\rm{deg}_{min}(\mf{S}^{*}) \leq \mu(\Omega(\mf{S}^{*})) \leq \rm{deg}_{max}(\mf{S}^{*})$,
%,\nonumber
%\end{eqnarray}
%and
%\begin{equation}\rm{inc}(\mf{L}^{*}) \leq \xi(\mathcal{T}(\mf{L}^{*})) \leq 2\rm{inc}(\mf{L}^{*}),\label{inc2}\end{equation}
we can control the degree of transversality between $\mathcal{L}(r)$ and $\mathcal{S}(s)$ by the following assumption.

\begin{assumption}\label{alg}
For all $p\in\mathrm N$, there exist $\kappa_{L},\kappa_{S}>0$ independent of $p$, %with $k_L\geq{\sqrt{r}}/{2}$,
with $\kappa_{S}\leq\delta_{2}$, $\kappa_{S}{\sqrt{r}}/\kappa_{L}\leq 1/24$, $(\kappa_{L} p)^{\delta_{1}}\leq \sqrt{p}$,
such that $\xi(\mathcal{T}(\mf{L}^{*}))={\sqrt r}/(\kappa_{L} p)^{\delta_{1}}$ and $\mu(\Omega({\mf{S^{*}}}))=\kappa_{S} p^{\delta_{1}}$ for all $p \in \N$.
%\item[(iii)] $\underline{\delta},\delta_{2}'$ with $\underline{\delta}\leq \delta_{1}+\frac{1}{2}$ such that $\Vert S_{\infty} \Vert_{1,v} \leq \delta_{2}' p^{\underline{\delta}}$.
\end{assumption}
Assumption \ref{alg} states that the maximum degree of $\mf{S^{*}}$ is $O(p^{\delta_1})$,
where $\delta_1 < \alpha_r$ by Assumptions \ref{eigenvalues}(i) and \ref{sparsity}(i).
%which means that tends to $0$ as $p\to\infty$.
More, the incoherence of $\mf{L}^{*}$ is assumed to scale to $O(p^{-\delta_1})$, in order to keep the product $\xi(\mathcal{T}(\mf{L}^{*}))\mu(\Omega(\mf{S}^{*}))$ proportional to $O(1)$, which is crucial to prove algebraic consistency (see Section \ref{Cons_both}). %Theorem \ref{thm_main}
This assumption resembles in nature the approximate factor model of \cite{chamberlain1983arbitrage}, because $\delta_1<\alpha_r$, such that the number of residual nonzeros will become negligible with respect to latent eigenvalues, and the manifold underlying $\mf{L}^{*}$
will be progressively easier to retrieve as $p\to\infty$.

We now define the function %$$f_{\delta_\epsilon}(p,n)=\frac{p^{2/(1+\delta_\epsilon)}}{n^{1/2}},$$
$$f_{\delta_\epsilon}(p,n)={\frac{\max\left(p^{{2}/{(1+\delta_\epsilon)}},\sqrt{\ln(p)}\right)}{n^{1/2}}},$$
with $\delta_\epsilon \in \R^{+}$ as defined in Assumption \ref{tails},
%$$f_{\delta_f}(p,n)=
%\bigl\{
%\begin{array}{rl}
%\frac{r^{2/(1+\delta_f)}}{n^{1/2}}, & \delta_f \in \R^{+}, \nonumber\\
%\frac{\ln(r)}{n}, & \delta_f=+\infty, \nonumber
%\end{array}
%\bigr.
%$$
%and
%$$f_{\delta_\epsilon}(p,n)=
%\left\{
%\begin{array}{rl}
%p^{\frac{2}{1+\delta_\epsilon}}\sqrt{\frac{\ln(p)}{n}}, & \delta_\epsilon \in \R^{+}, \nonumber\\
%\sqrt{\frac{\ln(p)}{n}}, & \delta_\epsilon=+\infty. \nonumber
%\end{array}
%\right.
%$$
and the scalar $\psi_{0}(\delta_\epsilon,p,n)={f_{\delta_\epsilon}(p,n)}/{\xi(\mathcal{T}(\mf{L}^{*}))}$
(hereafter, we will write $\psi_{0}$ to avoid notation burden).

\begin{assumption}\label{lowerbounds}
There exist $\delta_L,\delta_S>0$ such that
\begin{inparaenum}
\item[(i)] the minimum eigenvalue of $\mf{L}^{*}$, $\lambda_r(\mf{L}^{*})$, is greater than
$\delta_L{\psi_{0}}/{\xi^2(\mathcal{T}(\mf{L}^{*}))}$
\item[(ii)] the minimum absolute value of the non-zero off-diagonal entries of $\mf{S}^{*}$,
$\Vert\mf{S}^{*}\Vert_{\mathrm{min,off}}$, is greater than $\delta_S{\psi_{0}}/{\mu(\Omega(\mf{S}^{*}))}$.
\end{inparaenum}
\end{assumption}

Assumption \ref{lowerbounds} is crucial for identifiability, as it guarantees
that the solution pair of \eqref{obj} lies on the ``right" manifolds, i.e. that $\widehat{\mf{L}}_n\in\mathcal{L}(r)$ and
$\widehat{\mf{S}}_n\in\mathcal{S}(s)$ with high probability as $n\to\infty$.
We remark that $\psi_{0}$, which controls the error bound of the solution pair $\bigl(\widehat{\mf{L}}_n,\widehat{\mf{S}}_n\bigr)$,
directly depends on the moment condition imposed on the vector of residuals $\vf{\epsilon}$ by Assumption \ref{tails}.

\subsection{Main results}\label{main_res}

Suppose that model \eqref{mod} is a strict factor model with pervasive factors, i.e., for all $p\in\mathrm N$:
\begin{compactenum}
\item[(i)] \textbf{tails}: the elements of factor score vectors $\vf{f}_k$ and residual vectors $\vf{\epsilon}_k$, $k=1,\ldots,n$, present sub-exponential tails (meaning $\delta_f=\delta_\epsilon=+\infty$),
i.e., there exist\\ $b_{1},b_{2},c_{1},c_{2}>0$ independent of $p$ such that, for any $l>0$, $k \in \{1,\ldots,n\}$, $i \in \{1,\ldots,r\}$, $j \in \{1,\ldots,p\}$:
\begin{eqnarray}
\Pr(\vert{f}^{i}_{k}\vert>l) \leq \exp\{-(l/{b_{1}})^{c_1}\}, \qquad
\Pr(\vert{\epsilon}^{j}_{k}\vert>l) \leq \exp\{-(l/b_{2})^{c_2}\}; \nonumber
\end{eqnarray}
\item[(ii)] \textbf{latent eigenvalues}: %s.t. $\far{\Vert{\mf{S}^{*}}\Vert_{\infty}=\kappa_v}$,
%with $\mathrm{E}(\vf{\epsilon})=\vf{0}_p$, $\mathrm{V}(\vf{\epsilon})=\mf{I}_{p}$,
%$\lambda_p(\mf{S}^{*})>0$, $
%\mathrm{E}(\vf{f}\vf{\epsilon}')=\vf{0}_{r \times p}$
$\lambda_i(\mf{L}^{*}) \simeq O(p)$ for $i=1,\ldots,r$ (which means $\alpha_1=1,\ldots,\alpha_r=1$), $r=O(\ln(p))$ and $\Vert \vf{b}_j\Vert _{\infty}=O(1)$ for all $j=1,\ldots,p$;
\item[(iii)] \textbf{residual component}: $\mf{S}^{*}$ is diagonal, which implies that $\mu(\mf{S}^{*})=1$;
\item[(iv)] \textbf{incoherence}: $\mathrm{inc}(\mf{L}^{*})=(\widetilde{\kappa}_{L} r/p)^{1/3}$, with $\widetilde{\kappa}_{L} \in \R$ s.t. $(\widetilde{\kappa}_{L} r/p)^{1/3}\geq \sqrt{r/p}$,
and $\xi(\mathcal{T}(\mf{L}^{*}))={(\widetilde{\kappa}_{L} r/p)}^{1/3}$. Consequently, $\psi_{0}=\sqrt{\ln(p)/n}/
{(\widetilde{\kappa}_{L} r/p)}^{1/3}$;
%there exist $\widetilde{\kappa}_{L}>0$ such that $2\mathrm{inc}(\mf{L}^{*})={\sqrt r}/(\widetilde{\kappa}_{L})$,
%\far{with $\mathrm{inc}(\mf{L}^{*})=\max_{j=1,\ldots,p}\Vert \mathbb{P}_{{L}^{*}} \vf{e}_j \Vert$,
%where $\vf{e}_j$ is the $j$-th canonical basis vector, and $\mathbb{P}_{{L}^{*}}$ is the projection operator onto the row--column space of $\mf{L}^{*}$};
\item[(v)]  \textbf{minimum latent eigenvalue}: there exist $\delta_L>0$ such that the minimum eigenvalue of $\mf{L}^{*}$, $\lambda_r(\mf{L}^{*})$, is greater than $\delta_L{\psi_{0}}/{\xi^2(\mathcal{T}(\mf{L}^{*}))}$, i.e., greater than $\frac{\delta_L p\sqrt{\ln(p)/n}}{\widetilde{\kappa}_{L} r}$.
%$\delta_L\widetilde{\kappa}_{L}^{3}\sqrt{\ln(p)/n}/(r\sqrt{r})$.
%with $C_r=2\sqrt{2}\widetilde{\kappa}_{L}^{3/2}$,
\item[(vi)]  \textbf{sample size requirement}: since it must hold $\frac{\delta_L p\sqrt{\ln(p)/n}}{\widetilde{\kappa}_{L} r}=o(p)$
to cope with $\alpha_r=1$, it holds $\ln(p)/n \to 0$ as $n\to\infty$.
%by setting $r=O(\ln(p))$ we get
%$$\frac{\delta_L \sqrt{\ln(p)}}{\widetilde{\kappa}_{L} \sqrt{\ln(p)}}=o(\sqrt{n}),$$
%which is automatically satisfied as $n \to \infty$.
%which means $\sqrt{n} \geq \delta_L/\widetilde{\kappa}_{L}$, leading to $n=O(1)$.
\end{compactenum}

\begin{theorem}\label{thm_simple}
%Let us suppose $\widetilde{\kappa}_{L} \geq 54 \sqrt{r}$.
Let us set $\psi_{0}=\sqrt{\ln(p)/n}/(\widetilde{\kappa}_{L} r/p)^{1/3}$,
$\gamma \in \left[2{(\widetilde{\kappa}_{L} r/p)^{1/3}},{1}/{4}\right]$, $\rho_{0}=\gamma\psi_{0}$,
%with $\gamma \in [9\sqrt{r}/\widetilde{\kappa}_{L},1/6]$,
with $p \geq 512 \widetilde{\kappa}_{L} r$, $\psi=p\psi_{0}$, and $\rho=\rho_0$,
%$\rho=p^{\delta}\rho_{0}=p^{\delta} \gamma \psi_{0}$,
where $\psi$ and $\rho$ are the thresholds in \eqref{obj}.
Suppose that conditions (i)--(vi) hold.
Then, %there exists a positive real $\kappa$ independent of $p$ and $n$ such that,
for all $p \in \N$ such that $\ln(p)/n \to 0$, as $n\to\infty$ the pair of solutions defined in \eqref{obj} satisfies:
\begin{inparaenum}
\item[(i)] $\mathrm{P} \left(\frac{\Vert\widehat{\mf{L}}_n-\mf{L}^{*}\Vert_{2}}{{p^{4/3}}} \leq 5.3125 \frac{\sqrt{\ln(p)/n}}{(\widetilde{\kappa}_{L} r)^{1/3}}\right) \to 1$; %21.5/2/2
\item[(ii)] $\mathrm{P} \left(\Vert\widehat{\mf{S}}_n-\mf{S}^{*}\Vert_{\infty} \leq 5.3125 \frac{\sqrt{\ln(p)/n}}{(\widetilde{\kappa}_{L} r)^{1/3}}\right) \to 1$; %21.5/2/2
\item[(iii)] $\mathrm{P} (\mathrm{rk}(\widehat{\mf{L}}_n)=\mathrm{rk}(\mf{L}^{*})) \to 1$;
\item[(iv)] $\mathrm{P} (\mathrm{sgn}(\widehat{\mf{S}}_n)=\mathrm{sgn}(\mf{S}^{*})) \to 1$.
\end{inparaenum}
\end{theorem}

\begin{corollary}\label{corollsimple}
Under the assumptions of Theorem \ref{thm_simple},
%there exists a positive real $C$ independent of $p$ and $n$
%such that,
for all $p \in \N$, as $n\to\infty$ it holds:
%$\Vert \widehat{\mf{S}}_{\rm{ALCE}}-\mf{S}^{*}\Vert _{2}~\leq \phi_S $, $\Vert \widehat{{\mf{\Sigma}}}_{\rm{ALCE}}-{\mf{\Sigma}}^{*}\Vert _{2} \leq \phi$, $\Vert \widehat{\mf{S}}_{\rm{ALCE}}^{-1}-\mf{S}^{*-1}\Vert _{2} \leq \phi_S$,
%and $\Vert \widehat{{\mf{\Sigma}}}_{\rm{ALCE}}^{-1}-{\mf{\Sigma}}^{*-1}\Vert _{2} \leq \psi$, with probability approaching $1$.
\begin{inparaenum}
\item[(i)] $\mathrm{P} \left(\frac{\Vert \widehat{\mf{S}}_n-\mf{S}^{*}\Vert_{2}}{p{^{1/3}}}\leq 5.3125 \frac{\sqrt{\ln(p)/n}}{(\widetilde{\kappa}_{L} r)^{1/3}}\right) \to 1$;
\item[(ii)] $\mathrm{P} \left(\frac{\Vert \widehat{{\mf{\Sigma}}}_n-{\mf{\Sigma}}^{*}\Vert_{2}}{{p^{4/3}}}\leq 5.3125 \frac{\sqrt{\ln(p)/n}}{(\widetilde{\kappa}_{L} r)^{1/3}}\right) \to 1$;
\item[(iii)] $\mathrm{P} (\lambda_p(\widehat{\mf{S}}_n)>0) \to 1$;
\item[(iv)] $\mathrm{P} (\lambda_p(\widehat{\mf{\Sigma}}_n)>0) \to 1$.
%\item $\mathrm{P} (\Vert \widehat{\mf{S}}_{\rm{ALCE}}^{-1}-\mf{S}^{*-1}\Vert_{2} \leq C \phi_S) \to 1$;
%\item $\mathrm{P} (\Vert \widehat{{\mf{\Sigma}}}_{\rm{ALCE}}^{-1}-{\mf{\Sigma}}^{*-1}\Vert_{2} \leq C \phi) \to 1$.
\end{inparaenum}
In addition, supposing that {$\lambda_p(\mf{S}^{*})=O(p^{-\varepsilon})$ and $\lambda_p({\mf{\Sigma}}^{*})=O(p^{\-1-\varepsilon})$} for some $\varepsilon>0$,
%and $\sum_{i=1}^p\sum_{j=1}^p \mathbbm{1}(\mf{S}_{ij}^{*} \ne 0) \leq O(p^{2\delta})$,
the following statements hold for all $p \in \N$ as $n \to \infty$:
\begin{inparaenum}
\item[(v)] $\mathrm{P} \left(\frac{\Vert\widehat{\mf{S}}_n^{-1}-\mf{S}^{*-1}\Vert_{2}}{p^{1/3+2\varepsilon}} \leq 5.3125 \frac{\sqrt{\ln(p)/n}}{(\widetilde{\kappa}_{L} r)^{1/3}}\right)\to 1$;
\item[(vi)] $\mathrm{P} \left(\frac{\Vert\widehat{{\mf{\Sigma}}}_n^{-1}-{\mf{\Sigma}}^{*-1}\Vert_{2}}{p^{4/3+2\varepsilon}}
\leq 5.3125 \frac{\sqrt{\ln(p)/n}}{(\widetilde{\kappa}_{L} r)^{1/3}} \right) \to 1$.
%\item if $\lambda_{p}({\mf{S}^{*}})>\phi_{S}$, then  $\widehat{\mf{S}}_{\rm{ALCE}}$
%is positive definite;
%\item if $\lambda_{p}({{\mf{\Sigma}}^{*}})>\phi$, then $\widehat{{\mf{\Sigma}}}_{\rm{ALCE}}$ is positive definite;
%\item if $\lambda_{p}({\mf{S}^{*}}) \geq 2 \phi_{{S}}$, then $\widehat{\mf{S}}_{\rm{ALCE}}^{-1}$ is positive definite;
%\item if $\lambda_{p}({{\mf{\Sigma}}^{*}}) \geq 2 \phi$, then $\widehat{{\mf{\Sigma}}}_{\rm{ALCE}}^{-1}$ is positive definite.
\end{inparaenum}
\end{corollary}

Theorem \ref{thm_simple} and Corollary \ref{corollsimple} contain many powerful results.
First, problem \eqref{obj} is able to ensure that $\widehat{\mf{L}}_n$ is rank-consistent (part (iii) of Theorem \ref{thm_simple}),
and $\widehat{\mf{S}}_n$ is sparsistent (part (iv) of Theorem \ref{thm_simple}). This is particularly important in factor modelling,
because the latent rank is hard to identify in the presence of weak factors in high dimensions.
This identifiability result is obtained thanks to the control of the tangent spaces (via $\xi(\mathcal{T}(\mf{L}^{*}))$ and $\mu(\Omega(\mf{S}^{*}))$) %$\mathcal{T}(\mf{L}^{*})$ and $\Omega(\mf{S}^{*})$
to the low rank and the sparse matrix varieties containing $\mf{L}^{*}$ and $\mf{S}^{*}$ respectively, and is ensured under the mild condition $p \geq 512 \widetilde{\kappa}_{L} r$,
with $\widetilde{\kappa}_{L}$ possibly smaller than $1$ (as long as $(\widetilde{\kappa}_{L} r/p)^{-1/3} \geq \sqrt{r/p}$). Note that the corresponding condition for identifiability in \cite{chandrasekaran2012latent}
is in the best case $p \geq 54^2 \widetilde{\kappa}_{L} r$ with $\widetilde{\kappa}_{L} \geq 1$, which actually requires a super-large dimension.

Then, parts (i) and (ii) of  Theorem \ref{thm_simple}, as well as parts (i) and (ii) of Corollary \ref{corollsimple}, contain convergence results in spectral norm
for the covariance matrix estimates $\widehat{\mf{L}}_n$, $\widehat{\mf{S}}_n$ and $\widehat{\mf{\Sigma}}_n=\widehat{\mf{L}}_n+\widehat{\mf{S}}_n$,
which work under the mild condition $\ln(p)/n \to 0$ as $n \to \infty$ for all $p \in \N$.
Therefore, a situation is configured where the sample size requirement is almost negligible, and a high dimension $p$ is encouraged,
although identifiability still works for relatively small values of $p$, provided that the latent rank $r$ is very small.
Parts (iii) and (iv) ensure that $\widehat{\mf{S}}_n$ and $\widehat{\mf{\Sigma}}_n$ are positive definite, which means that their inverses can be
safely obtained when necessary (e.g., to estimate factor scores). Finally, we also obtain the convergence rates for the two inverses (parts (v) and (vi) of Corollary \ref{corollsimple}).

All these results are obtained with a little inconvenience: an additional term $p^{1/3}$ in the convergence rates in spectral norm,
which is due to the assumption $\xi(\mathcal{T}(\mf{L}^{*}))={(\widetilde{\kappa}_{L} r/p)}^{1/3}$. This is the price to pay to solve
a single-step problem like \eqref{obj} and obtain all these results at the same time.

In the following, we extend Theorem \ref{thm_simple} and Corollary \ref{corollsimple}, by generalizing the divergence rates of latent eigenvalues (allowing for weak factors),
by permitting a non-diagonal residual component (element-wise sparse), by tolerating that the distributions of factors and residuals do not possess all moments,
and by exploring what happens to the identifiability of underlying manifolds and to the convergence rates when $\xi(\mathcal{T}(\mf{L}^{*}))$ and $\mu(\mf{S}^{*})$
present a more general, less favorable shape.

\section{Consistency definitions}\label{Cons_both}

A distinctive feature of the nuclear norm plus $\ell_1$-norm penalization approach is to obtain a double type of consistency for covariance matrix estimators:
the usual parametric one (in a suitable norm) and the algebraic one, which is defined as follows.
\begin{definition}\label{def_alg}
A sequence of estimator pairs $\bigl(\widehat{\mf{L}}_n,\widehat{\mf{S}}_n\bigr)$, defined as in \eqref{obj_all} and indexed in $n$,
%with $\mf{S}, \mf{L} \in \R^{p\times p}$
is an algebraically consistent estimate of the low rank plus sparse decomposition \eqref{cov}
for the covariance matrix $\mf{\Sigma}^{*}$ if the following conditions hold with probability approaching $1$ as $n \to \infty$:
\begin{compactenum}
  \item the low rank estimate $\widehat{\mf{L}}_n$ is positive semidefinite and rank-consistent, i.e., $\mathrm{rk}(\widehat{\mf{L}}_n)=r$;%\mathrm{rk}(\mf{L}^{*})=r$;
  \item the residual estimate $\widehat{\mf{S}}_n$ is positive definite and sparsistent, i.e., $\mathrm{sgn}(\widehat{\mf{S}}_n)=\mathrm{sgn}(\mf{S}^{*})$;
  \item $\widehat{\mf{\Sigma}}_n=\widehat{\mf{L}}_n+\widehat{\mf{S}}_n$ is positive definite.
\end{compactenum}
\end{definition}

We now recall that
$\mathcal{P}(\mf{L},\mf{S})=\psi \Vert\mf{L}\Vert_{*}+\rho \Vert\mf{S}\Vert_{1},$
where $\Vert\mf{L}\Vert_{*}=\sum_{i=1}^p \lambda_i(\mf{L})$,
%is the nuclear norm of $\mf{L}$,
%i.e. the sum of the eigenvalues of $\mf{L^{*}}$,
$\Vert\mf{S}\Vert_{1}=\sum_{i=1}^{p} \sum_{j=1}^p \vert\mf{S}_{ij}\vert$,
$\psi$ and $\rho$ are non-negative threshold parameters.
From Section \ref{geometry}, we get that $\psi_{0}={f_{\delta_\epsilon}(p,n)}/{\xi(\mathcal{T}(\mf{L}^{*}))}$,
and we set $\rho_{0}=\gamma\psi_{0}$, %where $\gamma$ lies in the range of Proposition \ref{11},
$\psi=p^{\alpha_1}\psi_{0}$, and $\rho=\rho_0$.
%$\rho=p^{\delta}\rho_{0}=p^{\delta} \gamma \psi_{0}$,
%where $\psi$ and $\rho$ are the threshold parameters.
%observing that the low rank target $\mf{L}^{*}$ is such that $\Vert\mf{L}^{*}\Vert_{2}=O(p^{\alpha_1})$,
Then, we can rewrite $\mathcal{P}(\mf{L},\mf{S})$ as
%\mathcal{P}(\mf{L},\mf{S})&=&\psi\Vert\mf{L}\Vert_{*}+\rho \Vert\mf{S}\Vert_{1}\nonumber\\
$\mathcal{P}(\mf{L},\mf{S})=p^{\alpha_1}\psi_0{\Vert\mf{L}\Vert_{*}}+\rho_0{\Vert\mf{S}\Vert_{1}}$.
More, defining $\mathcal{P}_{\gamma}(\mf{L},\mf{S})=\psi^{-1}\mathcal{P}(\mf{L},\mf{S})$, we obtain
$\mathcal{P}_{\gamma}(\mf{L},\mf{S}) = {\Vert\mf{L}\Vert_{*}}+p^{-\alpha_1}\gamma\Vert\mf{S}\Vert_{1}$.
The dual norm of $\mathcal{P}_{\gamma}(\mf{L},\mf{S})$ can thus be defined as
$g_\gamma(\mf{L},\mf{S})=\max\left(\frac{\Vert\mf{L}\Vert_{2}}{p^{\alpha_1}},\frac{\Vert\mf{S}\Vert_{\infty}}{\gamma}\right)$.
Such definition is meaningful because, under Assumptions \ref{eigenvalues}(i) and \ref{sparsity}(i) with $\delta_{1}=0$,
we obtain that $\frac{\Vert\mf{L}\Vert_{2}}{p^{\alpha_1}}=O(1)$ and $\frac{\Vert\mf{S}\Vert_{\infty}}{\gamma}=O(1)$,
which renders the two quantities perfectly comparable.

Relying on these considerations, we focus on the solution of the following equivalent version of problem \eqref{obj}:
\begin{eqnarray}
%\bigl(\widehat{\mf{L}}_n,\widehat{\mf{S}}_n\bigr)=
%\arg\min_{{\mf{L},\mf{S},\\ \mf{L} \succeq 0, \mf{S} \succ 0,  \mf{L}+\mf{S} \succ 0}}
%%\phi(\mf{L},\mf{S})=
%\frac{1}{2}\ln \det (\mf{I}_p+\mf{\Delta}_{n}\mf{\Delta}_{n}')
%+\psi \Vert\mf{L}\Vert_{*}+\rho \Vert\mf{S}\Vert_{1}.
\bigl(\widehat{\mf{L}}_n,\widehat{\mf{S}}_n\bigr)=
\argmin_{\mf{L},\mf{S}}
%\phi(\mf{L},\mf{S})=
0.5\ln \det (\mf{I}_p+p^{-2\alpha_1}\mf{\Delta}_{n}\mf{\Delta}_{n}')
+\psi_0 \Vert\mf{L}\Vert_{*}+\rho_0 \frac{\Vert\mf{S}\Vert_{1}}{p^{\alpha_1}},
\end{eqnarray}
which may also be expressed as
\begin{eqnarray}
%\bigl(\widehat{\mf{L}}_n,\widehat{\mf{S}}_n\bigr)=
%\arg\min_{{\mf{L},\mf{S},\\ \mf{L} \succeq 0, \mf{S} \succ 0,  \mf{L}+\mf{S} \succ 0}}
%%\phi(\mf{L},\mf{S})=
%\frac{1}{2}\ln \det (\mf{I}_p+\mf{\Delta}_{n}\mf{\Delta}_{n}')
%+\psi \Vert\mf{L}\Vert_{*}+\rho \Vert\mf{S}\Vert_{1}.
\bigl(\widehat{\mf{L}}_n,\widehat{\mf{S}}_n\bigr)=
\argmin_{\mf{L},\mf{S}}
%\phi(\mf{L},\mf{S})=
0.5\ln \det (\mf{I}_p+p^{-2\alpha_1}\mf{\Delta}_{n}\mf{\Delta}_{n}')
+\Vert\mf{L}\Vert_{*}+\gamma \frac{\Vert\mf{S}\Vert_{1}}{p^{\alpha_1}}.\label{prob_ld_rescaled}
\end{eqnarray}
The consistency norm $g_\gamma$,
with which the direct sum $\mathcal{L}(r)\oplus\mathcal{S}(s)$ is naturally equipped (cfr. paragraph 3.3 in \cite{chandrasekaran2012latent}),
%which is the dual norm of the composite penalty $\psi_{0}\Vert\cdot\Vert_{*}+\rho_{0}\Vert\cdot\Vert_{1}$,
is defined as
%with $\gamma=\frac{\rho_{0}}{\psi_{0}}=\frac{\rho}{\psi}$:
%\begin{equation}\label{gg}
%g_\gamma({L},{S})=\max\bigl(\frac{\Vert{S}\Vert_{\infty}}{\gamma\,p^\delta},\frac{\Vert{L}\Vert_{2}}{p^\alpha}\bigr).
%\end{equation}
\begin{equation}
g_\gamma(\widehat{\mf{L}}_n-\mf{L}^{*},\widehat{\mf{S}}_n-\mf{S}^{*})=
\max\left(\frac{\Vert\widehat{\mf{L}}_n-{\mf{L}^{*}\Vert_{2}}}{p^{\alpha_1}},
\frac{\Vert\widehat{\mf{S}}_n-\mf{S}^{*}\Vert_{\infty}}{\gamma}\right).\label{ggamma}
% \, \Vert \mf{S}^{*} \Vert_{0,v}
\end{equation}
$g_\gamma$-consistency is our target, which also implies $\ell_{2}$-consistency for $\widehat{\mf{\Sigma}}_n=\widehat{\mf{L}}_n+\widehat{\mf{S}}_n$, as we later show.
For this reason, we define parametric consistency as follows.
\begin{definition}\label{par_cons} A sequence of estimator pairs $\bigl(\widehat{\mf{L}}_n,\widehat{\mf{S}}_n\bigr)$, defined as in \eqref{obj_all} and indexed in $n$,
is a parametrically consistent estimate of the low rank plus sparse decomposition \eqref{cov} for the covariance matrix $\mf{\Sigma}^{*}$ if the norm $g_\gamma(\widehat{\mf{L}}_n-\mf{L}^{*},\widehat{\mf{S}}_n-\mf{S}^{*})$ converges to $0$ with probability approaching $1$ as $n \to \infty$.\end{definition}
%Obviously, this $g_\gamma$-consistency implies consistency in $\ell_{2}$ norm.

\section{Algebraic consistency}\label{Cons_alg}

\subsection{Local convexity}\label{LC}

Let us reconsider our optimization problem
\begin{equation}\label{obj2}
\bigl(\widehat{\mf{L}}_n,\widehat{\mf{S}}_n\bigr)=
%\arg\min_{\mf{L},\mf{S}, \mf{L} \succeq 0, \mf{S} \succ 0,  \mf{L}+\mf{S} \succ 0}
\argmin_{\mf{L},\mf{S}\in\R^{p \times p},\mf{\Sigma}=\mf{L}+\mf{S}} \phi(\mf{L},\mf{S}),\nonumber
\end{equation}
where $\phi(\mf{L},\mf{S})=\phi_D(\mf{L},\mf{S})+\mathcal{P}(\mf{L},\mf{S})$,
%where $\phi(\mf{L},\mf{S})=\mathcal{L}^{(ld)}(\mf{\Sigma},\mf{\Sigma}_n)+\mathcal{P}(\mf{L},\mf{S})$,
$\phi_D(\mf{L},\mf{S})=\mathcal{L}^{(ld)}(\mf{\Sigma},\mf{\Sigma}_n)$,
%$=0.5\ln \det (\mf{I}_p+\mf{\Delta}_{n}\mf{\Delta}_{n}')$
%is the smooth component of $\phi(\mf{L},\mf{S})$,
$\mathcal{P}(\mf{L},\mf{S})=\psi \Vert\mf{L}\Vert_{*}+\rho \Vert\mf{S}\Vert_{1}$.
%is the non-smooth component of $\phi(\mf{L},\mf{S})$.
%We can disentangle the objective $\phi(\mf{L},\mf{S})$ as follows:
%where $\phi_D(\mf{L},\mf{S})=\mathcal{L}^{(ld)}(\mf{\Sigma},\mf{\Sigma}_n)$.
%In this section,
%%is the smooth component of $\phi(\mf{L},\mf{S})$ and $\mathcal{P}(\mf{L},\mf{S})$ is the non-smooth component of $\phi(\mf{L},\mf{S})$.
%we derive the conditions that guarantee the local convexity of $\phi_D(\mf{L},\mf{S})$,
%we calculate the first and the second derivative of $\phi_D(\mf{L},\mf{S})$,
%and we prove the Lipschitzianity of the gradient of $\phi_D(\mf{L},\mf{S})$.
%%These results are proved in this section.
%%We refer to \cite{bernardi2022log} for the proofs.
%%We refer to the Supplement for the proofs.

%\smallskip
%\far{da qui in poi, probabilmente va tolto per essere meno matematico e più discorsivo. Tutto in Supplement?}

In order to optimize $\phi_D(\mf{L},\mf{S})$, we need to prove that
$\phi_D(\mf{L},\mf{S})=0.5\ln \det (\mf{I}_p+\mf{\Delta}_{n}\mf{\Delta}_{n}')$, with $\mf{\Delta}_{n}=\mf{\Sigma}-\mf{\Sigma}_n$ and $\mf{\Sigma}=\mf{L}+\mf{S}$,
is convex within some range for $\mf{\Delta}_{n}$.
In the univariate context, the function $0.5\ln \det (1+x^2)$ is convex if and only if $\vert x \vert<{1}/{\sqrt{2}}$.
In the multivariate context, it is therefore reasonable to suppose that a similar condition on $\mf{\Delta}_{n}\mf{\Delta}_{n}'$ ensures local convexity.
%The proof requires to show the positive definiteness of the Hessian of $\phi_D(\mf{L},\mf{S})$ evaluated
%around as $\mf{\Sigma}_n$
%as $\mf{\Sigma}$ varies.
In other words, we aim to show that there exists a positive real $C$ such that, whenever $\Vert\mf{\Delta}_{n}\Vert<C$,
the function $0.5\ln\det(\mf{\Delta}_{n}\mf{\Delta}_{n}')$ is convex with high probability.

%%Changing variables in an obvious way, we have therefore proven the following.
%\begin{proposition}\label{conv_delta}
%For any $ \delta_\phi >0 $ the function $\ln\det\bigl(\delta_\phi^{-2}\mf{I}_p + \mf{\Delta}_{n} \mf{\Delta}_{n}'\bigr)$
%is convex on the closed  ball $ \mathcal{C_{\delta_\phi}}= \{ \mf{\Delta}_{n} | \mf{\Delta}_{n} \mathrm{~~is~ a~ real~~}
%p\times p \mathrm{~matrix~}, \Vert\mf{\Delta}_{n}\Vert_{2} \leq (3\delta_\phi p)^{-1} \}$.
%\end{proposition}
%
%\begin{proof}
%Proposition \ref{conv_delta} is a direct specialization of \cite{bernardi2022log}, Lemma 2.
%\end{proof}

%Changing variables in an obvious way, we have therefore proven the following.
\begin{proposition}\label{conv_delta}
For any $ \delta_\phi >0 $ the function $\ln\det\bigl(\delta_\phi^{-2}\mf{I}_p + \mf{\Delta}_{n} \mf{\Delta}_{n}'\bigr)$
is convex on the closed  ball $ \mathcal{C_{\delta_\phi}}= \{ \mf{\Delta}_{n} | \mf{\Delta}_{n} \mathrm{~~is~ a~ real~~}
p\times p \mathrm{~matrix~}, \Vert \mf{\Delta}_{n}\Vert_{2} \leq (3\delta_\phi p)^{-1} \}$.%\Vert\mf{L}^{*}\Vert_{2}^{-1}\Vert
\end{proposition}

\paragraph{Proof}%subsection*{Proof}%\begin{proof}
Proposition \ref{conv_delta} is a direct specialization of \cite{bernardi2022log}, Lemma 2. \qed
%\end{proof}

%In conclusion, even though the function $\ln\det(\mf{I}_p + \mf{A}) $ is always concave, Corollary
%\ref{conv_delta} shows that the matrix function $\ln\det\bigl(\delta_\phi^{-2}\mf{I}_p + \mf{A} \mf{A}'\bigr)$
Proposition \ref{conv_delta} evidences that the function $\ln\det\bigl(\delta_\phi^{-2}\mf{I}_p + \mf{\Delta}_{n} \mf{\Delta}_{n}'\bigr)$
can be made locally convex in arbitrary ball near $ 0 $, choosing a suitable $ \delta_\phi $ respecting the convexity radius
$\Vert \mf{\Delta}_{n} \Vert_{2} \leq (3\delta_\phi p)^{-1}$.
This proposition is necessary to study the random behaviour of $\phi_D(\mf{L},\mf{S})$ around $\mf{\Sigma}^{*}$ with respect to convexity,
because $\mf{\Delta}_{n}$ is a random matrix, and so the convexity condition $\Vert\mf{\Delta}_{n} \Vert_{2} \leq (3\delta_\phi p)^{-1}$
must be ensured for $\delta_\phi=1$. First and second derivative of $\phi_D(\mf{L},\mf{S})$, as well as its asymptotic behaviour,
are studied in the Supplement Section B.

\subsection{Fisher information}\label{fish_main}

The key to obtain the algebraic consistency of Definition \ref{def_alg}
is to control the algebraic features of the low rank and sparse matrix varieties $\mathcal{L}(r)$ and $\mathcal{S}(s)$ containing $\mf{L}^{*}$ and $\mf{S}^{*}$ respectively, as well as their tangent spaces $\mathcal{T}(\mf{L}^{*})$ and $\Omega(\mf{S}^{*})$,
because the low rank variety $\mathcal{L}(r)$ is locally curve, and so its recovery by using $\mathcal{L}^{(ld)}(\mf{\Sigma},\mf{\Sigma}_n)$ in problem \eqref{obj} may be very sensitive
to small perturbations in $\mf{\Sigma}_n$ (see \cite{chandrasekaran2012latent} for reference, as well as \cite{chen2013low} and \cite{hsu2011robust} for alternative approaches).
%around $\mf{\Sigma}^{*}$
%Our aim is that the two matrix varieties $\mathcal{L}(r)$ and $\mathcal{S}(s)$ can be recovered,
%without resorting to the need to control sample eigenvectors as done by eigenvalue shrinkage approaches (see Section \ref{motiv}).
%\far{Add Hsu et al. and Chen et al.}
%because we can recover the true latent eigenvectors and residual sparsity pattern.
%It follows that we need to control how impacts on ensuring algebraic consistency.
%At this stage, we can observe that $\Vert\mf{L}^{*}\Vert_{2}=O(p^{\alpha_1})$ under Assumption \ref{eigenvalues}(i).
Let us define the following reshaped version of $\phi_D(\mf {L},\mf{S})$:
\begin{equation}\widetilde{\phi}_{D}(\mf{L},\mf{S})=\widetilde{\mathcal{L}}^{(ld)}(\mf{\Sigma},\mf{\Sigma}_n)=0.5\ln \det (\mf{I}_p+p^{-2\alpha_1}\mf{\Delta}_{n}\mf{\Delta}_{n}'),\nonumber\end{equation}
which can also be expressed as
$\widetilde{\phi}_{D}(\widetilde{\mf{L}},\widetilde{\mf{S}})=0.5\ln \det (\mf{I}_p+\widetilde{\mf{\Delta}}_{n}\widetilde{\mf{\Delta}}_{n}')$,
where $\widetilde{\mf{\Delta}}_{n}=p^{-\alpha_1}\mf{\Delta}_{n}=\widetilde{\mf{\Sigma}}-\mf{\Sigma}_{n}$, $\widetilde{\mf{\Sigma}}=\widetilde{\mf{L}}+\widetilde{\mf{S}}$, with $\widetilde{\mf{L}}=p^{-\alpha_{1}}\mf{L}$ and $\widetilde{\mf{S}}=p^{-\alpha_{1}}\mf{S}$.
%\far{
%Let us define $\psi_{0}=\frac{1}{\xi(\mathcal{T}(\mf{L}^{*}))}\frac{1}{\sqrt{n}}=\frac{\psi}{p^{\alpha}}$,
%$\rho=\gamma p^{\delta} \psi_{0}$, and $\rho_{0}=\gamma\psi_{0}=\frac{\rho}{p^{\delta}}$,
%with $\gamma \in [9\xi(\mathcal{T}),\frac{1}{6\mu(\Omega)}]$.
We regard the function $-\widetilde{\mathcal{L}}^{(ld)}(\mf{\Sigma},\mf{\Sigma}_n)$
%negative reshaped smooth component of \eqref{obj2},
%$\widetilde{\phi}_{D}(\mf{L},\mf{S})$, as $$\widetilde{\phi}_{D}(\mf{L},\mf{S})=\widetilde{\mathcal{L}}^{(ld)}(\mf{\Sigma},\mf{\Sigma}_n)=\frac{1}{2}\ln \det (\mf{I}_p+p^{-2\alpha_1}\mf{\Delta}_{n}\mf{\Delta}_{n}'),$$
%with $\mf{\Delta}_n=\mf{\Sigma}-\mf{\Sigma}_{n}$, $\mf{\Sigma}=\mf{L}+\mf{S}$,
%defined in equation \eqref{reshaped},
as a nonlinear function of the squared sample covariance matrix,
which we need to maximize. $-\widetilde{\mathcal{L}}^{(ld)}(\mf{\Sigma},\mf{\Sigma}_n)$ is concave into
the convexity range of Proposition \ref{conv_delta}.
Therefore, to ensure that $-\widetilde{\mathcal{L}}^{(ld)}(\mf{\Sigma},\mf{\Sigma}_n)$ is concave for $\mf{\Sigma}=\widehat{\mf{\Sigma}}_n$,
once defined $\widehat{\mf{\Delta}}_{n}=\widehat{\mf{\Sigma}}_n-\mf{\Sigma}_n$, we need that
%$$\mathcal{P}(\Vert\mf{L}^{*}\Vert_{2}^{-1}\Vert\widehat{\mf{\Delta}}_{n}\Vert \geq (3p)^{-1}) \to 0$$
$p^{-\alpha_{1}}\Vert\widehat{\mf{\Delta}}_{n}\Vert \leq 1/(3p)$, which is dominated by the condition $\Vert\widehat{\mf{\Delta}}_{n}\Vert \leq 1/3$.

%because it is a quadratic function of $\mf{\Delta}_n$.
In this view, under its conditions, we can define
%the conditions of Lemma \ref{random_conv}
$\mathcal{I}^{*}$, the operator associating to each pair $(\widetilde{\mf{L}},\widetilde{\mf{S}})$
the Fisher information of $-\widetilde{\phi}_{D}(\widetilde{\mf{L}},\widetilde{\mf{S}})$, as
\begin{eqnarray}
\mathcal{I}^{*}(\widetilde{\mf{L}},\widetilde{\mf{S}})&=&-\mathrm{E}\left\{\frac{\partial^{2} [-\widetilde{\mathcal{L}}^{(ld)}(\widetilde{\mf{L}},\widetilde{\mf{S}})]}{\partial^{2}\widetilde{\mf{L}}}\right\}
=-\mathrm{E}\left\{\frac{\partial^{2} [-\widetilde{\mathcal{L}}^{(ld)}(\widetilde{\mf{L}},\widetilde{\mf{S}})]}{\partial^{2}\widetilde{\mf{S}}}\right\}.\nonumber
%-\mathrm{E}\left\{\frac{\partial^{2}}[-\widetilde{\mathcal{L}}^{(ld)}(\mf{L},\mf{S})]}{\partial^{2} \mf{S}}\right\}.\nonumber
\end{eqnarray}

In the Supplement (Section B.3), it is proved that,
imposing the condition $\Vert\widetilde{\mf{\Delta}}_{n}\Vert \leq \delta_{\widetilde{\Delta}}$,
we can estimate from below the expectation
\begin{equation}
\Bigg\Vert\mathrm{E}\left(\frac{\partial{\widetilde{\phi}'_{D}(\widetilde{\mf{L}},\widetilde{\mf{S}})}}{\partial{\widetilde{\mf{L}}}}\right)\Bigg\Vert_{2}\geq
\frac{(1-\delta_{\widetilde{\Delta}})^2}{(1+\delta_{\widetilde{\Delta}}^2)^2}.\label{minI}
\end{equation}
Analogously, imposing the condition $\Vert\widetilde{\mf{\Delta}}_{n}\Vert_{\infty} \leq \delta^{\infty}_{\widetilde{\Delta}}$, we can estimate from below the expectation
\begin{equation}
\Bigg\Vert\mathrm{E}\left(\frac{\partial{\widetilde{\phi}'_{D}(\widetilde{\mf{L}},\widetilde{\mf{S}})}}{\partial{\widetilde{\mf{L}}}}\right)\Bigg\Vert_{\infty}\geq
\frac{(1-\delta^{\infty}_{\widetilde{\Delta}})^2}{(1+\delta^{\infty 2}_{\widetilde{\Delta}})^2}.\nonumber
\end{equation}

Importantly, since $\Vert . \Vert_{\infty} \leq \Vert . \Vert_{2}$, it holds $\delta^{\infty}_{\widetilde{\Delta}}\leq\delta_{\widetilde{\Delta}}$.
Consequently, under the conditions of Proposition \ref{conv_delta}, %setting $\Vert\mf{L}^{*}\Vert_{2}=p$ we obtain
we can note that $\delta^{\infty}_{\widetilde{\Delta}}\leq\delta_{\widetilde{\Delta}} \leq 1/3p$,
so that %in the worst case
\begin{equation}
\frac{(1-\delta^{\infty}_{\widetilde{\Delta}})^2}{(1+\delta^{\infty 2}_{\widetilde{\Delta}})^2} \geq
\frac{(1-\delta_{\widetilde{\Delta}})^2}{(1+\delta_{\widetilde{\Delta}}^2)^2} \geq \frac{9p^2(9p^2+1-6p)}{(9p^2+1)^2} \geq \frac{9}{25}.\label{minimum_Fisher}
\end{equation}
This means that $\widetilde{\phi}_{D}(\widetilde{\mf{L}},\widetilde{\mf{S}})$ is strongly convex for all $p$ and $n$
in $\widetilde{\mf{\Sigma}}=p^{-\alpha_{1}}\widehat{\mf{\Sigma}}_n$,
whenever the convexity condition $\Vert p^{-\alpha_{1}}\widehat{\mf{\Sigma}}_n \Vert_{2} \leq 1/(3p)$ of Proposition \ref{conv_delta},
dominated by $\Vert\widehat{\mf{\Sigma}}_n\Vert_{2} \leq 1/3$, is satisfied. %%\textbf{More conditions inactive!}

\subsection{Manifold retrieval}

Let us define the following measure of transversality between two algebraic matrix varieties $\mathcal{T}_{1}$ and $\mathcal{T}_{2}$ as in  \cite{chandrasekaran2012latent}, equation (2.7):
\begin{equation}
\varrho(\mathcal{T}_{1},\mathcal{T}_{2})=
%\Vert\mathbb{P}_{\mathcal{T}_{1}}-\mathbb{P}_{\mathcal{T}_{2}}\Vert_{2 \rightarrow 2}=
\max_{\Vert \mf{N} \Vert_{2} \leq 1} \Vert \mathbb{P}_{\mathcal{T}_{1}}\mf{N}-\mathbb{P}_{\mathcal{T}_{2}}\mf{N}\Vert_{2},
\label{varrho_def}
\end{equation}
where $\mathbb{P}_{\mathcal T_{1}}$ and $\mathbb{P}_{\mathcal T_{2}}$ are the projection operators onto $\mathcal T_{1}$ and $\mathcal T_{2}$, respectively.
Hereafter, let $\Omega=\Omega(\mf{S}^{*})$ and
$\mathcal T=\mathcal{T}(\mf{L}^{*})$, where ${\Omega}$ is the space tangent to $\mathcal S=\mathcal S(s)$ (see \eqref{var:S}) at $\mf{S}^{*}$ and $\mathcal T$ is the space tangent to $\mathcal L=\mathcal L(r)$ (see \eqref{var:L}) at $\mf{L}^{*}$.
We define the manifold $\mathcal{T}'$ as the tangent space to a generic matrix $\mf{L}$:
\begin{eqnarray}
\mathcal{T}'=\mathcal{T}(\mf{L})=\{\mf{M}\in \R^{p \times p} \mid \mf{M}=\mf{U} \mf{Y}_1'+\mf{Y}_2 \mf{U}',
\mf{Y}_1, \mf{Y}_2 \in \R^{p \times r}, \mf{U}\in \R^{p\times r}, \mf{U}' \mf{U}=\mf{I}_r, \nonumber\\
\mf{U}' \mf{L}\mf{U} \in \R^{r \times r} \mbox{diagonal}\},\nonumber
\end{eqnarray}
such that $\varrho(\mathcal{T},\mathcal{T}')\leq \kappa_\mathcal{T}\xi(\mathcal{T})$ with $\kappa_\mathcal{T}\in(0,1/2]$,
and we set the Cartesian product $\mathcal{Y}={\Omega}\times\mathcal{T}'$.
In order to assess the effect of a deviation from $\mf{\Sigma}^{*}$ along the varieties $\mathcal{T'}$, $\mathcal{T'}^{\perp}$,
$\Omega$, $\Omega^{\perp}$ under the constraints $\mf{L}\in \mathcal{T}'$,
$\varrho(\mathcal{T}', \mathcal{T})\leq \kappa_\mathcal{T}\xi(\mathcal{T})$ with $\kappa_\mathcal{T}\in(0,1/2]$, and $\mf{S}\in \Omega$,
%it follows from \eqref{fish_bound} that we need to evaluate the quantity $\frac{5}{4}\Vert\mf{W}\Vert_{2}$ under each case.
we need to evaluate \eqref{minI} under each case.

%%under the constraint $\Vert \mf{W} \Vert \leq 1$
%%as $n \to \infty$.
%We then define the matrices
%%$\mf{\Delta}^{*}_{n}=\mf{\Sigma}^{*}-\mf{\Sigma}_{n}$, with $\mf{\Sigma}^{*}=\mf{L}^{*}+\mf{S}^{*}$,
%$\mf{H}_{\mathcal{T'}}=\mathbb{P}_{\mathcal{T'}}(\mf{L})$,
%$\mf{H}_{\mathcal{T'}^{\perp}}=\mathbb{P}_{\mathcal{T'}^{\perp}}(\mf{L})=\mf{L}-\mathbb{P}_{\mathcal{T'}}(\mf{L})$,
%$\mf{H}_{\Omega}=\mathbb{P}_{\Omega}(\mf{S})$, $\mf{H}_{\Omega^{\perp}}=\mathbb{P}_{\Omega^{\perp}}(\mf{S})=\mf{S}-\mathbb{P}_{\Omega}(\mf{S})$.
%%We now express the matrices of interest for this analysis.

Following \cite{chandrasekaran2012latent} and previous reasonings, we can define and compute the quantities
%we can observe that $\mf{\Sigma}=\mathbb{P}_{\mathcal{T'}}(\mf{\Sigma}^{*})=\mf{L}^{*}+\mathbb{P}_{\mathcal{T'}}(\mf{S}^{*})$, so that $\mf{W}=\mf{H}_{\mathcal{T'}^{\perp}}$,
\begin{eqnarray}&&\beta_{\mathcal{T}}=\max_{\mf{M}\in {\mathcal{T'}},\Vert \mf{M} \Vert_{\infty}=1}\Vert\mathcal{I}^{*}(\mf{M},\mf{S}^{*})\Vert_{\infty}=1
\qquad \mbox{and} \qquad
\beta_{\Omega}=\max_{\mf{M}\in {\Omega},\Vert \mf{M} \Vert_{2}=1}\Vert\mathcal{I}^{*}(\mf{L}^{*},\mf{M})\Vert_{2}=1.\nonumber\end{eqnarray}
Analogously, to study the behaviour of $\mathcal{I}^{*}(\mf{L},\mf{S})$ along the varieties $\mathcal{T'}$, $\mathcal{T'}^{\perp}$, $\Omega$, and $\Omega^{\perp}$,
following \cite{chandrasekaran2012latent} we define:
\begin{eqnarray}
&&\alpha_{\mathcal{T}}=\min_{\mf{M}\in {\mathcal{T'}},\Vert\mf{M}\Vert_{2}=1}\Vert\mathbb{P}_{\mathcal{T'}}\mathcal{I}^{*}\mathbb{P}_{\mathcal{T'}}(\mf{M})\Vert_2 \qquad \mbox{and} \qquad
\delta_{\mathcal{T}}=\max_{\mf{M}\in {\mathcal{T'}^{\perp}},\Vert \mf{M} \Vert_{2}=1}\Vert\mathbb{P}_{\mathcal{T'}^{\perp}}\mathcal{I}^{*}\mathbb{P}_{\mathcal{T'}^{\perp}}(\mf{M})\Vert_{2};\nonumber\\
&&\alpha_{\Omega}=\min_{\mf{M}\in {\Omega},\Vert \mf{M} \Vert_{\infty}=1}\Vert\mathbb{P}_{\Omega}\mathcal{I}^{*}\mathbb{P}_{\Omega}(\mf{M})\Vert_{\infty} \qquad \mbox{and} \qquad
\delta_{\Omega}=\max_{\mf{M}\in {\Omega^{\perp}},\Vert \mf{M} \Vert_{\infty}=1}\Vert\mathbb{P}_{\Omega^{\perp}}\mathcal{I}^{*}\mathbb{P}_{\Omega^{\perp}}(\mf{M})\Vert_{\infty}.\nonumber
\end{eqnarray}

Then, recalling that $\Vert\widetilde{\mf{\Delta}}_{n}\Vert_{\infty}\leq\Vert\widetilde{\mf{\Delta}}_{n}\Vert\leq\delta_{\widetilde{\Delta}}$ and exploiting \eqref{minI}, we can derive that:
\begin{compactenum}
\item[$\bullet$] setting $\widetilde{\mf{\Delta}}_{n}=p^{-\alpha_1}(\mf{\Sigma}_n-(\mf{M}+\mf{S}^{*}))$, we get
$\alpha_{\mathcal{T}}=\frac{(1-\delta_{\widetilde{\Delta}})^2}{(1+\delta_{\widetilde{\Delta}}^2)^2}$ and $\delta_{\mathcal{T}}=\beta_{\mathcal{T}}-\alpha_{\mathcal{T}}$;
%\item setting $\widetilde{\mf{\Delta}}_{n}=p^{-\alpha_1}(\mf{\Sigma}_n-(\mf{M}+\mf{S}^{*}))$,
%$$=1-\left(\frac{1}{\mathrm{E}(\lambda_1(\widetilde{\mf{\Delta}}_{n}^2))+1}\right)^2;$$
\item[$\bullet$] setting $\widetilde{\mf{\Delta}}_{n}=p^{-\alpha_1}(\mf{\Sigma}_n-(\mf{L}^{*}+\mf{M}))$, we get
$\alpha_{\Omega}=\frac{(1-\delta_{\widetilde{\Delta}})^2}{(1+\delta_{\widetilde{\Delta}}^2)^2}$ and
%\item setting $\widetilde{\mf{\Delta}}_{n}=p^{-\alpha_1}(\mf{\Sigma}_n-(\mf{L}^{*}+\mf{M}))$,
$\delta_{\Omega}=\beta_{\Omega}-\alpha_{\Omega}$.
%$$1-\left(\frac{1}{\mathrm{E}(\lambda_1(\widetilde{\mf{\Delta}}_{n}^2))+1}\right)^2.$$
\end{compactenum}

Since $\mathcal{I}^{*}(\mf{L},\mf{S})$ may be regarded as a map, we can consequently write
$\mathcal{A}^{\dag}\mathcal{I}^{*}\mathcal{A}(\mf{L},\mf{S})=\mathcal{I}^{*}(\mf{L}+\mf{S},\mf{L}+\mf{S})$;
$\mathbb{P}_{\mathcal{Y}}\mathcal{A}^{\dag}\mathcal{I}^{*}\mathcal{A}\mathbb{P}_{\mathcal{Y}}(\mf{L},\mf{S})=
\mathcal{I}^{*}(\mathbb{P}_{\Omega}\mf{L}+\mf{S},\mf{L}+\mathbb{P}_{\mathcal{T}'}\mf{S})$;
$\mathbb{P}_{\mathcal{Y}^\perp}\mathcal{A}^{\dag}\mathcal{I}^{*}\mathcal{A} \mathbb{P}_{\mathcal{Y}}(\mf{L},\mf{S})=\mathcal{I}^{*}(\mathbb{P}_{\Omega^{\perp}}\mf{L}+\mf{S},\mf{L}+\mathbb{P}_{\mathcal{T}'^{\perp}}\mf{S})$.
%The additional constraints are needed to ensure that the Hessian of $\frac{1}{2}\Vert \mf{\Sigma}_{n}-(\underline{\mf{L}}+\underline{\mf{S}})\Vert_{F}^2$ is positive definite,
%such that the optimum of \eqref{probtang} is unique.
These quantities need to be studied in order to ensure that the elements of the spaces ${\Omega}$ and $\mathcal{T}$ are individually identifiable
under the map $\mathcal{I}^{*}(\mf{L},\mf{S})$, and that the sensitivity of $\mathcal{I}^{*}(\mf{L},\mf{S})$ to small perturbations
%around $\mf{\Sigma}^{*}$ in the final solution of problem \eqref{obj}
has a limited impact on the recovery of the direct sum $\mathcal{L}(r)\oplus\mathcal{S}(s)$ by problem \eqref{obj}.
%(see \cite{chandrasekaran2012latent}, Section 3.4).
%\far{Invert $L$ and $S$!}

Therefore, to ensure algebraic consistency,
we need to explicitly control the quantities $\alpha_\Omega$, $\alpha_{\mathcal{T'}}$, $\delta_\Omega$, $\delta_{\mathcal{T'}}$, $\beta_\Omega$, $\beta_{\mathcal{T'}}$.
Analogously to \cite{chandrasekaran2012latent}, we define
$\alpha_{\mathcal{Y}}=\min(\alpha_\Omega,\alpha_{\mathcal{T'}})$, $\delta_{\mathcal{Y}}=\max(\delta_\Omega,\delta_{\mathcal{T'}})=1-\alpha_{\mathcal{Y}}$,
$\beta_{\mathcal{Y}}=\max(\beta_\Omega,\beta_{\mathcal{T'}})$, and we assume the following. %to ensure identifiability.
\begin{assumption}\label{ass_alg}
For some $\nu \in (0,{1}/{2}]$, it holds
${\alpha_{\mathcal{Y}}}^{-1}{\delta_{\mathcal{Y}}} \leq 1-2\nu$.
\end{assumption}
Combining $\gamma_{\mathcal{Y}}=1-\alpha_{\mathcal{Y}}$ and Assumption \ref{ass_alg}, we get
%$$\frac{\delta_{\mathcal{Y}}}{\alpha_{\mathcal{Y}}} \leq 1-2\nu$$
${\alpha_{\mathcal{Y}}}^{-1}{(1-\alpha_{\mathcal{Y}})} \leq 1-2\nu$,
which leads to the condition $\nu \geq \frac{2\alpha_{\mathcal{Y}}-1}{2\alpha_{\mathcal{Y}}}$.
At this stage, we can state the following proposition, which improves and extends Proposition 3.3 in \cite{chandrasekaran2012latent} to the case of the log-det loss $\widetilde{\mathcal{L}}^{(ld)}(\mf{\Sigma},\mf{\Sigma}_{n})$.

\begin{proposition}\label{11}
%Suppose that
Suppose that $\Vert\widehat{\mf{\Delta}}_{n}\Vert_{2} \leq 1/3$ and Assumption \ref{ass_alg} hold.
%Proposition \ref{conv_delta} holds for $\delta_{\phi}=1$ with $\mf{\Delta}_{n}=p^{-\alpha_1}\widehat{\mf{\Delta}}_{n}$.
%\far{and $\Vert\mf{L}^{*}\Vert_{2}^{-1}\Vert\widehat{\mf{\Delta}}_{n}\Vert \leq (3p)^{-1}$ with $\Vert\mf{L}^{*}\Vert_{2}<3p$}.
Let $\alpha_{\mathcal{Y}}\geq 0.70711$, $\kappa_\mathcal{T} \in (0,0.5]$,
$$\gamma \in \left[\frac{\frac{2(1-\kappa_\mathcal{T})}{(1+\kappa_\mathcal{T})}\xi(\mathcal{T}(\mf{L}^{*}))(1-\nu)}{\nu\alpha_{\mathcal{Y}}},
\frac{\nu\alpha_{\mathcal{Y}}}{4\mu(\Omega(\mf{S}^{*}))\beta_{\mathcal{Y}}(1-\nu)}\right],$$
$$\xi(\mathcal{T}(\mf{L}^{*})\mu(\Omega(\mf{S}^{*}))\leq\frac{1-\kappa_\mathcal{T}}
{8(1+\kappa_\mathcal{T})}\Biggl(\frac{\nu\alpha_{\mathcal{Y}}}{\beta_{\mathcal{Y}}(1-\nu)}\Biggr)^2,$$
with $\alpha_{\mathcal{Y}},\beta_{\mathcal{Y}},\gamma,\nu$ as previously defined.
%For all $p \in \N$, as $n \to \infty$, with probability $1-O(1/n^2)$ it holds:
%\far{STOP}
%Suppose that Assumptions \ref{alg} and \ref{ass_alg} hold,
%and that $p$ and $n$ are such that the convexity condition $\Vert \mf{\Delta}^{*}_n \Vert \leq {(3\delta_\phi p)^{-1}}$ is respected for $\delta_\phi=1$.
%Let $\frac{\sqrt{r}\kappa_{S}}{\kappa_{L}}\leq\frac{1}{6}\bigl(\frac{\nu\alpha_{\mathcal{Y}}}{\beta_{\mathcal{Y}}(2-\nu)}\bigr)^2$,
%and $\gamma \in [\frac{3\xi(\mathcal{T}(\mf{L}^{*}))(2-\nu)}{\nu\alpha_{\mathcal{Y}}},\frac{\nu\alpha_{\mathcal{Y}}}{2\mu(\Omega(\mf{S}^{*}))\beta_{\mathcal{Y}}(2-\nu)}]$,
%with $\alpha_{\mathcal{Y}},\beta_{\mathcal{Y}},\gamma,\nu$ as previously defined.
Then, %under Assumptions \ref{eigenvalues}-\ref{tails},
for all $(\mf{S},\mf{L})\in\mathcal{Y}$ such that $\mathcal{Y}={\Omega}\times\mathcal{T}'$ with $\varrho(\mathcal{T},\mathcal{T}')\leq \kappa_\mathcal{T}\xi(\mathcal{T})$,
the following statements hold:
%for all $p,n \in \N$: %with probability $1-O(1/n^2)$:
\begin{itemize}
%\item
%$g_{\gamma}(\mathbb{P}_{\mathcal{Y}}\mathcal{A}^{\dag}\mathcal{I}^{*}\mathcal{A}
%\mathbb{P}_{\mathcal{Y}}(\mf{S},\mf{L})) \geq {\alpha_{\mathcal{Y}}}g_{\gamma}(\mf{S},\mf{L})$;
%%And this implies that for all $(S,L)\in\mathcal{Y}$, \far{NON VEDO PERCHE'}
%%\begin{equation}
%%g_{\gamma}(\mathbb{P}_{\mathcal{Y}^\perp}\mathcal{A}^{\dag}\mathcal{A}
%%\mathbb{P}_{\mathcal{Y}}(S,L)\geq \frac{1}{2}g_{\gamma}(S,L).\nn
%%\end{equation}
%%\item The effect of elements in $\mathcal{Y}=\Omega \times \mathcal{T}'$ on the orthogonal complement
%%$\mc Y^\perp=\Omega^{\perp} \times \mathcal{T}'^{\perp}$  is bounded above as
%%\begin{equation}
%%\Vert\mathbb{P}_{\mathcal{Y}^\perp}\mathcal{A}^{\dag}\mathcal{A}\mathbb{P}_{\mathcal{Y}}
%%(\mathbb{P}_{\mathcal{Y}}\mathcal{A}^{\dag}\mathcal{A}\mathbb{P}_{\mathcal{Y}})^{-1}\Vert_{g_\gamma \rightarrow g_{\gamma}}\leq \frac{1}{2}\nn
%%\end{equation}
%%And that implies that for all $(S,L)\in\mathcal{Y}$, \far{NOTAZIONE NON CHIARA $g_{\gamma}$ PRENDE DUE ARGOMENTI MANCANO DELLE PARENTESI?}
%\item $g_{\gamma}(\mathbb{P}_{\mathcal{Y}^\perp}\mathcal{A}^{\dag}\mathcal{I}^{*}\mathcal{A}\mathbb{P}_{\mathcal{Y}}(\mf{S},\mf{L}))\leq
%g_{\gamma}(\mathbb{P}_{\mathcal{Y}^\perp}\mathcal{A}^{\dag}\mathcal{I}^{*}\mathcal{A}
%\mathbb{P}_{\mathcal{Y}}(\mf{S},\mf{L}))$.
%%%\frac{1}{2}g_{\gamma}(\mf{S},\mf{L})
\item[i)]
$\min_{(\mf{L},\mf{S})\in \mathcal{Y},\Vert \mf{L} \Vert_{2}=1,\Vert \mf{S} \Vert_{\infty}=\gamma} g_{\gamma}(\mathbb{P}_{\mathcal{Y}}\mathcal{A}^{\dag}\mathcal{I}^{*}\mathcal{A}
\mathbb{P}_{\mathcal{Y}}(\mf{S},\mf{L}))
\geq \alpha_{\mathcal{Y}}\Bigl(\frac{3}{2}-\alpha_{\mathcal{Y}}\Bigr)g_{\gamma}(\mf{S},\mf{L})$;
%$0.5{\alpha_{\mathcal{Y}}}g_{\gamma}(\mf{S},\mf{L})$;
%And this implies that for all $(S,L)\in\mathcal{Y}$, \far{NON VEDO PERCHE'}
%\begin{equation}
%g_{\gamma}(\mathbb{P}_{\mathcal{Y}^\perp}\mathcal{A}^{\dag}\mathcal{A}
%\mathbb{P}_{\mathcal{Y}}(S,L)\geq \frac{1}{2}g_{\gamma}(S,L).\nn
%\end{equation}
%\item The effect of elements in $\mathcal{Y}=\Omega \times \mathcal{T}'$ on the orthogonal complement
%$\mc Y^\perp=\Omega^{\perp} \times \mathcal{T}'^{\perp}$  is bounded above as
%\begin{equation}
%\Vert\mathbb{P}_{\mathcal{Y}^\perp}\mathcal{A}^{\dag}\mathcal{A}\mathbb{P}_{\mathcal{Y}}
%(\mathbb{P}_{\mathcal{Y}}\mathcal{A}^{\dag}\mathcal{A}\mathbb{P}_{\mathcal{Y}})^{-1}\Vert_{g_\gamma \rightarrow g_{\gamma}}\leq \frac{1}{2}\nn
%\end{equation}
%And that implies that for all $(S,L)\in\mathcal{Y}$, \far{NOTAZIONE NON CHIARA $g_{\gamma}$ PRENDE DUE ARGOMENTI MANCANO DELLE PARENTESI?}
\item[ii)] %if $\alpha_{\mathcal{Y}} \in [0.70711,1]$, then
$g_{\gamma}(\mathbb{P}_{\mathcal{Y}^\perp}\mathcal{A}^{\dag}\mathcal{I}^{*}\mathcal{A}\mathbb{P}_{\mathcal{Y}}(\mf{S},\mf{L}))
\leq(1-\nu)g_{\gamma}(\mathbb{P}_{\mathcal{Y}}\mathcal{A}^{\dag}\mathcal{I}^{*}\mathcal{A} \mathbb{P}_{\mathcal{Y}}(\mf{S},\mf{L}))$.
%%\frac{1}{2}g_{\gamma}(\mf{S},\mf{L})
\end{itemize}
\end{proposition}

\begin{remark}
Remarkably, the quantities $\alpha_{\mathcal{Y}},\beta_{\mathcal{Y}},\gamma_{\mathcal{Y}},\nu$ can now be explicitly computed, unlike in \cite{chandrasekaran2012latent}.
Comparing our identifiability condition with the one in \cite{chandrasekaran2012latent}, we can observe that, although we now need to impose the constraint $\alpha_{\mathcal{Y}}\geq 0.70711$,
our condition is much weaker. In particular, imposing the maximum value $\kappa_\mathcal{T}=1/2$ as they do, we get $\xi(\mathcal{T}(\mf{L}^{*}))\mu(\Omega(\mf{S}^{*}))\leq 1/24$ instead of
$\xi(\mathcal{T}(\mf{L}^{*}))\mu(\Omega(\mf{S}^{*}))\leq 1/54$.
Considering the value $\kappa_\mathcal{T}=1/4$, which is obtained as part (iii) of Corollary C.1 in the Supplement,
our identifiability condition reads as $\xi(\mathcal{T}(\mf{L}^{*}))\mu(\Omega(\mf{S}^{*}))\leq 3/40$.
\end{remark}

\begin{remark}
The condition $\alpha_{\mathcal{Y}}\geq 0.70711$ is explicitly assumed.
%ensured by inequality \eqref{minimum_Fisher},
%descending from the convexity condition $\Vert\widehat{\mf{\Delta}}_{n}\Vert \leq 1/3$.
However, inequality \eqref{minimum_Fisher} also shows that, as $p \to \infty$,
the actual lower bound for $\alpha_{\mathcal{Y}}$ gets increasingly larger, progressively approaching $1$.
%In principles, the condition $\alpha_{\mathcal{Y}}\in [0.70711,1]$ of part 2 ensures that Proposition \ref{11} holds for all $p,n\in\N$.
%However, from Lemma \ref{Lemma_ort}, it follows that this condition is more and more likely as $n\to\infty$ under the conditions of Lemma \ref{Lemma_cons}.
Moreover, the lower bound $\alpha_{\mathcal{Y}}(3/2-\alpha_{\mathcal{Y}})$ of part (i) cannot be smaller than $\alpha_{\mathcal{Y}}/2$ (which is the lower bound in \cite{chandrasekaran2012latent}),
since $\alpha_{\mathcal{Y}} \in (1/2,1]$ by Assumption \ref{ass_alg}, which leads to $\alpha_{\mathcal{Y}}(3/2-\alpha_{\mathcal{Y}})\geq 1/2$.
\end{remark}
%Note that it always holds $\alpha_{\mathcal{Y}}(3/2-\alpha_{\mathcal{Y}})\geq\frac{\alpha_{\mathcal{Y}}}{2}$, since $\alpha \in (1/2,1]$ by Assumption \ref{ass}.

%\begin{remark}
%Note that the upper bound for the product $\xi(\mathcal{T}(\mf{L}^{*}))\mu(\Omega(\mf{S}^{*}))$ is much larger than in \cite{chandrasekaran2012latent}.
%In their case, imposing $\kappa_\mathcal{T}=1/2$, we get $\xi(\mathcal{T}(\mf{L}^{*}))\mu(\Omega(\mf{S}^{*}))\leq 1/54$, while we get $\xi(\mathcal{T}(\mf{L}^{*}))\mu(\Omega(\mf{S}^{*}))\leq 1/24$.
%More, as Proposition 7 (see the Supplement) later shows, $\kappa_\mathcal{T}=1/4$ is always admissible under Assumptions \ref{eigenvalues}-\ref{ass_alg}, so that we go further up to $\xi(\mathcal{T}(\mf{L}^{*}))\mu(\Omega(\mf{S}^{*}))\leq 3/40$. Importantly, here the quantities $\alpha_{\mathcal{Y}},\beta_{\mathcal{Y}},\gamma,\nu$ can be explicitly computed, for any $\phi_D(\mf{L},\mf{S})$ such that \ldots
%\end{remark}

%%

\section{Parametric consistency}\label{Cons_par}

%\subsection{Consistency results}

%Let us define the following norm $g_\gamma$:%, with $\gamma=\frac{\rho_{0}}{\psi_{0}}=\frac{\rho}{\psi}$:
%%\begin{equation}\label{gg}
%%g_\gamma({L},{S})=\max\bigl(\frac{\Vert{S}\Vert_{\infty}}{\gamma\,p^\delta},\frac{\Vert{L}\Vert_{2}}{p^\alpha}\bigr).
%%\end{equation}
%\begin{equation}
%g_\gamma(\widehat{\mf{L}}_n-\mf{L}^{*},\widehat{\mf{S}}_n-\mf{S}^{*})=
%\max\left(\frac{\Vert\widehat{\mf{L}}_n-{\mf{L}^{*}\Vert_{2}}}{\Vert\mf{L}^{*}\Vert_{2}},
%\frac{\Vert\widehat{\mf{S}}_n-\mf{S}^{*}\Vert_{\infty}}{\gamma},\right),\label{ggamma}
%% \, \Vert \mf{S}^{*} \Vert_{0,v}
%\end{equation}
%where the range for $\gamma$ is given by Proposition \ref{11}, %$\psi=\frac{p^{\alpha}}{\xi(\mathcal{T}(\mf{L}^{*}))}\frac{1}{\sqrt{n}}$
%and set $\psi_{0}={f_{\delta_\epsilon}(p,n)}/{\xi(\mathcal{T}(\mf{L}^{*}))}$,
%$\rho_{0}=\gamma\psi_{0}$,
%$\psi=\Vert\mf{L}^{*}\Vert_{2}\psi_{0}$, and $\rho=\rho_0$,
%%$\rho=p^{\delta}\rho_{0}=p^{\delta} \gamma \psi_{0}$,
%where $\psi$ and $\rho$ are the thresholds in \eqref{obj}.
%The norm \eqref{ggamma} is the dual norm of the composite penalty $\psi_{0}\Vert\cdot\Vert_{*}+\rho_{0}\Vert\cdot\Vert_{1}$,
%with which the direct sum $\mathcal{L}(r)\oplus\mathcal{S}(s)$ is naturally equipped \far{(cfr. paragraph 3.3 in \cite{chandrasekaran2012})}.

%%

We now explicitly compare the rescaled heuristics of \cite{farne2020large,farne2024large}, which is defined as
\begin{equation}
\widetilde{\phi}^{(F)}(\mf{L},\mf{S})=\min_{\mf{L},\mf{S}} \frac{1}{2p^{\alpha_1}}
\Vert \mf{\Sigma}_{n} - (\mf{L}+\mf{S}) \Vert_{F}^2+\mathcal{P}_{\gamma}(\mf{L},\mf{S}),\label{fro}
\end{equation}
to the solution of problem \eqref{prob_ld_rescaled}.
%to the rescaled heuristics of this paper:
%\begin{equation}\widetilde{\phi}^{(ld)}(\mf{L},\mf{S})=\min_{\mf{L},\mf{S}}
%0.5 \ln \det (\mf{I}_p+p^{-2\alpha_1}(\mf{\Sigma}_{n} - (\mf{L}+\mf{S}))(\mf{\Sigma}_{n} - (\mf{L}+\mf{S}))')+\psi
%\Vert \mf{L} \Vert_{*} + \rho \Vert \mf{S} \Vert_{1}.\label{logdet}
%\end{equation}
We write $$\widetilde{\phi}^{(F)}(\mf{L},\mf{S})=\widetilde{\phi}^{(F)}_D(\mf{L},\mf{S})+\mathcal{P}(\mf{L},\mf{S}),
\qquad \mbox{with} \qquad \widetilde{\phi}^{(F)}_D(\mf{L},\mf{S})=0.5p^{-\alpha_1}\Vert \mf{\Sigma}_{n} - (\mf{L}+\mf{S}) \Vert_{F}^2,$$
and
$$\widetilde{\phi}^{(ld)}(\mf{L},\mf{S})=\widetilde{\phi}^{(ld)}_D(\mf{L},\mf{S})+\mathbb{P}_{\gamma}(\mf{L},\mf{S}),
\qquad \mbox{with} \qquad \widetilde{\phi}^{(ld)}_D(\mf{L},\mf{S})=0.5\ln \det (\mf{I}_p+p^{-2\alpha_1}(\mf{\Sigma}_{n} - (\mf{L}+\mf{S}))^2).$$
The Supplement (Section B.1) shows that
$\widetilde{\phi}_{D}^{'(ld)}(\mf{L},\mf{S})=(\mf{I}_p+p^{-2\alpha_1}\mf{\Delta}_n\mf{\Delta}_n')^{-1}p^{-\alpha_1}\mf{\Delta}_n,$
where $\mf{\Delta}_n=(\mf{\Sigma}_{n} - (\mf{L}+\mf{S}))$. We can easily derive
that $\widetilde{\phi}_{D}^{'(F)}(\mf{L},\mf{S})=p^{-\alpha_1}\mf{\Delta}_n$.
%$p^{-\alpha_1}(\mf{\Sigma}_{n} - (\mf{L}+\mf{S}))$.

Let us define the pair of solutions
$(\widehat{\mf{L}}_n^{(F)},\widehat{\mf{S}}_n^{(F)})=\arg\min_{\mf{L},\mf{S}}\widetilde{\phi}^{(F)}(\mf{L},\mf{S})$
and $(\widehat{\mf{L}}_n^{(ld)},\widehat{\mf{S}}_n^{(ld)})=\arg\min_{\mf{L},\mf{S}}\widetilde{\phi}^{(ld)}(\mf{L},\mf{S}),$
with $\widehat{\mf{\Sigma}}_n^{(F)}=\widehat{\mf{L}}_n^{(F)}+\widehat{\mf{S}}_n^{(F)}$ and
$\widehat{\mf{\Sigma}}_n^{(ld)}=\widehat{\mf{L}}_n^{(ld)}+\widehat{\mf{S}}_n^{(ld)}$.
The following important theorem holds
(see Section C in the Supplement for the proof).

\begin{theorem}\label{thm_comp}
%Let $\nu=\frac{1}{2}$. Then,
%Set $\nu={1}/{2}$ in Assumption \ref{ass_alg}.
%Under all the conditions of Theorem \ref{thm_main}, \far{it holds for any $p$ and $n$}:\\
Suppose that Assumption \ref{lowerbounds} holds and $\alpha_{\mathcal{Y}} \geq 0.77155$.
Under the conditions of Proposition \ref{11},
%assuming that $\Vert\widehat{\mf{\Delta}}_n\Vert_{2} \leq 1/3$,
%\far{assuming that $p^{-2\alpha_1}\Vert \widehat{\mf{\Delta}}_n \Vert< (3p)^{-1}$ with $\Vert\mf{L}^{*}\Vert_{2}<3p$},
%if $\alpha_{\mathcal{Y}}=1$ and $\nu={1}/{2}$,
%and the convexity condition $\Vert\mf{\Delta}_{n} \Vert_{2} \leq \frac{1}{3\delta_\phi p}$ is respected.
for all $p \in \N$ it holds with probability $1-O(1/n^2)$ as $n \to \infty$:
\begin{equation}
\frac{\max g_\gamma(\widehat{\mf{S}}^{(ld)}-\mf{S}^{*},\widehat{\mf{L}}^{(ld)}-\mf{L}^{*})}
{\max g_\gamma(\widehat{\mf{S}}^{(F)}-\mf{S}^{*},\widehat{\mf{L}}^{(F)}-\mf{L}^{*})} \leq 1.\nonumber
%\label{ineq_ld}
\end{equation}
\end{theorem}
Theorem \ref{thm_comp} states that the minimax error bound in $g_\gamma$-norm of $(\widehat{\mf{L}}_n^{(ld)},\widehat{\mf{S}}_n^{(ld)})$
is systematically not larger than the corresponding bound of $(\widehat{\mf{L}}_n^{(F)},\widehat{\mf{S}}_n^{(F)})$.
This holds because the gradient of $\widetilde{\phi}_D(\mf{L},\mf{S})$ exists if and only if the series whose limit is
$(\mf{I}_p+p^{-2\alpha_1}\mf{\Delta}_n\mf{\Delta}_n')^{-1}$ converges, and any norm of $(\mf{I}_p+p^{-2\alpha_1}\mf{\Delta}_n\mf{\Delta}_n')^{-1}$
is smaller than $1$ if Proposition \ref{conv_delta} holds. This implies that $g_\gamma(\widetilde{\phi}^{(ld)}_D(\mf{L},\mf{S}))\leq g_\gamma(\widetilde{\phi}^{(F)}_D(\mf{L},\mf{S}))$.
We may finally state the main results of the paper, showing algebraic and parametric consistency for $\widehat{\mf{L}}_n^{(ld)}$, $\widehat{\mf{S}}_n^{(ld)}$, and $\widehat{\mf{\Sigma}}_n^{(ld)}$.
\begin{theorem}\label{thm_main}
Suppose that Theorem \ref{thm_comp} holds, Assumptions \ref{tails}-\ref{alg} are met,
%$m_p\ll C(p^{\alpha})$ for $q=0$,
%\textbf{IF}:\\
%1) exactly
%Define $\psi_{0}=\frac{f_{\delta_\epsilon}(p,n)}{\xi(\mathcal{T}(\mf{L}^{*}))}$ and $\rho_{0}=\gamma\psi_{0}$.
%$\psi={p^{\alpha_{1}}}{\psi_{0}}$, and $\rho=\gamma p^{\delta}\psi_{0}$,
%\lambda=
%\bigl\{
%\begin{array}{rl}
% & \frac{1}{2} \leq \alpha \leq 1,\\
%\frac{1}{\xi(\mathcal{T})}\frac{p^{\alpha}}{\sqrt{n}} & 0\leq \alpha < \frac{1}{2},
%\end{array}
%\bigr.
%where $\gamma \in [9\xi(\mathcal{T}(\mf{L}^{*})),\frac{1}{6\mu(\Omega(\mf{S}^{*}))}]$.
%\textbf{Identification} conditions:\\
%2)$\mu(\Omega)\xi(\mathcal{T})\leq \frac{1}{54}$;\\
%3)\\
%%1) and 2) are needed for \textbf{identification}.\\ %(impose a bound on the curvature of T).\\
%4)
%5)
%6)
%$n\geq p$ (needed for ${\Sigma}_n$).\\ %%possibly different!
%In addition, suppose that \far{$n>\overline{\delta}_n p^{6\delta_{1}-2+2(\far{\alpha_{1}}-\far{\alpha_r})}$}
%<\overline{\delta}_n p^{6\delta_{1}}
%\far{for some $\overline{\delta}_n$ such that $\overline{\delta}_n>0$.}
and $\delta_{1} \leq {\alpha_r}/{3}$. %and $\frac{\sqrt{r}\kappa_{S}}{\kappa_{L}}\leq \frac{1}{54}$.
%$\frac{p^{2-2\alpha_r}}{n}=o(1)$ and
Then, there exists a positive real $\kappa$ independent of $p$ and $n$
such that, for all $p \in \N$ as $n\to\infty$, the pair of solutions defined in \eqref{obj} satisfies:
%\\\vspace{0.1cm}
%Moreover, the matrix losses for each component are bounded as follows:
\begin{inparaenum}
\item[(i)] $\mathrm{P} ({p^{-\alpha_{1}}} \Vert\widehat{\mf{L}}_n-\mf{L}^{*}\Vert_{2} \leq \kappa\psi_{0}) = 1 - O(1/n^2)$;
\item[(ii)] $\mathrm{P} (\Vert\widehat{\mf{S}}_n-\mf{S}^{*}\Vert_{\infty} \leq \kappa\rho_{0}) = 1 - O(1/n^2)$;
%\frac{1}{p^{\delta_{1}}}\
%\end{compactenum}
%\far{In addition, if there exists some $\underline{\delta}_n>0$ such that
%$n>\underline{\delta}_n p^{2(\far{\alpha_{1}}-\kappa_{0})}$}, then
%\begin{compactenum}
\item[(iii)] $\mathrm{P} (\mathrm{rk}(\widehat{\mf{L}}_n)=\mathrm{rk}(\mf{L}^{*})) = 1 - O(1/n^2)$;
%\end{inparaenum}
%Further, if $p^{1-\delta_1}=o(\sqrt{n})$, then: \begin{inparaenum}
\item[(iv)] $\mathrm{P} (\mathrm{sgn}(\widehat{\mf{S}}_n)=\mathrm{sgn}(\mf{S}^{*})) = 1 - O(1/n^2)$.
\end{inparaenum}
\end{theorem}

%Theorem \ref{thmALCE} is proved in Appendix \ref{pr_thmALCE}. %is proved in \cite{farne2020large}.
%The proof can be found in Appendix \ref{proofs}.
%In the end, since
%$\mf{\Delta}_n$ is random, we recall the random matrix theory result from \cite{farne2020large}, that prescribes that
%\begin{lemma}\label{lemmaprob}
%$$\Vert\mf{\Delta}_n \Vert \leq C\frac{p^{\alpha}}{\sqrt{n}},$$
%\end{lemma}
%which allows to set $g_\gamma(\mathcal{A}^{\dag}\mf{\Delta}_n)\leq\psi=C\frac{p^{\alpha}}{\xi(\mathcal{T})\sqrt{n}}$.
%It follows that $\widetilde{r}$ is bounded by $\frac{80}{9}\psi$.
%Finally, naming the solutions to problem \eqref{logdet} $\mf{\widehat{L}}^{(ld)}$, $\mf{\widehat{S}}^{(ld)}$, and $\mf{\widehat{\Sigma}}^{(ld)}=\mf{\widehat{L}}^{(ld)}+\mf{\widehat{S}}^{(ld)}$, and the solutions to problem \eqref{fro} $\mf{\widehat{L}}^{fro}$, $\mf{\widehat{S}}^{fro}$, and $\mf{\widehat{\Sigma}}^{fro}=\mf{\widehat{L}}^{fro}+\mf{\widehat{S}}^{fro}$,
%our main theorem can be stated.

%\subsection*{Proof of Corollary \ref{coroll_main}}\label{pr_corollALCE}
%Once defined $\phi_S=C s' \xi(\mathcal{T}) \psi$ and $\phi=C(s'\xi(\mathcal{T})+1)\psi$,
%where $s'=\Vert \mf{S}^{*} \Vert_{0,v}$ is the maximum number of non-zeros per row--column in $\mf{S}^{*}$,
%we can state the following Corollary from \cite{farne2020large} (proved in Appendix \ref{pr_corollALCE}).
\begin{corollary}\label{coroll_main}
Under all the assumptions and conditions of Theorem \ref{thm_main},
%there exists a positive real $C$ independent of $p$ and $n$
%such that,
the following statements hold for all $p \in \N$ as $n\to\infty$:
%$\Vert \widehat{\mf{S}}_n-\mf{S}^{*}\Vert _{2}~\leq \phi_S $, $\Vert \widehat{\mf{\Sigma}}_n-\mf{\Sigma}^{*}\Vert _{2} \leq \phi$, $\Vert \widehat{\mf{S}}_n^{-1}-\mf{S}^{*-1}\Vert _{2} \leq \phi_S$,
%and $\Vert \widehat{\mf{\Sigma}}_n^{-1}-\mf{\Sigma}^{*-1}\Vert _{2} \leq \psi$, with probability approaching $1$.
\begin{inparaenum}
\item[(i)] $\mathrm{P} ({p^{-\delta_{1}}}\Vert \widehat{\mf{S}}_n-\mf{S}^{*}\Vert_{2}
\leq \kappa f_{\delta_\epsilon}(p,n))\to 1$;
\item[(ii)] $\mathrm{P} ({p^{-(\alpha_{1}+\delta_{1})}}\Vert \widehat{\mf{\Sigma}}_n-\mf{\Sigma}^{*}\Vert _{2} \leq \kappa f_{\delta_\epsilon}(p,n))\to 1$;
\item[(iii)] $\mathrm{P} (\lambda_p(\widehat{\mf{S}}_n)>0) \to 1$;
\item[(iv)] $\mathrm{P} (\lambda_p(\widehat{\mf{\Sigma}}_n)>0) \to 1$.
%\item $\mathcal{P} (\Vert \widehat{\mf{S}}_n^{-1}-\mf{S}^{*-1}\Vert_{2} \leq C \phi_S) \to 1$;
%\item $\mathcal{P} (\Vert \widehat{\mf{\Sigma}}_n^{-1}-\mf{\Sigma}^{*-1}\Vert_{2} \leq C \phi) \to 1$.
\end{inparaenum}
%Further, supposing that {$\lambda_p(\mf{S}^{*})=O(p^{\alpha_1-1-\varepsilon})$ and $\lambda_p({\mf{\Sigma}}^{*})=O(p^{\alpha_1-1-\varepsilon})$}
%for some $\varepsilon>0$,
%%and $\sum_{i=1}^p\sum_{j=1}^p \mathbbm{1}(\mf{S}_{ij}^{*} \ne 0) \leq O(p^{2\delta})$,
%the following statements hold for all $p \in \N$ as $n \to \infty$:\\
%\begin{inparaenum}
%\item[5.] $\mathcal{P} \left({p^{-\delta_{1}}}{p^{-2(1-\alpha_1+\varepsilon)}}\Vert\widehat{\mf{S}}_n^{-1}-\mf{S}^{*-1}\Vert_{2}
%\leq \kappa f_{\delta_\epsilon}(p,n))\right)\to 1$;\\
%\item[6.] $\mathcal{P} \left({p^{-(\alpha_{1}+\delta_{1})}}{p^{-2(1-\alpha_1+\varepsilon)}}\Vert\widehat{{\mf{\Sigma}}}^{-1}-{\mf{\Sigma}}^{*-1}\Vert_{2}
%\leq \kappa f_{\delta_\epsilon}(p,n))\right)\to 1$.
%\end{inparaenum}
%%
Further, supposing that $\lambda_p(\mf{S}^{*})=O(1)$ and $\lambda_p(\mf{\Sigma}^{*})=O(1)$,
%and $\sum_{i=1}^p\sum_{j=1}^p \mathbbm{1}(\mf{S}_{ij}^{*} \ne 0) \leq O(p^{2\delta})$,
the following statements hold for all $p \in \N$  as $n \to \infty$:
\begin{inparaenum}
\item[(v)] $\mathrm{P} ({p^{-\delta_{1}}}\Vert\widehat{\mf{S}}_n^{-1}-\mf{S}^{*-1}\Vert_{2}
\leq \kappa f_{\delta_\epsilon}(p,n))\to 1$;
\item[(vi)] $\mathrm{P} ({p^{-(\alpha_{1}+\delta_{1})}}\Vert\widehat{\mf{\Sigma}}_n^{-1}-\mf{\Sigma}^{*-1}\Vert_{2}
\leq \kappa f_{\delta_\epsilon}(p,n))\to 1$.
\end{inparaenum}
%%\item if $\lambda_{p}({\mf{S}^{*}})>\phi_{S}$, then  $\widehat{\mf{S}}_n$
%%is positive definite;
%%\item if $\lambda_{p}({\mf{\Sigma}^{*}})>\phi$, then $\widehat{\mf{\Sigma}}_n$ is positive definite;
%%\item if $\lambda_{p}({\mf{S}^{*}}) \geq 2 \phi_{{S}}$, then $\widehat{\mf{S}}_n^{-1}$ is positive definite;
%%\item if $\lambda_{p}({\mf{\Sigma}^{*}}) \geq 2 \phi$, then $\widehat{\mf{\Sigma}}_n^{-1}$ is positive definite.
%\end{inparaenum}
\end{corollary}

Theorem \ref{thm_main} and Corollary \ref{coroll_main} establish the algebraic
and parametric consistency of the estimator pair in \eqref{obj}. Their proofs can be found in the Supplement Section C,
and rely on the consistency of the sample covariance matrix expressed in Lemma A.3 in the Supplement under Assumptions \ref{tails}-\ref{sparsity}.
The following remarks clarify the most relevant theoretical aspects.

%\far{Specificare che i TASSI, la latent rank/rsp recovery, e la positiva definitezza di $\widehat{\mf{S}}_n$ e $\widehat{\mf{\Sigma}}_n$
%sono di per sè non-asintotiche. Esprimere la condizione, e mostrare che essa vale con probabilità 1 per $p$ fissato ed $n$ che diverge, e vale con alta probabilità per $n$ piccolo.}

%\far{Note that the function $\mathcal{L}^{(ld)}(\mf{\Sigma},\mf{\Sigma}_n)$ is convex with probability $1$ as $n \to \infty$, i.e.,
%$\mathcal{P}(\Vert\mf{L}^{*}\Vert_{2}^{-1}\Vert\widehat{\mf{\Delta}}_{n}\Vert \geq (3p)^{-1}) \to 0$,
%because %there exists a value $\underbar{n}$ such that for any $n>\underbar{n}$
%$\mathcal{P}(\Vert\mf{L}^{*}\Vert_{2}^{-1}\Vert{\mf{\Delta}}_{n}\Vert >0) \to 0$ by Lemma \ref{Lemma_cons}
%and $\mathcal{P}(\Vert\mf{L}^{*}\Vert_{2}^{-1}\Vert\widehat{\mf{\Sigma}}_n-\mf{\Sigma}^{*}\Vert >0) \to 0$ by part 2 of Corollary \ref{coroll_main}
%for all $p\in\N$ as $n \to \infty$.}

\begin{remark}
Parts (i) and (ii) of  Theorem \ref{thm_main} establish the parametric consistency of \eqref{obj} according to Definition \ref{par_cons}.
Parts (i) and (ii) of Corollary \ref{coroll_main} contain subsequent convergence results in spectral norm.
%for the covariance matrix estimates $\widehat{\mf{L}}_n$, $\widehat{\mf{S}}_n$ and $\widehat{\mf{\Sigma}}_n=\widehat{\mf{L}}_n+\widehat{\mf{S}}_n$,
Parts (v) and (vi) of Corollary \ref{corollsimple} ensure the invertibility of $\widehat{\mf{S}}_n$ and $\widehat{\mf{\Sigma}}_n$
with probability tending to $1$ for all $p\in\N$ as $n\to\infty$.
\end{remark}

\begin{remark}
Parts (iii) and (iv) of Theorem \ref{thm_main} and Corollary \ref{coroll_main}, jointly considered,
ensure the algebraic consistency of \eqref{obj} according to Definition \ref{def_alg}. This result is established by adapting the results of \cite{chandrasekaran2012latent}
to take into account the random nature of the second derivative of the smooth loss $\mathcal{L}^{(ld)}(\mf{\Sigma},\mf{\Sigma}_{n})$
(see the Supplement Section B.2),
%(see Propositions \ref{11} and \ref{12} in the Supplement),
and by controlling the manifolds containing $\mf{L}^{*}$ and $\mf{S}^{*}$ as in Assumptions \ref{alg}-\ref{ass_alg},
causing the term $O(p^{\delta_1})$ to appear in the rates of $\widehat{\mf{S}}_n$ and $\widehat{\mf{\Sigma}}_n$ (parts (i) and (ii)  of Corollary \ref{coroll_main}).
The condition $\delta_{1} \leq {\alpha_r}/{3}$ is needed to ensure the compatibility of Assumptions \ref{eigenvalues}-\ref{sparsity} and \ref{alg}-\ref{lowerbounds} (see the proof), while the prevalence of the latent factor structure versus the residual one
is preserved by the condition $\delta_1<\alpha_r$, resulting from Assumptions \ref{eigenvalues} and \ref{sparsity}.
\end{remark}

%\begin{remark}
%The requirement $\delta_1>0$ of Assumption \ref{sparsity}(i) ensures that the equation
%$\xi(\mathcal{T}(\mf{L}^{*}))={\sqrt r}/(\kappa_{L} p)^{\delta_{1}}$ (Assumption \ref{alg})
%leads to an increasingly small identifiability error as $p\to\infty$. In this respect, a high dimension is rather a blessing than a curse for our method.
%Also, note that
%\end{remark}

\begin{remark}\label{rem_tails}
We note that Theorem \ref{thm_main} and Corollary \ref{coroll_main} hold as well for the solution pair of \eqref{obj_all}
with $\mathcal{L}(\mf{\Sigma},\mf{\Sigma}_{n})=\widetilde{\mathcal{L}}^{(F)}(\mf{\Sigma},\mf{\Sigma}_{n})$, that corresponds to ALCE estimator \citep{farne2020large,farne2024large},
since the Hessian of $\widetilde{\mathcal{L}}^{(F)}(\mf{\Sigma},\mf{\Sigma}_{n})$ is constant and equal to $\mf{I}_p \otimes \mf{I}_p$, i.e., $\widetilde{\mathcal{L}}^{(F)}(\mf{\Sigma},\mf{\Sigma}_{n})$ is globally convex.
In addition, in \cite{farne2020large,farne2024large} sub-exponential tails are assumed for factors and residuals, implying the existence of infinite moments, while Assumption \ref{tails} here only imposes the second moment to exist, although this leads to stronger requirements in terms of sample size to achieve consistency.
%{(also see Remark \ref{size_requirements})}.
In the case $\delta_\epsilon=\infty$, $f_{\delta_\epsilon}(p,n)$ reduces to $\sqrt{\ln(p)/{n}}$ into all rates,
thus obtaining the same rates as in \cite{farne2020large,farne2024large}
(see the proof of Theorem 1 in \cite{bickel2008covariance} to grasp technical reasons).
Theorem \ref{thm_simple} and Corollary \ref{corollsimple} constitute a specialized case of Theorem \ref{thm_main} and Corollary \ref{coroll_main}
under the sub-exponential tails condition (i) of Section \ref{main_res}.
\end{remark}

\begin{remark}\label{size_requirements}
Should we allow $p$ to diverge to infinity as well as $n$, the condition $f_{\delta_\epsilon}(p,n)\to 0$ would read as ${p^{2/(1+\delta_\epsilon)}}=o(\sqrt{n})$, i.e. ${p^{4/(1+\delta_\epsilon)}}=o(n)$,
which equals to ${p^2}/{n}=o(1)$ if $\delta_\epsilon=1$ (fourth moment existence),
to ${p^{4/3}}/{n}=o(1)$ if $\delta_\epsilon=2$ (sixth moment existence),
to ${p}/{n}=o(1)$ if $\delta_\epsilon=3$ (eight moment existence).
%The condition $\frac{p^{6\delta_1-2\alpha_r}{p^{4/(1+\delta_\epsilon)}}}{n} \to 0$ would be needed to ensure $\lambda_r(\mf{L}^{*})\simeq p^{\alpha_r}$ when $\delta_1 \in [\frac{\alpha_r}{3},\frac{1}{2}]$.
The condition $p^{1-\delta_1}f_{\delta_\epsilon}(p,n)=o(1)$, equivalent to ${p^{(1-\delta_1)+{4/(1+\delta_\epsilon)}}}=o(n)$,
would instead rise to ensure the recovery of the residual sparsity pattern, i.e. the sparsistency of $\widehat{\mf{S}}_n$. %(\far{cite Lam!})
%Conditions $\lambda_p(\mf{S}^{*})=O(p^{\alpha_1-1-\varepsilon})$ and $\lambda_p({\mf{\Sigma}}^{*})=O(p^{\alpha_1-1-\varepsilon})$
%for some $\varepsilon>0$ are needed to ensure compatibility with Assumption \ref{sparsity}(iv).
\end{remark}

%%BASTA!!!%%

%OK! Diagonal of $S$ unconstrained?

%This is the worst case scenario, with $\xi(T(\mf{L}^{*}))=2$. In the best case scenario, with $\xi(T(\mf{L}^{*}))=\sqrt{r/p}$,
%we get the rate $$\kappa(\kappa_{L} \sqrt{p\ln(p)/n}++\kappa_{LS}\sqrt{\ln(p)/rn}).$$

%Since $\gamma \geq 9 \sqrt{r}/\kappa_{L}$, $$=\kappa(\kappa_{L} p\sqrt{r\ln(p)/n}+\kappa_{LS}\sqrt{p\ln(p)/n}),$$
%with $\kappa_{LS} \geq 9$.
%This is the worst case scenario, with $\xi(T(\mf{L}^{*}))=2$. In the best case scenario, with $\xi(T(\mf{L}^{*}))=\sqrt{r/p}$,
%we get the rate $$\kappa(\kappa_{L} \sqrt{p\ln(p)/n}++\kappa_{LS}\sqrt{\ln(p)/rn}).$$

%\far{What happens with $\delta_1>0$?}
%\far{Move later, to the remarks of the main theorem?

\section{Computation}\label{practice}

%We refer to Algorithms 1 and 2 in the Supplement Section E to compute in practice
%$(\widehat{\mf{L}}_n^{(ld)},\widehat{\mf{S}}_n^{(ld)})$ and $(\widehat{\mf{L}}_n^{(F)},\widehat{\mf{S}}_n^{(F)})$,
%respectively.

For each threshold pair $(\psi,\rho)$ (we refer to Section E.2 in the Supplement for their initialisation), we can calculate
$\widehat{\mf{\Sigma}}^{(ld)}_{\rm{A}}(\psi,\rho)=\widehat{\mf{L}}^{(ld)}_{\rm{A}}(\psi,\rho)+\widehat{\mf{S}}^{(ld)}_{\rm{A}}(\psi,\rho)$ from Algorithm 1 (see Section E.2)
%$\widehat{\mf{\Sigma}}_n^{(F)}_{\rm{A}}(\psi,\rho)=\widehat{\mf{L}}_n^{(F)}_{\rm{A}}(\psi,\rho)+\widehat{\mf{S}}^{(F)}_{\rm{A}}(\psi,\rho)$,
and $\widehat{\mf{\Sigma}}^{(F)}_{\rm{A}}(\psi,\rho)=\widehat{\mf{L}}^{(F)}_{\rm{A}}(\psi,\rho)+\widehat{\mf{S}}^{(F)}_{\rm{A}}(\psi,\rho)$ from Algorithm 2 (see Section E.2).
Following \cite{farne2020large,farne2024large}, we also perform the unshrinkage of estimated latent eigenvalues,
as this operation improves the sample total loss as much as possible in the finite sample. We thus get
the UNALCE (UNshrunk ALCE) estimates as:
%\begin{eqnarray}
$$\widehat{\mf{L}}_{\rm{U}}=\widehat{\mf{U}}_{\rm{A}}(\widehat{\mf{\Lambda}}_{\rm{A}}+\psi \mf{I}_r)\widehat{\mf{U}}_{\rm{A}}', \; %\label{unshr1}\\
\mathrm{diag}(\widehat{\mf{S}}_{\rm{U}})=\mathrm{diag}(\widehat{\mf{\Sigma}}_{\rm{A}})-\mathrm{diag}(\widehat{\mf{L}}_{\rm{U}}), \; %\label{unshr2}\\
\mathrm{off-diag}(\widehat{\mf{S}}_{\rm{U}})=\mathrm{off-diag}(\widehat{\mf{S}}_{\rm{A}}),$$ %\label{unshr3}
%\end{eqnarray}
where $\psi>0$ is any chosen eigenvalue threshold parameter.
In this way, for each threshold pair $(\psi,\rho)$,
%in the bivariate grid $\psi_{init},\rho_{init}$,
we can derive the pairs of estimates
%\begin{eqnarray}
$\bigl(\widehat{\mf{L}}^{(ld)}_{\rm{U}}(\psi,\rho),\widehat{\mf{S}}^{(ld)}_{\rm{U}}(\psi,\rho)\bigr)$ or
%\nonumber\\
%$\bigl(\widehat{\mf{L}}^{(F)}_{\rm{U}}(\psi,\rho),\widehat{\mf{S}}^{(F)}_{\rm{U}}(\psi,\rho)\bigr)$,%\nonumber\\
$\bigl(\widehat{\mf{L}}^{(F)}_{\rm{U}}(\psi,\rho),\widehat{\mf{S}}^{(F)}_{\rm{U}}(\psi,\rho)\bigr)$,
%\nonumber
and as a consequence, the overall UNALCE estimates as
$\widehat{\mf{\Sigma}}^{(ld)}_{\rm{U}}(\psi,\rho)=\widehat{\mf{L}}^{(ld)}_{\rm{U}}(\psi,\rho)+\widehat{\mf{S}}^{(ld)}_{\rm{U}}(\psi,\rho)$ and
%$\widehat{\mf{\Sigma}}^{(F)}_{\rm{U}}(\psi,\rho)=\widehat{\mf{L}}^{(F)}_{\rm{U}}(\psi,\rho)+\widehat{\mf{S}}^{(F)}_{\rm{U}}(\psi,\rho)$,
$\widehat{\mf{\Sigma}}^{(F)}_{\rm{U}}(\psi,\rho)=\widehat{\mf{L}}^{(F)}_{\rm{U}}(\psi,\rho)+\widehat{\mf{S}}^{(F)}_{\rm{U}}(\psi,\rho)$.

Then, %for each $h=1,\ldots,100$,
given the latent variance proportions $\widehat{\theta}(\psi,\rho)_{\rm{A}}=(\mathrm{tr}(\widehat{\mf{\Sigma}}(\psi,\rho)_{\rm{A}}))^{-1}{\mathrm{tr}(\widehat{\mf{L}}(\psi,\rho)_{\rm{A}})}$ and
$\widehat{\theta}(\psi,\rho)_{\rm{U}}=(\mathrm{tr}(\widehat{\mf{\Sigma}}(\psi,\rho)_{\rm{U}}))^{-1}{\mathrm{tr}(\widehat{\mf{L}}(\psi,\rho)_{\rm{U}})}$,
%Then, for each $h=1,\ldots,100$
we can select the optimal threshold pairs $(\psi_{U},\rho_{U})$ and $(\psi_{A},\rho_{A})$
by minimizing the MC criteria
\begin{eqnarray}
MC(\psi,\rho)_{U}&=&\max
\left\{\frac{{\widehat{r}\Vert\widehat{\mf{L}}(\psi,\rho)_{\rm{U}}}\Vert_{2}}{\widehat{\theta}(\psi,\rho)_{\rm{U}}},
\frac{{\Vert\widehat{\mf{S}}(\psi,\rho)_{\rm{U}}}\Vert_{1,v}}{{\gamma}(1-\widehat{\theta}(\psi,\rho)_{\rm{U}})}\right\},\nonumber\\
MC(\psi,\rho)_{A}&=&\max\left\{\frac{{\widehat{r}\Vert\widehat{\mf{L}}(\psi,\rho)_{\rm{A}}}\Vert_{2}}{\widehat{\theta}(\psi,\rho)_{\rm{A}}},
\frac{{\Vert\widehat{\mf{S}}(\psi,\rho)_{\rm{A}}}\Vert_{1,v}}{{\gamma}(1-\widehat{\theta}(\psi,\rho)_{\rm{A}})}\right\},\nonumber
\end{eqnarray}
where ${\gamma}=\psi^{-1}{\rho}$ is the ratio between the sparsity and the latent eigenvalue threshold
(see \cite{farne2020large} for more details). In this way, we can select the optimal threshold pairs
$(\psi_{A},\rho_{A})=\arg \min_{(\psi,\rho)} MC(\psi,\rho)_{A}$ and
$(\psi_{U},\rho_{U})=\arg \min_{(\psi,\rho)} MC(\psi,\rho)_{U}$.
%This procedure is applied both for Algorithms \ref{alg_ld} and \ref{alg_fro_2}, after
%For Algorithms \ref{alg_ld} and \ref{alg_fro_2},
%$\vf{\psi}_{init}$ and $\vf{\rho}_{init}$. %(see Section \ref{real}).
%%POET

\section{Real data analysis}\label{real}

In this section, we compute $(\widehat{\mf{L}}^{(ld)}_{\rm{U}}(\psi_{U},\rho_{U}),\widehat{\mf{S}}^{(ld)}_{\rm{U}}(\psi_{U},\rho_{U}))$ and $(\widehat{\mf{L}}^{(F)}_{\rm{U}}(\psi_{U},\rho_{U}),\widehat{\mf{S}}^{(F)}_{\rm{U}}(\psi_{U},\rho_{U}))$
on a selection of $361$
macroeconomic indicators provided by the European Central Bank
for $364$ systematically important Euro Area banks.
The indicators, taken in logarithms, mainly are financial items in the banks' balance sheet,
reported at a high level of granularity. All data refer to Q4-2014.

\begin{figure}
\begin{center}
\includegraphics[width=3in]{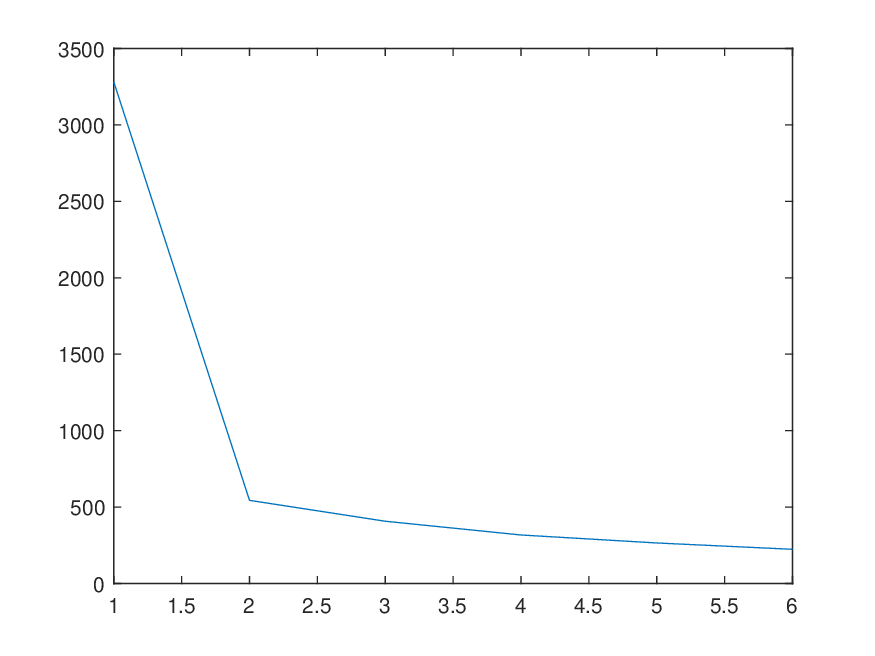}
\caption{ECB data: top six sample eigenvalues.}
\label{eigECB}
\end{center}
\end{figure}

%Algorithms \ref{alg_ld} and \ref{alg_fro_2} are applied on the input sample covariance matrix,
%setting the vectors of initial thresholds $\vf{\psi}_{init}$ and $\vf{\rho}_{init}$ with $n_{thr}=9$, $c^{\psi}_{thr}=1/2$, $c^{\rho}_{thr}=1$.
%%={i}/{p}$, with $i=0.05,0.1,0.2,0.\bar{3},0.5,1,2,5,10,20$,
%%and the vector of initial sparsity thresholds $\vf{\rho}_{init}$ as ${p}^{-1/2}\vf{\psi}_{init}$.
%%\footnote{Whenever one of the thresholds selected by the MC criterion \eqref{MC} lies in the grid extremes, it is advisable to shift the vector $\vf{\psi}_{init}$ to the left. In this case, the MC %%criterion selected the pair with positions $(9,2)$ for UNALCE-ld and $(5,3)$ for UNALCE-F.}
%\footnote{The grid positions of the optimal threshold pairs obtained by MC criterion are $(6,6)$ for UNALCE-LD and $(7,7)$ for UNALCE-F.}
Table \ref{ecb_data} reports estimation results.
The scree plot of sample eigenvalues (Figure \ref{eigECB}) highlights the presence of only one latent eigenvalue.
It follows that the estimated latent rank is $1$. The latent variance proportion is a bit smaller for $\widehat{\mf{L}}^{(ld)}_{\rm{U}}$ compared to $\widehat{\mf{L}}^{(F)}_{\rm{U}}$. More, $\widehat{\mf{S}}^{(ld)}_{\rm{U}}$ is a bit more selective than $\widehat{\mf{S}}^{(F)}_{\rm{U}}$
for residual nonzeros, and this results in a lower presence of non-zeros.

%\begin{figure}[H]
%\centering

\begin{table}
\caption{ECB data: estimation results for $(\widehat{\mf{L}}^{(ld)}_{\rm{U}},\widehat{\mf{S}}^{(ld)}_{\rm{U}})$ and $(\widehat{\mf{L}}^{(F)}_{\rm{U}},\widehat{\mf{S}}^{(F)}_{\rm{U}})$: $\widehat{r}$ is the estimated rank, $\widehat{\theta}$ is the estimated latent variance proportion, $\widehat{\rho}_{\widehat{\mf{S}}}$ is the estimated residual covariance proportion, $\widehat{\pi}_{\widehat{s}}$ is the estimated proportion of non-zeros, $\Vert \widehat{\mf{\Sigma}}-\mf{\Sigma}_{n} \Vert_{F}$ and $\Vert \widehat{\mf{\Sigma}}-\mf{\Sigma}_{n} \Vert_{2}$ are the Frobenius loss and the spectral loss from the sample covariance matrix $\mf{\Sigma}_{n}$, respectively.}
\label{ecb_data}
\begin{center}
\begin{tabular}{rrrrr}
\hline
            & $(\widehat{\mf{L}}^{(ld)}_{\rm{U}},\widehat{\mf{S}}^{(ld)}_{\rm{U}})$ & $(\widehat{\mf{L}}^{(F)}_{\rm{U}},\widehat{\mf{S}}^{(F)}_{\rm{U}})$\\
            \hline
            $\widehat{r}$ & 1 & 1\\
            $\widehat{\theta}$ & 0.2391 & 0.2473 \\
            $\widehat{\rho}_{\widehat{\mf{S}}}$  & 0.0058 & 0.0062\\
            $\widehat{\pi}_{\widehat{s}}$  & 0.0053 & 0.0060 \\
            $\Vert \widehat{\mf{\Sigma}}-\mf{\Sigma}_{n} \Vert_{F}$ & 1068.66 & 915.78 \\
            $\Vert \widehat{\mf{\Sigma}}-\mf{\Sigma}_{n} \Vert_{2}$ & 552.89 & 478.01 \\
\hline
\end{tabular}
\end{center}
\end{table}

From the estimates $\widehat{\mf{L}}$, $\widehat{\mf{S}}$, and $\widehat{\mf{\Sigma}}$,
we get for each variable $i=1,\ldots,p$ the estimated commonality as $(\widehat{\mf{\Sigma}}_{ii})^{-1}{\widehat{\mf{L}}_{ii}}$ and
the estimated idiosyncrasy as $(\widehat{\mf{\Sigma}}_{ii})^{-1}{\widehat{\mf{S}}_{ii}}$.
The estimated residual degree is obtained as $\mathrm{deg}_{\widehat{\mf{S}},i}=\sum_{j=1}^{p}\mathbbm{1}(\widehat{\mf{S}}_{ij} \ne 0).$
We obtain the spectral decomposition of $\widehat{\mf{L}}$ as $\widehat{\mf{U}}_{L}\widehat{\mf{D}}_{L}\widehat{\mf{U}}_{L}'$,
and we compute the vector of loadings $\widehat{\mf{U}}_{L}\widehat{\mf{D}}_{L}^{0.5}$.
%whose entries are plotted for
The estimated loadings are very similar for $\widehat{\mf{L}}^{(ld)}_{\rm{U}}$ and $\widehat{\mf{L}}^{(F)}_{\rm{U}}$,
although for $\widehat{\mf{L}}^{(ld)}_{\rm{U}}$ they are slightly more concentrated.
%in Figure \ref{fig:fac}.
They denote the contrast between loans and receivables and the rest of supervisory indicators,
as Table \ref{loadings} shows.

\begin{table}
\caption{Top three financial indicators by maximum and minimum loading value, as derived by $\widehat{\mf{U}}^{(ld)}_{L}\sqrt{\widehat{\mf{D}}_{L}^{(ld)}}$.\label{loadings}}
\tiny
\begin{center}
\begin{tabular}{rr}
\hline
Supervisory indicator & Loading value \\
\hline
Loans and advances - Central governments - Impaired assets & -0.0896\\
Available-for-sale assets - Other financial corporations - Carrying amount & -0.0855\\
Loans and advances - Central governments - Specific allowances for collectively assessed assets & -0.0839\\
Off-balance sheet exposures - Households - Nominal amount & 0.1143\\
Deposits - Current accounts / overnight deposits - Designated at fair value through profit or loss & 0.1163\\
Derivatives - Credit - Credit spread option - Notional amount - Sold & 0.1246
\end{tabular}
\end{center}
\end{table}

%\begin{figure}[t!]
%          \caption{UNALCE-ld and UNALCE-F - factor loadings}\label{fig:fac}
%                  \centering
%          \begin{tabular}{cc}
%              \includegraphics[width=.4\textwidth]{log-det_fac}&
%              \includegraphics[width=.4\textwidth]{fro_fac}
%              %{\footnotesize $\frac{\wh{\lambda}_j(\wt{\Sigma}(f_h))}p$, $j=1,2,3,4$.}&{\footnotesize %$\wh{\beta}(f_h)=\frac{\mathrm{tr}(\wh{L}(f_h))}{\mathrm{tr}(\wh{\Sigma}(f_h))}$.}\\
%          \end{tabular}
%\end{figure}

Table \ref{Comm_def} shows that the extracted factor is mainly connected to total assets, followed by
variables representing deposits and equity. Table \ref{deg_def} shows that the variables most connected with all the others are related to
loans and advances, debt securities and equity instruments held for trading.
Table \ref{Idio_def} shows that the most marginal variables \emph{wrt} the factor structure are related to financial assets for non-financial corporations
and derivatives for trading.

\begin{table}
\caption{\label{Comm_def}
Top six variables by estimated commonality, with respect to $\widehat{\mf{L}}^{(ld)}_{U}$.
This measure provides a ranking of the variables by systemic importance in determining the latent structure.}
%\centering
%\fbox{%
\tiny
\begin{center}
\begin{tabular}{rr}
\hline
Supervisory indicator & Commonality \\
\hline
Total assets & $0.5699$ \\
Deposits - Debt securities issued - Cumulative change in fair values due to changes in credit risk & $0.4203$\\
Financial assets held for trading - Equity instruments - At cost - Carrying amount & $0.4085$\\
\hline
\end{tabular}
\end{center}
\end{table}

\begin{table}
\caption{\label{deg_def}Top four variables by estimated degree, with respect to $\widehat{\mf{S}}^{(ld)}_{U}$.
This measure provides a ranking of the most connected variables with all the others, conditionally on the latent factor.}
%NFC stands for Non-Financial Corporations.
%\centering
%\fbox{%
\tiny
\begin{center}
\begin{tabular}{rr}
\hline
Supervisory indicator & Degree \\
\hline
Financial assets at fair value through profit or loss - Loans and advances - Changes due to credit risk & $6$ \\
Other debt securities issued - Held for trading & $5$\\
Loans and advances - Trade receivables - General governments & $5$\\
Financial assets held for trading - Equity instruments - Other financial corporations - Carrying amount & $5$\\
\hline
\end{tabular}
\end{center}
\end{table}

\begin{table}
\caption{\label{Idio_def}Top three variables by estimated idiosyncracy, with respect to $\widehat{\mf{S}}^{(ld)}_{U}$.
This measure provides a ranking of the variables by systemic irrelevance in determining the latent structure.}
%NFC stands for Non-Financial Corporations.
%\centering
%\fbox{%
\tiny
\begin{center}
\begin{tabular}{rr}
\hline
Supervisory indicator & Idiosyncracy \\
\hline
Financial assets at fair value through profit or loss - Non-financial corporations - Carrying amount & $1$ \\
Financial assets held for trading - Non-financial corporations - Carrying amount & $1$\\
Derivatives - Trading - Credit - Notional amount - Sold & $1$\\
\hline
\end{tabular}
\end{center}
\end{table}

%% log-det is a bit more parsimonious than Fro: less residual non-zeros and covariance,
%% smaller weights, less pervasive!

%% less extreme eigenvalues!
%% somehow better conditioned!

%% factor: + loans and advances, - the rest!
%% more comm: total assets!
%% more deg: Derivatives: Trading. Credit spread option. Sold.
%% more idio: Equity instruments. Of which: non-financial corporations. Credit amount.

%%

%% linearized algebraic problem!

\section{Conclusions}\label{concl}

In this paper, we study the estimation of large covariance matrices in high dimensions
under the low rank plus sparse assumption by minimizing a log-det heuristics augmented by a nuclear norm plus $\ell_1$-norm penalty.
In particular, we prove the local convexity and the Lipschitzianity of the proposed log-det heuristics,
which allows to solve the optimization problem via a proximal gradient algorithm. %We bound the curvature of the
%log-det heuristics under an appropriate random matrix theory framework. Then, %by adapting the results of \cite{chandrasekaran2012},
Then, we solve the low rank and sparse algebraic variety identification problem behind
the optimization problem, by showing that the log-det heuristics can be made locally convex and by controlling its Fisher information.
%and its gradient is a linear function.
Relying on these results, we prove
the algebraic and parametric consistency of the ensuing pair of low rank and sparse covariance matrix estimators.
We also prove that the same pair of estimators performs systematically not
worse than the corresponding estimator obtained by nuclear norm plus $\ell_1$-norm penalized Frobenius loss minimization.
A new solution algorithm, that also permits to control for the input threshold parameters, is proposed,
%A wide simulation study proves the validity of our theoretical results,
and an ECB supervisory data example shows the usefulness of our approach on a real dataset.

\bibliographystyle{agsm}
\bibliography{Bernardi_Farne_JASA}

\appendix

\section{Consistency of the sample covariance matrix}\label{Lemma_sample_2}

Let us consider the sample covariance matrix $\mf{\Sigma}_{n}=n^{-1}\sum_{k=1}^n \vf{x}_k\vf{x}_k'$,
assuming this form because, for the sake of simplicity, Assumption \ref{tails} imposes $\vf{x}$ to be mean centered.
Let us similarly define $\mf{\Sigma}_f=n^{-1}\sum_{k=1}^n \vf{f}_k\vf{f}_k'$, and
$\mf{\Sigma}_\epsilon=n^{-1}\sum_{k=1}^n \vf{\epsilon}_k\vf{\epsilon}_k'$.
%The following lemmas hold.

\begin{lemma}\label{lemma_prob_both}
Let $C(\delta_f)$ and $C(\delta_\epsilon)$ be two constants depending only on $\delta_f$ and $\delta_\epsilon$, respectively. Then, under Assumption \ref{tails} it holds:
\begin{eqnarray}
\Pr (\max_{i,j} \vert \mf{\Sigma}_{f,ij}-\mf{I}_r \vert \geq t)\leq r^2 M_{f} C(\delta_f)\frac{n^{-(1+\delta_f)/2}}{t^{1+\delta_f}};\nonumber\\%\label{prob_rate_f}
\Pr (\max_{i,j} \vert \mf{\Sigma}_{\epsilon,ij}-\mf{S}^{*} \vert \geq t)\leq p^2 M_{\epsilon} C(\delta_\epsilon)\frac{n^{-(1+\delta_\epsilon)/2}}{t^{1+\delta_\epsilon}}.\nonumber%\label{prob_rate_eps}
\end{eqnarray}
\end{lemma}

\paragraph{Proof}
%See section 2.3 in \cite{bickel2008covariance} (\far{argument more!})
This follows from \cite{bickel2008covariance}, inequality (24). \qed

\begin{lemma}\label{lemma_rate_both}
Under the assumptions of Lemma \ref{lemma_prob_both}, it holds
\begin{eqnarray}
\max_{i,j} \vert \mf{\Sigma}_{f,ij}-\mf{I}_r \vert=O_P\left(r^{\frac{2}{1+\delta_f}}\frac{1}{\sqrt{n}}\right);\label{rate_moment_f}\\
\max_{i,j} \vert \mf{\Sigma}_{\epsilon,ij}-\mf{S}^{*} \vert=O_P\left(p^{\frac{2}{1+\delta_\epsilon}}\frac{1}{\sqrt{n}}\right).\label{rate_moment_eps}
\end{eqnarray}
\end{lemma}

\paragraph{Proof}
%See section 2.3 in \cite{bickel2008covariance} (\far{argument more!})
This follows directly from Lemma \ref{lemma_prob_both} (also see \cite{bickel2008covariance}, Section 2.3). \qed

\begin{lemma}\label{Lemma_cons}
%Let ${\lambda}_r(\mf{\Sigma}_{n})$ be the $r-$th largest eigenvalue of the sample covariance matrix $\mf{\Sigma}_{n}=\frac{1}{n}\sum_{k=1}^n \vf{x}_k \vf{x}_k'$.
Under Assumptions \ref{tails}-\ref{sparsity}, it holds %\far{if $f_{\delta_\epsilon}(p,n) \to 0$
$$\frac{1}{p^{\alpha_{1}}}\Vert \mf{\Sigma}_{n}-\mf{\Sigma}^{*} \Vert_{2}=o_P(1)$$
for all $p \in \N$ with probability $1-O(1/n^2)$ as $n \to \infty$.
%${\lambda_r(\mf{\Sigma}_{n})} \simeq p^{\alpha_r}$.
%for some $C_{1}>0$.
\end{lemma}

\paragraph{Proof}

For any $t \geq 0$, we define $\mathcal{T}^{(H)}_{t}$, the hard-thresholding operator with parameter $t$, such that the $p \times p$
matrix $\mathcal{T}^{(H)}_{t}(\mf{M})$ has $(i,j)$ element $\mf{M}_{ij}$ if $\vert \mf{M}_{ij} \vert \geq t$, $0$ otherwise.

%On one hand, we note that, since $r+p-p=r \leq p$,
%dual Weyl inequality (see \cite{tao2011topics})
%can be applied, leading to
%\begin{equation}\label{bound_eig_left}
%\lambda_r(\mf{\Sigma}^{*})\geq\lambda_r(\mf{L}^{*})+\lambda_p(\mf{S}^{*}).
%\end{equation}
%From \eqref{bound_eig_left}, we can write
%\begin{equation}\label{bound_eig_left_{2}}
%\lambda_r(\mf{\Sigma}^{*})\succeq O(p^{\alpha_r})+O(p^{\delta_{1}})=O(p^{\alpha_r}),\nonumber
%\end{equation}
%because ${\lambda_r(\mf{L}^{*})} \simeq {p^{\alpha_r}}$ by Assumption \ref{eigenvalues}(i),
%and ${\lambda_p(\mf{S}^{*})}=O(p^{\delta_{1}})$ by Assumption \ref{sparsity}, with $\delta_{1}<\alpha_r$.
%
%On the other hand, Lidskii inequality (see \cite{tao2011topics}) leads to
%\begin{equation}\label{bound_eig_right}
%\lambda_r(\mf{\Sigma}^{*})\leq\lambda_r(\mf{L}^{*})+\sum_{j=1}^r\lambda_j(\mf{S}^{*}).
%\end{equation}
%From \eqref{bound_eig_right}, we can write
%\begin{equation}\label{bound_eig_right_{2}}
%\lambda_r(\mf{\Sigma}^{*})\preceq O(p^{\alpha_r})+O(r p^{\delta_{1}})=O(p^{\alpha_r}),\nonumber
%\end{equation}
%because ${\lambda_r(\mf{L}^{*})} \simeq p^{\alpha_r}$ by Assumption \ref{eigenvalues}(i),
%${\lambda_1(\mf{S}^{*})}=O(p^{\delta_{1}})$ with $\delta_{1}<\alpha_r$ by Assumption \ref{sparsity}{(i)},
%and $r$ is finite for all $p \in \N$ by Assumption \ref{eigenvalues}(ii).
%It follows that ${\lambda_r(\mf{\Sigma}^{*})} \simeq p^{\alpha_r}$.

Recalling that $\mf{\Sigma}_{n}=n^{-1}\sum_{k=1}^n \vf{x}_k \vf{x}_k^\top$ and $\vf{x}_k=\mf{B}\vf{f}_k+{\vf{\epsilon}}_k$,
where ${\vf{f}}_k$ and ${\vf{\epsilon}}_k$, $k\in\{1,\ldots,n\}$, are respectively the vectors of factor scores and residuals for each observation, we can decompose the error matrix ${\mf{E}}_n=\mf{\Sigma}_{n}-\mf{\Sigma}^*$ in four components as follows
(see \cite{fan2013large}):
$${\mf{E}}_n=\mf{\Sigma}_{n}-\mf{\Sigma}^*={\mf{D}}_{1}+{\mf{D}}_{2}+{\mf{D}}_3+{\mf{D}}_4,$$ where
%\begin{eqnarray}
${\mf{D}}_{1}={n}^{-1} \mf{B} \left(\sum_{k=1}^n \vf{f}_k \vf{f}_k^\top-\mf{I}_r\right)\mf{B}^\top$,%\nonumber\\
${\mf{D}}_{2}={n}^{-1} \sum_{k=1}^n \left({\vf{\epsilon}}_k {\vf{\epsilon}}_k^\top-\mf{S}^{*}\right)$,
%\nonumber\\
${\mf{D}}_3= {n}^{-1} \mf{B}\sum_{k=1}^n \vf{f}_k {\vf{\epsilon}}_k^\top$,
%\nonumber\\
${\mf{D}}_4={\mf{D}}_3^\top$.%\nonumber
%\end{eqnarray}

%%%%

%Following \cite{fan2013large}, we note that
%$$\Vert{\mf{D}}_{1}\Vert_{2} \leq \bigg\Vert\frac{1}{n} \sum_{k=1}^n \vf{f}_{ik}\vf{f}_{jk}-\mathrm{E}(\vf{f}_{ik}\vf{f}_{jk})\bigg\Vert_{2}\Vert{\mf{B}\mf{B}^\top}\Vert_{2}\leq rp^{\alpha_{1}}
%\bigg\Vert\frac{1}{n} \sum_{k=1}^n \vf{f}_{ik}\vf{f}_{jk}-\mathrm{E}(\vf{f}_{ik}\vf{f}_{jk})\bigg\Vert_{\infty},$$
%since $\mathrm{E}(\vf{f})=\vf{0}_r$ and $\mathrm{Var}(\vf{f})=\mf{I}_r$, and $$\bigg\Vert\frac{1}{n} \sum_{k=1}^n \vf{f}_{ik}\vf{f}_{jk}-\mathrm{E}(\vf{f}_{ik}\vf{f}_{jk})\bigg \Vert_{2} \leq r \bigg\Vert\frac{1}{n} \sum_{k=1}^n \vf{f}_{ik}\vf{f}_{jk}-\mathrm{E}(\vf{f}_{ik}\vf{f}_{jk})\bigg \Vert_{\infty}$$ by Assumption \ref{eigenvalues}(ii) and $\Vert{\mf{B}\mf{B}^\top}\Vert_{2}=O(p^{\alpha_{1}})$ by Assumption \ref{eigenvalues}(i).

Following \cite{fan2013large}, we note that
\begin{eqnarray}
\Vert{\mf{D}}_{1}\Vert_{2} &\leq& \bigg\Vert\frac{1}{n} \left(\sum_{k=1}^n \vf{f}_k \vf{f}_k^\top-\mf{I}_r\right)\bigg
\Vert_{2}\Vert{\mf{B}\mf{B}^\top}\Vert_{2}\nonumber\\
&\leq&  rp^{\alpha_{1}} {\mathrm{max}_{i,j \leq r} \biggl\vert  \frac{1}{n} \sum_{k=1}^n {f}^{i}_{k}{f}^{j}_{k}-\mathrm{E}[{f}^{i}_{k}{f}^{j}_{k}] \biggl\vert},\nonumber
\end{eqnarray}
since $\mathrm{E}[\vf{f}]=\vf{0}_r$ and $\mathrm{Var}[\vf{f}]=\mf{I}_r$ {by Assumption \ref{tails}}, $\Vert{\mf{B}\mf{B}^\top}\Vert_{2}=O(p^{\alpha_{1}})$ by Assumption \ref{eigenvalues}(i), and $$\bigg\Vert\frac{1}{n} \left(\sum_{k=1}^n \vf{f}_k \vf{f}_k^\top-\mf{I}_r\right)\bigg \Vert_{2} \leq r\bigg\Vert\frac{1}{n} \left(\sum_{k=1}^n \vf{f}_k \vf{f}_k^\top-\mf{I}_r\right)\bigg \Vert_{\infty}.$$
Let us define %$$f_{\delta_f}(r,n)=\frac{r^{2/(1+\delta_f)}}{n^{1/2}},$$
$$f_{\delta_f}(r,n)={\frac{\max\left(r^{{2}/{(1+\delta_f)}},\sqrt{\ln(r)}\right)}{n^{1/2}}},$$
with $\delta_f \in \R^{+}$ as defined in Assumption \ref{tails},
analogously to $f_{\delta_\epsilon}(p,n)$.
%$$f_{\delta_f}(r,n)=
%\left\{
%\begin{array}{rl}
%\frac{r^{2(1+\delta_f)}}{n^{1/2}}, & \delta_f \in \R^{+}, \nonumber\\
%\frac{1}{\sqrt{n}}, & \delta_f=+\infty, \nonumber
%\end{array}
%\right.
%$$
%analogously to $f_{\delta_\epsilon}(p,n)$.
Now, from \eqref{rate_moment_f} in Lemma \ref{lemma_rate_both}, it holds
\begin{equation}
\Vert{\mf{D}}_{1}\Vert_{\infty}=\bigg\Vert\frac{1}{n} \left(\sum_{k=1}^n \vf{f}_k \vf{f}_k^\top-\mf{I}_r\right)\bigg
\Vert_{\infty} \leq f_{\delta_f}(r,n).
\label{Lemma4_{1}}
\end{equation}
%and $$\bigg\Vert\frac{1}{n} \left(\sum_{k=1}^n \vf{f}_k \vf{f}_k^\top-\mf{I}_r\right)\bigg \Vert_{2} - \left(\sum_{k=1}^n \vf{f}_k \vf{f}_k^\top-\mf{I}_r\right)\bigg \Vert_{\infty},$$
%
%or, more extensively,
%$$\Pr(\mathrm{max}_{i,j \leq r} \biggl\vert \frac{1}{n} \sum_{k=1}^n {f}^{i}_{k}{f}^{j}_{k}-\mathrm{E}[{f}^{i}_{k}{f}^{j}_{k}])$$
%with probability $1-O(1/n^2)$ (${\tilde{C}}$ is a real positive constant).
Therefore, %applying \eqref{rate_moment_f} and recalling Assumption \ref{sparsity}(ii),
we obtain
\begin{equation}\Vert{\mf{D}}_{1}\Vert_{2} \leq rp^{\alpha_{1}}f_{\delta_f}(r,n)=
{\tilde{C}}_r \frac{p^{\alpha_{1}}}{\sqrt{n}}=O_P\left(\frac{p^{\alpha_{1}}}{\sqrt{n}}\right)\label{bound1}\end{equation}
{with} $\tilde{C}_r = r \sqrt{n} f_{\delta_f}(r,n)=r^{1+2/(1+\delta_f)}$, because Lemma \ref{lemma_prob_both} ensures that %\eqref{prob_rate_f} in
\begin{equation}\Pr (\Vert{\mf{D}}_{1}\Vert_{\infty} \geq t) \leq r^2 M_{f} C(\delta_f)\frac{n^{-(1+\delta_f)/2}}{t^{1+\delta_f}},\label{prob_f}
\end{equation}
where the \emph{rhs} tends to $0$ for all $p \in \N$ as $n \to \infty$ because $r$ is finite and independent of $p$ by Assumption \ref{eigenvalues}(ii).
%because Assumption \ref{pr} prescribes that $r=\delta_3 \ln(p)$ and $\ln(p)=o(n)$.

%%

Then, we note that the diagonal elements of the matrix ${\mf{S}^{*}}$
are bounded by a finite constant, due to Assumption \ref{sparsity}(ii).
Under Assumption \ref{tails}, %and the condition $f_{\delta_\epsilon}(p,n) \to 0$, as $n \to \infty$,
%%p inactive if \delta=\infty%%
%%discuss samplign theory, double!%%
%$p^2 M_{f} C(\delta_\epsilon)\frac{n^{-(1+\delta_\epsilon)/2}}{t^{1+\delta_\epsilon}}=o(1)$,
%which leads $to \frac{p^2}{n^{-(1+\delta_\epsilon)/2}}=o(1)$,
%{and the condition $\ln(p)/n \to 0$},
%(12) in \cite{bickel2008covariance} thus holds for the matrix ${\mf{S}^{*}}$, leading to:
we obtain from \eqref{rate_moment_eps} in Lemma \ref{lemma_rate_both} that
\begin{equation}
\Vert{\mf{D}}_{2}\Vert_{\infty}=\mathrm{max}_{i,j \leq p}
\biggl\vert  \frac{1}{n} \sum_{k=1}^n {\epsilon}_{i,k}{\epsilon}_{j,k}-\mathrm{E}({\epsilon}^{i}_{k}{\epsilon}^{j}_{k})\biggl
\vert \leq \tilde{C}_2 f_{\delta_\epsilon}(p,n),\label{Lemma4_{2}}
\end{equation}
%$\Vert{\mf{D}}_{2}\Vert_{\infty}=o_P(1)$,
and Lemma \ref{lemma_prob_both} ensures that %from \eqref{prob_rate_eps} in
\begin{equation}\Pr (\Vert{\mf{D}}_{2}\Vert_{\infty} \geq t) \leq p^2 M_{\epsilon} C(\delta_\epsilon)\frac{n^{-(1+\delta_\epsilon)/2}}{t^{1+\delta_\epsilon}},
\label{prob_eps}\end{equation}
where the r.h.s tends to $0$ for all $p \in \N$ as $n \to \infty$.
%if $f_{\delta_\epsilon}(p,n) \to 0$

%\begin{equation}
%\Vert{\mf{D}}_{2}\Vert_{\infty}=\mathrm{max}_{i,j \leq p}
%\biggl\vert  \frac{1}{n} \sum_{k=1}^n {\epsilon}^{i}_{k}{\epsilon}^{j}_{k}-\mathrm{E}({\epsilon}^{i}_{k}{\epsilon}^{j}_{k})\biggl
%\vert \leq \tilde{C}_2 \sqrt{\frac{\ln(p)}{{n}}},\label{Lemma4_{2}}
%\end{equation}
%that holds with probability $1-O(1/n^2)$.

Now, by the triangular inequality we can write
\begin{equation}
\Vert {\mf{D}}_{2}\Vert_{2} \leq \Vert {\mf{D}}^{(1)}_{2} \Vert_{2}+\Vert {\mf{D}}^{(2)}_{2} \Vert_{2},\label{bound_pre}
\end{equation}
where $\mf{D}^{(1)}_{2}=\mathcal{T}^{(H)}_{\Vert \mf{S}^{*} \Vert_{\mathrm{\mathrm{min,off}}}}(\mf{D}_{2})$,
with $\mathcal{T}^{(H)}$ hard-thresholding operator of parameter $\Vert \mf{S}^{*} \Vert_{\mathrm{\mathrm{min,off}}}$,
and $\mf{D}^{(2)}_{2}=\mf{D}_{2}-\mf{D}^{(1)}_{2}$.
Since by Assumption \ref{sparsity}(i) $\Vert \mf{S}^{*} \Vert_{0,v} = O(p^{\delta_{1}})$,
%$\Vert\mf{S}^{*}\Vert_{\mathrm{\mathrm{min,off}}}\Vert \leq \Vert \mf{S}^{*} \Vert_{\infty}=O(1)$,
%and $\Vert \mf{S}^{*} \Vert_{\infty}=O(p^{\kappa_{0}})$,
it follows from \eqref{Lemma4_{2}} that
%, if $f_{\delta_\epsilon}(p,n) \to 0$, as $n \to \infty$,
\begin{equation}
\Vert {\mf{D}}^{(1)}_{2} \Vert_{2}\leq \Vert {\mf{D}}^{(1)}_{2} \Vert_{0,v} \Vert {\mf{D}}^{(1)}_{2} \Vert_{\infty} \leq %O(p^{\delta_{1}})O\left(\sqrt{\frac{\ln(p)}{{n}}}\right);\label{Thr_top}\\
\tilde{C}_2 \delta_2 p^{\delta_{1}}f_{\delta_\epsilon}(p,n).\label{Thr_top}\\
\end{equation}
Similarly, it follows from \eqref{Lemma4_{2}} that
%, as $n \to \infty$,
\begin{equation}
\Vert {\mf{D}}^{(2)}_{2} \Vert_{2} \leq \Vert {\mf{D}_{2}}^{(2)} \Vert_{0,v} \Vert {\mf{D}_{2}}^{(2)} \Vert_{\infty} < {p}\Vert\mf{S}^{*}\Vert_{\mathrm{\mathrm{min,off}}}.\label{Thr_down}
%< O(p^{1-\delta_{1}})O(\Vert\mf{S}^{*}\Vert_{\infty})=O(p^{1-\delta_{1}})O(1).\label{Thr_down}
\end{equation}
Imposing the condition ${p}\Vert\mf{S}^{*}\Vert_{\mathrm{\mathrm{min,off}}}=o(p^{\delta_{1}}f_{\delta_\epsilon}(p,n))$, we get
${p}^{1-\delta_1}\Vert\mf{S}^{*}\Vert_{\mathrm{\mathrm{min,off}}}=o(f_{\delta_\epsilon}(p,n))$, which leads to
${p}^{1-\delta_1}\Vert\mf{S}^{*}\Vert_{\mathrm{\mathrm{min,off}}}=o(1)$
for all $p \in \N$ as $n \to \infty$. Therefore, combining \eqref{Thr_top} and \eqref{Thr_down},
%since $\Vert {\mf{D}}^{(2)}_{2} \Vert_{2} = O(p^{1-\delta_{1}})o(1)=o(1)$,
by Assumption \ref{sparsity}(iii) it follows from \eqref{bound_pre} that, for all $p \in \N$ as $n \to \infty$,
%\begin{equation}\Vert{\mf{D}}_{2}\Vert_{\infty}=
%\biggl\vert  \frac{1}{n} \sum_{k=1}^n \vf{\epsilon}_{ik}\vf{\epsilon}_{jk}-\mathrm{ E}(\vf{\epsilon}_{ik}\vf{\epsilon}_{jk})\biggl
%\vert \leq C' \Vert \mf{S}^{*} \Vert_{\infty} \sqrt{\frac{\ln(p)}{{n}}},
%\end{equation}
%from which it descends
\begin{equation}
\Vert {\mf{D}}_{2} \Vert_{2} \leq %\Vert {\mf{D}}_{2} \Vert_{0,v} \Vert{\mf{D}}_{2}\Vert_{\infty}=
\tilde{C}_2 \delta_2 p^{\delta_{1}}f_{\delta_\epsilon}(p,n)=O_P\left(\frac{p^{\delta_{1}}}{\sqrt{n}}\right).\label{bound2}
\end{equation}
%where we used the fact that $$\mathrm{P}(\Vert {\mf{D}}_{2} \Vert_{0,v} {\leq} \Vert\mf{S}^{*}\Vert_{0,v}) \to 1$$
%as $n \to \infty$.

At this stage, we consider $\Vert {\mf{D}}_{3} \Vert_{2}=\bigg\Vert n^{-1} \sum_{k=1}^n \vf{f}_k\vf{\epsilon}_k^\top \bigg\Vert_{2}$.
We first observe that\\
$\bigg\Vert n^{-1} \sum_{k=1}^n \vf{f}_k\vf{\epsilon}_k^\top \bigg\Vert_{2} \leq \bigg\Vert n^{-1} \sum_{k=1}^n \vf{f}_k\vf{\epsilon}_k^\top \bigg\Vert_{F}$.
We then write
\begin{eqnarray}\bigg\Vert n^{-1} \sum_{k=1}^n \vf{f}_k\vf{\epsilon}_k^\top \bigg\Vert_{F}=\sqrt{\sum_{i=1}^r\sum_{j=1}^p \widehat{\mathrm{Cov}}(f_{i},\epsilon_{j})^2} \leq \sqrt{\sum_{i=1}^r\sum_{j=1}^p {\widehat{\mathrm{V}}(f_{i})}{\widehat{\mathrm{V}}(\epsilon_{j})}}.\label{ineq_D3_1}\end{eqnarray}
We note that
$\sqrt{\sum_{i=1}^r\sum_{j=1}^p {\widehat{\mathrm{V}}(f_{i})}{\widehat{\mathrm{V}}(\epsilon_{j})}}$ converges to
$\sqrt{\sum_{i=1}^r\sum_{j=1}^p {{\mathrm{V}}(f_{i})}{{\mathrm{V}}(\epsilon_{j})}}$ for each $i \in \{1,\ldots,r\}$ and $j \in \{1,\ldots,p\}$
for all $p \in \N$ as $n \to \infty$ by (\ref{Lemma4_{1}}) and (\ref{Lemma4_{2}}).
Therefore, we can write
\begin{eqnarray}
\sqrt{\sum_{i=1}^r\sum_{j=1}^p {{\mathrm{V}}(f_{i})}{{\mathrm{V}}(\epsilon_{j})}}&\leq&\sqrt{r \sum_{j=1}^p \mf{S}_{jj}^{*}}
\leq\sqrt{r o_P(p^{\alpha_{1}})}=o_P(p^{\alpha_{1}/2}),\label{ineq_D3_2}\end{eqnarray}
by Cauchy--Schwarz inequality and Assumptions \ref{eigenvalues}(ii), \ref{sparsity}(iv) and \ref{tails}.
From \eqref{ineq_D3_1} and \eqref{ineq_D3_2}, it follows that \begin{equation}\bigg\Vert n^{-1} \sum_{k=1}^n \vf{f}_k\vf{\epsilon}_k^\top \bigg\Vert_{2}=o_P(p^{\alpha_{1}/2})\label{ineq_3}\end{equation} for all $p \in \N$ as $n \to \infty$.
Consequently, for all $p \in \N$ as $n \to \infty$, from \eqref{ineq_3} we obtain with probability $1-O(1/n^2)$ the following claim:
\begin{eqnarray}
\Vert{\mf{D}}_3\Vert_{2} &\leq& \bigg\Vert \frac{1}{n} \sum_{k=1}^n \vf{f}_k\vf{\epsilon}_k^\top\bigg\Vert_{2}\Vert{\mf{B}}\Vert_{2}
%C'\left(p^{\frac{\delta_{1}}{2}}f_{\delta_\epsilon}(p,n)\right)
=o_P(p^{\frac{\alpha_{1}}{2}})O\left(p^{\frac{\alpha_{1}}{2}}\right)=o_P(p^{\alpha_{1}})
%C' p^{\frac{\alpha_{1}}{2}+\frac{\delta_{1}}{2}}f_{\delta_\epsilon}(p,n),
\label{bound3}\end{eqnarray}
because $\Vert{\mf{B}}\Vert=O(p^{{\alpha_{1}}/{2}})$ by Assumption \ref{eigenvalues}(i).

Putting \eqref{bound1}, \eqref{bound2}, and \eqref{bound3} together, the following inequality is proved for all $p \in \N$ as $n \to \infty$ with probability $1-O(1/n^2)$:
\begin{equation}\Vert{{\mf{\Sigma}}}_{n}-\mf{\Sigma}^{*}\Vert_{2}\leq C' \frac{p^{\alpha_{1}}}{\sqrt{n}},\label{bound_sample_2}\end{equation}
because $\delta_{1} < \alpha_r \leq \alpha_{1}$ by Assumptions \ref{eigenvalues}{(i)} and \ref{sparsity}{(i)},
and because of \eqref{prob_f} and \eqref{prob_eps} for all $p \in \N$ as $n \to \infty$.
It follows that
$$\Vert \mf{\Sigma}_{n}-\mf{\Sigma}^{*} \Vert_{2} = O_P\left(\frac{p^{\alpha_{1}}}{\sqrt{n}}\right),$$
for all $p \in \N$ as $n \to \infty$ with probability $1-O(1/n^2)$,
%as $p$ and $n$ diverge to infinity,
which proves the lemma. \qed
%In fact, if $\delta > \alpha$, the condition of Theorem \ref{thm1} $\lambda_r(\mf{L}^{*})> {(C_{2} \psi)}/{\xi^2(T)}$ would result in $\lambda_r(\mf{L}^{*})>C_{2} p^{\delta}$, thus violating Assumption \ref{spikyPOETnew} under Assumption \ref{pr}.

%In other words, the bound \eqref{boundtop} means
%$\Vert\mf{E}_n\Vert_{2} \rightarrow 0 \Leftrightarrow {p^{\alpha_{1}}}/{\sqrt{n}} \rightarrow 0$.
%Exploiting the basic property $\Vert.\Vert_{\infty}\leq \Vert.\Vert_{2}$
%and the minimum for $\gamma$ in the range of Theorem \ref{thm1}, we can also write
%$\Vert\mf{E}_n\Vert_{\infty} \rightarrow 0 \Leftrightarrow \xi(\mathcal{T}){p^{\alpha}}/{\sqrt{n}} \rightarrow 0$.

%\section{Lemma \ref{Lemma_bmax}}\label{Lemma_sample_max}

\begin{lemma}\label{Lemma_bmax}
Under Assumptions \ref{tails}, \ref{eigenvalues}(ii), and \ref{sparsity}(ii),
it holds
$$\Vert \mf{\Sigma}_{n}-\mf{\Sigma}^{*} \Vert_{\infty} = o_P(1)$$
for all $p \in \N$ with probability $1-O(1/n^2)$ as $n \to \infty$.
\end{lemma}

\paragraph{Proof}
Under Assumptions \ref{tails} and \ref{eigenvalues}{(ii)}, we get with probability $1-O(1/n^2)$:
\begin{equation}
\Vert{\mf{D}}_{1}\Vert_{\infty} \leq \bigg\Vert\frac{1}{n} \left(\sum_{k=1}^n \vf{f}_k \vf{f}_k^\top-\mf{I}_r\right)\bigg\Vert_{\infty}\Vert{\mf{B}\mf{B}^\top}\Vert_{\infty}\leq C'f_{\delta_f}(r,n),\label{bmax1}
\end{equation}
because the r.h.s tends to $0$ for all $p \in \N$ as $n \to \infty$ by definition,
and $$\Vert{\mf{B}\mf{B}^\top}\Vert_{\infty}\leq \left(\max_{j\in\{1,\ldots,p\}} \Vert \vf{b}_j \Vert\right)^2 \leq
r^2 \Vert \mf{B} \Vert_{\infty}^2=O(1).$$

Under Assumptions \ref{tails} and \ref{sparsity}{(ii)}, \eqref{Lemma4_{2}} ensures that, with probability $1-O(1/n^2)$,
\begin{equation}
\Vert{\mf{D}}_{2}\Vert_{\infty}=\mathrm{max}_{i,j \leq p}
\biggl \vert \frac{1}{n} \sum_{k=1}^n {\epsilon}^{i}_{k}{\epsilon}^{j}_{k}-\mathrm{E}({\epsilon}^{i}_{k}{\epsilon}^{j}_{k})
\biggl \vert \leq C' f_{\delta_\epsilon}(p,n),\label{bmax2}
\end{equation}
where the r.h.s tends to $0$ for all $p \in \N$ as $n \to \infty$.

Under Assumptions \ref{tails}, \ref{eigenvalues}{(ii)}, and \ref{sparsity}{(ii)}, from \eqref{bmax1} and \eqref{bmax2} we get, with probability $1-O(1/n^2)$,
\begin{equation}
\Vert{\mf{D}}_{3}\Vert_{\infty}=\biggl\Vert \frac{1}{n} \sum_{k=1}^n \vf{f}_k\vf{\epsilon}_k^\top\biggl\Vert_{\infty} \leq C' f_{\delta_\epsilon}(p,n),\label{bmax3}
\end{equation}
where the r.h.s tends to $0$ for all $p \in \N$ as $n \to \infty$.

Putting together \eqref{bmax1}, \eqref{bmax2}, \eqref{bmax3}, with probability $1-O(1/n^2)$ for all $p \in \N$ as $n \to \infty$, we obtain
\begin{eqnarray}
\Vert \mf{\Sigma}_{n}-\mf{\Sigma}^{*} \Vert_{\infty}&=&O_P(\max(f_{\delta_f}(p,n),f_{\delta_\epsilon}(p,n)))\nonumber\\
&=&O_P(f_{\delta_\epsilon}(p,n))),
\label{bound_sample_max}
\end{eqnarray}
which proves the lemma. \qed

\section{Mathematical and probabilistic guarantees}\label{math_anal}

%Following \cite{nesterov2013gradient},
%we can apply the proximal gradient algorithm to solve the optimization problem \eqref{obj}.

The mathematical properties of $\mathcal{L}^{(ld)}(\mf{\Sigma},\mf{\Sigma}_n)$ have been extensively studied in \cite{bernardi2022log}.
This section extends the recalled mathematical analysis results related to
$\mathcal{L}^{(ld)}(\mf{\Sigma},\mf{\Sigma}_n)=0.5\ln \det (\mf{I}_p+\mf{\Delta}_{n}\mf{\Delta}_{n}')$.
To this end, we recall the first and the second derivatives of
$\mathcal{L}^{(ld)}(\mf{\Sigma},\mf{\Sigma}_n)$.

\subsection{First derivative}
%Therefore,

Let us denote the spectral decomposition of the non-definite symmetric random matrix
$\mf{\Delta}_n$ as $\mf{U}_{\Delta} \mf{\Lambda}_{\Delta} \mf{U}_{\Delta}'$, so that the spectral decomposition of the semi-definite matrix
$\mf{\Delta}_n \mf{\Delta}_n'$ results to be $\mf{U}_{\Delta} \mf{\Lambda}_{\Delta}^2 \mf{U}_{\Delta}'$.
We study the first derivative of $\phi_D(\mf{L},\mf{S})$.
%$=0.5\ln \det (\mf{I}_p+\mf{\Delta}_n\mf{\Delta}_n')$ \emph{wrt} (``with respect to'')
%$\mf{L}$ and $\mf{S}$, where $\mf{\Delta}_{n}=\mf{\Sigma}-\mf{\Sigma}_n$ and $\mf{\Sigma}=\mf{L}+\mf{S}$.
%The following proposition then holds.

%%%%%%%%%%%%%%%%%%%%%%%%%%%%%%%%%%%%%%%%%%%%%%%%%%%%

\begin{proposition}\label{first_der}
%\begin{inparaenum}
%Suppose that $\Vert\mf{\Sigma}^{*}\Vert_{2}<3p$.
Define the matrix $\mf{D}_{\Delta}=(\mf{I}_p-(\mf{\Lambda}_{\Delta}^{-2}+\mf{I}_p)^{-1})\mf{\Lambda}_{\Delta}$.
Assume the conditions of Proposition \ref{conv_delta} (main paper) with $\delta_\phi=1$. Then, the following equations hold:
\begin{eqnarray}
%Then, $\phi_{D}'(\mf{L},\mf{S})$ exists.
\phi'_{D}(\mf{L},\mf{S})=\frac{\partial \phi_D(\mf{L},\mf{S})}{\partial \mf{L}}=
\frac{\partial \phi_D(\mf{L},\mf{S})}{\partial \mf{S}}=(\mf{I}_p+\mf{\Delta}_{n}\mf{\Delta}_{n}')^{-1}\mf{\Delta}_{n}\label{grad}\\
(\mf{I}_p+\mf{\Delta}_{n}\mf{\Delta}_{n}')^{-1}\mf{\Delta}_{n}=\mf{U}_{\Delta}\mf{D}_{\Delta}\mf{U}_{\Delta}'.\label{corr}
%\mf{\Delta}_{n}^{-2} 2 \mf{\Delta}_{n} = \mf{\Delta}_{n}^{-1}.
\end{eqnarray}
%\end{inparaenum}
\end{proposition}

%Proposition \ref{first_der} gives us ground to wonder about the nature of the gradient in \eqref{grad}, i.e. $(\mf{I}_p+\mf{\Delta}_{n}\mf{\Delta}_{n}')^{-1}\mf{\Delta}_{n}$.

\paragraph{Proof}%\begin{proof}
The proof of equation \eqref{grad} is the proof of equation (6) in \cite{bernardi2022log},
provided that $(\mf{I}_p+\mf{\Delta}_{n}\mf{\Delta}_{n}')^{-1}$ exists.
Since it holds
\begin{eqnarray}
\phi'(\mf{L},\mf{S})^{(ld)}&=&
(\mf{I}_p+\mf{\Delta}_{n}\mf{\Delta}^{'}_n)^{-1}\mf{\Delta}_{n}=\sum_{j=0}^{\infty}(-\mf{\Delta}_{n}\mf{\Delta}^{'}_n)^{j}\mf{\Delta}_{n}=
\sum_{j=0}^{\infty}(-1)^{j}(\mf{\Delta}_{n})^{2j+1},\label{phi'_series1}
\end{eqnarray}
we can observe that series \eqref{phi'_series1} converges if and only if $\lambda_1(\mf{\Delta}_{n})^{2}<1$.
At this stage, since Proposition \ref{conv_delta} with $\delta_\phi=1$ ensures that $\lambda_1(\mf{\Delta}_{n})\leq (3p)^{-1}$,
the existence of $\phi'(\mf{L},\mf{S})^{(ld)}$ is ensured. %$\Vert\mf{\Sigma}^{*}\Vert_{2}<3p$ ensures

To prove equation \eqref{corr}, recalling Woodbury formula, we write
\begin{eqnarray}
(\mf{I}_p+\mf{\Delta}_n\mf{\Delta}_n')^{-1}&=&(\mf{I}_p+\mf{U}_{\Delta} \mf{\Lambda}_{\Delta}^2 \mf{U}_{\Delta}')^{-1}
=\mf{I}_p-\mf{U}_{\Delta}(\mf{\Lambda}_{\Delta}^{-2}+\mf{U}_{\Delta}'\mf{U}_{\Delta})^{-1}\mf{U}_{\Delta}'\nonumber\\
&=&\mf{I}_p-\mf{U}_{\Delta}(\mf{\Lambda}_{\Delta}^{-2}+\mf{I}_p)^{-1}\mf{U}_{\Delta}',\nonumber
\end{eqnarray}
because $\mf{U}_{\Delta}$ is orthogonal.
Therefore, we obtain
\begin{eqnarray}
(\mf{I}_p+\mf{\Delta}_n\mf{\Delta}_n')^{-1}\mf{\Delta}_n&=&(\mf{I}_p-\mf{U}_{\Delta}(\mf{\Lambda}_{\Delta}^{-2}+\mf{I}_p)^{-1}\mf{U}_{\Delta}')\mf{\Delta}_n\nonumber\\
&=&\mf{\Delta}_n-\mf{U}_{\Delta}(\mf{\Lambda}_{\Delta}^{-2}+\mf{I}_p)^{-1}\mf{U}_{\Delta}'\mf{\Delta}_n.\nonumber
%&=&\mf{U}_{\Delta} \mf{\Lambda}_{\Delta} \mf{U}_{\Delta}'-\mf{U}_{\Delta}
%(\mf{\Lambda}_{\Delta}^{-2}+\mf{I}_p)^{-1}\mf{\Lambda}_{\Delta} \mf{U}_{\Delta}'.\nonumber
\end{eqnarray}
Going on, since $\mf{U}_{\Delta}$ is orthogonal, we can write
\begin{eqnarray}
\mf{\Delta}_n-\mf{U}_{\Delta}(\mf{\Lambda}_{\Delta}^{-2}+\mf{I}_p)^{-1}\mf{U}_{\Delta}'\mf{\Delta}_n
&=&\mf{U}_{\Delta}\mf{\Lambda}_{\Delta}\mf{U}_{\Delta}'-\mf{U}_{\Delta}(\mf{\Lambda}_{\Delta}^{-2}+\mf{I}_p)^{-1}
\mf{U}_{\Delta}'\mf{U}_{\Delta}\mf{\Lambda}_{\Delta}\mf{U}_{\Delta}'\nonumber\\
&=&\mf{U}_{\Delta} [\mf{\Lambda}_{\Delta}-(\mf{\Lambda}_{\Delta}^{-2}+\mf{I}_p)^{-1}\mf{\Lambda}_{\Delta}] \mf{U}_{\Delta}'\nonumber
%&=&\mf{U}_{\Delta} [(\mf{I}_p-(\mf{\Lambda}_{\Delta}^{-2}+\mf{I}_p)^{-1})\mf{\Lambda}_{\Delta}] \mf{U}_{\Delta}',\nonumber
\end{eqnarray}
which proves the thesis because $\mf{D}_{\Delta}=(\mf{I}_p-(\mf{\Lambda}_{\Delta}^{-2}+\mf{I}_p)^{-1})\mf{\Lambda}_{\Delta}$.
%$\mf{U}_{\Delta}(\mf{\Lambda}_{\Delta}^{-2}+\mf{I}_p)^{-1}\mf{U}_{\Delta}')\mf{\Delta}_n$\\
%$\mf{U}_{\Delta}(\mf{\Lambda}_{\Delta}^{-2}+\mf{I}_p)^{-1}\mf{U}_{\Delta}')\mf{U}_{\Delta} \mf{\Lambda}_{\Delta} \mf{U}_{\Delta}'=$\\
%$\mf{U}_{\Delta}(\mf{\Lambda}_{\Delta}^{-2}+\mf{I}_p)^{-1}\mf{\Lambda}_{\Delta} \mf{U}_{\Delta}'=$\\
%$\mf{U}_{\Delta} \mf{\Lambda}_{\Delta} \mf{U}_{\Delta}'-\mf{U}_{\Delta}(\mf{\Lambda}_{\Delta}^{-2}+\mf{I}_p)^{-1}\mf{\Lambda}_{\Delta} \mf{U}_{\Delta}'=$\\
%$\mf{U}_{\Delta} [\mf{\Lambda}_{\Delta}-(\mf{\Lambda}_{\Delta}^{-2}+\mf{I}_p)^{-1}\mf{\Lambda}_{\Delta}] \mf{U}_{\Delta}'$.\\
%\end{proof}
\qed

Proposition \ref{first_der} provides a direct link between the first derivative of %the smooth component in \eqref{obj},
$\phi_D(\mf{L},\mf{S})=0.5\ln \det \varphi(\mf{\Sigma})$ (where $\varphi(\mf{\Sigma})=(\mf{I}_p+\mf{\Delta}_{n}\mf{\Delta}_{n}')$)
and the matrix $\mf{D}_{\Delta}$, which is a $p \times p$ diagonal matrix
whose $i$-th element, $i=1,\ldots,p$, is $\mf{D}_{ii}=\mf{\Lambda}_{\Delta,ii}\bigl(1-{(1+\mf{\Lambda}_{\Delta,ii}^{-2})^{-1}}\bigr)$.
This result is of extreme interest, because the matrix $\mf{D}_{\Delta}$ performs a sample eigenvalue correction: it shrinks down very large or very small sample eigenvalues,
i.e. those most affecting matrix inversion. This is why minimizing $\phi_D(\mf{L},\mf{S})=0.5\ln \det \varphi(\mf{\Sigma})$
leads to eigenvalue-regularized covariance matrix estimates. This is a main difference between $\mathcal{L}^{(F)}(\mf{\Sigma},\mf{\Sigma}_n)=0.5\Vert \mf{\Sigma}_n-(\mf{L}+\mf{S})\Vert_{F}^2=
0.5\sum_{i=1}^p \lambda^2_i(\mf{\Delta}_n)$, whose gradient is just $\mf{\Delta}_n$,
and $\mathcal{L}^{(ld)}(\mf{\Sigma},\mf{\Sigma}_n)=0.5\ln \det (\mf{I}_p+\mf{\Delta}_{n}\mf{\Delta}_{n}')$,
%%in place of \frac{1}{2}\Vert \mf{\Sigma}_n-\mf{\Sigma)\Vert_{F}$.
because the gradient of the smooth loss $\mathcal{L}^{(ld)}(\mf{\Sigma},\mf{\Sigma}_n)$ is designed to avoid
rank deficiency by intrinsically regularizing sample eigenvalues, unlike $\mathcal{L}^{(F)}(\mf{\Sigma},\mf{\Sigma}_n)$.
%\begin{remark}\label{alg_null}
%$[\mf{\Lambda}_{\Delta}-(\mf{\Lambda}_{\Delta}^{-2}+\mf{I}_p)^{-1}\mf{\Lambda}_{\Delta}]=$
%$(\mf{\Lambda}_{\Delta}^{-2}+\mf{I}_p)^{-1}$, which is a diagonal matrix whose $i$-th diagonal element is
%$1/(1+\lambda_i^{-2})$. Therefore, the $i$-th diagonal element of the diagonal matrix $\mf{\Lambda}_{\Delta}-(\mf{\Lambda}_{\Delta}^{-2}+\mf{I}_p)^{-1}\mf{\Lambda}_{\Delta}$ is
%$\lambda_i-\frac{\lambda_i}{1+\lambda_i^{-2}}$.
%This expression can be intended as an eigenvalue regularization one, as it tends to shrink towards $0$
%very large or very small eigenvalues.
%Fisher information matrix $\mathcal{I}^{*}(\mf{L},\mf{S})$ can also be computed as $\phi_D'(\mf{L},\mf{S})\otimes\phi_D'(\mf{L},\mf{S})$.
%Defining the matrix $\mf{\Psi}_{\Delta}=\mf{U}_{\Delta} [\mf{\Lambda}_{\Delta}-(\mf{\Lambda}_{\Delta}^{-2}+\mf{I}_p)^{-1}\mf{\Lambda}_{\Delta}] \mf{U}_{\Delta}$,
%and recalling Proposition \ref{corr_der}, we can write $\mathcal{I}^{*}(\mf{L},\mf{S})=\mf{\Psi}_{\Delta}\otimes\mf{\Psi}_{\Delta}$. Qualitatively,

%$\mf{U}_{\Delta} [\mf{\Lambda}_{\Delta}-(\mf{\Lambda}_{\Delta}^{-2}+\mf{I}_p)^{-1}\mf{\Lambda}_{\Delta}] \mf{U}_{\Delta}'$
%$\mf{U}_{\Delta}[\mf{\Lambda}_{\Delta}-(\mf{\Lambda}_{\Delta}^{-2}+\mf{I}_p)^{-1}\mf{\Lambda}_{\Delta}]\mf{U}_{\Delta}'=\Psi_{F}$.
%Therefore, we can write $\mathcal{I}(\mf{L},\mf{S})=\Psi_{F} \bigotimes \Psi_{F}$.
%Eigenvalues towards zero, like the no error situation!
%\end{remark}

%%
\subsection{Second derivative}

%From Proposition \ref{second_der}, it follows that if  %we get
Let us define the $p^2 \times p^2$ Hessian of $\phi_D(\mf{L},\mf{S})$ for $i,j,h,k \in \{1,\ldots,p\}$ as
\begin{eqnarray}
\bigl(\mathrm{Hess} \phi_D(\mf{L},\mf{S})\bigr)_{i j h k}&=&\frac{\partial^{2}}{\partial \sigma_{ij}\partial \sigma_{hk}}0.5\ln\det\varphi(\mf{\Sigma})
= \bigl(0.5\mathrm{Hess} \ln\det \varphi(\mf{\Sigma})\bigr)_{i j h k}.\nonumber
\end{eqnarray}
%$=\frac{\partial^{2}}{\partial \sigma_{ij}\partial \sigma_{hk}}\frac{1}{2} \phi_D(\mf{L},\mf{S})$.
We study the behaviour of $\bigl(\mathrm{Hess} \phi_D(\mf{L},\mf{S})\bigr)_{i j h k}$ by the following relevant proposition,
which holds whenever the random matrix $\mf{\Delta}_{n}$ is identically zero. %for $\mathrm{Hess} \phi_D(\mf{L},\mf{S})$
\begin{proposition}\label{second_der}
Suppose that $\mf{\Delta}_{n}=\mf{\Sigma}-\mf{\Sigma}_n = \mf{0}_{p \times p}$. Then, under the conditions of Proposition \ref{first_der}, it holds
\begin{equation}
\mathrm{Hess}\phi_D(\mf{L},\mf{S}) = \mf{I}_p \otimes \mf{I}_p. \label{fishone_bis}
\end{equation}
\end{proposition}

\paragraph{Proof}%\begin{proof}
Equation \eqref{fishone_bis} follows from $\mathrm{Hess}\phi_D(\mf{L},\mf{S}) = 0.5\mathrm{Hess} \ln\det \varphi(\mf{\Sigma})$
and the equation $0.5\mathrm{Hess} \ln\det \varphi(\mf{\Sigma})=\mf{I}_p \otimes \mf{I}_p$,
proved as equation (11) in \cite{bernardi2022log}. \qed

%because
%\begin{equation}
%$$\frac{\partial^{2}}{\partial \sigma_{ij}\partial \sigma_{hk}}\frac{1}{2} \ln\det\varphi(\mf{\Sigma}) =  \bigl(\frac{1}{2}\mathrm{Hess}\ln\det \varphi(\mf{\Sigma})\bigr)_{i j h k}
%= \bigl(\mf{I}_p \otimes \mf{I}_p\bigr) _{i j h k}.$$
%\end{equation}
%\end{proof}
Proposition \ref{second_der} will be crucial to study the random behaviour of the Hessian of $\phi_D(\mf{L},\mf{S})$ around $\mf{\Sigma}^{*}$,
which must be controlled to ensure solution stability for problem \eqref{obj}.
To this purpose, we observe that
$$\frac{\partial^{2}{\phi_{D}(\mf{L},\mf{S})}}{\partial^{2}{\mf{L}}}%=\frac{\partial{\phi'_{D}(\mf{L},\mf{S})}}{\partial{\mf{L}}} \qquad \mbox{and} \qquad
=\frac{\partial{\phi'_{D}(\mf{L},\mf{S})}}{\partial{\mf{L}}}=\frac{\partial{\phi'_{D}(\mf{L},\mf{S})}}{\partial{\mf{S}}}=\frac{\partial^{2}{\phi_{D}(\mf{L},\mf{S})}}{\partial^{2}{\mf{S}}}.$$

\subsection{Fisher information}\label{fish_bis}

In order to control the sensitivity of $\mathcal{I}^{*}(\widetilde{\mf{L}},\widetilde{\mf{S}})$ around $(\widetilde{\mf{L}}^{*},\widetilde{\mf{S}}^{*})$,
we set $\widetilde{\mf{\Sigma}}=\widetilde{\mf{\Sigma}}^{*}+\widetilde{\mf{W}}$, where $\widetilde{\mf{\Sigma}}$ is a generic $p \times p$ matrix s.t. $\widetilde{\mf{\Sigma}}=\widetilde{\mf{L}}+\widetilde{\mf{S}}$,
and we study the effect of its deviation $\widetilde{\mf{W}}$ from the target $\widetilde{\mf{\Sigma}}^{*}=\widetilde{\mf{L}}^{*}+\widetilde{\mf{S}}^{*}$
onto $\mathcal{I}^{*}(\widetilde{\mf{L}},\widetilde{\mf{S}})$. %This equals to postulate $\mf{\Sigma}_{n}=\mf{L}^{*}+\mf{S}^{*}+\mf{W}$.
We first recall that $\mf{\Delta}^{*}_n=\mf{\Sigma}^{*}-\mf{\Sigma}_{n}$, so that
$$\mf{\Delta}_n=\mf{\Sigma}-\mf{\Sigma}_{n}=\mf{\Sigma}^{*}-\mf{\Sigma}_{n}+\mf{W}=\mf{\Delta}^{*}_n+\mf{W}.$$
We need to ensure that the conditions of Proposition \ref{conv_delta} are satisfied for $\delta_\phi=1$ in the solution $\widehat{\mf{\Sigma}}_n$ of problem \eqref{obj},
i.e. to ensure that $-\widetilde{\mathcal{L}}^{(ld)}(\mf{\Sigma},\mf{\Sigma}_n)$ is concave for $\mf{\Sigma}=\widehat{\mf{\Sigma}}_n$.
%with probability approaching 1 as $n \to \infty$.
In other words, once defined $\widehat{\mf{\Delta}}_{n}=\widehat{\mf{\Sigma}}_n-\mf{\Sigma}_n$, we need that
%$$\mathcal{P}(\Vert\mf{L}^{*}\Vert_{2}^{-1}\Vert\widehat{\mf{\Delta}}_{n}\Vert \geq (3p)^{-1}) \to 0$$
$p^{-\alpha_{1}}\Vert\widehat{\mf{\Delta}}_{n}\Vert \leq 1/(3p)$, which is dominated by the condition $\Vert\widehat{\mf{\Delta}}_{n}\Vert \leq 1/3$.

Let us define %$\widetilde{\mf{\Delta}}_{n}=p^{-\alpha_1}\mf{\Delta} _{n}$,
the spectral decomposition of $\widetilde{\mf{\Delta}}_{n}$ as $\widetilde{\mf{U}}_{\Delta}\widetilde{\mf{\Lambda}}_{\Delta}\widetilde{\mf{U}}_{\Delta}$.
%$p^{-\alpha_1}\widetilde{\mf{\Lambda}}_{\Delta}=\mf{\Lambda}_{\Delta}$ and
For a generic pair $(\widetilde{\mf{L}},\widetilde{\mf{S}})$, we observe from \eqref{grad} and the definition of $\mathcal{I}^{*}(\widetilde{\mf{L}},\widetilde{\mf{S}})$ that
\begin{eqnarray}
\Vert\mathcal{I}^{*}(\widetilde{\mf{L}},\widetilde{\mf{S}})\Vert&=&\Bigg\Vert\mathrm{E}\left(\frac{\partial{\widetilde{\phi}'_{D}(\widetilde{\mf{L}},\widetilde{\mf{S}})}}{\partial{\widetilde{\mf{L}}}}\right)\Bigg\Vert
%=\Bigg\Vert\mathrm{E}\left(\frac{\partial{\phi'_{D}(\mf{L},\mf{S})}}{\partial{\mf{S}}}\right)\Bigg\Vert=
%=\Bigg\Vert\mathrm{E}(\mf{I}_p\otimes\mf{I}_p)\Bigg\Vert-
=\Bigg\Vert\mathrm{E}\left(\frac{\partial(\mf{I}_p+\widetilde{\mf{\Delta}}_{n}\widetilde{\mf{\Delta}}_{n}')^{-1}\widetilde{\mf{\Delta}}_{n}}{\partial{\widetilde{\mf{L}}}}\right)\Bigg\Vert.\nonumber
%&=&1-\Bigg\Vert\mathrm{E}\left(\frac{\partial{\widetilde{\mf{U}}_{\Delta}[\mf{I}_p+(\widetilde{\mf{\Lambda}}_{\Delta}^{-2}+\mf{I}_p)^{-1}]\widetilde{\mf{\Lambda}}_{\Delta}\widetilde{\mf{U}}_{\Delta}'}}{\partial{\widetilde{\mf{L}}}}\right)\Bigg\Vert.
%&=&1-(1+1/\lambda_1(\mf{\Delta}_n\mf{\Delta}_n'))^{-2}=
%1-(\frac{\lambda_1(\mf{\Delta}_n\mf{\Delta}_n')}{\lambda_1(\mf{\Delta}_n\mf{\Delta}_n')+1})^2\label{fish_bound}
\end{eqnarray}
Therefore, we directly compute
\begin{equation}
\frac{\partial^{2}\widetilde{\phi}_{D}(\widetilde{\mf{L}},\widetilde{\mf{S}})}{\partial^{2}{\widetilde{\mf{L}}}}
=\frac{\partial(\mf{I}_p+\widetilde{\mf{\Delta}}_{n}\widetilde{\mf{\Delta}}_{n}')^{-1}\widetilde{\mf{\Delta}}_{n}}{\partial{\widetilde{\mf{L}}}}.\nonumber
\end{equation}
We observe that $\widetilde{\mf{\Delta}}_{n}'=\widetilde{\mf{\Delta}}_{n}$,
$\frac{\partial\widetilde{\mf{\Delta}}_{n}}{\partial{\widetilde{\mf{L}}}}=\mf{I}_p$,
$\frac{\partial(\mf{I}_p+\widetilde{\mf{\Delta}}^{2}_{n})}{\partial{\widetilde{\mf{L}}}}=2\widetilde{\mf{\Delta}}_{n}$,
which leads to
\begin{eqnarray}
&&\frac{\partial[(\mf{I}_p+\widetilde{\mf{\Delta}}^{2}_{n})^{-1}\widetilde{\mf{\Delta}}_{n}]}{\partial{\widetilde{\mf{L}}}}
=-(\mf{I}_p+\widetilde{\mf{\Delta}}^{2}_{n})^{-1}(2\widetilde{\mf{\Delta}}_{n})(\mf{I}_p+\widetilde{\mf{\Delta}}^{2}_{n})^{-1}
=-(\mf{I}_p+\widetilde{\mf{\Delta}}^{2}_{n})^{-2}(2\widetilde{\mf{\Delta}}_{n}),\nonumber
\end{eqnarray}
because all the multiplicands commute.
Finally,
\begin{eqnarray}
\frac{\partial^{2}\widetilde{\phi}_{D}(\widetilde{\mf{L}},\widetilde{\mf{S}})}{\partial^{2}{\widetilde{\mf{L}}}}&=&(\mf{I}_p+\widetilde{\mf{\Delta}}^{2}_{n})^{-1}-(\mf{I}_p+
\widetilde{\mf{\Delta}}^{2}_{n})^{-2}(2\widetilde{\mf{\Delta}}_{n})
=(\mf{I}_p+\widetilde{\mf{\Delta}}^{2}_{n})^{-1}(\mf{I}_p-2(\mf{I}_p+\widetilde{\mf{\Delta}}^{2}_{n})^{-1}\widetilde{\mf{\Delta}}_{n}).\nonumber
\end{eqnarray}
Since the condition $\Vert\widetilde{\mf{\Delta}}_{n}\Vert_{2}<1$ is necessary to ensure the existence of $(\mf{I}_p+\widetilde{\mf{\Delta}}^{2}_{n})^{-1}$,
it follows that
\begin{eqnarray}
\frac{\partial^{2}\widetilde{\phi}_{D}(\widetilde{\mf{L}},\widetilde{\mf{S}})}{\partial^{2}{\widetilde{\mf{L}}}}&=&(\mf{I}_p+\widetilde{\mf{\Delta}}^{2}_{n})^{-1}-
(\mf{I}_p+\widetilde{\mf{\Delta}}^{2}_{n})^{-2}(2\widetilde{\mf{\Delta}}_{n})\nonumber\\
&=&(\mf{I}_p+\widetilde{\mf{\Delta}}^{2}_{n})^{-2}(\mf{I}_p+\widetilde{\mf{\Delta}}^{2}_{n}-2\widetilde{\mf{\Delta}}_{n})
=(\mf{I}_p+\widetilde{\mf{\Delta}}^{2}_{n})^{-2}(\mf{I}_p-\widetilde{\mf{\Delta}}_{n})^2.\nonumber
\end{eqnarray}
The last expression allows the expected value $\mathrm{E}$ to preserve the sign, because $\mf{I}_p-\widetilde{\mf{\Delta}}_{n}$ is positive definite since $\Vert\widetilde{\mf{\Delta}}_{n}\Vert_{2}<1$.
Imposing the condition $\Vert\widetilde{\mf{\Delta}}_{n}\Vert \leq \delta_{\widetilde{\Delta}}$,
%with $\delta_{\widetilde{\Delta}} \leq 1/3$,
and observing that $(\mf{I}_p-\widetilde{\mf{\Delta}}_{n})^2 \geq (1-\delta_{\widetilde{\Delta}})^2$,
the estimate from below of Section \ref{fish_main} follows.

\subsection{Asymptotic behaviour}

In this section, we describe the asymptotic behaviour of $\phi_D(\mf{L}^{*},\mf{S}^{*})$
%\frac{1}{2}\ln \det (\mf{I}_p+\mf{\Delta}_{n}\mf{\Delta}_{n}')$ a
as $\mf{\Sigma}_n$ lies around $\mf{\Sigma}^{*}=\mf{L}^{*}+\mf{S}^{*}$,
in terms of local convexity, first and second derivative.
%i.e., with $\mf{\Delta}_{n}=\mf{\Sigma}^{*}-\mf{\Sigma}_n$.
To this purpose, it is essential to identify the conditions ensuring the convergence of the sample covariance matrix $\mf{\Sigma}_n$
to the true $\mf{\Sigma}^{*}$ for all $p \in \N$ as $n \to \infty$. Such convergence is established in Lemma \ref{Lemma_cons},
%reported and proved in the Appendix,
which is the key result behind all the following lemmas.
% the convergence of $\mf{\Delta}_{n}$ to $\mf{0}$, that is,

%\far{Discuti per bene le conditioni di Lemma \ref{Lemma_cons} e il suo senso!}
Lemma \ref{Lemma_cons} essentially states that, with probability $1-O(1/n^2)$,
\begin{equation}\frac{1}{p^{\alpha_{1}}}\Vert \mf{\Delta}^{*}_{n} \Vert_{2}=o_P(1).\label{eq_Lemma_cons}\end{equation}
The proof of Lemma \ref{Lemma_cons} is based on controlling the moments of factor and residual distributions by Assumption \ref{tails},
latent eigenvalues and factor loadings by Assumption \ref{eigenvalues}, the residual sparsity pattern by Assumption \ref{sparsity}.
A lemma from \cite{bickel2008covariance} is the key to allow for a limited number of moments as prescribed by Assumption \ref{tails}.
%at the price of imposing $f_{\delta_\epsilon}(p,n) \to 0$ as $n \to \infty$.
The validity of these conditions is a requirement for all the following lemmas, whose proof is reported in the Supplement.

Since $\widetilde{\phi}_{D}(\mf {L},\mf{S})$ is a stochastic function of $\mf{L}$ and $\mf{S}$,
as $\mf{\Sigma}_n$ is a random matrix, we need to ensure that the condition of Proposition \ref{conv_delta} is respected.
For this purpose, it is necessary to control the probability %$\mathcal{P}(\Vert \mf{\Delta}_{n}\mf{\Delta}_{n}' \Vert \geq \frac{1}{3\delta p})$. Since $\mf{\Delta}_{n}$ is symmetric, the same probability can be expressed as
%$\mathcal{P}(\Vert \mf{\Delta}_{n}^2 \Vert \geq \frac{1}{3\delta p} )$, which is equivalent to prove
$\mathcal{P}(p^{-\alpha_1}\Vert \mf{\Delta}^{*}_{n} \Vert \geq (3p)^{-1})$, %for $\delta_\phi=1$
where $\mf{\Delta}^{*}_{n}=\mf{\Sigma}^{*}-\mf{\Sigma}_n$.

\begin{lemma}\label{random_conv}
Under Assumptions \ref{tails}-\ref{sparsity},
%if $\Vert\mf{\Delta}^{*}_{n}\Vert_{2}\leq \Vert\mf{L}^{*}\Vert_{2}/(3p)$,
the function $\widetilde{\phi}_{D}(\mf {L},\mf{S})$
is convex for $\mf{\Delta}^{*}_n=\mf{\Sigma}^{*}-\mf{\Sigma}_n$
with probability tending to $1$ for all $p \in \N$ as $n \to \infty$.
%if and only if $\delta_\phi\succeq\frac{\sqrt{n}}{p^{\alpha_1+1}}$ as $n \to \infty$.
\end{lemma}

\paragraph{{Proof}}%\label{random_conv_proof}

Lemma \ref{Lemma_cons} shows that the claim ${p^{-\alpha_1}} \Vert \mf{\Delta}^{*}_n \Vert_{2} \preceq C n^{-1/2}$ holds for some $C>0$
under Assumptions \ref{tails}-\ref{sparsity} for all $p \in \N$ with probability $1-O(1/n^2)$ as $n \to \infty$.
%if $f_{\delta_\epsilon}(p,n) \to 0$  as $p,n \to \infty.$
Then, setting $\delta_\phi=1$, the convexity condition of Proposition \ref{conv_delta}
is always satisfied under Assumption \ref{eigenvalues} if it holds the condition $\Vert\mf{\Delta}^{*}_{n}\Vert_{2} \leq 1/3$.
Since such conditions is automatically satisfied for all $p \in \N$ with probability $1-O(1/n^2)$ as $n \to \infty$
by Lemma \ref{Lemma_cons}, the thesis follows. \qed

%which is dominated by the condition $\Vert \mf{\Delta}^{*}_{n}\Vert_{2} \leq 1/3$,
%is asymptotically satisfied under the same conditions.

%For any $ \delta_\phi >0 $ the function $\ln\det\bigl(\delta_\phi^{-2}\mf{I}_p + \mf{\Delta}_{n} \mf{\Delta}_{n}'\bigr)$
%is convex on the closed  ball $ \mathcal{C_{\delta_\phi}}= \{ \mf{\Delta}_{n} | \mf{\Delta}_{n} \mathrm{~~is~ a~ real~~}
%p\times p \mathrm{~matrix~}, \Vert \mf{\Delta}_{n}\Vert_{2} \leq (3\delta_\phi p)^{-1} \}$.%\Vert\mf{L}^{*}\Vert_{2}^{-1}\Vert

%\far{Can this consideration be extended to a general function? Stein's loss?}

Lemma \ref{random_conv} ensures for all $p \in \N$ the convexity of $\widetilde{\phi}_{D}(\mf{L},\mf{S})$ with probability tending to $1$ as $n \to \infty$, and is a necessary condition to derive the asymptotic error rates of the solution pair of problem \eqref{obj} (main paper)
for all $p \in \N$ as $n \to \infty$.

\begin{lemma}\label{random:first}
%Under Assumptions \ref{tails}-\ref{sparsity}, for all $p \in \N$, it holds with probability tending to $1$ as $n \to \infty$:
Under the conditions of Lemma \ref{random_conv}, for all $p \in \N$,
%\far{the conditions of Lemma \ref{Lemma_cons}},
it holds with probability tending to $1$ as $n \to \infty$:
\begin{equation}
\Big\Vert\frac{\partial \widetilde{\phi}_{D}(\mf{L}^{*},\mf{S}^{*})}{\partial \mf{L}^{*}}\Big\Vert_{2}=
\Big\Vert\frac{\partial \widetilde{\phi}_{D}(\mf{L}^{*},\mf{S}^{*})}{\partial \mf{S}^{*}}\Big\Vert_{2}
\stackrel{n\to\infty}\longrightarrow \mf{0}_{p \times p}.\label{grad_0}
%\mf{\Delta}_{n}^{-2} 2 \mf{\Delta}_{n} = \mf{\Delta}_{n}^{-1}.
\end{equation}
\end{lemma}
Lemma \ref{random:first} shows that the gradient of $\widetilde{\phi}_{D}(\mf{L},\mf{S})$ \emph{wrt} $\mf{L}^{*}$ and $\mf{S}^{*}$ converges to $0$ in probability if Lemma \ref{random_conv} holds.
Its proof relies on Proposition \ref{first_der}.

\paragraph{{Proof}}%\label{random_first_proof}

The proof of Lemma \ref{Lemma_cons} has highlighted that
\begin{equation}\frac{1}{p^{\alpha_{1}}}\Vert \mf{\Delta}^{*}_{n} \Vert_{2} \stackrel{n\to\infty}\longrightarrow 0,\label{Lemma_input}\end{equation}
with probability $1-O(1/n^2)$ for all $p \in \N$ as $n \to \infty$ under all its conditions, where $\mf{\Delta}^{*}_{n}=\mf{\Sigma}^{*}-\mf{\Sigma}_n$.
We recall from Proposition \ref{first_der} that
\begin{eqnarray}
\widetilde{\phi}_{D}'(\mf{L}^{*},\mf{S}^{*})^{(ld)}&=&
\frac{\partial \widetilde{\phi}_{D}(\mf{L}^{*},\mf{S}^{*})}{\partial \mf{L}^{*}}=\frac{\partial \widetilde{\phi}_{D}(\mf{L}^{*},\mf{S}^{*})}
{\partial \mf{S}^{*}}\nonumber\\
&=&(\mf{I}_p+p^{-2\alpha_1}\mf{\Delta}^{*}_{n}\mf{\Delta}^{*'}_{n})^{-1}p^{-\alpha_1}\mf{\Delta}^{*}_{n}.\nonumber
\end{eqnarray}
We can observe that $\mf{\Delta}^{*2}_{n}=\mf{\Delta}^{*}_{n}\mf{\Delta}^{*'}_{n}$ (due to symmetry) and that
all the eigenvalues of $\mf{\Delta}^{*2}_{n}$ are strictly smaller than $1$ by the convexity condition of Lemma \ref{random_conv},
$\Vert\mf{\Delta}^{*}_{n}\Vert_{2} \leq 1/3$, which is implicitly satisfied for all $p \in \N$ as $n \to \infty$.
%i.e., $\Vert\mf{L}^{*}\Vert_{2}^{-1}\Vert\mf{\Delta}^{*}_{n}\Vert_{2} \leq (3p)^{-1}$,
%\far{since $\Vert\mf{L}^{*}\Vert_{2}<3p$}.
Then, it holds
\begin{eqnarray}
\phi'(\mf{L}^{*},\mf{S}^{*})^{(ld)}&=&
(\mf{I}_p+p^{-2\alpha_1}\mf{\Delta}^{*2}_{n})^{-1}p^{-\alpha_1}\mf{\Delta}^{*}_{n}\nonumber\\
&=&\sum_{j=0}^{\infty}(p^{-2\alpha_1}\mf{\Delta}^{*2}_{n})^{j}p^{-\alpha_1}\mf{\Delta}^{*}_{n}\nonumber\\
&=&\sum_{j=0}^{\infty}(-1)^{j}(p^{-\alpha_1}\mf{\Delta}^{*}_{n})^{2j+1}.\label{phi'_series}
\end{eqnarray}
Equation \eqref{phi'_series} also implies that
\begin{eqnarray}
\Vert(\mf{I}_p+p^{-2\alpha_1}\mf{\Delta}^{*2}_{n})^{-1}\Vert&=&
\Bigg\Vert\sum_{j=0}^{\infty}(-1)^{j}(p^{-\alpha_1}\mf{\Delta}^{*}_{n})^{2j}\Bigg\Vert\nonumber\\
&=&\sum_{j=0}^{\infty}(-1)^{j}p^{-2\alpha_1}\Vert(\mf{\Delta}^{*}_{n})^2\Vert^{j}\nonumber\\
&=&\frac{1}{1+p^{-2\alpha_1}\Vert(\mf{\Delta}^{*}_{n})^2\Vert} \leq 1\label{norm1}
\end{eqnarray}
for any norm. For spectral norm, we get from \eqref{norm1} that $$\Vert(\mf{I}_p+p^{-2\alpha_1}\mf{\Delta}^{*}_{n}\mf{\Delta}^{*'}_n)^{-1}\Vert_{2}=\frac{1}{1+p^{-2\alpha_1}\lambda_{1}(\mf{\Delta}^{*}_{n})^2},$$
implying that the minimum for $\Vert(\mf{I}_p+p^{-2\alpha_1}\mf{\Delta}^{*}_{n}\mf{\Delta}^{*'}_n)^{-1}\Vert_{2}$, under the convexity condition of Lemma \ref{random_conv}, is
$$\frac{1}{(1+1/9)}=\frac{9}{10},$$
while the maximum is $1$, and is reached asymptotically for all $p \in \N$ as $n\to\infty$ by \eqref{Lemma_input}.
Finally, by triangular inequality, we observe that
\begin{eqnarray}\label{D}
\Vert(\mf{I}_p+p^{-2\alpha_1}\mf{\Delta}^{*}_{n}\mf{\Delta}^{*'}_{n})^{-1}p^{-\alpha_1}\mf{\Delta}^{*}_{n}\Vert_{2}\leq\nonumber\\
\Vert(\mf{I}_p+p^{-2\alpha_1}\mf{\Delta}^{*}_{n}\mf{\Delta}^{*'}_{n})^{-1}\Vert_{2}\Vert p^{-\alpha_1}\mf{\Delta}^{*}_{n}\Vert_{2}.\nonumber
\end{eqnarray}
Putting together \eqref{Lemma_input} and \eqref{norm1}, the thesis follows. \qed

\begin{lemma}\label{random:second}
%More, if $ \mf{\Sigma} = \mf{\Sigma}_n$, we get
Under the conditions of Lemma \ref{random_conv}, for all $p \in \N$,
it holds with probability tending to $1$ as $n \to \infty$:
%\begin{equation}
%  \label{kron0}
%  \frac{1}{p^{\alpha_1}}\bigl(\frac{1}{2}\mathrm{Hess}\ln\det
%                     \varphi(\mf{\Sigma}^{*})\bigr)_{i j h k} \stackrel{n\to\infty}\longrightarrow \delta_{jk}\otimes
%                   \delta_{i h} = \bigl(\mf{I}_p \otimes \mf{I}_p\bigr) _{i j h k},
%\end{equation}
%that is,
%
\begin{equation}
\mathrm{Hess}\widetilde{\phi}_{D}(\mf{L}^{*},\mf{S}^{*})\stackrel{n\to\infty}\longrightarrow \mf{I}_p \otimes \mf{I}_p. \label{fishone}
\end{equation}
\end{lemma}

\paragraph{{Proof}}%\label{random_second_proof}

%We consider equation \eqref{fishone} in the main paper. %\eqref{kron0} and
We consider the thesis $\mathrm{Hess}\widetilde{\phi}_{D}(\mf{L}^{*},\mf{S}^{*})\stackrel{n\to\infty}\longrightarrow \mf{I}_p \otimes \mf{I}_p$.
%We start by equation \eqref{lips_top}) where we set, $\epsilon=1$,
%and $\mf{H}=\mf{\Delta}_{n}=\mf{\Sigma}_n-\mf{\Sigma}^{*}$. Furthermore we note that
We observe that, for all $p \in \N$, $p^{-\alpha_1}\Vert \mf{\Delta}^{*}_{n} \Vert_{2} \stackrel{n\to\infty}\longrightarrow 0$
with probability $1-O(1/n^2)$ as $n \to \infty$ by \eqref{Lemma_input}, under the conditions of Lemma \ref{Lemma_cons}.
%the conditions of Lemma \ref{Lemma_cons},
%such that $\frac{1}{2p^{\alpha_{1}}}\Vert F(\mf{\Delta}_n)-F(\mf{0})\Vert \stackrel{n\to\infty}\longrightarrow \mf{0}_{p \times p}$
%under those assumptions,
Consequently, under such conditions, Proposition \ref{second_der} holds for all $p \in \N$ with probability $1-O(1/n^2)$ as $n \to \infty$.
The thesis then follows. \qed
%\eqref{kron0_bis} and
%which means that
%with $F(\mf{0})=\mf{I}_p \otimes \mf{I}_p$ by \ref{fishone}.
%$$\frac{1}{2p^{\alpha_{1}}}\Vert\mathrm{Hess}\ln\det\varphi(\mf{\Sigma}_n)-\mathrm{Hess}\ln\det\varphi(\mf{\Sigma}^{*})\Vert
%\stackrel{n\to\infty}\longrightarrow \mf{0}_{p \times p}.$$ Therefore the thesis follows.

%\section{Proofs to integrate in Section \ref{Cons_both}}

%\paragraph{Proof of Proposition \ref{12}}

%%

Lemma \ref{random:second} establishes that the Hessian of $\widetilde{\phi}_{D}(\mf{L}^{*},\mf{S}^{*})$
converges to $\mf{I}_p \otimes \mf{I}_p$ in probability if Lemma \ref{random_conv} holds.
Its proof relies on Proposition \ref{second_der}.

%\far{Studiare il comportamento FINITO di $\mathrm{Hess}\phi_D(\mf{L}^{*},\mf{S}^{*})$ per $p,n$ fissati, in particolare provando che il suo autovalore più piccolo è strettamente positivo
%con alta probabilità.}

%%
%\far{START}
\begin{lemma}\label{lemma_strong}
Under the conditions of Lemma \ref{random_conv}, for all $p \in \N$, $\widetilde{\phi}_{D}(\mf{L}^{*},\mf{S}^{*})$ is strongly convex with probability tending to $1$ as $n \to \infty$.
\end{lemma}

\paragraph{{Proof}}%\label{random_second_proof}

Lemma \ref{random_conv} ensures that $\widetilde{\phi}_{D}(\mf{L}^{*},\mf{S}^{*})$ is locally convex around $\mf{\Sigma}^{*}$
with probability tending to $1$ as $n \to \infty$. Lemma \ref{random:second} ensures that the Hessian of $\widetilde{\phi}_{D}(\mf{L}^{*},\mf{S}^{*})$ is positive definite with probability tending to $1$ as $n \to \infty$, because its limit is $\mf{I}_p \otimes \mf{I}_p$, which is positive definite as all its eigenvalues are equal to 1.
Putting together these two facts, the thesis then follows.

Lemma \ref{lemma_strong} ensures the strong convexity of $\widetilde{\phi}_{D}(\mf{L}^{*},\mf{S}^{*})$ for all $p\in\N$ as $n\to\infty$,
which is necessary to establish a direct comparison with \cite{10.1214/12-STS400}. \qed

%Generally speaking, the six quantities $\alpha_{\Omega}$, $\alpha_{\mathcal{T}}$, $\beta_{\Omega}$, $\beta_{\mathcal{T}}$,
%$\delta_{\Omega}$, $\delta_{\mathcal{T}}$ all lie within the interval $[0,1]$. However,
%Relevantly, the following result holds with probability $1-O(1/n^2)$ as $n \to \infty$.
\begin{lemma}\label{Lemma_ort}
Under the conditions of Lemma \ref{random_conv},
the following statements hold with probability tending to 1 for all $p \in \N$ as $n \to \infty$:
%$1-O(1/n^2)$
$\alpha_{\Omega} \to 1$, $\alpha_{\mathcal{T}} \to 1$, %$\beta_{\Omega} \to 1$, $\beta_{\mathcal{T}} \to 1$,
$\delta_{\Omega} \to 0$, $\delta_{\mathcal{T}} \to 0$, $\nu \to 1/2$.
%$\mf{W}=\mf{\Delta}_n$, and $\mf{\Delta}_n \to 0$ under the conditions of Lemma \ref{Lemma_cons} for all $p \in \N$, as $n \to \infty$.
%for all $p \in \N$ with probability $1-O(1/n^2)$ as $n \to \infty$.
\end{lemma}
%\begin{proof}
%The statement follows by the definitions of $\alpha_{\Omega}$, $\alpha_{\mathcal{T}}$, $\beta_{\Omega}$, $\beta_{\mathcal{T}}$,
%$\delta_{\Omega}$, $\delta_{\mathcal{T}}$, and because, with probability $1-O(1/n^2)$, for all $p \in \N$ as $n \to \infty$ it holds
%$\frac{1}{p^{\alpha_{1}}}\Vert \mf{W} \Vert_{2} \to 0$ and $\Vert \mf{W} \Vert_{\infty} \to 0$
%under the conditions of Lemma \ref{Lemma_cons} (and, as a consequence, of Lemma \ref{Lemma_bmax}).
%\end{proof}
%Also note that  under the conditions of Lemma \ref{Lemma_ort}.
%\far{END}

\paragraph{{Proof}}%\label{random_second_proof}
%which translates to
%$$\Vert \mathcal{I}^{*}(\mf{L},\mf{S}) \Vert \leq \Vert \mf{I}_p \otimes \mf{I}_p \Vert - \Vert \mf{W}\otimes\mf{W} \Vert=
%1-\Vert \mf{W}\otimes\mf{W} \Vert$$ %plus the squared bias $\mathrm{E}(F(\mf{\Delta}_n) \otimes F(\mf{\Delta}_n))$.
%as $n\to\infty$.
Under the conditions of Lemma \ref{random_conv}, Lemmas \ref{random:first} and \ref{random:second} ensure that
the Hessian of $\widetilde{\mathcal{L}}^{(ld)}(\mf{\Sigma},\mf{\Sigma}_{n})$ approaches the Hessian of
$\widetilde{\mathcal{L}}^{(F)}(\mf{\Sigma},\mf{\Sigma}_{n})=0.5p^{-\alpha_1}
\Vert \mf{\Sigma}_{n}-(\mf{L}+\mf{S})\Vert_{F}^2$ for all $p \in \N$ as $n \to \infty$.
This proves the lemma, because for $\widetilde{\mathcal{L}}^{(F)}(\mf{\Sigma},\mf{\Sigma}_{n})$ it holds
$\alpha_{\Omega} = 1$, $\alpha_{\mathcal{T}} = 1$, $\beta_{\Omega} = 1$, $\beta_{\mathcal{T}} = 1$,
$\delta_{\Omega} = 0$, $\delta_{\mathcal{T}} = 0$, $\nu = 1/2$ for any $p$ and $n$ (cf. \cite{luo2011high}). \qed

\begin{remark}
Under the conditions of Lemma \ref{Lemma_ort},
the identifiability condition of Proposition \ref{11} becomes $\xi(\mathcal{T}(\mf{L}^{*}))\mu(\Omega(\mf{S}^{*}))\leq\frac{1-\kappa_\mathcal{T}}{8(1+\kappa_\mathcal{T})}$
for all $p \in \N$ as $n \to \infty$, which also holds for $\widetilde{\mathcal{L}}^{(F)}(\mf{\Sigma},\mf{\Sigma}_{n})$ for all $p,n \in \N$.
Moreover, %the approximation error of $\mf{\Sigma}_n$ becomes negligible, so that we can assume to be in the noiseless setting, under which $\kappa_\mathcal{T}=0$.
%under the conditions of Lemma \ref{random_conv},
the identifiability condition progressively becomes $\xi(\mathcal{T}(\mf{L}^{*}))\mu(\Omega(\mf{S}^{*}))\leq 1/8$,
because $\kappa_\mathcal{T} \to 0$ as $n \to \infty$.
When $\xi(\mathcal{T}(\mf{L}^{*}))$ attains its minimum $\sqrt{r/p}$ and $\mf{S}^{*}$ is diagonal, the identifiability condition simply becomes $p \geq 64 r$.
\end{remark}

\section{Proofs of theorems and corollaries}\label{proofs}

%Let us define $$f_{\delta_f}(r,n)=
%\bigl\{
%\begin{array}{rl}
%\frac{r^{2/(1+\delta_f)}}{n^{1/2}}, & \delta_f \in \R^{+}, \nonumber\\
%\frac{\ln(r)}{n}, & \delta_f=+\infty, \nonumber
%\end{array}
%\bigr.
%$$
%and
%$$f_{\delta_\epsilon}(p,n)=
%\bigl\{
%\begin{array}{rl}
%\frac{p^{2/(1+\delta_\epsilon)}}{n^{1/2}}, & \delta_\epsilon \in \R^{+}, \nonumber\\
%\frac{\ln(p)}{n}, & \delta_\epsilon=+\infty, \nonumber
%\end{array}
%\bigr.
%$$
%and
%$$f_{\delta_x}(p,n)=
%\bigl\{
%\begin{array}{rl}
%\frac{p^{2/(1+\delta_x)}}{n^{1/2}}, & \delta_x \in \R^{+}, \nonumber\\
%\frac{\ln(p)}{n}, & \delta_x=+\infty, \nonumber
%\end{array}
%\bigr.
%$$
%with $\delta_x=\max{(\delta_f,\delta_\epsilon)}$.

%Consequently, setting $M_x=\max(M_f,M_\epsilon)$ and $\delta_x=\max{(\delta_f,\delta_\epsilon)}$, the following corollary holds.
%\begin{corollary}\label{coroll_x}
%Under the assumptions of Lemma \ref{lemma_rate_both}, it holds
%\begin{eqnarray}
%\max_{i,j} \vert \mf{\Sigma}_{n,ij}-\mf{\Sigma}_{ij}^{*} \vert=O_P\bigl(\frac{p^{2/(1+\delta_x)}}{n^{1/2}}\bigr);\label{rate_moment_x}\\
%\Pr (\max_{i,j} \vert \mf{\Sigma}_{\epsilon,ij}-\mf{S}^{*} \vert \leq t)\leq p^2 M_x C(\delta_x)\frac{n^{-(1+\delta_x)/2}}{t^{1+\delta_x}}.\label{prob_rate_x}
%\end{eqnarray}
%\end{corollary}

\subsection*{{Proof of Proposition \ref{11}}}%\label{random_conv_proof}

The condition $\Vert\widehat{\mf{\Delta}}_{n}\Vert \leq 1/3$ implies by
Proposition \ref{conv_delta} that $\widetilde{\phi}_{D}(\mf{L},\mf{S})$ is locally convex in $\widehat{\mf{\Sigma}}_{n}$,
and by inequality \eqref{minimum_Fisher} (main paper) that $\widetilde{\phi}_{D}(\mf{L},\mf{S})$ is strongly convex.
%with $\alpha_{\mathcal{Y}}\geq 0.81$.
Then, we observe that
\begin{eqnarray}
\xi(\mathcal{T}') &\leq& \frac{\xi(\mathcal{T})+\varrho(\mathcal{T},\mathcal{T}')}{1-\varrho(\mathcal{T},\mathcal{T}')}
=\frac{\xi(\mathcal{T})+\kappa_\mathcal{T}\xi(\mathcal{T})}{1-\kappa_\mathcal{T}\xi(\mathcal{T})}\nonumber\\
&=&\frac{1+\kappa_\mathcal{T}}{1-\kappa_\mathcal{T}}\xi(\mathcal{T}).\label{rho_ineq}
\end{eqnarray}

We note that the condition
\begin{equation}
\gamma \in \left[\frac{2\xi(\mathcal{T'})\beta_{\mathcal{Y}}(1-\nu)}{\nu\alpha_{\mathcal{Y}}},
\frac{\nu\alpha_{\mathcal{Y}}}{4\mu(\Omega)\beta_{\mathcal{Y}}(1-\nu)}\right]\label{gamma_range}
\end{equation} becomes by \eqref{rho_ineq}
$$\gamma \in \left[\frac{\frac{2(1+\kappa_\mathcal{T})}{(1-\kappa_\mathcal{T})}\xi(\mathcal{T}(\mf{L}^{*}))\beta_{\mathcal{Y}}(1-\nu)}{\nu\alpha_{\mathcal{Y}}},
\frac{\nu\alpha_{\mathcal{Y}}}{4\mu(\Omega(\mf{S}^{*}))\beta_{\mathcal{Y}}(1-\nu)}\right].$$
For this to be meaningful, it must hold
$$\frac{\frac{2(1+\kappa_\mathcal{T})}{(1-\kappa_\mathcal{T})}\xi(\mathcal{T}(\mf{L}^{*}))\beta_{\mathcal{Y}}(1-\nu)}{\nu\alpha_{\mathcal{Y}}} \leq \frac{\nu\alpha_{\mathcal{Y}}}{4\mu(\Omega(\mf{S}^{*}))\beta_{\mathcal{Y}}(1-\nu)},$$ which is equivalent to
$$\xi(\mathcal{T}(\mf{L}^{*}))\mu(\Omega(\mf{S}^{*}))\leq\frac{1-\kappa_\mathcal{T}}{8(1+\kappa_\mathcal{T})}\Bigl(\frac{\nu\alpha_{\mathcal{Y}}}{\beta_{\mathcal{Y}}(1-\nu)}\Bigr)^2.$$

From \eqref{gamma_range}, it follows that
$$\max\left(\frac{\xi(\mathcal{T}(\mf{L}^{*}))}{\gamma},2\mu(\Omega(\mf{S}^{*}))\gamma\right)\leq\frac{\nu\alpha_{\mathcal{Y}}}{2\beta_{\mathcal{Y}}(1-\nu)}.$$
Therefore, following the same passages of the proof of Proposition 3.3 in \cite{chandrasekaran2012latent}, we get by Assumption \ref{ass_alg} that
\begin{eqnarray}
\min_{(\mf{L},\mf{S})\in \mathcal{Y},\Vert \mf{L} \Vert_{2}=1,\Vert \mf{S} \Vert_{\infty}=\gamma} g_{\gamma}(\mathbb{P}_{\mathcal{Y}}\mathcal{A}^{\dag}\mathcal{I}^{*}\mathcal{A}
\mathbb{P}_{\mathcal{Y}}(\mf{S},\mf{L}))\nonumber\\
\geq \alpha_{\mathcal{Y}}-\beta_{\mathcal{Y}}\frac{\nu\alpha_{\mathcal{Y}}}{2\beta_{\mathcal{Y}}(1-\nu)}=
\alpha_{\mathcal{Y}}-\frac{\nu\alpha_{\mathcal{Y}}}{2(1-\nu)}.\nonumber
\end{eqnarray}

Recalling the lower bound $\nu\geq\frac{2\alpha_{\mathcal{Y}}-1}{2\alpha_{\mathcal{Y}}}$,
we can derive the consequent upper bound for $\frac{\nu}{1-\nu}$ as $\frac{2\alpha_{\mathcal{Y}}-1}{2\alpha_{\mathcal{Y}}}{2\alpha_{\mathcal{Y}}}=2\alpha_{\mathcal{Y}}-1$,
because $$\frac{1}{1-\nu} \geq \frac{2\alpha_{\mathcal{Y}}-2\alpha_{\mathcal{Y}}+1}{2\alpha_{\mathcal{Y}}}=\frac{1}{2\alpha_{\mathcal{Y}}}.$$
As a consequence, the lower bound $\alpha_{\mathcal{Y}}-\frac{\nu\alpha_{\mathcal{Y}}}{2(1-\nu)}$ by \eqref{gamma_range} becomes

\begin{eqnarray}
\alpha_{\mathcal{Y}}-\frac{\nu\alpha_{\mathcal{Y}}}{2(1-\nu)}&=&\frac{(2\alpha_{\mathcal{Y}}-(2\alpha_{\mathcal{Y}}-1)\alpha_{\mathcal{Y}})}{2}\nonumber\\
&=&\frac{(3\alpha_{\mathcal{Y}}-2\alpha_{\mathcal{Y}}^2)}{2}=\alpha_{\mathcal{Y}}\Bigl(\frac{3}{2}-\alpha_{\mathcal{Y}}\Bigr).\nonumber
\end{eqnarray}

Then, part (i) follows:
\begin{eqnarray}
\min_{(\mf{L},\mf{S})\in \mathcal{Y},\Vert \mf{L} \Vert_{2}=1,\Vert \mf{S} \Vert_{\infty}=\gamma} g_{\gamma}(\mathbb{P}_{\mathcal{Y}}\mathcal{A}^{\dag}\mathcal{I}^{*}\mathcal{A}
\mathbb{P}_{\mathcal{Y}}(\mf{S},\mf{L}))\geq \alpha_{\mathcal{Y}}\Bigl(\frac{3}{2}-\alpha_{\mathcal{Y}}\Bigr).\nonumber
\end{eqnarray}

Following the proof of Proposition 3.3 in \cite{chandrasekaran2012latent}, we get
$$\frac{g_{\gamma}(\mathbb{P}_{\mathcal{Y}^\perp}\mathcal{A}^{\dag}\mathcal{I}^{*}\mathcal{A}\mathbb{P}_{\mathcal{Y}}(\mf{S},\mf{L}))}
{g_{\gamma}(\mathbb{P}_{\mathcal{Y}}\mathcal{A}^{\dag}\mathcal{I}^{*}\mathcal{A} \mathbb{P}_{\mathcal{Y}}(\mf{S},\mf{L}))}\leq
\frac{\delta_{\mathcal{Y}}+0.5\nu\alpha_{\mathcal{Y}}}{\alpha_{\mathcal{Y}}-0.5\nu\alpha_{\mathcal{Y}}}.$$
%\leq 1-\nu.$

Then, once recalled that $\delta_{\mathcal{Y}}=1-\alpha_{\mathcal{Y}}$, we can observe that the thesis of part (ii)
$$\frac{(1-\alpha_{\mathcal{Y}})+0.5\nu\alpha_{\mathcal{Y}}}{\alpha_{\mathcal{Y}}-0.5\nu\alpha_{\mathcal{Y}}} \leq 1-\nu$$
is dominated by the thesis
\begin{equation}\frac{(1-\alpha_{\mathcal{Y}})+\nu/(2\alpha_{\mathcal{Y}})}{\alpha_{\mathcal{Y}}-\nu/(2\alpha_{\mathcal{Y}})}\leq\frac{1}{2\alpha_{\mathcal{Y}}}.\label{eq_ort}\end{equation}
We further develop \eqref{eq_ort} as follows:
\begin{eqnarray}
\frac{2\alpha_{\mathcal{Y}}(1-\alpha_{\mathcal{Y}})+\nu}{2\alpha_{\mathcal{Y}}^2-\nu} &\leq& 1\nonumber\\
{2\alpha_{\mathcal{Y}}(1-\alpha_{\mathcal{Y}})+\nu} &\leq& {2\alpha_{\mathcal{Y}}^2-\nu}\nonumber\\
4\alpha_{\mathcal{Y}}^2-2\alpha_{\mathcal{Y}}-2\nu &\geq& 0.\nonumber
\end{eqnarray}
In the worst admissible case $\nu=\frac{2\alpha_{\mathcal{Y}}-1}{2\alpha_{\mathcal{Y}}}$,
we obtain the inequality $$4\alpha_{\mathcal{Y}}^2-2\alpha_{\mathcal{Y}}-2\frac{2\alpha_{\mathcal{Y}}-1}{2\alpha_{\mathcal{Y}}} \geq 0.$$
This leads to the inequality $4\alpha_{\mathcal{Y}}^3-2\alpha_{\mathcal{Y}}^2-2\alpha_{\mathcal{Y}}+1 \geq 0$,
which is verified within the range $\alpha_{\mathcal{Y}}\in (1/2,1]$ (prescribed by Assumption \ref{ass_alg})
as long as $\alpha_{\mathcal{Y}} \geq 0.70711$ (corresponding to $\nu \geq 0.29290$).
Since such condition is explicitly assumed, part (ii) then follows.

%$$\frac{(2\alpha_{\mathcal{Y}}-1)(\alpha_{\mathcal{Y}}-1)+1}{(2\alpha_{\mathcal{Y}}-1)(1-\alpha_{\mathcal{Y}})+1}\leq\frac{1}{2\alpha_{\mathcal{Y}}}$$
%$$2\alpha_{\mathcal{Y}}{(2\alpha_{\mathcal{Y}}-1)(\alpha_{\mathcal{Y}}-1)+1}\leq {(2\alpha_{\mathcal{Y}}-1)(1-\alpha_{\mathcal{Y}})+1}.$$
%This leads to
%$$2\alpha_{\mathcal{Y}}(\alpha_{\mathcal{Y}}-1)\leq(1-\alpha_{\mathcal{Y}})$$
%$$2\alpha_{\mathcal{Y}}^2-\alpha_{\mathcal{Y}}-1 \leq 0,$$
%which is always verified as long as $\alpha_{\mathcal{Y}} \in [-1/2,1]$, well beyond the range $[-1/2,1]$ of $\alpha_{\mathcal{Y}}$.

%\far{Report and explain (5.1) page 1955 in \cite{chandrasekaran2012latent}}

\subsection*{{Proof of Theorem \ref{thm_comp}}}

First of all, we observe that Lemmas \ref{random_conv}-\ref{Lemma_ort} apply,
such that $\widetilde{\phi}(\mf{L},\mf{S})^{(ld)}$ is strongly convex, and for all $p\in\N$ as $n \to \infty$,
$\alpha_{\mathcal{Y}} \to 1$, $\nu \to 1/2$, and $\kappa_\mathcal{T} \to 0$.

%%%%%
Then, we consider the solution of the following algebraic problem:
%for $\theta_h=\frac{h\pi}{M_T}$, $\vert h \vert \leq M_T$
\begin{equation}\label{probtang}
\bigl(\widehat{\mf{L}}_{\mathcal{T}'},\widehat{\mf{S}}_{\Omega}\bigr)=\arg\!\!\!\!\!
\min_{\underline{\mf{L}} \in \mathcal{T}',\underline{\mf{S}} \in \Omega}
\widetilde{\phi}_{D}(\underline{\mf{L}},\underline{\mf{S}})
+\psi_{0} \Vert \underline{\mf{L}} \Vert_{*}+\rho_{0} \Vert \underline{\mf{S}}\Vert_{1},
\end{equation}
where $\widetilde{\phi}_{D}(\underline{\mf{L}},\underline{\mf{S}})=0.5\ln \det (\mf{I}_p+p^{-2\alpha_1}\underline{\mf{\Delta}}_{n}\underline{\mf{\Delta}}_{n}')$,
with $\underline{\mf{\Delta}}_{n}=\mf{\Sigma}_n-(\underline{\mf{L}}+\underline{\mf{S}})$.
Problem \eqref{probtang} is an equivalent version of minimization \eqref{obj} (main paper)
rescaled to cope with Lemma \ref{Lemma_cons}.
%by $\Vert\mf{L}^{*}\Vert_{2}$.
In Proposition \ref{11} (main paper), we have bounded the degree of transversality between the low rank and the sparse variety,
to control the impact of $\mathcal{I}^{*}(\mf{L},\mf{S})$ on matrix variety identification.
%with high probability.
Here, we bound the error norm \eqref{ggamma} (main paper) for the solution pair $(\widehat{\mf{L}},\widehat{\mf{S}})$
defined in \eqref{prob_ld_rescaled} (main paper). To reach that goal, following \cite{chandrasekaran2012latent},
we need before to bound the error norm of the solution pair \eqref{probtang}.
%by the following proposition (proved in the Supplement).

\begin{proposition}\label{12}
Let $\varrho(\mathcal{T}',\mathcal T)\leq \kappa_{\mathcal{T}}\xi(\mathcal{T})$ and define
$\mf{\Delta}^{*}_n=\mf{\Sigma}_n-\mf{\Sigma}^{*}$, $\mf{C}_{\mathcal{T}'}=\mathbb{P}_{\mathcal{T}'^{\perp}}(\mf{L}^{*})$,
$$r_{\gamma}(\psi_0)=\frac{4}{\alpha_{\mathcal{Y}}(3-2\alpha_{\mathcal{Y}})}[g_{\gamma}(\mathcal{A}^{\dag}\mf{\Delta}^{*}_n)
+g_{\gamma}(\mathcal{A}^{\dag}\mathcal{I}^{*}\mf{C}_{\mathcal{T}'})+\psi_{0}],$$ and
$$\widetilde{r}=\max\left\{r_{\gamma}(\psi_0),\frac{\Vert \mf{C}_{\mathcal{T}'}\Vert_{2}}{p^{\alpha_1}}\right\}.$$
Then, if $\Vert\mf{\Delta}^{\mathcal{Y}}_n\Vert_{2} \leq 1/3$,
under the conditions of Proposition \ref{11}
%\far{if $\Vert\mf{L}^{*}\Vert_{2}^{-1}\Vert\mf{\Delta}^{\mathcal{Y}}_{n}\Vert_{2} \leq (3p)^{-1}$ and $\Vert\mf{L}^{*}\Vert_{2}<3p$},
%and Lemma \ref{random_conv},
the solution of problem \eqref{probtang}
$(\widehat{\mf{L}}_{\mathcal{T}'},\widehat{\mf{S}}_{\Omega})$ satisfies
\begin{equation}
g_\gamma(\widehat{\mf{L}}_{\mathcal{T}'}-\mf{L}^{*},\widehat{\mf{S}}_{\Omega}-\mf{S}^{*})\leq 2\widetilde{r}.\nonumber
%\label{ggamma_ineq}
\end{equation}
%with high probability as $n \to \infty$.
\end{proposition}

\paragraph{Proof}

Let us set $\mathcal{Y}=\Omega\oplus\mathcal{T'}$ and define $\mf{\Delta}^{\mathcal{Y}}_n=\mf{\Sigma}_n-(\widehat{\mf{L}}_{\mathcal{T}'}+\widehat{\mf{S}}_{\Omega})$.
We observe that $\widetilde{\phi}_{D}(\mf{L},\mf{S})$
is convex within the convexity range $p^{-\alpha_1}\Vert \mf{\Delta}^{\mathcal{Y}}_n \Vert_{2} \leq (3p)^{-1}$ by Proposition \ref{conv_delta},
which is dominated by the condition $\Vert\mf{\Delta}^{\mathcal{Y}}_n\Vert_{2} \leq 1/3$.
%From the tangent space constraints $\mf{S}\in\Omega$ and $\mf{L}\in\mathcal{T}'$, we know that the optimum $\bigl(\widehat{\mf{L}}_{\mathcal{T}'},\widehat{\mf{S}}_{\Omega}\bigr)$ the objective $\widetilde{\phi}_{D}(\mf{L},\mf{S})$,
It follows that $\bigl(\widehat{\mf{L}}_{\mathcal{T}'},\widehat{\mf{S}}_{\Omega}\bigr)$ is also the unique solution of \eqref{probtang}.
Following \cite{clarke1990optimization}, we know that the optimum $\bigl(\widehat{\mf{L}}_{\mathcal{T}'},\widehat{\mf{S}}_{\Omega}\bigr)$
satisfies for two Lagrangian multipliers in the spaces $\mathcal{T'}^{\perp}$ and $\Omega^{\perp}$,
$\mf{Q}_{\mathcal{T'}^{\perp}} \in \mathcal{T'}^{\perp}$ and $\mf{Q}_{\Omega^{\perp}}\in\Omega^{\perp}$,
the following conditions: $$\widetilde{\phi}'_{D}(\mf{L},\mf{S})^{(ld)}+\mf{Q}_{\mathcal{T'}^{\perp}}\in -\psi_{0}\partial\Vert \widehat{\mf{L}}_{\mathcal{T}'}\Vert_{*},$$
$$\widetilde{\phi}'_{D}(\mf{L},\mf{S})^{(ld)}+\mf{Q}_{\Omega^{\perp}}\in -\gamma\psi_{0}\partial\Vert \widehat{\mf{S}}_{\Omega}\Vert_{1},$$
where $\partial\Vert \widehat{\mf{L}}_{\mathcal{T}'}\Vert_{*}$ and
$\partial\Vert \widehat{\mf{S}}_{\Omega}\Vert_{1}$ denote the sub-differentials of
$\Vert \widehat{\mf{L}}_{\mathcal{T}'}\Vert_{*}$ and $\Vert \widehat{\mf{S}}_{\Omega}\Vert_{1}$ (see \cite{watson1992characterization}).

%In the proof of Proposition 12 by \cite{luo2011high}, the error matrix
%$\widehat{\mf{S}}_{\Omega}+\widehat{\mf{L}}_{\mathcal{T}'}-\mf{\Sigma}_n$
%is decomposed with respect to the Cartesian sum $\mathcal{Y}=\mathcal{T'}+\Omega$.
Recalling equation \eqref{grad} in the main paper, we can write
$$\phi'(\widehat{\mf{L}}_{\mathcal{T}'},\widehat{\mf{S}}_{\Omega})^{(ld)}=
(\mf{I}_p+\mf{\Delta}^{\mathcal{Y}}_n\mf{\Delta}^{\mathcal{Y}}_n)^{-1}
\mf{\Delta}^{\mathcal{Y}}_n.$$
%where $\mf{\Delta}^{\mathcal{Y}}_n=\widehat{\mf{S}}_{\Omega}+\widehat{\mf{L}}_{\mathcal{T}'}-\mf{\Sigma}_n$.
%Lemma \ref{lemmatop} guarantees that
Following the same steps as in the proof of Lemma \ref{random:first},
since imposing the convexity condition $\Vert{\mf{\Delta}}^{\mathcal{Y}}_n\Vert_{2} \leq 1/3$
all the eigenvalues of $\mf{\Delta}^{\mathcal{Y}}_n$ are strictly smaller than $1$,
%\far{because $\Vert\mf{L}^{*}\Vert_{2}^{-1}\Vert\mf{\Delta}^{\mathcal{Y}}_{n}\Vert_{2} \leq (3p)^{-1}$ and $\Vert\mf{L}^{*}\Vert_{2}<3p$},
%we can define the spectral decomposition of
%the positive semi-definite matrix $\mf{\Delta}^{\mathcal{Y}}_n\mf{\Delta}^{*}_n$
%as $\mf{U}_{\mf{\Delta}}^{\mathcal{Y}}(\mf{\Lambda}_{\mf{\Delta}}^{\mathcal{Y}})^2\mf{U}_{\mf{\Delta}}^{\mathcal{Y}'}$,
%and recalling \eqref{max_logdet}
we can show that
$$
\Vert \phi'(\widehat{\mf{L}}_{\mathcal{T}'},\widehat{\mf{S}}_{\Omega})^{(ld)}\Vert_{2}\leq
(1+\lambda_1(\mf{\Lambda}_{\mf{\Delta}}^{\mathcal{Y}})^2)^{-1}\Vert\mf{\Delta}^{\mathcal{Y}}_n\Vert_{2},
%\Vert(\mf{I}_p+\mf{\Delta}^{\mathcal{Y}}_n\mf{\Delta}^{*}_n)^{-1}\Vert\Vert\mf{\Delta}^{\mathcal{Y}}_n\Vert,
$$
because
\begin{eqnarray}\Vert (\mf{I}_p+\mf{\Delta}^{\mathcal{Y}}_n\mf{\Delta}^{\mathcal{Y}}_n)^{-1}\Vert_{2}=
\Vert(\mf{\Lambda}_{\mf{\Delta}}^{\mathcal{Y}})^2+\mf{I}_p)^{-1}\Vert_{2}
=\frac{1}{1+\lambda_1(\mf{\Lambda}_{\mf{\Delta}}^{\mathcal{Y}})^2},\label{inv_bound_one}
\end{eqnarray}
%and $\lambda_p(\mf{\Lambda}_{\mf{\Delta}}^{\mathcal{Y}})^2$ tends to $0$ with probability tending to 1 under the conditions of Lemma \ref{Lemma_cons},
so that in general $\Vert(\mf{I}_p+\mf{\Delta}^{\mathcal{Y}}_n\mf{\Delta}^{\mathcal{Y}'}_n)^{-1}\Vert_{2} \leq 1$.
Since it holds $$\Vert(\mf{I}_p+\mf{\Delta}^{\mathcal{Y}}_n\mf{\Delta}^{\mathcal{Y}'}_n)^{-1}\Vert_{\infty}\leq
\Vert(\mf{I}_p+\mf{\Delta}^{\mathcal{Y}}_n\mf{\Delta}^{\mathcal{Y}'}_n)^{-1}\Vert_{2} \leq 1,$$
%if follows that $g_\gamma((\mf{I}_p+\mf{\Delta}^{\mathcal{Y}}_n\mf{\Delta}^{*}_n)^{-1}) \leq \gamma$
it follows that
$$
g_\gamma(\phi'(\widehat{\mf{L}}_{\mathcal{T}'},\widehat{\mf{S}}_{\Omega})^{(ld)})\leq
g_\gamma(\phi'(\widehat{\mf{L}}_{\mathcal{T}'},\widehat{\mf{S}}_{\Omega})^{(F)}).
$$
%As a consequence,
%%and that under the conditions of Lemma \ref{Lemma_cons}
%%$\Vert (\mf{I}_p+\mf{\Delta}_{n}\mf{\Delta}_{n}')^{-1}\Vert$ tends to $1$.
%%Consequently, the following inequality holds:
%$\Vert \widetilde{\phi}'(\widehat{\mf{L}}_{\mathcal{T}'},\widehat{\mf{S}}_{\Omega})^{(ld)}\Vert \leq \Vert \mf{\Delta}^{\mathcal{Y}}_n \Vert$.\\
%Since $\widetilde{\phi}'(\widehat{\mf{L}}_{\mathcal{T}'},\widehat{\mf{S}}_{\Omega})^{(F)}=\mf{\Delta}^{\mathcal{Y}}_n$,
%it immediately follows that
%\begin{equation}
%\Vert \widetilde{\phi}'(\widehat{\mf{L}}_{\mathcal{T}'},\widehat{\mf{S}}_{\Omega})^{(ld)} \Vert \leq
%\Vert \widetilde{\phi}'(\widehat{\mf{L}}_{\mathcal{T}'},\widehat{\mf{S}}_{\Omega})^{(F)} \Vert.
%\label{ineq_grad}
%\end{equation}

%%

Let us observe that, since $\widehat{\mf{S}}_{\Omega}\in\Omega$ and $\widehat{\mf{L}}_{\mathcal{T}'}\in\mathcal{T}'$,
$$\mathbb{P}_{{\Omega}}(\widetilde{\phi}'(\widehat{\mf{L}}_{\mathcal{T}'},\widehat{\mf{S}}_{\Omega})^{(ld)})=\mf{Z}_{\Omega} \; \mbox{and} \;
\mathbb{P}_{\mathcal{T'}}(\widetilde{\phi}'(\widehat{\mf{L}}_{\mathcal{T}'},\widehat{\mf{S}}_{\Omega})^{(ld)})=\mf{Z}_{\mathcal{T'}},$$
where $\mf{Z}_{\mathcal{T'}}=-\psi_{0} \widetilde{\mf{U}}_L \widetilde{\mf{U}}_L'$, $\mf{Z}_{\Omega}=-\psi_{0}\gamma\;\mathrm{sgn}(\mf{S}^{*})$, and $\mf{Z}=\bigl(\mf{Z}_{\mathcal{T'}},\mf{Z}_{\Omega}\bigr)$.
The two conditions on the projected gradient are essential to ensure the optimality of $(\widehat{\mf{L}}_{\mathcal{T}'},\widehat{\mf{S}}_{\Omega})$.
%We can write %$\gamma=\psi/\rho$,
%$\mathbb{P}_{{\Omega}}(\mf{\Delta}^{\mathcal{Y}}_n)=\mf{Z}_{\Omega}$ and
%$\mathbb{P}_{\mathcal{T'}}(\mf{\Delta}^{\mathcal{Y}}_n)=\mf{Z}_{\mathcal{T'}}$.
%Being $\widehat{\mf{S}}_{\Omega}\in\Omega$ and $\widehat{\mf{L}}_{\mathcal{T}'}$ the optimal solution pair,
It follows that $\Vert\mf{Z}_{\Omega}\Vert_{\infty}=\psi_{0}\gamma$,
$\Vert\mf{Z}_{\mathcal{T'}}\Vert_{2} \leq 2\psi_0$, and therefore, $g_\gamma(\mf{Z})\leq 2\psi_{0}$.
%%

%At this stage, we can write $$g_\gamma((\mf{I}_p+\mf{\Delta}^{\mathcal{Y}}_n\mf{\Delta}^{'\mathcal{Y}}_n)^{-1}\mf{\Delta}^{\mathcal{Y}}_n)= $$g_\gamma((\mf{I}_p+\widetilde{\mf{\Delta}}^{\mathcal{Y}}_n\mf{\Delta}^{'\mathcal{Y}}_n)^{-1}\widetilde{\mf{\Delta}}^{\mathcal{Y}}_n)= %$$\frac{2}{1+\lambda_p^{2}(\mf{\Delta}^{\mathcal{Y}}_n)}\psi_0.$$

%because, by Proposition \ref{conv_delta},
%\begin{eqnarray}
%g_\gamma((\mf{I}_p+\mf{\Delta}^{\mathcal{Y}}_n\mf{\Delta}^{'\mathcal{Y}}_n)^{-1})&=&
%%\max\bigl({\Vert(\mf{I}_p+\mf{\Delta}^{\mathcal{Y}}_n\mf{\Delta}^{'\mathcal{Y}}_n)^{-1}\Vert_{\infty},
%%\Vert(\mf{I}_p+\mf{\Delta}^{\mathcal{Y}}_n\mf{\Delta}^{'\mathcal{Y}}_n)^{-1}\Vert_{2}}\bigr)\nonumber\\
%\gamma^{-1}\Vert(\mf{I}_p+\mf{\Delta}^{\mathcal{Y}}_n\mf{\Delta}^{'\mathcal{Y}}_n)^{-1}\Vert_{2} \leq
%\frac{24}{\xi(\mathcal{T}(\mf{L}^{*}))}\Vert\sum_{j=0}^{\infty}(\mf{\Delta}^{\mathcal{Y}}_n)^{2j}\Vert \nonumber\\
%&\leq& 24\sqrt{\frac{p}{r}} \sum_{j=0}^{\infty}\Vert(\mf{\Delta}^{\mathcal{Y}}_n)\Vert^{2j} =
%24\sqrt{\frac{p}{r}} \frac{1}{1-\Vert(\mf{\Delta}^{\mathcal{Y}}_n)\Vert^{2}}
%%\sum_{j=0}^{\infty}\Vert\mf{L}^{*}\Vert_{2}^{-1}\Vert(\mf{\Delta}^{\mathcal{Y}}_n)\Vert^{2j}=\frac{\Vert\mf{L}^{*}\Vert_{2}^{-1}}{1-\Vert(\mf{\Delta}^{\mathcal{Y}}_n)\Vert^{2}}
%%=\frac{\Vert\mf{L}^{*}\Vert_{2}^{-1}}{1-\Vert\mf{L}^{*}\Vert_{2}^2/9p}
%%\Vert\mf{L}^{*}\Vert_{2}^{-1}\Vert\mf{\Delta}^{\mathcal{Y}}_{n}\Vert_{2} \leq (3p)^{-1}.\label{ineq_ggamma}
%\end{eqnarray}
%\Vert\mf{L}^{*}\Vert_{2}^{-1}\Vert\mf{\Delta}_{n}\Vert_{2} \leq (3\delta_{\phi} p)^{-1}

%Recalling $\mathcal{Y}=\Omega\oplus\mathcal{T'}$,
Let us define $\mf{\Delta}_L=\widehat{\mf{L}}_{\mathcal{T}'}-\mf{L}^{*}$ and $\mf{\Delta}_S=\widehat{\mf{S}}_{\Omega}-\mf{S}^{*}$.
Then, following the proof of Proposition 12 in \cite{luo2011high}, we can note that,
since $\Vert \mf{\Delta}^{\mathcal{Y}}_n \Vert_{2}<1$, under the convexity condition, it holds
%{because $\Vert \mf{L}^{*} \Vert_{2}<3p$},
%\begin{equation}\mathbb{P}_{\mathcal{Y}}\mathcal{A}^{\dag}(\mf{I}_p+\mf{\Delta}^{\mathcal{Y}}_n\mf{\Delta}^{'\mathcal{Y}}_n)^{-1}\mf{\Delta}^{\mathcal{Y}}_n=
%(\mf{I}_p+\mf{\Delta}^{\mathcal{Y}}_n\mf{\Delta}^{'\mathcal{Y}}_n)^{-1}\mf{Z},\label{eq_1_ld}\end{equation}
%and
%$$\mf{\Delta}^{\mathcal{Y}}_n=\widehat{\mf{S}}_{\Omega}+\widehat{\mf{L}}_{\mathcal{T}'}-\mf{\Sigma}_n=\mf{\Delta}^{*}_n+
%\mathcal{A}\mathbb{P}_{\mathcal{Y}}(\mf{\Delta}_L,\mf{\Delta}_S)-\mf{C}_{\mathcal{T}'},\label{eq_1_ld}$$
%This translates into
\begin{eqnarray}(\mf{I}_p+\mf{\Delta}^{\mathcal{Y}}_n\mf{\Delta}^{'\mathcal{Y}}_n)^{-1}\mf{\Delta}^{\mathcal{Y}}_n=
(\mf{I}_p+\mf{\Delta}^{\mathcal{Y}}_n\mf{\Delta}^{'\mathcal{Y}}_n)^{-1}
(\mathcal{I}^{*}\mathcal{A}\mathbb{P}_{\mathcal{Y}}(\mf{\Delta}_{L},\mf{\Delta}_{S})-\mf{\Delta}^{*}_n-\mathcal{I}^{*}\mf{C}_{\mathcal{T}'}),\label{eq_2_ld}
\end{eqnarray}
and we can apply Brouwer's fixed point theorem, to seek for the fixed point of the function
\begin{eqnarray}
F(\mf{\Delta}_L,\mf{\Delta}_S)=(\mf{\Delta}_L,\mf{\Delta}_S)
-(\mathbb{P}_{\mathcal{Y}^\perp}\mathcal{A}^{\dag}\mathcal{I}^{*}\mathcal{A}
\mathbb{P}_{\mathcal{Y}})^{-1}\mathbb{P}_{\mathcal{Y}}\mathcal{A}^{\dag}\nonumber\\
\{(\mf{I}_p+\mf{\Delta}^{\mathcal{Y}}_n\mf{\Delta}^{'\mathcal{Y}}_n)^{-1}
(\mathcal{I}^{*}\mathcal{A}\mathbb{P}_{\mathcal{Y}}(\mf{\Delta}_{L},\mf{\Delta}_{S})-\mf{\Delta}^{*}_n-\mathcal{I}^{*}\mf{C}_{\mathcal{T}'})\}.\nonumber
%\label{brouwer}
\end{eqnarray}
We know that $\mathbb{P}_{\mathcal{Y}}\bigl(\mf{\Delta}_L,\mf{\Delta}_S\bigr)$ is a fixed point of
$F(\mf{\Delta}_L,\mf{\Delta}_S)$, and it is unique because it satisfies the two optimality conditions on projected gradients.
%due to \eqref{eq_1_ld} and \eqref{eq_2_ld}.
Then, relying on Proposition \ref{11} (part (i)) and the inequalities $g_\gamma(\mf{Z})\leq 2\psi_0$ and \eqref{inv_bound_one}, since
$${\alpha_{\mathcal{Y}}^{-1}\Bigl(\frac{3}{2}-\alpha_{\mathcal{Y}}\Bigr)}^{-1}=\frac{2}{\alpha_{\mathcal{Y}}(3-2\alpha_{\mathcal{Y}})},$$
we can define $\mf{M}_{\mathcal{Y}}=\mathcal{I}^{*}\mathcal{A}\mathbb{P}_{\mathcal{Y}}(\mf{\Delta}_{L},\mf{\Delta}_{S})-\mf{\Delta}^{*}_n-\mathcal{I}^{*}\mf{C}_{\mathcal{T}'}$
and $\mf{M}_{\mathbb{P}_{\mathcal{Y}}}=\mathbb{P}_{\mathcal{Y}}(\mf{\Delta}_{L},\mf{\Delta}_{S})-\mathcal{A}^{\dag}\mf{\Delta}^{*}_n
-\mathcal{A}^{\dag}\mathcal{I}^{*}\mf{C}_{\mathcal{T}'}$,
and we can write
\begin{eqnarray}
&&g_\gamma\{F(\mf{\Delta}_L,\mf{\Delta}_S)\}\nonumber\\
&\leq&\frac{2g_\gamma(\mathbb{P}_{\mathcal{Y}}\mathcal{A}^{\dag}\{(\mf{I}_p+\mf{\Delta}^{\mathcal{Y}}_n\mf{\Delta}^{'\mathcal{Y}}_n)^{-1}
(\mf{M}_{\mathcal{Y}})\}}{\alpha_{\mathcal{Y}}(3-2\alpha_{\mathcal{Y}})}
\nonumber\\
&\leq&\frac{4g_{\gamma}(\mathcal{A}^{\dag}(\mf{I}_p+\mf{\Delta}^{\mathcal{Y}}_n\mf{\Delta}^{'\mathcal{Y}}_n)^{-1}
(\mf{M}_{\mathcal{Y}}))}{\alpha_{\mathcal{Y}}(3-2\alpha_{\mathcal{Y}})}
,\nonumber\\%\leq\tilde{r}.\nonumber\\
&\leq&\frac{4\Vert\mathcal{A}^{\dag}(\mf{I}_p+\mf{\Delta}^{\mathcal{Y}}_n\mf{\Delta}^{'\mathcal{Y}}_n)^{-1}\Vert_{2}}{\alpha_{\mathcal{Y}}(3-2\alpha_{\mathcal{Y}})}%
%g_{\gamma}(\mathcal{A}^{\dag}(\mf{I}_p+\mf{\Delta}^{\mathcal{Y}}_n\mf{\Delta}^{'\mathcal{Y}}_n)^{-1})
g_{\gamma}(\mf{M}_{\mathbb{P}_{\mathcal{Y}}})\label{final_g}\\
&\leq&\frac{4}{\alpha_{\mathcal{Y}}(3-2\alpha_{\mathcal{Y}})}%\frac{1}{1+\lambda_p(\mf{\Lambda}_{\mf{\Delta}}^{\mathcal{Y}})^2}
g_{\gamma}(\mf{M}_{\mathbb{P}_{\mathcal{Y}}})\leq\tilde{r}.\nonumber
%\end{eqnarray}
%&\leq&g_\gamma\bigl\{(\mf{I}_p+\mf{\Delta}^{\mathcal{Y}}_n\mf{\Delta}^{'\mathcal{Y}}_n)^{-1}\bigr\}\tilde{r}\leq\tilde{r}
\end{eqnarray}
%Since $3.558494\leq\frac{4}{\alpha_{\mathcal{Y}}(3-2\alpha_{\mathcal{Y}})}\leq 4$
%because $0.77155\leq \alpha_{\mathcal{Y}} \leq 1$,
Finally, it is enough to observe that
\begin{eqnarray}
g_\gamma(\widehat{\mf{L}}_{\mathcal{T}'}-\mf{L}^{*},\widehat{\mf{S}}_{\Omega}-\mf{S}^{*})\leq
g_\gamma(F(\mf{\Delta}_L,\mf{\Delta}_S))+p^{-\alpha_1}{\Vert\mf{C}_{\mathcal{T}'}\Vert_{2}},
\label{cs_ineq}
\end{eqnarray}
from which the thesis follows.
\qed

%%
%
% \section*{Appendix}\label{appendix}

% 1) Differential of the log-det (smooth part of \ref{logdet}).\\

% \noindent
% 2) Showing the conditions ensuring convexity of the log-det (smooth part of \ref{logdet}).

% \noindent
% 3) Local convexity of the log-det.

% \noindent
% 4) Derivation of the Hessian matrix.

%%

In order to enforce that the low rank and sparse solutions recovered via problem \eqref{prob_ld_rescaled} (main paper)
possess respectively the true latent rank $r$ and the true sparsity pattern $\mathrm{sgn}(\mf{S}^{*})$,
following \cite{luo2011high} we define the matrix class
\begin{eqnarray}
\mathcal{M}=\{(\mf{L},\mf{S}) \mid \mathrm{rk}(\mf{L})\leq r, \mf{S} \in \Omega(\mf{S}^{*}),
\Vert\mathbb{P}_{\mathcal{T}^{\perp}}(\mf{L}-\mf{L}^{*})\Vert_{2}\leq \xi(\mathcal{T}(\mf{L}^{*}))\psi_0,\nonumber\\
g_\gamma(\mathcal{A}\mathcal{I}^{*}\mathcal{A}^{\dag}(\mf{L}-\mf{L}^{*},\mf{S}-\mf{S}^{*}))\leq 11\psi_0 \}\nonumber
\end{eqnarray}
and we solve the constrained problem
\begin{eqnarray}
(\widehat{\mf{L}}_{\mathcal{M}},\widehat{\mf{S}}_{\mathcal{M}})=\arg
\min_{(\mf{L},\mf{S})\in\mathcal{M}}0.5\ln \det (\mf{I}_p+p^{-2\alpha_1}\mf{\Delta}_{n}\mf{\Delta}_{n}')
+\psi_{0} \Vert \mf{L} \Vert_{*}+\rho_{0} \Vert \mf{S} \Vert_{1},\label{probtang1}
\end{eqnarray}
where $\mf{\Delta}_{n}=\mf{\Sigma}_{n}-(\mf{L}+\mf{S})$.

%Following Proposition \ref{11} and \cite{chandrasekaran2012latent},
%$\alpha_{\mathcal{Y}}=\beta_{\mathcal{Y}}=1$, $\delta_{\mathcal{Y}}=0$.
The solution pair of problem \eqref{probtang1} possesses a number of interesting properties, summarized in the following corollary, adapted from \cite{chandrasekaran2012latent}.
\begin{corollary}\label{algconsprop}
Consider any $(\mf{L},\mf{S})\in\mathcal{M}$. %Recall that $\mf{\Delta}_L={\mf{L}}-\mf{L}^{*}$, $\mf{\Delta}_S={\mf{S}}-\mf{S}^{*}$.
Suppose that $\gamma$ is in the range of Proposition \ref{11}, and let
$$\delta_{max}=\max\left(1,\frac{\nu\alpha_{\mathcal{Y}}(1+\kappa_{\mathcal{T}})}{2(1-\nu)\beta_{\mathcal{Y}}(1-\kappa_{\mathcal{T}})}\right),$$
%$D=\max\left(1,\frac{\nu\alpha_{\mathcal{Y}}(1+\kappa_{\mathcal{T}})}{2(1-\nu)\beta_{\mathcal{Y}}(1-\kappa_{\mathcal{T}})}\right)$,
$$C_2=\frac{48}{\alpha_{\mathcal{Y}}}+\frac{1}{\Vert \mathcal{I}^{*} (\mf{L,\mf{S}})\Vert_{2}}=48\alpha_{\mathcal{Y}}^{-1}+(\mathrm{E}(\lambda_1(\widetilde{\mf{\Delta}}_{n}^2))+1)^2,$$
$$C_3=\frac{\nu}{4(1-\nu)}C_2^2 \delta_{max}, \qquad C_4=C_2+\frac{\alpha_{\mathcal{Y}}C_2^2(1-\nu)}{4(2-\nu)}.$$
Suppose that Assumption \ref{lowerbounds} holds, i.e., %the minimum eigenvalue of $\mf{L}^{*}$ is such that
$$\lambda_r(\mf{L}^{*})>\frac{\delta_L{\psi_{0}}}{\xi^2(T)} \qquad \mbox{and} \qquad \Vert \mf{S}^{*} \Vert_{\mathrm{min,off}}>\frac{\delta_S\psi_{0}}{\mu(\Omega)},$$
with $\delta_L=\max{(C_3,C_4)}$ and $\delta_S=\frac{\nu\alpha_{\mathcal{Y}}}{\beta_{\mathcal{Y}}2(1-\nu)}$.\\
Then, setting $\mf{C}_{\widetilde{\mathcal T}}=\mathbb{P}_{\widetilde{\mathcal T}^{\perp}}(\mf{L}^{*})$ and $\widetilde{\mathcal T}=\mathcal T'$, we get:
\begin{inparaenum}
  \item[(i)] $\mathrm{rk}(\mf{L})=r$, i.e., $\mf{L}$ is a smooth point of the low rank variety $\mathcal{L}(r)$;
  \item[(ii)] $\Vert\mathbb{P}_{\mathcal{T}^{\perp}}(\mf{L}-\mf{L}^{*})\Vert_{2}\leq \frac{\xi(T)\psi}{19\delta_{max}\Vert \mathcal{I}^{*} (\mf{L,\mf{S}})\Vert_{2}}$;
  \item[(iii)] $\varrho(\mathcal{T},\mathcal{T}')\leq\xi(T)/4$;
  \item[(iv)] $g_\gamma(\mathcal{A}^{\dag}\mathcal{I}^{*}\mf{C}_{\widetilde{\mathcal T}})\leq \frac{\psi_{0}\nu}{4(1-\nu)}$;
  \item[(v)] $p^{-\alpha_1}{\Vert \mf{C}_{\widetilde{\mathcal T}} \Vert_{2}}\leq \frac{4(2-\nu)}{\alpha_{\mathcal{Y}}(1-\nu)}\psi_{0}$;
  \item[(vi)] $\mathrm{sgn}(\mf{S})=\mathrm{sgn}(\mf{S}^{*})$.
\end{inparaenum}
\end{corollary}
\paragraph{Proof}
The proof is analogous to the proof of Corollary 3.4 in the Supplement of \cite{chandrasekaran2012latent}.
\qed

Let us define the tangent space to $\mathcal L(r)$ in a generic $\widetilde{\mf{L}} \ne \mf{L}^{*}$:
\begin{eqnarray}
\widetilde{\mathcal T}(\widetilde{\mf{L}})=\{\mf{M}\in \R^{p \times p} \mid \mf{M}=\mf{U} \mf{Y}_{1}'+\mf{Y}_{2} \mf{U}' \mid  \mf{Y}_{1}, \mf{Y}_{2} \in \R^{p \times r},\nonumber\\
\mf{U}\in \R^{p\times r}, \mf{U}' \mf{U}=\mf{I}_r, \mf{U}' \widetilde{\mf{L}}\mf{U} \in \R^{r \times r} \mbox{diagonal}, \widetilde{\mf{L}} \in \mathcal L(r)\}.\nonumber
\end{eqnarray}
Problem \eqref{probtang1} is nonconvex in nature, due to the rank constraint. The following problem (which is equal to \eqref{probtang})
may be thought of as the convex version of problem \eqref{probtang1}:
%Consider the solution pair %define the following problem:
\begin{eqnarray}
(\widehat{\mf{L}}_{\widetilde{\mathcal{T}}},\widehat{\mf{S}}_{{\Omega}})=\arg
\min_{{\underline{\mf{L}} \in \widetilde{\mathcal T},\underline{\mf{S}} \in \Omega}}0.5\ln \det (\mf{I}_p+p^{-2\alpha_1}\underline{\mf{\Delta}}_{n}\underline{\mf{\Delta}}_{n}')
+\psi_{0} \Vert\underline{\mf{L}} \Vert_{*}+\rho_{0} \Vert\underline{\mf{S}} \Vert_{1},\label{probtang2}
\end{eqnarray}
where $\underline{\mf{\Delta}}_{n}=\mf{\Sigma}_n-(\underline{\mf{L}}+\underline{\mf{S}})$.
The following proposition proves that problems \eqref{probtang1} and \eqref{probtang2} are equivalent
under specific conditions on $g_\gamma(\mathcal{A}^{\dag}\mf{\Delta}^{*}_n)$ and $\alpha_{\mathcal{Y}}$.
\begin{proposition}\label{plus}
%Let $\gamma$ be in the range of Proposition \ref{11} and suppose that the minimum eigenvalue of $\mf{L}^{*}$ is such that $\lambda_r(\mf{L}^{*})>\delta_L{\psi_{0}}/{\xi^2(T)}$ and $\Vert \mf{S}^{*} \Vert_{\mathrm{min,off}}>\delta_S{\psi_{0}}/{\mu(\Omega)}$.
Suppose that Corollary \ref{algconsprop} and Proposition \ref{12} hold true.
Suppose also that $g_\gamma(\mathcal{A}^{\dag}\mf{\Delta}^{*}_n)\leq \frac{\psi_{0}\nu}{4(1-\nu)}$, with $\mf{\Delta}^{*}_n=\mf{\Sigma}_n-\mf{\Sigma}^{*}$.
%if $\alpha_{\mathcal{Y}}\geq 0.77155$,
Then, we have that $$(\widehat{\mf{L}}_{\widetilde{\mathcal{T}}},\widehat{\mf{S}}_{{\Omega}})=(\widehat{\mf{L}}_{\mathcal{M}},\widehat{\mf{S}}_{\mathcal{M}}).$$
\end{proposition}

\paragraph{Proof}
We need to prove that the solution pair $(\widehat{\mf{L}}_{\widetilde{\mathcal{T}}},\widehat{\mf{S}}_{{\Omega}})$ satisfies
the constraints of set $\mathcal{M}$. First, passing by the nonconvex constraint $\mathrm{rk}(\mf{L})\leq r$
to the convex constraint $\underline{\mf{L}} \in \widetilde{\mathcal T}$, since Corollary \ref{algconsprop} holds true,
under Proposition \ref{11} $\widehat{\mf{L}}_{\mathcal{M}}$ is a smooth point of the low rank variety $\mathcal{L}(r)$,
which implies that the solutions of problems \eqref{probtang1} and \eqref{probtang2} are unique,
because the Hessian of $\widetilde{\phi}_{D}^{(ld)}(\underline{\mf{L}},\underline{\mf{S}})$ is positive definite
under the constraints ${\underline{\mf{L}} \in \widetilde{\mathcal T},\underline{\mf{S}} \in \Omega}$.
Analogously, the constraint $\Vert\mathbb{P}_{\mathcal{T}^{\perp}}(\underline{\mf{L}}-\mf{L}^{*})\Vert_{2}\leq \xi(\mathcal{T}(\mf{L}^{*}))\psi_0$
is implied by part (ii) of Corollary \ref{algconsprop}.
Consequently, it remains to prove that
$$g_\gamma(\mathcal{A}\mathcal{I}^{*}\mathcal{A}^{\dag}(\mf{L}-\mf{L}^{*},\mf{S}-\mf{S}^{*}))\leq 11\psi_0.$$

First, we recall that $g_\gamma(\mathcal{A}^{\dag}\mathcal{I}^{*}\mathcal{A}(\mf{\Delta}_{L},\mf{\Delta}_{S}))\leq 2\widetilde{r}$ by Proposition \ref{12}.
%where
%$$\widetilde{r}=\max\left\{\frac{4}{\alpha_{\mathcal{Y}}(3-2\alpha_{\mathcal{Y}})}
%[g_{\gamma}(\mathcal{A}^{\dag}\mf{\Delta}^{*}_n)+g_{\gamma}(\mathcal{A}^{\dag}\mathcal{I}^{*}\mf{C}_{\mathcal{T}'})+\psi_{0}],
%\frac{\Vert \mf{C}_{\mathcal{T}'}\Vert_{2}}{p^{\alpha_1}}\right\}.$$
Then, we observe that the following inequality holds:
%$g_\gamma(\widehat{\mf{S}}_{\Omega}-\mf{S}^{*},\widehat{\mf{L}}_{\mathcal{T}'}-\mf{L}^{*})$
\begin{eqnarray}
g_\gamma(\mathcal{A}^{\dag}\mathcal{I}^{*}\mathcal{A}(\mf{\Delta}_{L},\mf{\Delta}_{S}))\leq
g_\gamma(\mathbb{P}_{\mathcal{Y}}\mathcal{A}^{\dag}\mathcal{I}^{*}\mathcal{A}\mathbb{P}_{\mathcal{Y}}(\mf{\Delta}_{L},\mf{\Delta}_{S}))
+g_\gamma(\mathbb{P}_{\mathcal{Y}^\perp}\mathcal{A}^{\dag}\mathcal{I}^{*}\mathcal{A} \mathbb{P}_{\mathcal{Y}}(\mf{\Delta}_{L},\mf{\Delta}_{S}))
+g_\gamma(\mathcal{A}^{\dag}\mathcal{I}^{*}\mf{C}_{\widetilde{\mathcal T}}).\label{ineq_ggamma}
\end{eqnarray}

By Proposition \ref{11} (part (ii)), we get
\begin{eqnarray}g_\gamma(\mathbb{P}_{\mathcal{Y}^\perp}\mathcal{A}^{\dag}\mathcal{I}^{*}\mathcal{A} \mathbb{P}_{\mathcal{Y}}(\mf{\Delta}_{L},\mf{\Delta}_{S}))\leq
(1-\nu)g_\gamma(\mathbb{P}_{\mathcal{Y}}\mathcal{A}^{\dag}\mathcal{I}^{*}\mathcal{A}\mathbb{P}_{\mathcal{Y}}(\mf{\Delta}_{L},\mf{\Delta}_{S})).\nonumber
\end{eqnarray}

By Proposition \ref{12}, we get
\begin{eqnarray}g_\gamma(\mathbb{P}_{\mathcal{Y}}\mathcal{A}^{\dag}\mathcal{I}^{*}\mathcal{A}\mathbb{P}_{\mathcal{Y}}(\mf{\Delta}_{L},\mf{\Delta}_{S}))
\leq \left\{\frac{4[g_{\gamma}(\mathcal{A}^{\dag}\mf{\Delta}^{*}_n)+g_{\gamma}(\mathcal{A}^{\dag}\mathcal{I}^{*}\mf{C}_{\mathcal{T}'})+\psi_{0}]}
{\alpha_{\mathcal{Y}}(3-2\alpha_{\mathcal{Y}})}\right\}.\nonumber\end{eqnarray}

By Corollary \ref{algconsprop}, we get
\begin{eqnarray}[g_{\gamma}(\mathcal{A}^{\dag}\mf{\Delta}^{*}_n)+g_{\gamma}(\mathcal{A}^{\dag}\mathcal{I}^{*}\mf{C}_{\mathcal{T}'})+\psi_{0}]\leq
\left[\frac{\psi_{0}\nu}{4(1-\nu)}+\frac{\psi_{0}\nu}{4(1-\nu)}+\psi_{0}\right]\leq
\leq \psi_0+\frac{\psi_{0}\nu}{2(1-\nu)}.\label{bound_up_alg}
\end{eqnarray}
%and also that $\Vert \mf{C}_{\widetilde{\mathcal T}} \Vert_{2} \leq \frac{8(2-\nu)}{2\alpha_{\mathcal{Y}}(1-\nu)}\psi_{0}$.

Summing up, from inequality \eqref{ineq_ggamma} we get
\begin{eqnarray}g_\gamma(\mathcal{A}^{\dag}\mathcal{I}^{*}\mathcal{A}(\mf{\Delta}_{L},\mf{\Delta}_{S})) &\leq&
\frac{4}{\alpha_{\mathcal{Y}}}\left(\psi_0+\frac{\psi_{0}\nu}{2(1-\nu)}\right)
+(1-\nu)\frac{4}{\alpha_{\mathcal{Y}}}\left(\psi_0+\frac{\psi_{0}\nu}{2(1-\nu)}\right)
+\frac{\psi_{0}\nu}{4(1-\nu)}\nonumber\\
&\leq&\frac{8}{\alpha_{\mathcal{Y}}}\left(\psi_0+\frac{\psi_{0}\nu}{2(1-\nu)}\right)
-\nu \frac{4}{\alpha_{\mathcal{Y}}}\left(\psi_0+\frac{\psi_{0}\nu}{2(1-\nu)}\right)
+\frac{\psi_{0}\nu}{4(1-\nu)}.\nonumber\end{eqnarray}

Since it must be $$g_\gamma(\mathcal{A}^{\dag}\mathcal{I}^{*}\mathcal{A}(\mf{\Delta}_{L},\mf{\Delta}_{S}))\leq 11 \psi_{0}$$
to ensure algebraic consistency, and considering that $\frac{\nu}{(1-\nu)}\leq 2\alpha_{\mathcal{Y}}-1$ (see the proof of Proposition \ref{11}),
we get $\alpha_{\mathcal{Y}}\geq 0.77155$, which is assumed by Theorem \ref{thm_comp}.
%Since the convexity condition $\Vert\mf{\Delta}^{\mathcal{Y}}_n\Vert_{2} \leq 1/3$
%ensures that $\alpha_{\mathcal{Y}} \geq 0.81$ by inequality \eqref{minimum_Fisher} (main paper),
The thesis then follows. \qed

Finally, we need to prove that
$(\widehat{\mf{L}}_{{\mathcal{T}'}},\widehat{\mf{S}}_{{\Omega}})$ (see problem \eqref{probtang}) is also a \emph{global} solution,
i.e. is the unique solution of the unconstrained problem \eqref{prob_ld_rescaled} (main paper).

\begin{proposition}\label{13}
Under the conditions of Corollary \ref{algconsprop} and Propositions \ref{12}-\ref{plus},
$(\widehat{\mf{L}}_{{\mathcal{T}'}},\widehat{\mf{S}}_{{\Omega}})$ is also the unique solution of problem \eqref{prob_ld_rescaled} (main paper).
\end{proposition}

\paragraph{Proof}

Let us restrict the analysis to $\mathcal{Y}_{\mathcal{M}}=\mathcal{T}_{\mathcal{M}}\oplus\Omega$.
%denote the spectral decomposition of $\widehat{\mf{L}}_{\mathcal{M}}$ as $\widetilde{\mf{U}}_L\widetilde{\mf{\Lambda}}_L\widetilde{\mf{U}}_L'$.
By Proposition \ref{plus}, we can define $\mf{Z}_{\mathcal{T}_{\mathcal{M}}}=-\psi_{0} \widetilde{\mf{U}}_L \widetilde{\mf{U}}_L'$, $\mf{Z}_{\Omega}=-\psi_{0}\gamma\;\mathrm{sgn}(\mf{S}^{*})$, $\mf{Z}_{\mathcal{M}}=\bigl(\mf{Z}_{\mathcal{T}_{\mathcal{M}}},\mf{Z}_{\Omega}\bigr)$,
and observe that
$$\mathbb{P}_{\mathcal{T}_{\mathcal{M}}}(\mathcal{A}^{\dag}\phi'(\widehat{\mf{L}}_{\mathcal{T}'},\widehat{\mf{S}}_{\Omega})^{(ld)})=\mf{Z}_{\mathcal{T}_{\mathcal{M}}} \qquad \mbox{and} \qquad \mathbb{P}_{\Omega}(\mathcal{A}^{\dag}\phi'(\widehat{\mf{L}}_{\mathcal{T}'},\widehat{\mf{S}}_{\Omega})^{(ld)})=\mf{Z}_{\Omega}.$$
Wrapping up, it holds that $$\mathbb{P}_{\mathcal{Y}}(\mathcal{A}^{\dag}\phi'(\widehat{\mf{L}}_{\mathcal{T}'},\widehat{\mf{S}}_{\Omega})^{(ld)})=\mf{Z}_{\mathcal{M}},$$
%which are proved because we can restrict $\mathcal{Y}=\mathcal{T}_{\mathcal{M}}\oplus\Omega$.
and consequently, $g_\gamma(\mathcal{A}^{\dag}\phi'(\widehat{\mf{L}}_{\mathcal{T}'},\widehat{\mf{S}}_{\Omega})^{(ld)})=\psi_0$.
Therefore, we have proved the first two optimality conditions ensuring the global optimality of $(\widehat{\mf{L}}_{{\mathcal{T}'}},\widehat{\mf{S}}_{{\Omega}})$.

Then, to complete the proof, we need to prove that the other two optimality conditions hold, namely,
%Proving the optimality conditions \cite{chandrasekaran2011rank}
$$\Vert \mathbb{P}_{\mathcal{T}_{\mathcal{M}}^{\perp}}(\mathcal{A}^{\dag}\phi'(\widehat{\mf{L}}_{\mathcal{T}'},\widehat{\mf{S}}_{\Omega})^{(ld)}) \Vert_{2} < \psi_{0} \qquad \mbox{and} \qquad \Vert \mathbb{P}_{{\Omega}^{\perp}}(\mathcal{A}^{\dag}\phi'(\widehat{\mf{L}}_{\mathcal{T}'},\widehat{\mf{S}}_{\Omega})^{(ld)}) \Vert_{\infty} < \psi_{0} \gamma,$$
which is equivalent to prove that
$$g_\gamma(\mathbb{P}_{\mathcal{Y}^{\perp}}(\mathcal{A}^{\dag}\phi'(\widehat{\mf{L}}_{\mathcal{T}'},\widehat{\mf{S}}_{\Omega})^{(ld)})<\psi_0.$$

By Proposition \ref{11} (part (ii)), we get
\begin{eqnarray}g_\gamma(\mathbb{P}_{\mathcal{Y}^\perp}\mathcal{A}^{\dag}\mathcal{I}^{*}\mathcal{A} \mathbb{P}_{\mathcal{Y}}(\mf{\Delta}_{L},\mf{\Delta}_{S}))\leq
(1-\nu)g_\gamma(\mathbb{P}_{\mathcal{Y}}\mathcal{A}^{\dag}\mathcal{I}^{*}\mathcal{A}\mathbb{P}_{\mathcal{Y}}(\mf{\Delta}_{L},\mf{\Delta}_{S})).\end{eqnarray}

It follows from \eqref{bound_up_alg} that
\begin{eqnarray}
g_\gamma(\mathbb{P}_{\mathcal{Y}^\perp}(\mf{\Delta}_{L},\mf{\Delta}_{S}))) \leq (1-\nu) \left(\psi_0 + \frac{\nu}{2(1-\nu)}\psi_0\right),
= (1-\nu)\psi_0 + 0.5\nu \psi_0 \leq \psi_0-0.5\nu\psi_0 < \psi_0 \nonumber
\end{eqnarray}
%g_\gamma(\mathbb{P}_{\mathcal{Y}^\perp}\mathcal{A}^{\dag}\mathcal{I}^{*}\mathcal{A} \leq
%(1-\nu) (\psi_0 + \frac{\nu}{2(1-\nu)}\psi_0) \leq (1-\nu)\psi_0 + 0.5\nu \psi_0 \leq \psi_0-0.5\nu < \psi_0$$
because $\nu \in (0,1/2]$. The thesis then follows. \qed

%%%%%%%%%%%%%%%%%%%%%%%%%%%%%%%%%%%%%%%%%%%%%
%\textbf{Equivalence of variety and tangent-space constrained problems. Corollary to add.}

%\textbf{START}

Proposition \ref{13} also implies that the condition $\Vert \widehat{\mf{\Delta}}_n \Vert \leq 1/3$ of Proposition \ref{11}, assumed in Theorem \ref{thm_comp},
%with $\Vert\mf{L}^{*}\Vert_{2}<3p$
is equivalent to the convexity condition $\Vert \widehat{\mf{\Delta}}^{\mathcal{Y}}_n\Vert \leq 1/3$, which is assumed in Proposition \ref{12}, and is thus unnecessary.
%the same condition involving $\widehat{\mf{\Delta}}^{\mathcal{Y}}_n$ instead of $\widehat{\mf{\Delta}}_n$, which motivates the requirement of Proposition \ref{11}.}

%%

Finally, we can write:
\begin{eqnarray}g_\gamma(\widehat{\mf{S}}_{\Omega}-\mf{S}^{*},\widehat{\mf{L}}_{\mathcal{T}'}-\mf{L}^{*}) &\leq& g_\gamma(\mathcal{A}^{\dag}\mathcal{I}^{*}\mathcal{A}(\mf{\Delta}_{L},\mf{\Delta}_{S}))
+p^{-\alpha_1}\Vert \mf{C}_{\mathcal{T}'}\Vert_{2},\label{ggamma_ineq}
\end{eqnarray}
and following Proposition \ref{12} and Corollary \ref{algconsprop}, we can derive that %in the worst admissible scenario ($p=4$) that
%\begin{eqnarray}g_\gamma(\widehat{\mf{S}}_{\Omega}-\mf{S}^{*},\widehat{\mf{L}}_{\mathcal{T}'}-\mf{L}^{*}) &\leq&
%g_\gamma(\mathcal{A}^{\dag}\mathcal{I}^{*}\mathcal{A}(\mf{\Delta}_{L},\mf{\Delta}_{S}))+ p^{-\alpha_1}{\Vert \mf{C}_{\mathcal{T}'}\Vert_{2}}\nonumber\\
%&\leq& 11\psi_{0}+\frac{4(2-\nu)}{\alpha_{\mathcal{Y}}(1-\nu)}\psi_{0}\nonumber\\
%&\leq& (11+12.82668)\psi_{0}=23.82668\psi_{0}.\nonumber
%\end{eqnarray}
in the best case scenario, i.e. $\alpha_{\mathcal{Y}}=1$, $\delta_{\mathcal{Y}}=1-\alpha_{\mathcal{Y}}=0$, $\nu=1/2$, we get
\begin{eqnarray}
g_\gamma(\mathcal{A}^{\dag}\mathcal{I}^{*}\mathcal{A}(\mf{\Delta}_{L},\mf{\Delta}_{S})) &\leq&
\frac{4}{\alpha_{\mathcal{Y}}}\left(\psi_0+\frac{\psi_{0}\nu}{2(1-\nu)}\right)
+(1-\nu) \frac{4}{\alpha_{\mathcal{Y}}}\left(\psi_0+\frac{\psi_{0}\nu}{2(1-\nu)}\right)
+\frac{\psi_{0}\nu}{4(1-\nu)}\nonumber\\
&=&(6+3+0.25)\psi_{0}=9.25\psi_{0},\qquad\mbox{and}\nonumber\\
p^{-\alpha_1} {\Vert \mf{C}_{\widetilde{\mathcal T}'} \Vert_{2}}
&\leq& \frac{8(2-\nu)}{2\alpha_{\mathcal{Y}}(1-\nu)}\psi_{0}=12\psi_{0}.\nonumber
\end{eqnarray}
Consequently, from \eqref{ggamma_ineq}, we get
$$g_\gamma(\widehat{\mf{S}}_{\Omega}-\mf{S}^{*},\widehat{\mf{L}}_{\mathcal{T}'}-\mf{L}^{*}) \leq 21.25\psi_{0}.$$

In the very end, let us define $\underline{\mf{\Delta}}_{n}=\mf{\Sigma}_n-(\underline{\mf{L}}+\underline{\mf{S}})$,
$\underline{\varphi}(\mf{\Sigma})=\mf{I}_p+p^{-2\alpha_1}\underline{\mf{\Delta}}_{n}\underline{\mf{\Delta}}_{n}'$,
$$\mf{\Delta}^{(ld)}_L=\widehat{\mf{L}}_{\mathcal{T}'}-\mf{L}^{*}, \mf{\Delta}^{(ld)}_S=\widehat{\mf{S}}_{\Omega}-\mf{S}^{*},
\qquad \mbox{and} \qquad \mf{\Delta}^{(F)}_L=\widehat{\mf{L}}^{(F)}_{\mathcal{T}'}-\mf{L}^{*}, \mf{\Delta}^{(F)}_S=\widehat{\mf{S}}^{(F)}_{\Omega}-\mf{S}^{*},$$
where
\begin{equation}%\label{probtangbis}
\bigl(\widehat{\mf{L}}^{(ld)}_{\mathcal{T}'},\widehat{\mf{S}}^{(ld)}_{\Omega}\bigr)=\arg\!\!\!\!\!
\min_{\underline{\mf{L}} \in \mathcal{T}',\underline{\mf{S}} \in \Omega}
0.5\ln \det (\underline{\varphi}(\mf{\Sigma}))
+\mathcal{P}_{\gamma}(\mf{L},\mf{S})\nonumber
\end{equation}
and
\begin{equation}%\label{probtangF}
\bigl(\widehat{\mf{L}}^{(F)}_{\mathcal{T}'},\widehat{\mf{S}}^{(F)}_{\Omega}\bigr)=\arg\!\!\!\!\!
\min_{\underline{\mf{L}} \in \mathcal{T}',\underline{\mf{S}} \in \Omega}
\frac{1}{2p^{\alpha_1}}\Vert \underline{\mf{\Delta}}_{n} \Vert_{F}
+\mathcal{P}_{\gamma}(\mf{L},\mf{S}),\nonumber
\end{equation}
with $\mathcal{P}_{\gamma}(\mf{L},\mf{S})=\psi^{-1}\mathcal{P}(\mf{L},\mf{S})$ and $\mathcal{P}(\mf{L},\mf{S})=\psi \Vert\mf{L}\Vert_{*}+\rho \Vert\mf{S}\Vert_{1}$.
%From inequality \eqref{ineq_grad}, it follows that
%\begin{eqnarray}
%%g_\gamma\{F(\mf{\Delta}^{(ld)}_L,\mf{\Delta}^{(ld)}_S)\}&\leq& g_\gamma\{(\mf{I}_p+\mf{\Delta}^{\mathcal{Y}}_n\mf{\Delta}^{'\mathcal{Y}}_n)^{-1}\}g_\gamma\{F(\mf{\Delta}^{(F)}_L,\mf{\Delta}^{(F)}_S)\}
%%\nonumber\\
%g_\gamma\{F(\mf{\Delta}^{(ld)}_L,\mf{\Delta}^{(ld)}_S)\} \leq g_\gamma\{F(\mf{\Delta}^{(F)}_L,\mf{\Delta}^{(F)}_S)\},\label{ineq_F}
%\end{eqnarray}
%due to \eqref{ineq_ggamma}.
Defining $\mf{Z}_{0}=\mathbb{P}_{\mathcal{Y}}(\mf{\Delta}^{\mathcal{Y}}_n)$, we can derive
\begin{eqnarray}\mathbb{P}_{\mathcal{Y}}(\mathcal{A}^{\dag}\phi'(\widehat{\mf{L}}_{\mathcal{T}'},\widehat{\mf{S}}_{\Omega})^{(ld)})
=\mathbb{P}_{\mathcal{Y}}\left(\sum_{j=0}^{\infty}(-1)^{j}(\mf{\Delta}^{\mathcal{Y}}_n)^{2j+1}\right)
=\sum_{j=0}^{\infty}(-1)^{j}\mathbb{P}_{\mathcal{Y}}(\mf{\Delta}^{\mathcal{Y}}_n)^{2j+1}=\sum_{j=0}^{\infty}(-1)^{j}\mf{Z}_{0}^{2j+1}.\nonumber
\end{eqnarray}
We can note that
\begin{eqnarray}
g_\gamma(\mathbb{P}_{\mathcal{Y}}(\mathcal{A}^{\dag}\phi'(\widehat{\mf{L}}_{\mathcal{T}'},\widehat{\mf{S}}_{\Omega}))^{(ld)})
&=&g_\gamma\left(\sum_{j=0}^{\infty}(-1)^{j}\mf{Z}_{0}^{2j+1}\right)
=\sum_{j=0}^{\infty}(-1)^{j}g_\gamma(\mf{Z}_{0}^{2j+1})\nonumber\\
&=&\sum_{j=0}^{\infty}(-1)^{j}g_\gamma(\mf{Z}_{0})^{2j+1}
=\frac{g_\gamma(\mf{Z}_{0})}{1+g_\gamma(\mf{Z}_{0}^{2})}.\nonumber
\end{eqnarray}
%\sum_{j=0}^{\infty}(-1)^{j}\psi_0^{2j+1}
Importantly, since $g_\gamma(\mathbb{P}_{\mathcal{Y}}(\mathcal{A}^{\dag}\phi'(\widehat{\mf{L}}_{\mathcal{T}'},\widehat{\mf{S}}_{\Omega})^{(F)}))={g_\gamma(\mf{Z}_{0})}$,
it follows that
\begin{equation}
g_\gamma(\mathbb{P}_{\mathcal{Y}}(\mathcal{A}^{\dag}\phi'(\widehat{\mf{L}}_{\mathcal{T}'},\widehat{\mf{S}}_{\Omega})^{(ld)}))\leq
g_\gamma(\mathbb{P}_{\mathcal{Y}}(\mathcal{A}^{\dag}\phi'(\widehat{\mf{L}}_{\mathcal{T}'},\widehat{\mf{S}}_{\Omega})^{(F)})),
\end{equation}
the equality holding if and only if $\mf{Z}_{0}=\mf{0}$, which happens only in the limit as $n \to \infty$.
Then, by \eqref{final_g}, we can write
\begin{equation}
%\max g_\gamma(\widehat{\mf{S}}^{(ld)}_{\Omega}-\mf{S}^{*},\widehat{\mf{L}}^{(ld)}_{\mathcal{T}'}-\mf{L}^{*})&=&
g_\gamma(F(\mf{\Delta}^{(ld)}_L,\mf{\Delta}^{(ld)}_S))\leq g_\gamma(F(\mf{\Delta}^{(F)}_L,\mf{\Delta}^{(F)}_S)),\label{ineqF}\\
%&=& \max g_\gamma(\widehat{\mf{S}}^{(F)}_{\Omega}-\mf{S}^{*},\widehat{\mf{L}}^{(F)}_{\mathcal{T}'}-\mf{L}^{*}),\nonumber
\end{equation}
%Under the conditions of Proposition \ref{13},
and by \eqref{cs_ineq} and \eqref{ineqF} it follows that
\begin{equation}
\frac{\max g_\gamma(\widehat{\mf{S}}^{(ld)}-\mf{S}^{*},\widehat{\mf{L}}^{(ld)}-\mf{L}^{*})}
{\max g_\gamma(\widehat{\mf{S}}^{(F)}-\mf{S}^{*},\widehat{\mf{L}}^{(F)}-\mf{L}^{*})} \leq 1,\nonumber
\end{equation}
%due to Proposition \ref{13} (part 3), t
which finally proves the theorem.

\subsection*{{Proof of Theorem \ref{thm_main}}}

%\far{First of all, we observe that, due to Lemma \ref{Lemma_cons}, as $n \to \infty$ local convexity of $\phi_D(\mf{L},\mf{S})$ and existence of $\phi_D'(\mf{L},\mf{S})$ are ensured around $(\mf{L}^{*},\mf{S}^{*})$.}

First, we need to ensure that under Assumption \ref{alg},
Assumptions \ref{lowerbounds}(i) and \ref{eigenvalues}(i) are compatible, i.e. that
$$\lambda_r(\mf{L}^{*}) > \delta_L \frac{\psi_{0}}{\xi^2(T(\mf{L}^{*}))} \geq \delta_L \left(\frac{\sqrt{r}}{\kappa_{L}^{\delta_{1}}}\right)^3 {p^{3\delta_{1}}}f_{\delta_\epsilon}(p,n)$$
under $\lambda_r(\mf{L}^{*})\simeq p^{\alpha_r}$,
%for some $\kappa_r>0$, which is inactive
%holds under $\frac{p^{6\delta_{1}-2\alpha_r}}{n} \to 0$.
%This condition actually reduces to the condition $\frac{p^{2-2\alpha_r}}{n} \to 0$,
%unless $6\delta_{1}-2\alpha_r>2-2\alpha_r$, that never holds since
which holds true if $\delta_1\leq\alpha_r/3$
for all $p \in \N$ as $n \to \infty$. %due to Lemma \ref{Lemma_cons}.
%, because $f_{\delta_\epsilon}(p,n)\to 0$ as $n \to \infty$, since Lemma \ref{random_conv} holds true.
Assumptions \ref{lowerbounds}(ii) and \ref{sparsity}(i) are instead always compatible, as
\footnotesize
$$0< 54\delta_S f_{\delta_\epsilon}(p,n)\leq \frac{\delta_S f_{\delta_\epsilon}(p,n)}{\xi(\mathcal{T}(\mf{L}^{*}))\mu(\Omega(\mf{S}^{*}))}$$
$$=\delta_S \frac{\psi_{0}}{\mu(\Omega(\mf{S}^{*}))}<\Vert\mf{S}^{*}\Vert_{\mathrm{\mathrm{min,off}}}<\Vert\mf{S}^{*}\Vert_{\infty}=O(1),$$
\normalsize
that is always verified for all $p \in \N$ as $n \to \infty$. %because $f_{\delta_\epsilon}(p,n)\to 0$.
{Similarly, Assumptions \ref{lowerbounds}(ii) and \ref{sparsity}(iii)
are always compatible, since $54\delta_S f_{\delta_\epsilon}(p,n)<\Vert\mf{S}^{*}\Vert_{\mathrm{\mathrm{min,off}}}=o(1/p^{1-\delta_1})$
%is always verified as $p^{1-\delta_1}\Vert\mf{S}^{*}\Vert_{\mathrm{\mathrm{min,off}}}=o(1)$.
is always true for all $p \in \N$ as $n \to \infty$.
%and ${f_{\delta_\epsilon}(p,n)}=o(1/p^{1-\delta_1})$, which leads to the condition $p^{1-\delta_1}f_{\delta_\epsilon}(p,n)=o(1)$.}

%\paragraph{Proof}
At this stage, we note that under Assumptions \ref{tails}-\ref{sparsity},
%there exists $C>0$, such that,
for all $p \in \N$ as $n \to \infty$
it holds:
\begin{eqnarray}
g_\gamma(\mathcal{A}^{\dag}\mf{\Delta}_n)&=&g_\gamma(\mf{\Sigma}_n-\mf{\Sigma}^{*},\mf{\Sigma}_n-\mf{\Sigma}^{*})\nonumber\\
&\leq& \max\left({{\gamma}^{-1}\Vert \mf{\Sigma}_n-\mf{\Sigma}^{*} \Vert_{\infty}},
{{\Vert \mf{\Sigma}^{*} \Vert^{-1}_{2}}\Vert \mf{\Sigma}_n-\mf{\Sigma}^{*} \Vert_{2}}\right)\nonumber\\
&\leq& \max\left({{\gamma}^{-1}\Vert \mf{\Sigma}_n-\mf{\Sigma}^{*} \Vert_{\infty}},
{p^{-\alpha_{1}}\Vert {\mf{\Sigma}_n-\mf{\Sigma}^{*} \Vert_{2}}}\right)\nonumber
%&\leq \max\bigl(
%C\frac{\Vert \mf{\Sigma}_n-\mf{\Sigma}^{*} \Vert_{2}}{\gamma p^{\delta}},
%C\frac{\Vert \mf{\Sigma}_n-\mf{\Sigma}^{*} \Vert_{2}}{p^{\alpha}}
%\bigr)\leq\nonumber\\
\end{eqnarray}
leading to
\begin{equation}
g_\gamma(\mathcal{A}^{\dag}\mf{\Delta}_n)\leq C'({2\xi(\mathcal{T})})^{-1}{f_{\delta_\epsilon}(p,n)}.\label{g_gamma_bound}
\end{equation}
%\leq\frac{1}{18\xi(\mathcal{T})}f_{\delta_\epsilon}(p,n),\nonumber
%\leq C'\frac{1}{\xi(\mathcal{T})}{f_{\delta_\epsilon}(p,n)}\nonumber.
This result descends with probability $1-O(1/n^2)$ from \eqref{bound_sample_max} (main paper),
%the inequality $\Vert \mf{\Sigma}_n-\mf{\Sigma}^{*} \Vert_{\infty}\leq\Vert \mf{\Sigma}_n-\mf{\Sigma}^{*} \Vert_{2}$,
%and because
%$\Vert \mf{\Sigma}^{*} \Vert_{0,v}\geq \Vert \mf{\Sigma}^{*} \Vert_{\infty} \Vert \mf{\Sigma}^{*} \Vert_{2} \geq O(p^{\alpha_{1}})$
%O(p^{\delta_{1}})+\Vert \mf{L}^{*} \Vert_{0,v} \geq O(p^{\delta_{1}})+\Vert \mf{L}^{*} \Vert_{2}
%with $\Vert \mf{\Sigma}^{*} \Vert_{\infty} \geq \Vert \mf{S}^{*} \Vert_{\infty}=O(1)$
%\geq\frac{\Vert \mf{\Sigma}^{*} \Vert_{2}}{\Vert \mf{\Sigma}^{*} \Vert_{\infty}}\geq O(p^{\alpha_{1}})$,
%as $\Vert \mf{S}^{*} \Vert_{\infty}\leq \Vert \mf{\Sigma}^{*} \Vert_{\infty}\leq \Vert \mf{L}^{*} \Vert_{\infty}
%+\Vert \mf{S}^{*} \Vert_{\infty}=O(1)$,
%and $\Vert \mf{S}^{*} \Vert_{\infty}=O(1)$ and $\Vert \mf{L}^{*} \Vert_{\infty}=O(1)$
%under Assumptions \ref{eigenvalues}(i) and \ref{sparsity}(ii),
%, and \ref{loadings}(i),
%from the condition $\kappa_\mathcal{T}\leq 1/4$ (part 2 of Proposition \ref{13}),
from the range
$$\gamma \in \left[\frac{\frac{2(1+\kappa_\mathcal{T})}{1-\kappa_\mathcal{T}}\xi(\mathcal{T}(\mf{L}^{*}))(1-\nu)}{\nu\alpha_{\mathcal{Y}}},
\frac{\nu\alpha_{\mathcal{Y}}}{4\mu(\Omega(\mf{S}^{*}))\beta_{\mathcal{Y}}(1-\nu)}\right]$$ of Proposition \ref{11},
where the minimum for $\gamma$, $\frac{2\xi(\mathcal{T}(\mf{L}^{*}))(1-\nu)}{\nu\alpha_{\mathcal{Y}}}$,
is attained for $\alpha_{\mathcal{Y}}=1$, $\nu={1}/{2}$, and $\kappa_{\mathcal{T}}=0$,
and from \eqref{bound_sample_2} (main paper), under Assumptions \ref{tails}-\ref{sparsity}.
%that requires Assumption \ref{alg}.

%due to Proposition \ref{11} and \ref{alg}, it holds %and Assumptions \ref{sparsity}(i)
%$\min \gamma \Vert \mf{\Sigma}^{*} \Vert_{0,v} \leq 9 \delta_{2} \xi(\mathcal{T}) p^{\delta} \leq 9 \delta_{2} \kappa_{L}$
%since $\gamma \in [9\xi(\mathcal{T}),\frac{1}{6\mu(\Omega)}]$ and
%$\Vert\mf{S}^{*}\Vert_{0,v}\leq\Vert\mf{\Sigma}^{*}\Vert_{0,v}\leq\Vert\mf{L}^{*} \Vert_{0,v}+\Vert\mf{S}^{*} \Vert_{0,v}$,
%and $\xi(\mathcal{T})\leq \kappa_{L} p^{-\delta}$.

Since we have set
$\psi_{0}=\frac{f_{\delta_\epsilon}(p,n)}{\xi(\mathcal{T}(\mf{L}^{*}))}$,
part (iii) of Corollary \ref{algconsprop} can be written as
$$g_\gamma(\mathcal{A}^{\dag}\mf{\Delta}_n) \leq \frac{\psi_{0}}{4}\leq\frac{(\kappa_{L} p)^{\delta_{1}}}{4\sqrt{r}}
f_{\delta_\epsilon}(p,n)$$ by Assumption \ref{alg}.
%where we have set $C'\leq \frac{1}{2}$.
Therefore, setting $C={\sqrt{\kappa_{L}}}/(4\sqrt{r})$ and $C'={1}/{2}$ in \eqref{g_gamma_bound},
under Assumptions \ref{tails}-\ref{ass_alg} %and \ref{lowerbounds}
Proposition \ref{13} ensures that
the solution $(\widehat{\mf{L}}_n,\widehat{\mf{S}}_n)$ of problem \eqref{prob_ld_rescaled} (main paper) satisfies
\begin{equation}
g_\gamma(\widehat{\mf{L}}_n-\mf{L}^{*},\widehat{\mf{S}}_n-\mf{S}^{*})\leq \kappa p^{\delta_{1}}f_{\delta_\epsilon}(p,n),%\leq C24.18437\psi_{0}
\nonumber
\end{equation}
with $\kappa\geq(21.25\kappa_{L}^{\delta_1})/(4\sqrt{r})$,
%$\leq\leq(23.82668\kappa_{L}^{\delta_1})/(4\sqrt{r})$,
where the minimum is verified as $n \to \infty$ when $\alpha_{\mathcal{Y}}=1$ and $\nu=1/2$.
%and the maximum is for the worst admissible case $p=4$, $\alpha_{\mathcal{Y}} = 0.8287277$, $\nu=0.3966655$ (by Assumption \ref{ass_alg}).
Recalling the definition of $g_\gamma$ in \eqref{ggamma} (main paper), it thus holds with probability $1-O(1/n^2)$:
\begin{eqnarray}
\Vert\widehat{\mf{L}}_n-\mf{L}^{*}\Vert_{2}&\leq& \kappa p^{\alpha_{1}}\psi_{0} = O(p^{\alpha_{1}+\delta_{1}}f_{\delta_\epsilon}(p,n)),\nonumber\\
\Vert\widehat{\mf{S}}_n-\mf{S}^{*}\Vert_{\infty}&\leq& \kappa\gamma \psi_{0} = O(f_{\delta_\epsilon}(p,n)).\nonumber
\end{eqnarray}
This proves parts (i) and (ii) of Theorem \ref{thm_main}.
Finally, Proposition \ref{13} ensures that
\begin{eqnarray}
\mathcal{P}(\mathrm{rk}(\widehat{\mf{L}}_n)=r) = 1-O(1/n^2), \nonumber\\
\mathcal{P}(\mathrm{sgn}(\widehat{\mf{S}}_n)=\mathrm{sgn}(\mf{S}^{*})) = 1-O(1/n^2). \nonumber
\end{eqnarray}
This proves parts (iii) and (iv) of Theorem \ref{thm_main}.

\subsection*{{Proof of Corollary \ref{coroll_main}}}

%\paragraph{Proof}

Suppose that all the assumptions and conditions of Theorem \ref{thm_main} hold, and recall that the pair $(\widehat{\mf{L}}_n,\widehat{\mf{S}}_n)$
is the solution of \eqref{prob_ld_rescaled}. Consequently part (i) holds true because of Theorem \ref{thm_main} part (ii) and Assumption \ref{sparsity}(i), as
\begin{eqnarray}
\Vert\widehat{\mf{S}}_n-\mf{S}^{*}\Vert_{2} &\leq& \Vert \widehat{\mf{S}}_n-\mf{S}^{*}\Vert_{0,v} \Vert\widehat{\mf{S}}_n-\mf{S}^{*}\Vert_{\infty}
\leq \kappa \Vert \mf{S}^{*}\Vert_{0,v} f_{\delta_\epsilon}(p,n)
\leq \kappa \delta_2 {p^{\delta_{1}}} f_{\delta_\epsilon}(p,n),\nonumber
%\leq C {p^{\delta_1}}f_{\delta_\epsilon}(p,n),\nonumber
\end{eqnarray}
where we used arguments analogous to \eqref{Thr_top}, \eqref{Thr_down}, and \eqref{bound2}.
%the fact that $\mathcal{P}(\Vert \widehat{\mf{S}}_n-\mf{S}^{*}\Vert_{0,v} = \Vert\mf{S}^{*}\Vert_{0,v}) \to 1$ as $n \to \infty$.

Part (ii) holds true under Assumptions \ref{eigenvalues}(i) and \ref{sparsity}(i) because
\begin{eqnarray}
\Vert\widehat{\mf{\Sigma}}_n-\mf{\Sigma}^{*}\Vert_{2}&\leq& \Vert\widehat{\mf{L}}_n-\mf{L}^{*}\Vert_{2}+\Vert\widehat{\mf{S}}_n-\mf{S}^{*}\Vert_{2}
\leq \kappa p^{\alpha_{1}+\delta_{1}} f_{\delta_\epsilon}(p,n) + \kappa \delta_2 p^{\delta_1} f_{\delta_\epsilon}(p,n).\nonumber
%\leq C\kappa p^{\max\bigl\{\alpha_{1}+\delta_{1},\delta_{1}\bigr\}}f_{\delta_\epsilon}(p,n)
%\leq C {p^{\alpha_{1}+ \delta_{1}}}f_{\delta_\epsilon}(p,n).\nonumber
\end{eqnarray}
%because $\delta_{1}<\alpha_{1}$ by Assumption \ref{sparsity}(i).

Then, Proposition \ref{13} ensures part (iii) of the Corollary, as $\widehat{\mf{S}}_n\succ 0$ because $\widehat{\mf{S}}_n \in \mathcal{S}(s)$
for all $p \in \N$ as $n \to \infty$.
Part (iv) of the Corollary descends by Proposition \ref{13},
because $\mathcal{P}(\mathrm{rk}(\widehat{\mf{L}}_n)=r) \to 1$ for all $p \in \N$ as $n \to \infty$ (part (iii) of Theorem \ref{thm_main}),
and $\widehat{\mf{\Sigma}}_n \succ 0$ because $$\lambda_p(\widehat{\mf{\Sigma}}_n)\geq\lambda_p(\widehat{\mf{L}}_n)+\lambda_p(\widehat{\mf{S}}_n)>0+\lambda_p(\widehat{\mf{S}}_n)>0,$$
by dual Lidksii inequality and part (iii) of the Corollary.
%as $n \to \infty$

Part (v) of the Corollary holds because
$$\Vert \widehat{\mf{S}}_n^{-1}-\mf{S}^{*-1} \Vert_{2} \leq \Vert \widehat{\mf{S}}_n-\mf{S}^{*} \Vert_{2}
\frac{1}{\lambda_p(\mf{S}^{*})} \frac{1}{\lambda_p(\widehat{\mf{S}}_n)},$$
$\lambda_p(\mf{S}^{*})=O(1)$ by assumption, and $\lambda_p(\widehat{\mf{S}}_n)$ tends to $\lambda_p(\mf{S}^{*})$
for all $p \in \N$ as $n \to \infty$.

Analogously, part (vi) holds because
$$\Vert \widehat{\mf{\Sigma}}_n^{-1}-\mf{\Sigma}^{*-1} \Vert_{2} \leq \Vert \widehat{\mf{\Sigma}}_n-\mf{\Sigma}^{*} \Vert_{2} \frac{1}{\lambda_p(\mf{\Sigma}^{*})}\frac{1}{\lambda_p(\widehat{\mf{\Sigma}}_n)},$$
$\lambda_p(\mf{\Sigma}^{*})=O(1)$ by assumption, and $\lambda_p(\widehat{\mf{\Sigma}}_n)$
tends to $\lambda_p(\mf{\Sigma}^{*})$ for all $p \in \N$ as $n \to \infty$.

\section{Comparison with \cite{agarwal2012noisy}}

By basic algebra, since $\mathrm{rk}(\widehat{\mf{L}}_n)=\mathrm{rk}(\mf{L}^{*})=r$ with probability $1-O(1/n^2)$ as $n \to \infty$
(part 3 of Theorem \ref{thm_main}), once assumed $\alpha_1=1$ as in \cite{agarwal2012noisy},
we can derive that
$$\Vert\widehat{\mf{L}}_n-\mf{L}^{*}\Vert_{F} \leq \sqrt{r} \Vert\widehat{\mf{L}}_n-\mf{L}^{*}\Vert_{2} \leq \kappa p\sqrt{r} \psi_0.$$
Similarly, we can derive
$$\Vert\widehat{\mf{S}}_n-\mf{S}^{*}\Vert_{F} \leq \sqrt{s} \Vert\widehat{\mf{S}}_n-\mf{S}^{*}\Vert_{\infty}=\kappa \sqrt{sp} \rho_{0}.$$
Consequently, we obtain
\begin{eqnarray}
\Vert \widehat{{\mf{\Sigma}}}_n-{\mf{\Sigma}}^{*}\Vert_{F}&\leq&
\Vert\widehat{{\mf{L}}}_n-{\mf{L}}^{*}\Vert_{F}+\Vert\widehat{{\mf{S}}}_n-{\mf{S}}^{*}\Vert_{F}
\leq \kappa(\sqrt{rp} \psi_0+\sqrt{s}\rho_{0}).\label{rateF}
\end{eqnarray}

%\begin{remark}\label{agarwal_comparison}

Assume now that factors and residuals possess all moments, i.e., set $\delta_f=\delta_\epsilon=+\infty$.
Then, we get $\psi_{0}=\sqrt{\ln(p)/n}/{\xi(\mathcal{T}(\mf{L}^{*}))}$.
Since $\rho_0=\gamma\psi_0$, the \emph{lhs} of \eqref{rateF} becomes \begin{equation}\frac{\kappa p\sqrt{\frac{r\ln(p)}{n}}}{\xi(\mathcal{T}(\mf{L}^{*}))}+\frac{\kappa\gamma\sqrt{\frac{rs\ln(p)}{n}}}{\xi(\mathcal{T}(\mf{L}^{*}))}.
\label{rate_agarwal}\end{equation}

%%$$\sqrt{s}\sqrt{p}\gamma\kappa_{L}\sqrt{\ln(p)/n}/\sqrt{r}).$$
%%Basic scenario: \xi(\mathcal{T}(\mf{L}^{*}))=\kappa_{L}\sqrt{\ln(p)/n}/\sqrt{r}

We now consider the standard scenario with $\delta_1=0$ and $\nu=1/2$.
Then, according to Proposition \ref{11},
we get ${\xi(\mathcal{T}(\mf{L}^{*}))}=\sqrt{r}/\kappa_{L}=O(1)$ and $\gamma\leq 1/(4\mu(\Omega(\mf{S}^{*})))=1/4=O(1)$.
%1/(6\kappa_s\Vert\mf{S}^{*}\Vert_{\infty})=(6\kappa_s)^{-1}=O(1)$.
If $r=O(\ln(p))$, rate \eqref{rate_agarwal} becomes \begin{equation}\kappa\kappa_{L} p \sqrt{\frac{\ln(p)}{n}}+\frac{\kappa\kappa_{L}}{6}\sqrt{\frac{s\ln(p)}{n}},\label{rate_agarwal2}\end{equation}
%with $\kappa_{LS} \geq 9$.
%Therefore, %considering the sum of squared losses
%$\Vert\widehat{{\mf{L}}}_n-{\mf{L}}^{*}\Vert^{2}_{F}+\Vert\widehat{{\mf{S}}}_n-{\mf{S}}^{*}\Vert^{2}_{F}$ as in \cite{agarwal2012noisy},
which is equivalent to
\begin{equation}O\left(\frac{p\ln(p)}{n}\right)+O\left(\frac{s\ln(p)}{n}\right).\label{rate_agarwal3}\end{equation}
Considering instead the opposite scenario $\delta_1=1/2$,
rate \eqref{rate_agarwal3} worsens by a factor ${O(p)}$,
becoming \begin{equation}O\left(\frac{p^2\ln(p)}{n}\right)+O\left(\frac{s\ln(p)}{n}\right),\label{rate_agarwal4}\end{equation}
since ${\xi(\mathcal{T}(\mf{L}^{*}))}=\sqrt{r/p}$.
%%Logical: stronger conditions!
Rate \eqref{rate_agarwal4} exactly corresponds to the rate of Corollary 3 in \cite{agarwal2012noisy}
%This rate corresponds exactly to the error rate of Corollary 3 in \cite{agarwal2012noisy},
in our case of perfect identifiability,
because under the condition $\alpha_1=1$
it follows that $\Vert \mf{\Sigma}^{*} \Vert_{2}= O(p)$ from Assumption \ref{eigenvalues}(i) and Assumption \ref{sparsity}(i),
and it follows that $\Vert \mf{\Sigma}^{*} \Vert_{\infty}= O(1)$
from Assumption \ref{eigenvalues}(ii) and Assumption \ref{sparsity}(ii). %This completes the analogy.
%where $\alpha=0$ (as we have perfect identifiability and $s=1$,
%because $\mf{S}^{*}=\mf{I}_p$, and $\Vert\mf{S}^{*}\Vert_{0,v}=1$).

This perfect correspondence cannot surprise, because
inequality \eqref{minimum_Fisher} ensures that %the rescaled smooth component of \eqref{logdet},
$\widetilde{\phi}_{D}(\mf{L},\mf{S})$ is strongly convex for all $p,n \in \N$ as Proposition \ref{11} holds.
This implies that, under the conditions of Proposition \ref{11}, $\widetilde{\phi}_{D}(\mf{L},\mf{S})$
satisfies the Restricted Strong Convexity as defined in \cite{10.1214/12-STS400},
%for all $p \in \N$ as $n \to \infty$,
%under certain conditions (see \far{Proof ENRICO}),
so that the error rates in Frobenius norm of \cite{agarwal2012noisy} hold for $(\widehat{\mf{L}}_n,\widehat{\mf{S}}_n)$.
As a consequence,  rate \eqref{rate_agarwal4} does not only hold for all $p \in \N$ as $n \to \infty$,
as implied by Lemma \ref{lemma_strong},
%because $\frac{1}{2}\ln \det (\mf{I}_p+\mf{\Delta}_{n}\mf{\Delta}_{n}')$ is asymptotically convex due to Lemma \ref{random:second}, but also means
but also for finite $p$ and $n$,
provided that $\widetilde{\phi}_{D}(\mf{L},\mf{S})$ is strongly convex,
which is ensured as Proposition \ref{11} holds.

%\far{START NEW PAPER}
%\textbf{Section non-asymptotic}\\
%Riemannian geometry\\ %Identifiability\\
%Lemma Enrico Hessian\\
%Lemma Enrico RSC\\
%Non-asymptotic claim (Theorem and Corollary), with probability of the conditioning event!\\
%Paragone con Agarwal.\\
%\far{END}
%$\sqrt{r} \leq \kappa_{L} \leq \sqrt{p}$!

\section{Computational aspects}\label{computation}

\subsection{Lipschitz-continuity}\label{LipC}

As explained in \cite{nesterov2013gradient},
a problem like \eqref{obj2} can be numerically solved by applying a proximal gradient method. %(see Section \ref{algorithm}).
To implement it, we need to prove the Lipschitzianity of the gradient of the smooth component $\phi_D(\mf{L},\mf{S})$.
Let us recall the two-argument matrix function
$$\phi(\mf{L},\mf{S})=\phi_D(\mf{L},\mf{S})+\mathcal{P}(\mf{L},\mf{S}),$$
where
$\phi_D(\mf{L},\mf{S})=0.5\ln\det(\varphi(\mf{\Sigma}))$, with $\varphi(\mf{\Sigma})=(\mf{I}_p+\mf{\Delta}_{n}\mf{\Delta}_{n}')$, $\mf{\Delta}_n=\mf{\Sigma}-\mf{\Sigma}_n$, $\mf{\Sigma}=\mf{L}+\mf{S}$,
and $\mathcal{P}(\mf{L},\mf{S})=\psi \Vert \mf{L} \Vert_{*}+\rho \Vert \mf{S} \Vert_{1}$.
The gradient of $\phi_D(\mf{L},\mf{S})$ with respect to $\mf{L}$ and to $\mf{S}$ is the same, and corresponds to
$\varphi(\mf{\Sigma})^{-1}\mf{\Delta}_n$ under the conditions of Proposition \ref{first_der}.

We now define the $2p^2$-dimensional vectorized gradient %$\mathrm{vec}\bigl(\frac{\partial \phi_D(\mf{L},\mf{S})}{\partial (\mf{L},\mf{S})}\bigr)$,
\begin{eqnarray}
\vf{g}(\mf{\Sigma})&=&\mathrm{vec}\left(\frac{\partial \phi_D(\mf{L},\mf{S})}{\partial (\mf{L},\mf{S})}\right)
=\{\mathrm{vec}\bigl[\varphi(\mf{\Sigma})^{-1}\mf{\Delta}_n\bigr]'\mathrm{vec}\bigl[\varphi(\mf{\Sigma})^{-1}\mf{\Delta}_n\bigr]'\}'.\nonumber
\end{eqnarray}
%We define two $p \times p$ matrices $\mf{\Sigma}_2=\mf{L}_2+\mf{S}_2$, $\mf{\Sigma}_1=\mf{\Sigma}_2+\epsilon \mf{H}$, $\epsilon>0$,
%such that $\mf{\Delta}_{1,n}=\mf{\Sigma}_1-\mf{\Sigma}_n$, $\mf{\Delta}_{2,n}=\mf{\Sigma}_2-\mf{\Sigma}_n$, \far{and $\mf{H}$ is a perturbation matrix}.
%%belong to the convexity ball of $\phi_D(\mf{L},\mf{S})$, i.e $\Vert\mf{\Delta}_{1,n}\Vert \leq \frac{1}{\delta 3p}$ and
%%$\Vert\mf{\Delta}_{2,n}\Vert \leq \frac{1}{\delta 3p}$ for a specific $\delta>0$.
%We set the difference vector $\vf{d}(\mf{\Sigma}_1,\mf{\Sigma}_2)=\vf{g}(\mf{\Sigma}_1)-\vf{g}(\mf{\Sigma}_2)$, which is
%a $2p^2$-dimensional vector, composed of two identical components of $p^2$ elements, stacked one below the other.
%\begin{equation}
%\vf{d}(\mf{\Sigma}_1,\mf{\Sigma}_2)=\{\mathrm{vec}(\varphi(\mf{\Sigma}_1)^{-1}\mf{\Delta}_{1,n})'-\mathrm{vec}(\varphi(\mf{\Sigma}_2)^{-1}\mf{\Delta}_{2,n})'\}'.\label{diff_def}
%\end{equation}
The following proposition controls the effect on the vectorized gradient of a perturbation matrix $\mf{H}$.
%by means of the difference $\mf{\Delta}_{1,n}-\mf{\Delta}_{2,n}$.
\begin{proposition}\label{lips_diff}
Define two $p \times p$ matrices $\mf{\Sigma}_2=\mf{L}_2+\mf{S}_2$, $\mf{\Sigma}_1=\mf{\Sigma}_2+\mf{H}$, $\epsilon>0$,
such that $\mf{\Delta}_{1,n}=\mf{\Sigma}_1-\mf{\Sigma}_n$, $\mf{\Delta}_{2,n}=\mf{\Sigma}_2-\mf{\Sigma}_n$, and $\mf{H}=\mf{\Delta}_{1,n}-\mf{\Delta}_{2,n}$ is a perturbation matrix.
%belong to the convexity ball of $\phi_D(\mf{L},\mf{S})$, i.e $\Vert\mf{\Delta}_{1,n}\Vert \leq \frac{1}{\delta 3p}$ and
%$\Vert\mf{\Delta}_{2,n}\Vert \leq \frac{1}{\delta 3p}$ for a specific $\delta>0$.
Set the difference vector $\vf{d}(\mf{\Sigma}_1,\mf{\Sigma}_2)=\vf{g}(\mf{\Sigma}_1)-\vf{g}(\mf{\Sigma}_2)$.
Then, it holds:
$$
\Vert\vf{d}(\mf{\Sigma}_1,\mf{\Sigma}_2)\Vert_F %\leq 2 \Vert (\mf{I}_p+\mf{\Delta}_{1,n}\mf{\Delta}_{1,n}')^{-1}
%\mf{\Delta}_{1,n}-(\mf{I}_p+\mf{\Delta}_{2,n}\mf{\Delta}_{2,n}')^{-1}\mf{\Delta}_{2,n} \Vert_{F}\nonumber\\
\leq \frac{5}{2}\Vert\mf{H}\Vert_{F}.\label{lips_const}
$$
\end{proposition}

\paragraph{Proof}

From the definition of $\vf{g}(\mf{\Sigma})$ and $\varphi(\mf{\Sigma})$, it follows that
%\begin{equation}
%\Vert\vf{d}(\mf{\Sigma}_1,\mf{\Sigma}_2)\Vert_2\leq 2 \Vert\mathrm{vec}((\mf{I}_p+\mf{\Delta}_{1,n}\mf{\Delta}_{1,n}')^{-1}
%\mf{\Delta}_{1,n}-(\mf{I}_p+\mf{\Delta}_{2,n}\mf{\Delta}_{2,n}')^{-1}\mf{\Delta}_{2,n})'\Vert_2,\nonumber
%\end{equation}
%or
\begin{eqnarray}
\Vert\vf{d}(\mf{\Sigma}_1,\mf{\Sigma}_2)\Vert_F\leq 2 \Vert (\mf{I}_p+\mf{\Delta}_{1,n}\mf{\Delta}_{1,n}')^{-1}\mf{\Delta}_{1,n}
-(\mf{I}_p+\mf{\Delta}_{2,n}\mf{\Delta}_{2,n}')^{-1}\mf{\Delta}_{2,n} \Vert_{F}.\label{diff_eq}
\end{eqnarray}
Therefore, to prove the thesis, we need to bound the \emph{rhs} of \eqref{diff_eq}.
To this purpose, we study the matrix $(\mf{I}_p+\mf{\Delta}_{1,n}\mf{\Delta}_{1,n}')^{-1}\mf{\Delta}_{1,n}-(\mf{I}_p+\mf{\Delta}_{2,n}\mf{\Delta}_{2,n}')^{-1}\mf{\Delta}_{2,n}$,
that is equal to
\begin{eqnarray}
(\mf{I}_p + (\mf{\Delta}_{2,n} + \epsilon \mf{H})(\mf{\Delta}_{2,n} +
\epsilon \mf{H})')^{-1}(\mf{\Delta}_{2,n} + \epsilon \mf{H})
- (\mf{I}_p + \mf{\Delta}_{2,n}\mf{\Delta}_{2,n}')^{-1}\mf{\Delta}_{2,n}.\nonumber
\end{eqnarray}
Then, we recall from \cite{bernardi2022log}, Lemmas 3 and 4, the Lipschitzianity of the smooth function
$$\phi_D(\mf{L},\mf{S})=0.5\ln\det(\mf{I}_p+\mf{\Delta}_n\mf{\Delta}_n'),$$
and of its gradient function, $$\frac{\partial \phi_D(\mf{L},\mf{S})}{\partial \mf{L}}=\frac{\partial \phi_D(\mf{L},\mf{S})}{\partial \mf{S}}=F(\mf{\Delta}_n)=(\mf{I}_p+\mf{\Delta}_n\mf{\Delta}_n')^{-1}\mf{\Delta}_n,$$
respectively. %(see \ref{grad}).
In particular, Lemma 4 in \cite{bernardi2022log} shows that the gradient function is Lipschitz continuous with Lipschitz constant equal to ${5}/{4}$:
%$\frac{\partial \mathcal{L}(\mf{\Sigma},\mf{\Sigma}_n)}{\partial \mf{L}}=\frac{\partial \mathcal{L}(\mf{\Sigma},\mf{\Sigma}_n)}
%{\partial \mf{S}}=F(\mf{\Delta}_n)$
\begin{equation}
  %\label{lips_top}
  \Vert F(\mf{\Delta}_n + \epsilon \mf{H}) - F(\mf{\Delta}_n)\Vert_{2} \leq
  \frac{5}{4}\epsilon \Vert\mf{H}\Vert_{2} + O(\epsilon^{2}),\label{lemma:lipschitz_first}
\end{equation}
with
$$F(\mf{\Delta}_n + \epsilon \mf{H}) =  (\mf{I}_p+(\mf{\Delta}_{n} + \epsilon \mf{H})(\mf{\Delta}_{n} + \epsilon \mf{H})')^{-1}(\mf{\Delta}_n + \epsilon \mf{H}),$$
for any $ \epsilon > 0 $.
%such that $\Vert\mf{\epsilon H}\Vert_{2}\leq\frac{1}{6 \delta p}$.
%\begin{equation}
%  \label{eq:24}
%  |\ln\det\varphi(\mf{\Sigma}_{1}) -  \ln\det\varphi(\mf{\Sigma}_{2})| \leq \Vert\mf{\Sigma}_{1} - \mf{\Sigma}_{2}\Vert_{2}
%\end{equation}
Therefore, from inequality \eqref{diff_eq}, noticing that inequality \eqref{lemma:lipschitz_first} holds for the Frobenius norm as well and setting $\epsilon=1$, the thesis follows. \qed

%Note that the condition $\Vert\mf{H}\Vert_{2}\leq\frac{1}{6 \delta p}$ is derived
%by triangular inequality to respect the convexity region of Lemma \ref{lemmaconv}.

Proposition \ref{lips_diff} enables to state the proximal gradient algorithm to solve problem \eqref{obj}. %in Section \ref{algorithm}.

\subsection{Solution algorithm}\label{algorithm}

%At this stage, it is sufficient to observe that the vectorized
%$\delta_{vec} \phi_D=vec(\frac{\delta \frac{1}{2}\ln \det (\mf{I}_p+\mf{\Delta}_{n}\mf{\Delta}_{n}')}{\partial \mf{L}})$
%allows to display that the $2p^2$-dimensional vectorized $\frac{\delta \frac{1}{2}\ln \det (\mf{I}_p+\mf{\Delta}_{n}\mf{\Delta}_{n}')}{\partial \mf{L}};\frac{\delta \frac{1}{2}\ln \det (\mf{I}_p+\mf{\Delta}_{n}\mf{\Delta}_{n}')}{\partial \mf{S}}$,
%which can be written as $\vf{\delta_{vec}}=[\delta_{vec} \phi_D; \delta_{vec} \phi_D]'$.
%Therefore, according to Proposition \ref{prop_lips}, the Euclidean norm of
%$\vf{\delta_{vec}}$ is majorized by $\delfta_{vec}(\mf{L},\mf{S})$

For any $t \geq 0$, we define:
\begin{itemize}
\item $\mathcal{T}^{(S)}_{t}$, the soft-thresholding operator with parameter $t$, such that the $p \times p$ matrix $\mathcal{T}^{(S)}_{t}(\mf{M})$ has $(i,j)$ element $\mathrm{sgn}(\mf{M}_{ij}) \max(\vert \mf{M}_{ij} \vert - t, 0)$;
\item $\mathcal{T}^{(SVT)}_{t}$, the singular value thresholding operator with parameter $t$, such that the $p \times p$  matrix $\mathcal{T}^{(SVT)}_{t}(\mf{M})$ is equal to $\mf{U}_{M}\mathcal{T}^{(S)}_{t}(\mf{\Lambda}_M)\mf{U}_{M}'$, where $\mf{U}_{M}\mf{\Lambda}_M\mf{U}_{M}'$ is the spectral decomposition of $\mf{M}$.
\end{itemize}

Exploiting the results of Section \ref{LipC}, we provide a solution algorithm for problem \eqref{obj}.
Following \cite{luo2011high}, \cite{nesterov2013gradient} and the supplement of \cite{farne2020large},
we set the relevant step-size to $\ell={5}/{2}$ from Proposition \ref{lips_diff}, %\eqref{lips_const},
and we derive Algorithm \ref{alg_ld}.
\begin{algorithm}
\caption{Pseudocode to solve problem \eqref{obj} given any input covariance matrix $\mf{\Sigma}_{n}$.}\label{alg_ld}
\begin{enumerate}
\item {Set} $(\mf{L}_{0},\mf{S}_{0})=\frac{1}{2\mathrm{tr}(\mf{\Sigma}_{n})}(\mathrm{diag}(\mf{\Sigma}_{n}),\mathrm{diag}(\mf{\Sigma}_{n}))$, $\eta_{0}=1$.
\item {Initialize} $\mf{Y}_{0}=\mf{L}_{0}$ and $\mf{Z}_{0}=\mf{S}_{0}$. Set $t=1$.
\item For $t\geq 1$, {repeat}:
\begin{description}
\item[(i)] {calculate} $\mf{\Delta}_{t,n}=\mf{Y}_{t-1}+\mf{Z}_{t-1}-\mf{\Sigma}_{n}$;
\item[(ii)] {compute}
\begin{eqnarray}
&&\frac{\partial 0.5\ln \det \bigl(\mf{I}_p+\mf{\Delta}_{t,n}\mf{\Delta}_{t,n}'\bigr)}{\partial \mf{Y}_{t-1}}
=\frac{\partial 0.5\ln \det \bigl(\mf{I}_p+\mf{\Delta}_{t,n}\mf{\Delta}_{t,n}'\bigr)}{\partial \mf{Z}_{t-1}}
=(\mf{I}_p+\mf{\Delta}_{t,n}\mf{\Delta}_{t,n}')^{-1}\mf{\Delta}_{t,n}.\nonumber
\end{eqnarray}
\item[(iii)] {define} $\mf{E}_{Y,t}=\mf{Y}_{t-1}-(2/5)(\mf{I}_p+\mf{\Delta}_{t,n}\mf{\Delta}_{t,n}')^{-1}\mf{\Delta}_{t,n}$, and {set} $\mf{L}_{t}=\mathcal{T}^{(SVT)}_{\psi}(\mf{E}_{Y,t})=\widehat{\mf{U}}\widehat{\mf{D}}_\psi \widehat{\mf{U}}^\top$;
\item[(iv)] {define} $\mf{E}_{Z,t}=\mf{Z}_{t-1}-(2/5)(\mf{I}_p+\mf{\Delta}_{t,n}\mf{\Delta}_{t,n}')^{-1}\mf{\Delta}_{t,n}$, and {set} $\mf{S}_{t}=\mathcal{T}^{(S)}_\rho(\mf{E}_{Z,t})$;
\item[(v)] {set} $(\mf{Y}_{t},\mf{Z}_{t})=(\mf{L}_{t},\mf{S}_{t})+\bigl\{\frac{{\eta_{t-1}-1}}{{\eta_{t}}}\bigr\}\{(\mf{L}_{t},\mf{S}_{t})-(\mf{L}_{t-1},\mf{S}_{t-1})\}$,
    where $\eta_{t}={0.5+0.5\sqrt{1+4 \eta_{t-1}^2}}$;
\item[(vi)] {stop} if the convergence criterion $\frac{\Vert\mf{L}_{t}-\mf{L}_{t-1}\Vert_F}{{1+\Vert \mf{L}_{t-1}\Vert_F}}
    +\frac{\Vert\mf{S}_{t}-\mf{S}_{t-1}\Vert_F}{{1+\Vert \mf{S}_{t-1}\Vert_F}} \leq \varepsilon$.
\end{description}
\item {Set} $\widehat{\mf{L}}^{(ld)}_{\rm{A}}=\mathrm{tr}(\mf{\Sigma}_{n})\mf{Y}_{t}$ and $\widehat{\mf{S}}^{(ld)}_{\rm{A}}=\mathrm{tr}(\mf{\Sigma}_{n})\mf{Z}_{t}$.
\end{enumerate}
%\end{algorithm}
\end{algorithm}
In Algorithm \ref{alg_ld}, we first rescale, at step 1, by the trace of the input $\mf{\Sigma}_{n}$, and we then restore the original scale at step 4.
This approach differs from the original application in \cite{farne2020large}. It has the advantage to set the threshold grid
and to perform threshold selection in a controllable way. By Algorithm \ref{alg_ld}, we derive $(\widehat{\mf{L}}^{(ld)}_{\rm{A}},\widehat{\mf{S}}^{(ld)}_{\rm{A}})$,
where the superscript $\rm{A}$ stands for ALCE (ALgebraic Covariance Estimator).
We then present Algorithm \ref{alg_fro_2}, which is the analog of Algorithm \ref{alg_ld} for problem \eqref{obj_all} with $\mathcal{L}(\mf{\Sigma},\mf{\Sigma}_{n})=\mathcal{L}^{(F)}(\mf{\Sigma},\mf{\Sigma}_{n})$.
%with a standardized trace of the input.
%In this paper we calculate both the pair $(\widehat{\mf{L}}_n^{(F)}_{\rm{A}},\widehat{\mf{S}}^{(F)}_{\rm{A}})$
%by Algorithm \ref{alg_fro}, and the pair $(\widehat{\mf{L}}_n^{(F)}_{\rm{A}},\widehat{\mf{S}}^{(F)}_{\rm{A}})$
%by Algorithm \ref{alg_fro_2}, that guarantees, exactly as Algorithm \ref{alg_ld}, to select the thresholds $\psi$ and $\rho$
%in a controllable way.
\begin{algorithm}
\caption{Pseudocode to solve problem \eqref{obj_all} with $\mathcal{L}(\mf{\Sigma},\mf{\Sigma}_{n})=\mathcal{L}^{(F)}(\mf{\Sigma},\mf{\Sigma}_{n})$.}\label{alg_fro_2}
\begin{enumerate}
\item {Set} $(\mf{L}_{0},\mf{S}_{0})=\frac{1}{2\mathrm{tr}(\mf{\Sigma}_{n})}
    (\mathrm{diag}(\mf{\Sigma}_{n}),\mathrm{diag}(\mf{\Sigma}_{n}))$, $\eta_{0}=1$.
\item {Initialize} $\mf{Y}_{0}=\mf{L}_{0}$ and $\mf{Z}_{0}=\mf{S}_{0}$. Set $t=1$.
\item For $t\geq 1$, {repeat}:
\begin{description}
\item[(i)] {compute}
\begin{eqnarray}
\frac{\partial 0.5\Vert\mf{Y}_{t-1}+\mf{Z}_{t-1}-\mf{\Sigma}_{n}\Vert^2_{F}}{\partial \mf{Y}_{t-1}}
=\frac{\partial 0.5\Vert\mf{Y}_{t-1}+\mf{Z}_{t-1}-\mf{\Sigma}_{n}\Vert^2_{F}}{\partial \mf{Z}_{t-1}}
=\mf{Y}_{t-1}+\mf{Z}_{t-1}-\mf{\Sigma}_{n}.\nonumber
\end{eqnarray}
\item[(ii)] {define} $\mf{E}_{Y,t}=\mf{Y}_{t-1}- 0.5(\mf{Y}_{t-1}+\mf{Z}_{t-1}-\mf{\Sigma}_{n})$ and {set} $\mf{L}_{t}=\mathcal{T}^{(SVT)}_\psi(\mf{E}_{Y,t})=\widehat{\mf{U}}\widehat{\mf{D}}_\psi \widehat{\mf{U}}^\top$;
\item[(iii)] {define} $\mf{E}_{Z,t}=\mf{Z}_{t-1}- 0.5(\mf{Y}_{t-1}+\mf{Z}_{t-1}-{\mf{\Sigma}}_n)$ and {set} $\mf{S}_{t}=\mathcal{T}^{(S)}_\rho(\mf{E}_{Z,t})$;
\item[(iv)] {set} $(\mf{Y}_{t},\mf{Z}_{t})=(\mf{L}_{t},\mf{S}_{t})+\bigl\{\frac{{\eta_{t-1}-1}}{{\eta_{t}}}\bigr\}\{(\mf{L}_{t},\mf{S}_{t})-(\mf{L}_{t-1},\mf{S}_{t-1})\}$
where $\eta_{t}={0.5+0.5\sqrt{1+4 \eta_{t-1}^2}}$;
\item[(v)] {stop} if the convergence criterion $\frac{\Vert\mf{L}_{t}-\mf{L}_{t-1}\Vert_F}{{1+\Vert \mf{L}_{t-1}\Vert_F}}
    +\frac{\Vert\mf{S}_{t}-\mf{S}_{t-1}\Vert_F}{{1+\Vert \mf{S}_{t-1}\Vert_F}} \leq \varepsilon$.
\end{description}
\item {Set} $\widehat{\mf{L}}^{(F)}_{\rm{A}}=\mathrm{tr}(\mf{\Sigma}_{n})\mf{Y}_{t}$ and $\widehat{\mf{S}}^{(F)}_{\rm{A}}=\mathrm{tr}(\mf{\Sigma}_{n})\mf{Z}_{t}$.
\end{enumerate}
%\end{algorithm}
\end{algorithm}

Algorithms \ref{alg_ld} and \ref{alg_fro_2}, unlike the algorithm in \cite{farne2020large,farne2024large},
enable us to define initializers in a reasoned and effective way.
In fact, starting by the formula $\psi_{0}={f_{\delta_\epsilon}(p,n)}/{\xi(\mathcal{T}(\mf{L}^{*}))}$ of Assumption \ref{lowerbounds},
we set $\delta_\epsilon=\infty$ and $\delta_1={1}/{2}$, so that the vector of initial eigenvalue thresholds $\vf{\psi}_{init}$ can be set as
proportional to $p^{-1}$ by Assumption \ref{alg}, and the vector of initial sparsity thresholds $\vf{\rho}_{init}$ can be set as proportional to $(p\sqrt{p})^{-1}$, because $\mu(\Omega(\mf{S}^{*}))=\Vert\mf{S}^{*}\Vert_{0,v}=\sqrt{p}$ is the maximum tolerated order of the residual degree under Assumption \ref{sparsity}.
This is due to the fact that we rescale by the trace of the input in both algorithms (see step 1).
\footnote{More precisely, we set $\vf{\psi}_{init}$ as the vector of $n_{thr}$ equi-spaced points from $\frac{1}{c^{\psi}_{thr}\ln(p)p}$ to $\frac{\ln(p)}{c^{\psi}_{thr}p}$,
and $\vf{\rho}_{init}$ as $\frac{0.5c^{\rho}_{thr}}{\sqrt{p}}\vf{\psi}_{init}$, where $c^{\psi}_{thr}$ and $c^{\rho}_{thr}$ are user-based,
and must be set in order to obtain optimal thresholds by MC criterion far from the grid boundaries.}
%%where $c_{\rho}=\mathrm{geomean}(\frac{1}{p},\frac{1}{\sqrt{p}})$, because
%these two values represent the maximum and the minimum redistribution factors of the maximum row total across
%the single row elements of the residual component.
%Algorithm \ref{alg_fro}
%does not allow for the same threshold grid setup, and since the eigenvalues of a matrix
%are not a linear function of the dimension, a reliable threshold grid cannot be easily determined in that case.
%[1/(20*p) 1/(10*p) 1/(5*p) 1/(2*p) 1/p 2/p 5/p 10/p 20/p]
In the end, for each threshold pair $(\psi,\rho)$ we can calculate
$\widehat{\mf{\Sigma}}^{(ld)}_{\rm{A}}(\psi,\rho)=\widehat{\mf{L}}^{(ld)}_{\rm{A}}(\psi,\rho)+\widehat{\mf{S}}^{(ld)}_{\rm{A}}(\psi,\rho)$
%$\widehat{\mf{\Sigma}}_n^{(F)}_{\rm{A}}(\psi,\rho)=\widehat{\mf{L}}_n^{(F)}_{\rm{A}}(\psi,\rho)+\widehat{\mf{S}}^{(F)}_{\rm{A}}(\psi,\rho)$,
and $\widehat{\mf{\Sigma}}^{(F)}_{\rm{A}}(\psi,\rho)=\widehat{\mf{L}}^{(F)}_{\rm{A}}(\psi,\rho)+\widehat{\mf{S}}^{(F)}_{\rm{A}}(\psi,\rho)$.

Following \cite{farne2020large,farne2024large}, we also perform the unshrinkage of estimated latent eigenvalues,
as this operation improves the sample total loss as much as possible in the finite sample. We thus get
the UNALCE (UNshrunk ALCE) estimates as:
\begin{eqnarray}
\widehat{\mf{L}}_{\rm{U}}=\widehat{\mf{U}}_{\rm{A}}(\widehat{\mf{\Lambda}}_{\rm{A}}+\psi \mf{I}_r)\widehat{\mf{U}}_{\rm{A}}',\label{unshr1}\\
\mathrm{diag}(\widehat{\mf{S}}_{\rm{U}})=\mathrm{diag}(\widehat{\mf{\Sigma}}_{\rm{A}})-\mathrm{diag}(\widehat{\mf{L}}_{\rm{U}}),\label{unshr2}\\
\mathrm{off-diag}(\widehat{\mf{S}}_{\rm{U}})=\mathrm{off-diag}(\widehat{\mf{S}}_{\rm{A}}),\label{unshr3}
\end{eqnarray}
where $\psi>0$ is any chosen eigenvalue threshold parameter.
%\mbox{and} $\widehat{\mf{S}}_{\rm{U}}$ \mbox{such that}

By setting $\widehat{r}_{A}=\mathrm{rk}(\widehat{\mf{L}}_{\rm{A}})$ and defining the spectral decomposition of $\widehat{\mf{L}}_{\rm{A}}$ as
$\widehat{\mf{L}}_{\rm{A}}=\widehat{\mf{U}}_{\rm{A}}\widehat{\mf{D}}_{\rm{A}}\widehat{\mf{U}}_{\rm{A}}'$,
with $\widehat{\mf{U}}_{\rm{A}}$ $p \times \widehat{r}_A$ matrix such that $\widehat{\mf{U}}_{\rm{A}}'\widehat{\mf{U}}_{\rm{A}}=\mf{I}_{\widehat{r}_A}$,
and $\widehat{\mf{D}}_{\rm{A}}$ $\widehat{r}_A \times \widehat{r}_A$ diagonal matrix,
it can be proved \citep{farne2020large} that it holds
\begin{equation}
\bigl(\widehat{\mf{L}}_{\rm{U}},\widehat{\mf{S}}_{\rm{U}}\bigr)=\arg\min_{\mf{L} \in \widehat{\mathcal{L}}(\widehat{r}_{A}), \mf{S} \in \widehat{\mathcal{S}}_{diag}}\frac 12 \Vert{\mf{\Sigma}}_{n}-(\mf{L}+\mf{S})\Vert_{2},\nonumber%\label{min2}
\end{equation}
where %$\widehat{\mathcal{S}}_{diag}$ is the following set of matrices:
\begin{eqnarray}
\widehat{\mathcal{L}}(\widehat{r}_{A}) =  \{\mf{L} \mid \mf{L} \succeq 0, {\mf{L}}={\widehat{\mf{U}}_{\rm{A}}\mf{D}\widehat{\mf{U}}_{\rm{A}}'}, \mf{D} \in \R^{r \times r} \mathrm{diagonal}\},\nonumber\\
\widehat{\mathcal{S}}_{diag}=\{\mf{S} \in \R^{p \times p} \mid \mathrm{diag}(\mf{L})+\mathrm{diag}(\mf{S})=\mathrm{diag}(\widehat{\mf{\Sigma}}_{\rm{A}}),\nonumber\\
\mathrm{off-diag}(\mf{S})=\mathrm{off-diag}(\widehat{\mf{S}}_{\rm{A}}),
\mf{L}\in\widehat{\mathcal{L}}(\widehat{r}_{A})\}.\nonumber
\end{eqnarray}

%%as in \cite{farne2020}

%Le soglie rappresentano \mu e \xi!

%\far{From \ref{diff_eq}, we note that the Lipschitz constant of the problem is $2\frac{5}{4}=\frac{5}{2}$,
%such that the Algorithm needs un update of $\frac{2}{5}$. See \cite{luo2011high}. Divide by $\mathrm{tr}(\mf{\Sigma}_{n})$.}

For this reason, we calculate $\widehat{\mf{L}}_{\rm{U}}$ and $\widehat{\mf{S}}_{\rm{U}}$ as in \eqref{unshr1}, \eqref{unshr2} and \eqref{unshr3}
by Algorithm \ref{alg_ld} or \ref{alg_fro_2},
and we obtain, for each threshold pair $(\psi,\rho)$,
%in the bivariate grid $\psi_{init},\rho_{init}$,
the pairs of estimates
%\begin{eqnarray}
$\bigl(\widehat{\mf{L}}^{(ld)}_{\rm{U}}(\psi,\rho),\widehat{\mf{S}}^{(ld)}_{\rm{U}}(\psi,\rho)\bigr)$ or
%\nonumber\\
%$\bigl(\widehat{\mf{L}}^{(F)}_{\rm{U}}(\psi,\rho),\widehat{\mf{S}}^{(F)}_{\rm{U}}(\psi,\rho)\bigr)$,%\nonumber\\
$\bigl(\widehat{\mf{L}}^{(F)}_{\rm{U}}(\psi,\rho),\widehat{\mf{S}}^{(F)}_{\rm{U}}(\psi,\rho)\bigr)$.
%\nonumber
As a consequence, we can derive the overall UNALCE estimates as
$\widehat{\mf{\Sigma}}^{(ld)}_{\rm{U}}(\psi,\rho)=\widehat{\mf{L}}^{(ld)}_{\rm{U}}(\psi,\rho)+\widehat{\mf{S}}^{(ld)}_{\rm{U}}(\psi,\rho)$ and
%$\widehat{\mf{\Sigma}}^{(F)}_{\rm{U}}(\psi,\rho)=\widehat{\mf{L}}^{(F)}_{\rm{U}}(\psi,\rho)+\widehat{\mf{S}}^{(F)}_{\rm{U}}(\psi,\rho)$,
$\widehat{\mf{\Sigma}}^{(F)}_{\rm{U}}(\psi,\rho)=\widehat{\mf{L}}^{(F)}_{\rm{U}}(\psi,\rho)+\widehat{\mf{S}}^{(F)}_{\rm{U}}(\psi,\rho)$.
%\nonumber
%\end{eqnarray}
%Once we have calculated UNALCE over a grid of values for $\psi$ and $\rho$,
%obtaining the solution pairs
%$(\widehat{\mf{L}}(\psi,\rho)_{h,A},\widehat{\mf{S}}(\psi,\rho)_{h,A})$ and $(\widehat{\mf{L}}(\psi,\rho)_{h,U},\widehat{\mf{S}}(\psi,\rho)_{h,U})$, the overall solutions
%$\widehat{\mf{\Sigma}}(\psi,\rho)_{h,A}=\widehat{\mf{L}}(\psi,\rho)_{h,A}+\widehat{\mf{S}}(\psi,\rho)_{h,A}$ and %$\widehat{\mf{\Sigma}}(\psi,\rho)_{h,U}=\widehat{\mf{L}}(\psi,\rho)_{h,U}+\widehat{\mf{S}}(\psi,\rho)_{h,U}$, and the latent variance proportions $\widehat{\theta}(\psi,\rho)_{h,A}=\frac{\mathrm{tr}(\widehat{\mf{L}}(\psi,\rho)_{h,A})}{\mathrm{tr}(\widehat{\mf{\Sigma}}(\psi,\rho)_{h,A})}$ and
%$\widehat{\theta}(\psi,\rho)_{h,U}=\frac{\mathrm{tr}(\widehat{\mf{L}}(\psi,\rho)_{h,U})}{\mathrm{tr}(\widehat{\mf{\Sigma}}(\psi,\rho)_{h,U})}$.

Then, %for each $h=1,\ldots,100$,
given the latent variance proportions $\widehat{\theta}(\psi,\rho)_{\rm{A}}=(\mathrm{tr}(\widehat{\mf{\Sigma}}(\psi,\rho)_{\rm{A}}))^{-1}{\mathrm{tr}(\widehat{\mf{L}}(\psi,\rho)_{\rm{A}})}$ and
$\widehat{\theta}(\psi,\rho)_{\rm{U}}=(\mathrm{tr}(\widehat{\mf{\Sigma}}(\psi,\rho)_{\rm{U}}))^{-1}{\mathrm{tr}(\widehat{\mf{L}}(\psi,\rho)_{\rm{U}})}$,
%Then, for each $h=1,\ldots,100$
we can select the optimal threshold pairs $(\psi_{U},\rho_{U})$ and $(\psi_{A},\rho_{A})$
by minimizing the MC criteria
\begin{eqnarray}
MC(\psi,\rho)_{U}&=&\max
\left\{\frac{{\widehat{r}\Vert\widehat{\mf{L}}(\psi,\rho)_{\rm{U}}}\Vert_{2}}{\widehat{\theta}(\psi,\rho)_{\rm{U}}},
\frac{{\Vert\widehat{\mf{S}}(\psi,\rho)_{\rm{U}}}\Vert_{1,v}}{{\gamma}(1-\widehat{\theta}(\psi,\rho)_{\rm{U}})}\right\},\nonumber\\
MC(\psi,\rho)_{A}&=&\max\left\{\frac{{\widehat{r}\Vert\widehat{\mf{L}}(\psi,\rho)_{\rm{A}}}\Vert_{2}}{\widehat{\theta}(\psi,\rho)_{\rm{A}}},
\frac{{\Vert\widehat{\mf{S}}(\psi,\rho)_{\rm{A}}}\Vert_{1,v}}{{\gamma}(1-\widehat{\theta}(\psi,\rho)_{\rm{A}})}\right\},\nonumber
\end{eqnarray}
where ${\gamma}=\psi^{-1}{\rho}$ is the ratio between the sparsity and the latent eigenvalue threshold
(see \cite{farne2020large} for more details). In this way, we can select the optimal threshold pairs
$(\psi_{A},\rho_{A})=\arg \min_{(\psi,\rho)} MC(\psi,\rho)_{A}$ and
$(\psi_{U},\rho_{U})=\arg \min_{(\psi,\rho)} MC(\psi,\rho)_{U}$.
This procedure is applied both for Algorithms \ref{alg_ld} and \ref{alg_fro_2}, after
%For Algorithms \ref{alg_ld} and \ref{alg_fro_2},
the possible initial threshold pairs are obtained by the Cartesian product of the initial vectors
$\vf{\psi}_{init}$ and $\vf{\rho}_{init}$. %(see Section \ref{real}).
%%POET

\section{Simulation study}\label{sim}

%Parameters
%Settings
%Metrics
%Results
%Comments

% parameters
% scenarios
% model input alg-s POET UNALCE ALCE! (both)!
% all metrics
% 2+2(+2(appl)) 4 3 3

\subsection{Simulation settings}

%In this section, we test the validity of Theorems \ref{bartlett_opt} and \ref{thompson_opt}
%on some data simulated for that purpose. We report our main simulation parameters:
In this section, we test the theoretical results of previous sections
on some data simulated with this purpose.
Hereafter, we report the key simulation parameters:
\begin{enumerate}
\item the dimension $p$ and the sample size $n$;
\item the rank $r$ and the condition number $c=\mathrm{cond}(\mf{L}^{*})=\lambda_{1}(\mf{L}^{*})/\lambda_{r}(\mf{L}^{*})$
of the low rank component $\mf{L}^{*}$;% in spectral norm;
\item the trace of $\mf{L}^{*}$, $\tau \theta p$, where $\tau$ is a magnitude parameter
and $\theta=\mathrm{tr}(\mf{L}^{*})/\mathrm{tr}(\mf{\Sigma}^{*})$ is the proportion of variance explained by $\mf{L}^{*}$;
\item the number of off-diagonal non-zeros $s$ in the sparse component $\mf{S}^{*}$;
\item the minimum latent eigenvalue $\lambda_r(\mf{L}^{*})$;
\item the minimum nonzero off-diagonal residual entry in absolute value $\Vert\mf{S}^{*}\Vert_{\mathrm{min,off}}$;
\item the proportion of non-zeros over the number of off-diagonal elements, $\pi_s=\frac{2s}{p(p-1)}$ ;
\item the proportion of (absolute) residual covariance $\rho_{\mf{S}^{*}}=\frac{\sum_{i=1}^p\sum_{j \ne i}{\vert \mf{S}_{ij}^*\vert}}{\sum_{i=1}^p\sum_{j \ne i}\vert\mf{\Sigma}_{ij}^*\vert}$;
\item $N=100$ replicates for each setting.
\end{enumerate}
%Essentially, the low rank component is simulated by setting $r$ equispaced eigenvalues with sum $\tau \theta p$ and by deriving an orthogonal $r$-dimensional basis by Gram--Schmidt algorithm from a $p$-dimensional random matrix.
%The residual variances are simulated by a $p$-dimensional Dirichlet distribution with sum $(1-\tau) \theta p$, and then matched to the previously simulated diagonal elements of the low rank component $\mf{L}^{*}$ according to their relative magnitude.
%The off-diagonal elements are first simulated entry-wise as
%$\mf{S}^{*}_{ij} \sim \mbox{sgn}(\mf{L}^{*}_{ij})\mathrm{Unif}(0,\sqrt{\mf{S}_{ii}^{*}\mf{S}_{jj}^{*}})$, $i,j=1,\ldots,p$, $i \ne j$,
%%by a uniform distribution with parameters $0$ and the product of the correspondent standard deviations.
%where `$\mathrm{Unif}$' stands for `uniform distribution'.
%The smallest $p(p-1)/2-s$ off-diagonal elements in absolute value are then set to $0$.
The detailed simulation algorithm is reported in \cite{farne2016large}. We consider two distributional schemes:
\begin{enumerate}
  \item $f_i \sim MVN(\vf{0},\mf{I}_r)$, $\epsilon_i \sim MVN(\vf{0},\mf{S}^{*})$, with MVN multivariate normal;
  \item $f_i \sim MVT(\vf{0},\mf{I}_r,df=5)$, $\epsilon_i \sim MVT(\vf{0},\mf{S}^{*},df=5)$,
  with MVT multivariate Student's $t$. %with $5$ degrees of freedom.
\end{enumerate}

The main parameters of simulated settings are reported in Table \ref{sett}. %and \ref{specond}.
Settings 1-5 all present the same condition number $c=2$ of the low rank component, and a similar degree of residual sparsity in terms of proportion of nonzeros $\pi_s$ and residual covariance proportion $\rho_{\mf{S}^{*}}$.
However, Settings 1-3 have $p/n=0.1$, Setting 4 has $p/n=1$ and Setting 5 has $p/n=2$.
Settings 1-3 also differ by latent rank (following the progression $r=2-4-6$), latent variance proportion ($\theta=0.7-0.8-0.8$) and magnitude of latent eigenvalues ($\tau=200-200-300$).
In this way, Settings from 1 to 5 are increasingly challenging for what concerns the recovery of latent rank, residual sparsity pattern,
covariance matrix and its low rank and sparse components.

\begin{table}
  \tabcolsep=0pt
  \caption{\label{sett} Simulated settings: parameters.}
  %\centering
  %%\hline
  %\fbox{
  \begin{tabular*}{\columnwidth}{@{\extracolsep{\fill}}lcrcrrrrrrrrrrr@{}}
    \hline
    % after \\: %\hline or \cline{col1-col2} \cline{col3-col4} ...
    Setting  &  $p$ &  $n$ &  $p/n$ &  $r$ &  $\theta$ & $\tau$ & $c$ &  $\pi_s$ &  $\rho_{\mf{S}^{*}}$\\
    % &  $\pi^{T}_s$ &  $\rho^{N}_{\mf{T}^{*}}$\\ %&  {spikiness}&  {sparsity}\\% &  $\VertL\Vert$\\
    \hline
    1 & $100$ & $1000$ & $0.1$ & $2$ & $0.7$ & $200$ & $2$ & $0.0438$ & $0.0121$\\ %& {low} &{high}\\% & $23.33$\\
    %\hline
    %  & $p$ & $n$ & $r$ & $\tau$ & $\theta$ & $c$ & $s$ & $\rho_{\widehat{\mf{S}}}$\\
    %%\hline
    %%\hline
    % after \\: %\hline or \cline{col1-col2} \cline{col3-col4} ...
    %& $p$ & $n$ & $r$ & $\tau$ & $\theta$ & $c$ & $s$ & $\rho_{\widehat{\mf{S}}}$\\
    %%\hline
    %\textbf{2} & $100$ & $1000$ & $0.1$ & $4$ & $0.7$ & $4$ & $0.0677$ & $0.0048$ & {middle} & {middle}\\% & $28$\\
    2 & $100$ & $1000$ & $0.1$ & $4$ & $0.8$ & $200$ & $2$ & $0.0503$ & $0.0194$\\ %& {high} & {low}\\% & $128$\\
    3 & $100$ & $1000$ & $0.1$ & $6$ & $0.8$ & $300$ &$2$ & $0.0446$ & $0.0087$\\ %& {high} & {low}\\% & $128$\\
    %%\hline
    4 & $150$ & $150$ & $1$ & $5$ & $0.8$ & $200$ & $2$ & $0.0497$ & $0.0082$\\ %& {middle} & {middle}\\% & $32$\\
    5 & $200$ & $100$ & $2$ & $6$ & $0.8$ & $200$ & $2$ & $0.0503$ & $0.0092$\\ %& {middle} & {middle}\\% & $35.56$\\
    %\matt{5} & $500$ & $50$ & $10$ & $3$ & $0.8$ & $2$ & $0.0000$ & $0.0000$ & {high} & {high}\\%
    %\matt{6} & $1000$ & $20$ & $50$ & $2$ & $0.8$ & $2$ & $0.0001$ & $0.0001$ & {high} & {high}\\%
    \hline%yes!
  \end{tabular*}
%}
\end{table}

\subsection{Simulation metrics}

For each scenario and distributional scheme, we simulate $N=100$ replicates from model \eqref{mod} (main paper),
thus getting $100$ instances of the input sample covariance matrix $\mf{\Sigma}_n$.
We set up a grid of $n_{thr}=9$ thresholds as explained in Section \ref{algorithm}.
We define $\psi_{ord}$ and $\rho_{ord}$ as the ordinal place occupied by the optimal threshold pairs selected by MC criterion.
Tuning parameters $c^{\psi}_{thr}$ and $c^{\rho}_{thr}$ are chosen in order to obtain average $\psi_{ord}$ and $\rho_{ord}$
over the $N$ replicates far from $1$ and $n_{thr}$.
\begin{itemize}
\item We apply Algorithm \ref{alg_ld} to each generated $\mf{\Sigma}_n$ to get the %pair of estimates \eqref{logdet},
%that we call
ALCE-ld pair: $\bigl(\widehat{\mf{L}}^{(ld)}_{A},\widehat{\mf{S}}^{(ld)}_{A}\bigr)$.
Then, we apply the unshrinkage steps in \eqref{unshr1}, \eqref{unshr2}, \eqref{unshr3},
and we get the UNALCE-ld pair of estimates $\bigl(\widehat{\mf{L}}^{(ld)}_{U},\widehat{\mf{S}}^{(ld)}_{U}\bigr)$.
\item We apply Algorithm \ref{alg_fro_2} to each generated $\mf{\Sigma}_n$ to get the %pair of estimates \eqref{alg_fro_2},
%that we call
ALCE-F pair: $\bigl(\widehat{\mf{L}}^{(F)}_{A},\widehat{\mf{S}}^{(F)}_{A}\bigr)$.
Then, we apply the unshrinkage steps in \eqref{unshr1}, \eqref{unshr2}, \eqref{unshr3},
and we get the UNALCE-F pair: $\bigl(\widehat{\mf{L}}^{(F)}_{U},\widehat{\mf{S}}^{(F)}_{U}\bigr)$.
%\item we get the POET pair of estimates $\bigl(\widehat{\mf{L}}_{P},\widehat{\mf{S}}_{P}\bigr)$,
%where $\widehat{\mf{L}}_{h,P}$ is the covariance matrix of the first $\widehat{r}_{A}$ principal components of $\mf{\Sigma}_n$,
%$\widehat{\mf{S}}_{h,P}$ is obtained by thresholding its orthogonal complement as in \cite{fan2013large}, and $\widehat{\mf{\Sigma}}_{h,P}=\widehat{\mf{L}}_{h,P}+\widehat{\mf{S}}_{h,P}$ is the overall POET estimate.
\end{itemize}

Let us denote the generic low rank estimate as $\widehat{\mf{L}}$, the generic sparse estimate as $\widehat{\mf{S}}$,
and the generic covariance matrix estimate $\widehat{\mf{\Sigma}}=\widehat{\mf{L}}+\widehat{\mf{S}}$.
The performance metrics to assess the quality of estimates are
%The error norms used are:% the following:
%\mathrm{{Loss}}=  \Vert\widehat{\mf{L}}-\mf{L}^*\Vert_{F} + \Vert\widehat{\mf{S}}-\mf{S}^*\Vert_{F},\label{Loss}\\
the Frobenius total loss $TLF = \Vert\widehat{\mf{\Sigma}}-\mf{\Sigma}^{*}\Vert_{F}$;
%\label{TLF}
the spectral total loss $TL2 = \Vert\widehat{\mf{\Sigma}}-\mf{\Sigma}^{*}\Vert_{2}$;
%\label{TL2}
the spectral low rank loss $LL2 = \Vert\widehat{\mf{L}}-\mf{L}^{*}\Vert_{2}$;
%,\label{LL2}
the sparse maximum loss $SLM= \Vert\widehat{\mf{S}}-\mf{S}^{*}\Vert_{\infty}$.
%\label{SLM}
%\end{itemize}
%The proportion of wrongly recovered latent ranks is
%Their estimation performance is measured by %the mean square error, defined for $\widehat{\theta}$ as
%%\begin{equation}\mathrm{MSE}(\widehat{\theta}) = \frac{1}{N} \sum_{h=1}^{N} (\widehat{\theta}_h-\theta)^2,\label{MSE}\end{equation}
%%where $\widehat{\theta}_h$ is the estimate of $\theta$ on the $h$-th replicate. We also compute
%the estimation bias for each parameter, defined as
%$\mathrm{bias}(\widehat{\theta}) = \widehat{\theta}_{mean}-\theta$,%\label{bias_theta}
%$\mathrm{bias}(\widehat{\rho}_{\widehat{\mf{S}}}) = \widehat{\rho}_{\widehat{\mf{S}},{mean}}-\rho_{\mf{S}^{*}}$,%\label{bias_rho}\\
%$\mathrm{bias}(\widehat{\pi}_{\widehat{s}}) = \widehat{\pi}_{\widehat{s},{mean}}-\pi_s$,
%%\label{bias_pi}
%%\end{eqnarray}
%where $\widehat{\theta}_{mean}$, $\widehat{\rho}_{\widehat{\mf{S}},{mean}}$ and $\widehat{\pi}_{\widehat{s},{mean}}$
%are the mean estimates of $\theta$, $\rho_{\mf{S}^{*}}$ and $\pi_s$ over the $N$ replicates.

The performance in terms of eigen-structure recovery is measured for $\mf{\Sigma}^{*}$ by $\lambda(\widehat{\mf{\Sigma}})$,
which is defined as the Euclidean distance between the estimated and true eigenvalues of $\mf{\Sigma}^{*}$:
\begin{equation}\label{eig}
\lambda(\widehat{\mf{\Sigma}})=\sqrt{\sum_{i=1}^p (\widehat{\lambda}_i(\widehat{\mf{\Sigma}})-\lambda_i(\mf{\Sigma^{*}}))^2}.
\end{equation}
Measure \eqref{eig} is similarly defined for $\mf{L}^{*}$ and $\mf{S}^{*}$ as $\lambda(\widehat{\mf{L}})$ and $\lambda(\widehat{\mf{S}})$,
respectively. All these measures are then averaged over the $N$ replicates.

The estimated latent rank $\widehat{r}$,
%$err(\widehat{r})=\frac{1}{N}\sum_{k=1}^K\mathbbm{1}(\widehat{r}_{k}=r)$.
the estimated proportion of latent variance $\widehat{\theta}$,
of residual covariance $\widehat{\rho}_{\widehat{\mf{S}}}$, and of residual non-zeros $\widehat{\pi}_{\widehat{s}}$ are also computed.
Averages and standard deviations over the $N$ replicates of these quantities, denoted by $M$ and $SD$, are derived.

In the end, we calculate the following metrics for sparsity pattern recovery:
\begin{itemize}
\item $poserr=(pos)^{-1}{pw}$, where $pw$ are the positive elements mislabeled as zero or negative and $pos$ is the number of positive elements;
\item $negerr=(neg)^{-1}{nw}$, where $nw$ are the negative elements mislabeled as zero or positive and $neg$ is the number of negative elements;
\item $zeroerr=(nn)^{-1}{zw}$, where $zw$ are the zero elements mislabeled as positive or negative and $nn$ is the number of null elements.
\end{itemize}
The averages of these quantities over the $N$ replicates are reported.

Tables \ref{tab:normal1}--\ref{tab:t5} report simulation results for Settings 1-5 under the normal case
and under Student's t case, respectively. Hereafter, we report the general observable patterns.
\begin{itemize}
\item Concerning performance metrics (for matrix components and eigenvalues), UNALCE is generally doing quite better than ALCE, and log-det estimates are generally doing slightly better than Frobenius ones.
\item UNALCE latent variance proportion is systematically closer to the true $\theta$ compared to ALCE, and log-det estimates are doing slightly better than Frobenius ones.
\item Importantly, under the most challenging scenarios (Settings 4-5 with Student's t), log-det estimates are less prone than Frobenius ones
    to latent rank overestimation.
\item As $p/n$ decreases and $r$ increases, recovering the residual sparsity pattern becomes more difficult, and this effect is amplified when the distribution of factors and residuals possess heavy tails (like Student's t with $5$ degrees of freedom).
\end{itemize}

\bibliographystyle{agsm}
\bibliography{Bernardi_Farne_JASA}

\newpage

% Table generated by Excel2LaTeX from sheet '1 norm Latex'
\begin{table}[htbp]
\tiny
  \tabcolsep=0pt
  \caption{Simulation results: Setting 1 under the standard normal case. Estimates with $c^{\psi}_{thr}=20$ and $c^{\rho}_{thr}=3$, $n_{thr}=9$.}
  \begin{tabular*}{\columnwidth}{@{\extracolsep{\fill}}lcrcrrrrrrrrrrr@{}}
    \hline
          & \multicolumn{1}{l}{UNALCE-LD} & \multicolumn{1}{l}{ALCE-LD} & \multicolumn{1}{l}{UNALCE-F} & \multicolumn{1}{l}{ALCE-F} \\
    \hline
    $M(\widehat{r})$ & 2     & 2     & 2     & 2 \\
    $SD(\widehat{r})$ & 0     & 0     & 0     & 0 \\
    $M(\widehat{\theta})$ & 0.7030 & 0.6994 & 0.7038 & 0.6996 \\
    $SD(\widehat{\theta})$ & 0.0077 & 0.0079 & 0.0078 & 0.0078 \\
    $M(\widehat{\pi}_s)$ & 0.0381 & 0.0349 & 0.0400 & 0.0398 \\
    $SD(\widehat{\pi}_s)$ & 0.0044 & 0.0063 & 0.0043 & 0.0048 \\
    $M(\widehat{\rho}_{\widehat{\mf{S}}})$ & 0.0053 & 0.0058 & 0.0057 & 0.0057 \\
    $SD(\widehat{\rho}_{\widehat{\mf{S}}})$ & 0.0008 & 0.0010 & 0.0008 & 0.0010 \\
    $M(LL2)$ & 5.5101 & 5.5331 & 5.5162 & 5.5005 \\
    $SD(LL2)$ & 2.0695 & 1.9893 & 2.0908 & 2.0015 \\
    $M(SLM)$ & 0.4057 & 0.4007 & 0.4023 & 0.3867 \\
    $SD(SLM)$ & 0.0873 & 0.0850 & 0.0907 & 0.0869 \\
    $M(TL2)$ & 5.4774 & 5.5562 & 5.467 & 5.5178 \\
    $SD(TL2)$ & 1.9114 & 1.8586 & 1.9242 & 1.8790 \\
    $M(TLF)$ & 6.8729 & 6.9855 & 6.8408 & 6.8898 \\
    $SD(TLF)$ & 1.6368 & 1.5942 & 1.6505 & 1.6067 \\
    $M(\lambda(\widehat{\mf{L}}))$ & 4.0375 & 4.0744 & 4.0448 & 4.0461 \\
    $SD(\lambda(\widehat{\mf{L}}))$ & 2.3986 & 2.3296 & 2.4223 & 2.3383 \\
    $M(\lambda(\widehat{\mf{S}}))$ & 1.3268 & 1.3632 & 1.2873 & 1.2521 \\
    $SD(\lambda(\widehat{\mf{S}}))$ & 0.1800  & 0.2249 & 0.1854 & 0.2133 \\
    $M(\lambda(\widehat{\mf{\Sigma}}))$ & 4.1541 & 4.2763 & 4.1351 & 4.2115 \\
    $SD(\lambda(\widehat{\mf{\Sigma}}))$ & 2.2107 & 2.1819 & 2.2202 & 2.1939 \\
    $M(poserr)$ & 0.1649 & 0.1649 & 0.1649 & 0.1649 \\
    $SD(poserr)$ & 0.0803 & 0.0803 & 0.0803 & 0.0803 \\
    $M(negerr)$ & 0.1835 & 0.1835 & 0.1835 & 0.1835 \\
    $SD(negerr)$ & 0.0798 & 0.0798 & 0.0798 & 0.0798 \\
    $M(nnerr)$ & 0.0037 & 0.0037 & 0.0037 & 0.0037 \\
    $SD(nnerr)$ & 0.0017 & 0.0017 & 0.0017 & 0.0017 \\
    $M(\psi_{ord})$ & 6.66  & 7.35  & 7.72  & 7.88 \\
    $M(\rho_{ord})$ & 5.57  & 6.09  & 6.45  & 6.52 \\
    \hline
  \end{tabular*}%
  \label{tab:normal1}%
\end{table}%

\newpage

% Table generated by Excel2LaTeX from sheet '1 t nu=5 Latex'
\begin{table}[htbp]
\tiny
  \tabcolsep=0pt
  \caption{Simulation results: Setting 1 under Student't case with $5$ dof. Estimates with $c^{\psi}_{thr}=10$ and $c^{\rho}_{thr}=5$, $n_{thr}=9$.}
  \begin{tabular*}{\columnwidth}{@{\extracolsep{\fill}}lcrcrrrrrrrrrrr@{}}
    \hline
          & \multicolumn{1}{l}{UNALCE-LD} & \multicolumn{1}{l}{ALCE-LD} & \multicolumn{1}{l}{UNALCE-F} & \multicolumn{1}{l}{ALCE-F} \\
    \hline
    $M(\widehat{r})$ & 2.02  & 2.04  & 2.06  & 2.07 \\
    $SD(\widehat{r})$ & 0.1407 & 0.1969 & 0.2387 & 0.2564 \\
    $M(\widehat{\theta})$ & 0.7008 & 0.6960 & 0.7031 & 0.6967 \\
    $SD(\widehat{\theta})$ & 0.0167 & 0.0165 & 0.0156 & 0.0163 \\
    $M(\widehat{\pi}_s)$ & 0.0181 & 0.0127 & 0.0223 & 0.0179 \\
    $SD(\widehat{\pi}_s)$ & 0.0055 & 0.0082 & 0.0074 & 0.0100 \\
    $M(\widehat{\rho}_{\widehat{\mf{S}}})$ & 0.0027 & 0.0020 & 0.0033 & 0.0027 \\
    $SD(\widehat{\rho}_{\widehat{\mf{S}}})$ & 0.0008 & 0.0013 & 0.0010 & 0.0015 \\
    $M(LL2)$ & 62.5681 & 61.6074 & 62.8119 & 61.6974 \\
    $SD(LL2)$ & 14.6297 & 14.5853 & 14.5984 & 14.585 \\
    $M(SLM)$ & 3.6095 & 3.6299 & 3.5902 & 3.6481 \\
    $SD(SLM)$ & 0.7387 & 0.7590 & 0.7247 & 0.7511 \\
    $M(TL2)$ & 62.4306 & 61.2988 & 62.7981 & 61.5594 \\
    $SD(TL2)$ & 14.6298 & 14.5441 & 14.5972 & 14.617 \\
    $M(TLF)$ & 70.7890 & 69.3244 & 71.2676 & 69.6512 \\
    $SD(TLF)$ & 13.0491 & 12.9854 & 12.9943 & 13.0421 \\
    $M(\lambda(\widehat{\mf{L}}))$ & 68.6955 & 66.2495 & 66.3665 & 64.3276 \\
    $SD(\lambda(\widehat{\mf{L}}))$ & 16.416 & 18.7969 & 21.2332 & 21.9519 \\
    $M(\lambda(\widehat{\mf{S}}))$ & 6.938 & 6.9961 & 6.9209 & 7.0647 \\
    $SD(\lambda(\widehat{\mf{L}}))$ & 0.9668 & 0.9646 & 0.8919 & 0.9864 \\
    $M(\lambda(\widehat{\mf{\Sigma}}))$ & 70.3783 & 68.8897 & 70.8665 & 69.2295 \\
    $SD(\lambda(\widehat{\mf{\Sigma}}))$ & 13.094 & 13.0331 & 13.0374 & 13.0901 \\
    $M(poserr)$ & 0.6315 & 0.6315 & 0.6315 & 0.6315 \\
    $SD(poserr)$ & 0.1891 & 0.1891 & 0.1891 & 0.1891 \\
    $M(negerr)$ & 0.6348 & 0.6348 & 0.6348 & 0.6348 \\
    $SD(negerr)$ & 0.1960 & 0.1960 & 0.1960 & 0.1960 \\
    $M(nnerr)$ & 0.0019 & 0.0019 & 0.0019 & 0.0019 \\
    $SD(nnerr)$ & 0.0021 & 0.0021 & 0.0021 & 0.0021 \\
    $M(\psi_{ord})$ & 6.16  & 6.67  & 6.92  & 7.39 \\
    $M(\rho_{ord})$ & 3.49  & 4.42  & 3.79  & 4.44 \\
    \hline
  \end{tabular*}%
  \label{tab:t1}%
\end{table}%

\newpage

%% Table generated by Excel2LaTeX from sheet '1 norm Latex'
%\begin{table}[htbp]
%  \centering
%  \caption{Add caption}
%    \begin{tabular}{lrrrr}
%    $(c^{\psi}_{thr}$ & 20    &       &       &  \\
%    $(c^{\rho}_{thr}$ & 3     &       &       &  \\
%    num_thr & 9     &       &       &  \\
%    \end{tabular}%
%  \label{tab:addlabel}%
%\end{table}%

% Table generated by Excel2LaTeX from sheet '2 norm Latex'
\begin{table}[htbp]
\tiny
  \tabcolsep=0pt
  \caption{Simulation results: Setting 2 under the standard normal case. Estimates with $c^{\psi}_{thr}=10$ and $c^{\rho}_{thr}=3$, $n_{thr}=9$.}
   \begin{tabular*}{\columnwidth}{@{\extracolsep{\fill}}lcrcrrrrrrrrrrr@{}}
    \hline
          & \multicolumn{1}{l}{UNALCE-LD} & \multicolumn{1}{l}{ALCE-LD} & \multicolumn{1}{l}{UNALCE-F} & \multicolumn{1}{l}{ALCE-F} \\
    \hline
    $M(\widehat{r})$ & 4     & 4     & 4     & 4 \\
    $SD(\widehat{r})$ & 0     & 0     & 0     & 0 \\
    $M(\widehat{\theta})$ & 0.6992 & 0.6901 & 0.7011 & 0.6885 \\
    $SD(\widehat{\theta})$ & 0.0059 & 0.0069 & 0.0058 & 0.0061 \\
    $M(\widehat{\pi}_s)$ & 0.0257 & 0.0273 & 0.0296 & 0.0187 \\
    $SD(\widehat{\pi}_s)$ & 0.0064 & 0.0166 & 0.0066 & 0.0095 \\
    $M(\widehat{\rho}_{\widehat{\mf{S}}})$ & 0.0056 & 0.0059 & 0.0064 & 0.0042 \\
    $SD(\widehat{\rho}_{\widehat{\mf{S}}})$ & 0.0013 & 0.0034 & 0.0015 & 0.0021 \\
    $M(LL2)$ & 4.6827 & 4.8314 & 4.7075 & 4.8029 \\
    $SD(LL2)$ & 1.0521 & 1.0023 & 1.0699 & 1.0077 \\
    $M(SLM)$ & 0.4142 & 0.5033 & 0.4051 & 0.4384 \\
    $SD(SLM)$ & 0.0875 & 0.1093 & 0.0828 & 0.0781 \\
    $M(TL2)$ & 4.6939 & 4.9692 & 4.6519 & 5.0217 \\
    $SD(TL2)$ & 1.0489 & 1.1567 & 1.0492 & 1.1393 \\
    $M(TLF)$ & 7.3782 & 7.6453 & 7.283 & 7.7881 \\
    $SD(TLF)$ & 1.0363 & 1.1444 & 1.0413 & 1.0157 \\
    $M(\lambda(\widehat{\mf{L}}))$ & 3.2652 & 3.4522 & 3.2683 & 3.4915 \\
    $SD(\lambda(\widehat{\mf{L}}))$ & 1.2983 & 1.2616 & 1.3209 & 1.2277 \\
    $M(\lambda(\widehat{\mf{S}}))$ & 1.7429 & 1.8722 & 1.6814 & 1.8922 \\
    $SD(\lambda(\widehat{\mf{S}}))$ & 0.2290 & 0.2839 & 0.2325 & 0.2751 \\
    $M(\lambda(\widehat{\mf{\Sigma}}))$ & 3.5286 & 3.9763 & 3.4463 & 4.0914 \\
    $SD(\lambda(\widehat{\mf{\Sigma}}))$ & 1.1745 & 1.3008 & 1.173 & 1.2208 \\
    $M(poserr)$ & 0.6037 & 0.6037 & 0.6037 & 0.6037 \\
    $SD(poserr)$ & 0.1728 & 0.1728 & 0.1728 & 0.1728 \\
    $M(negerr)$ & 0.6624 & 0.6624 & 0.6624 & 0.6624 \\
    $SD(negerr)$ & 0.1752 & 0.1752 & 0.1752 & 0.1752 \\
    $M(nnerr)$ & 0.0003 & 0.0003 & 0.0003 & 0.0003 \\
    $SD(nnerr)$ & 0.0010 & 0.0010 & 0.0010 & 0.0010 \\
    $M(\psi_{ord})$ & 4.51  & 5.05  & 5.28  & 6.59 \\
    $M(\rho_{ord})$ & 4.50   & 4.88  & 4.99  & 6.59\\
    \hline
%    $(c^{\psi}_{thr}$ & 10    &       &       &  \\
%    $(c^{\\rho}_{thr}$ & 3     &       &       &  \\
%    num_thr & 9     &       &       &  \\
  \end{tabular*}%
  \label{tab:normal2}%
\end{table}%

\newpage

% Table generated by Excel2LaTeX from sheet '2 t nu=5 bis Latex'
\begin{table}[htbp]
\tiny
  \tabcolsep=0pt
  \caption{Simulation results: Setting 2 under Student't case with $5$ dof. Estimates with $c^{\psi}_{thr}=10$ and $c^{\rho}_{thr}=3$, $n_{thr}=9$.}
   \begin{tabular*}{\columnwidth}{@{\extracolsep{\fill}}lcrcrrrrrrrrrrr@{}}
   \hline
          & \multicolumn{1}{l}{UNALCE-LD} & \multicolumn{1}{l}{ALCE-LD} & \multicolumn{1}{l}{UNALCE-F} & \multicolumn{1}{l}{ALCE-F} \\
   \hline
    $M(\widehat{r})$ & 4.06  & 4.05  & 4.08  & 4.07 \\
    $SD(\widehat{r})$ & 0.2387 & 0.219 & 0.2727 & 0.2564 \\
    $M(\widehat{\theta})$ & 0.6925 & 0.6792 & 0.6956 & 0.6814 \\
    $SD(\widehat{\theta})$ & 0.0125 & 0.0139 & 0.0118 & 0.0127 \\
    $M(\widehat{\pi}_s)$ & 0.0264 & 0.0179 & 0.0278 & 0.0249 \\
    $SD(\widehat{\pi}_s)$ & 0.0084 & 0.0090 & 0.0075 & 0.0092 \\
    $M(\widehat{\rho}_{\widehat{\mf{S}}})$ & 0.0064 & 0.0045 & 0.0067 & 0.0061 \\
    $SD(\widehat{\rho}_{\widehat{\mf{S}}})$ & 0.0019 & 0.0024 & 0.0016 & 0.0021 \\
    $M(LL2)$ & 33.5096 & 32.2751 & 33.7628 & 32.5335 \\
    $SD(LL2)$ & 5.5716 & 5.5653 & 5.5306 & 5.5491 \\
    $M(SLM)$ & 3.2487 & 3.3099 & 3.2223 & 3.3481 \\
    $SD(SLM)$ & 0.5365 & 0.5935 & 0.5086 & 0.5289 \\
    $M(TL2)$ & 33.7685 & 32.4463 & 34.0292 & 32.8041 \\
    $SD(TL2)$ & 5.5645 & 5.5736 & 5.5446 & 5.5902 \\
    $M(TLF)$ & 49.8538 & 47.4765 & 50.3261 & 48.1149 \\
    $SD(TLF)$ & 5.5594 & 5.5649 & 5.5082 & 5.6008 \\
    $M(\lambda(\widehat{\mf{L}}))$ & 44.6805 & 42.8866 & 44.2589 & 42.5162 \\
    $SD(\lambda(\widehat{\mf{L}}))$ & 12.6154 & 11.3429 & 14.1889 & 12.9488 \\
    $M(\lambda(\widehat{\mf{S}}))$ & 7.1923 & 7.3725 & 7.1221 & 7.5071 \\
    $SD(\lambda(\widehat{\mf{S}}))$ & 0.8926 & 1.1069 & 0.8537 & 0.9494 \\
    $M(\lambda(\widehat{\mf{\Sigma}}))$ & 48.9854 & 46.533 & 49.4696 & 47.2089 \\
    $SD(\lambda(\widehat{\mf{\Sigma}}))$ & 5.6432 & 5.6573 & 5.5921 & 5.6987 \\
    $M(poserr)$ & 0.5513 & 0.5513 & 0.5513 & 0.5513 \\
    $SD(poserr)$ & 0.1424 & 0.1424 & 0.1424 & 0.1424 \\
    $M(negerr)$ & 0.5741 & 0.5741 & 0.5741 & 0.5741 \\
    $SD(negerr)$ & 0.1458 & 0.1458 & 0.1458 & 0.1458 \\
    $M(nnerr)$ & 0.0030 & 0.0030 & 0.0030 & 0.0030 \\
    $SD(nnerr)$ & 0.0024 & 0.0024 & 0.0024 & 0.0024 \\
    $M(\psi_{ord})$ & 6.00     & 7.08  & 7.16  & 7.53 \\
    $M(\rho_{ord})$ & 5.18  & 6.42  & 6.08  & 6.55 \\
    \hline
%    $(c^{\psi}_{thr}$ & 10    &       &       &  \\
%    $(c^{\\rho}_{thr}$ & 3     &       &       &  \\
%    num_thr & 9     &       &       &  \\
    \end{tabular*}%
  \label{tab:t2}%
\end{table}%

%\end{document}

\newpage

% Table generated by Excel2LaTeX from sheet '3 norm Latex'
\begin{table}[htbp]
\tiny
  \tabcolsep=0pt
  \caption{Simulation results: Setting 3 under the standard normal case. Estimates with $c^{\psi}_{thr}=15$ and $c^{\rho}_{thr}=3$, $n_{thr}=9$.}
   \begin{tabular*}{\columnwidth}{@{\extracolsep{\fill}}lcrcrrrrrrrrrrr@{}}
   \hline
   & \multicolumn{1}{l}{UNALCE-LD} & \multicolumn{1}{l}{ALCE-LD} & \multicolumn{1}{l}{UNALCE-F} & \multicolumn{1}{l}{ALCE-F} \\
   \hline
    $M(\widehat{r})$ & 6     & 6     & 6     & 6 \\
    $SD(\widehat{r})$ & 0     & 0     & 0     & 0 \\
    $M(\widehat{\theta})$ & 0.7915 & 0.7814 & 0.7913 & 0.7811 \\
    $SD(\widehat{\theta})$ & 0.0055 & 0.0052 & 0.0047 & 0.0052 \\
    $M(\widehat{\pi}_s)$ & 0.0058 & 0.0043 & 0.0081 & 0.0057 \\
    $SD(\widehat{\pi}_s)$ & 0.0027 & 0.0025 & 0.0026 & 0.0021 \\
    $M(\widehat{\rho}_{\widehat{\mf{S}}})$ & 0.0011 & 0.0009 & 0.0015 & 0.0011 \\
    $SD(\widehat{\rho}_{\widehat{\mf{S}}})$ & 0.0005 & 0.0004 & 0.0004 & 0.0004 \\
    $M(LL2)$ & 6.2100  & 6.5112 & 6.1783 & 6.4838 \\
    $SD(LL2)$ & 1.0767 & 1.0896 & 1.0718 & 1.1447 \\
    $M(SLM)$ & 0.7156 & 0.7114 & 0.6757 & 0.6882 \\
    $SD(SLM)$ & 0.2611 & 0.2466 & 0.2693 & 0.2390 \\
    $M(TL2)$ & 6.2745 & 6.7120 & 6.2318 & 6.6743 \\
    $SD(TL2)$ & 1.1247 & 1.1440 & 1.1168 & 1.1941 \\
    $M(TLF)$ & 10.2945 & 10.6816 & 10.2239 & 10.6401 \\
    $SD(TLF)$ & 1.184 & 1.2334 & 1.1909 & 1.2822 \\
    $M(\lambda(\widehat{\mf{L}}))$ & 4.384 & 4.9884 & 4.3349 & 4.9805 \\
    $SD(\lambda(\widehat{\mf{L}}))$ & 1.2621 & 1.3306 & 1.2673 & 1.3296 \\
    $M(\lambda(\widehat{\mf{S}}))$ & 1.3075 & 1.3798 & 1.2315 & 1.3538 \\
    $SD(\lambda(\widehat{\mf{S}}))$ & 0.299 & 0.2829 & 0.3035 & 0.2521 \\
    $M(\lambda(\widehat{\mf{\Sigma}}))$ & 4.6536 & 5.4080 & 4.5684 & 5.3788 \\
    $SD(\lambda(\widehat{\mf{\Sigma}}))$ & 1.2388 & 1.3348 & 1.2273 & 1.3347 \\
    $M(poserr)$ & 0.8560 & 0.8560 & 0.8560 & 0.8560 \\
    $SD(poserr)$ & 0.0547 & 0.0547 & 0.0547 & 0.0547 \\
    $M(negerr)$ & 0.8928 & 0.8928 & 0.8928 & 0.8928 \\
    $SD(negerr)$ & 0.0434 & 0.0434 & 0.0434 & 0.0434 \\
    $M(nnerr)$ & 0.0001 & 0.0001 & 0.0001 & 0.0001 \\
    $SD(nnerr)$ & 0.0001 & 0.0001 & 0.0001 & 0.0001 \\
    $M(\psi_{ord})$ & 4.64  & 5.48  & 5.49  & 6.63 \\
    $M(\rho_{ord})$ & 6.62  & 7.46  & 7.05  & 7.89\\
    \hline
%    $(c^{\psi}_{thr}$ & 15    &       &       &  \\
%    $(c^{\\rho}_{thr}$ & 3     &       &       &  \\
%    num_thr & 9     &       &       &  \\
    \end{tabular*}%
  \label{tab:normal3}%
\end{table}%

\newpage

% Table generated by Excel2LaTeX from sheet '3 t nu=5 ter Latex'
\begin{table}[htbp]
\tiny
  \tabcolsep=0pt
  \caption{Simulation results: Setting 3 under Student't case with $5$ dof.
   Estimates with $c^{\psi}_{thr}=15$ and $c^{\rho}_{thr}=5$, $n_{thr}=9$.}
   \begin{tabular*}{\columnwidth}{@{\extracolsep{\fill}}lcrcrrrrrrrrrrr@{}}
   \hline
          & \multicolumn{1}{l}{UNALCE-LD} & \multicolumn{1}{l}{ALCE-LD} & \multicolumn{1}{l}{UNALCE-F} & \multicolumn{1}{l}{ALCE-F} \\
   \hline
    $M(\widehat{r})$ & 6.02  & 6.02  & 6.05  & 6.05 \\
    $SD(\widehat{r})$ & 0.1407 & 0.1407 & 0.2190 & 0.2190 \\
    $M(\widehat{\theta})$ & 0.7870 & 0.7773 & 0.7892 & 0.7789 \\
    $SD(\widehat{\theta})$ & 0.0134 & 0.0134 & 0.0109 & 0.0125 \\
    $M(\widehat{\pi}_s)$ & 0.0059 & 0.0047 & 0.0077 & 0.0068 \\
    $SD(\widehat{\pi}_s)$ & 0.0025 & 0.0022 & 0.0022 & 0.0025 \\
    $M(\widehat{\rho}_{\widehat{\mf{S}}})$ & 0.0012 & 0.0010 & 0.0015 & 0.0014 \\
    $SD(\widehat{\rho}_{\widehat{\mf{S}}})$ & 0.0005 & 0.0005 & 0.0004 & 0.0006 \\
    $M(LL2)$ & 40.7173 & 39.6084 & 40.9535 & 39.7755 \\
    $SD(LL2)$ & 5.2848 & 5.3406 & 5.2900  & 5.3559 \\
    $M(SLM)$ & 7.2486 & 7.3550 & 7.3072 & 7.3599 \\
    $SD(SLM)$ & 1.6003 & 1.6359 & 1.5770 & 1.6148 \\
    $M(TL2)$ & 41.3662 & 40.2975 & 41.6264 & 40.4842 \\
    $SD(TL2)$ & 5.2810 & 5.3375 & 5.2995 & 5.3505 \\
    $M(TLF)$ & 69.5908 & 67.2182 & 70.1709 & 67.6205 \\
    $SD(TLF)$ & 6.2298 & 6.2917 & 6.2164 & 6.3602 \\
    $M(\lambda(\widehat{\mf{L}}))$ & 64.8503 & 62.3077 & 64.088 & 61.3725 \\
    $SD(\lambda(\widehat{\mf{L}}))$ & 11.1999 & 10.9645 & 14.4942 & 14.0328 \\
    $M(\lambda(\widehat{\mf{S}}))$ & 10.5325 & 10.8211 & 10.5807 & 10.823 \\
    $SD(\lambda(\widehat{\mf{S}}))$ & 1.6456 & 1.6547 & 1.5606 & 1.6226 \\
    $M(\lambda(\widehat{\mf{\Sigma}}))$ & 68.2597 & 65.8355 & 68.8552 & 66.2533 \\
    $SD(\lambda(\widehat{\mf{\Sigma}}))$ & 6.2225 & 6.3009 & 6.2048 & 6.3472 \\
    $M(poserr)$ & 0.8443 & 0.8443 & 0.8443 & 0.8443 \\
    $SD(poserr)$ & 0.0571 & 0.0571 & 0.0571 & 0.0571 \\
    $M(negerr)$ & 0.8720 & 0.8720 & 0.8720 & 0.8720 \\
    $SD(negerr)$ & 0.0450 & 0.0450 & 0.045 & 0.045 \\
    $M(nnerr)$ & 0.0004 & 0.0004 & 0.0004 & 0.0004 \\
    $SD(nnerr)$ & 0.0005 & 0.0005 & 0.0005 & 0.0005 \\
    $M(\psi_{ord})$ & 5.58  & 6.40   & 6.44  & 7.24 \\
    $M(\rho_{ord})$ & 7.28  & 7.78  & 7.88  & 8.11 \\
    \hline
%    $(c^{\psi}_{thr}$ & 15    &       &       &  \\
%    $(c^{\\rho}_{thr}$ & 3     &       &       &  \\
%    num_thr & 9     &       &       &  \\
    \end{tabular*}%
  \label{tab:t3}%
\end{table}%

\newpage

% Table generated by Excel2LaTeX from sheet '4 norm Latex'
\begin{table}[htbp]
\tiny
  \tabcolsep=0pt
  \caption{Simulation results: Setting 4 under the standard normal case. Estimates with $c^{\psi}_{thr}=10$ and $c^{\rho}_{thr}=10$, $n_{thr}=9$.}
   \begin{tabular*}{\columnwidth}{@{\extracolsep{\fill}}lcrcrrrrrrrrrrr@{}}
   \hline
          & \multicolumn{1}{l}{UNALCE-LD} & \multicolumn{1}{l}{ALCE-LD} & \multicolumn{1}{l}{UNALCE-F} & \multicolumn{1}{l}{ALCE-F} \\
   \hline
    $M(\widehat{r})$ & 5     & 5     & 5     & 5.04 \\
    $SD(\widehat{r})$ & 0     & 0     & 0     & 0.1969 \\
    $M(\widehat{\theta})$ & 0.7883 & 0.7811 & 0.7933 & 0.7826 \\
    $SD(\widehat{\theta})$ & 0.0108 & 0.0112 & 0.0093 & 0.0098 \\
    $M(\widehat{\pi}_s)$ & 0.0053 & 0.0049 & 0.0061 & 0.0057 \\
    $SD(\widehat{\pi}_s)$ & 0.0040 & 0.0040 & 0.0035 & 0.0035 \\
    $M(\widehat{\rho}_{\widehat{\mf{S}}})$ & 0.0011 & 0.0010 & 0.0013 & 0.0012 \\
    $SD(\widehat{\rho}_{\widehat{\mf{S}}})$ & 0.0007 & 0.0007 & 0.0006 & 0.0006 \\
    $M(LL2)$ & 11.4366 & 11.4455 & 11.4472 & 11.4347 \\
    $SD(LL2)$ & 2.5940 & 2.5067 & 2.6278 & 2.4995 \\
    $M(SLM)$ & 0.5742 & 0.5679 & 0.5749 & 0.5715 \\
    $SD(SLM)$ & 0.2199 & 0.2198 & 0.2230 & 0.2186 \\
    $M(TL2)$ & 11.4482 & 11.5183 & 11.4366 & 11.4993 \\
    $SD(TL2)$ & 2.5324 & 2.4376 & 2.5526 & 2.4433 \\
    $M(TLF)$ & 17.4280 & 17.5020 & 17.4079 & 17.4807 \\
    $SD(TLF)$ & 2.4610 & 2.4233 & 2.4675 & 2.4265 \\
    $M(\lambda(\widehat{\mf{L}}))$ & 7.9636 & 8.1587 & 7.9092 & 8.1413 \\
    $SD(\lambda(\widehat{\mf{L}}))$ & 3.2225 & 3.0890 & 3.2720 & 3.1062 \\
    $M(\lambda(\widehat{\mf{S}}))$ & 0.9853 & 1.0170 & 0.9624 & 1.0039 \\
    $SD(\lambda(\widehat{\mf{S}}))$ & 0.2082 & 0.2070 & 0.2280 & 0.2063 \\
    $M(\lambda(\widehat{\mf{\Sigma}}))$ & 7.9122 & 8.1680 & 7.8320 & 8.1099 \\
    $SD(\lambda(\widehat{\mf{\Sigma}}))$ & 3.1273 & 3.0079 & 3.1666 & 3.0239 \\
    $M(poserr)$ & 0.9053 & 0.9053 & 0.9053 & 0.9053 \\
    $SD(poserr)$ & 0.0538 & 0.0538 & 0.0538 & 0.0538 \\
    $M(negerr)$ & 0.9084 & 0.9084 & 0.9084 & 0.9084 \\
    $SD(negerr)$ & 0.0539 & 0.0539 & 0.0539 & 0.0539 \\
    $M(nnerr)$ & 0.0009 & 0.0009 & 0.0009 & 0.0009 \\
    $SD(nnerr)$ & 0.0009 & 0.0009 & 0.0009 & 0.0009 \\
    $M(\psi_{ord})$ & 6.63  & 7.20   & 7.50   & 8.15 \\
    $M(\rho_{ord})$ & 2.87  & 3.16  & 3.05  & 3.18\\
    \hline
%    $(c^{\psi}_{thr}$ & 10    &       &       &  \\
%    $(c^{\\rho}_{thr}$ & 10    &       &       &  \\
%    num_thr & 9     &       &       &  \\
    \end{tabular*}%
  \label{tab:normal4}%
\end{table}%

\newpage

% Table generated by Excel2LaTeX from sheet '4 t nu=5 bis Latex'
\begin{table}[htbp]
\tiny
  \tabcolsep=0pt
  \caption{Simulation results: Setting 4 under Student't case with $5$ dof.
   Estimates with $c^{\psi}_{thr}=3$ and $c^{\rho}_{thr}=5$, $n_{thr}=9$.}
   \begin{tabular*}{\columnwidth}{@{\extracolsep{\fill}}lcrcrrrrrrrrrrr@{}}
   \hline
          & \multicolumn{1}{l}{UNALCE-LD} & \multicolumn{1}{l}{ALCE-LD} & \multicolumn{1}{l}{UNALCE-F} & \multicolumn{1}{l}{ALCE-F} \\
   \hline
    $M(\widehat{r})$ & 5.09  & 5.11  & 5.19  & 5.22 \\
    $SD(\widehat{r})$ & 0.2876 & 0.3145 & 0.3943 & 0.4163 \\
    $M(\widehat{\theta})$ & 0.7712 & 0.7579 & 0.7829 & 0.7639 \\
    $SD(\widehat{\theta})$ & 0.0307 & 0.0329 & 0.0248 & 0.0286 \\
    $M(\widehat{\pi}_s)$ & 0.0027 & 0.0018 & 0.0032 & 0.003 \\
    $SD(\widehat{\pi}_s)$ &  0.0022 & 0.0022 & 0.0029 & 0.0029 \\
    $M(\widehat{\rho}_{\widehat{\mf{S}}})$ & 0.0009 & 0.0005 & 0.0009 & 0.0008 \\
    $SD(\widehat{\rho}_{\widehat{\mf{S}}})$ & 0.0007 & 0.0007 & 0.0008 & 0.0007 \\
    $M(LL2)$ & 42.6011 & 41.0968 & 43.203 & 41.544 \\
    $SD(LL2)$ & 12.3140 & 12.3856 & 12.4008 & 12.3978 \\
    $M(SLM)$ & 3.2916 & 3.2528 & 3.1464 & 3.1858 \\
    $SD(SLM)$ & 1.3842 & 1.3556 & 1.2107 & 1.1764 \\
    $M(TL2)$ & 42.8417 & 41.3303 & 43.4106 & 41.7779 \\
    $SD(TL2)$ & 12.3264 & 12.4000  & 12.4053 & 12.3957 \\
    $M(TLF)$ & 57.1260 & 54.4483 & 58.1335 & 55.2066 \\
    $SD(TLF)$ & 13.8736 & 13.9476 & 13.9696 & 13.8934 \\
    $M(\lambda(\widehat{\mf{L}}))$ & 47.5626 & 43.4467 & 43.3635 & 39.3319 \\
    $SD(\lambda(\widehat{\mf{L}}))$ & 20.7962 & 20.5631 & 24.8669 & 24.6726 \\
    $M(\lambda(\widehat{\mf{S}}))$ & 5.6675 & 5.7639 & 5.2878 & 5.6197 \\
    $SD(\lambda(\widehat{\mf{S}}))$ & 1.8795 & 1.8550 & 1.5878 & 1.6010 \\
    $M(\lambda(\widehat{\mf{\Sigma}}))$ & 52.4152 & 49.4936 & 53.5001 & 50.3303 \\
    $SD(\lambda(\widehat{\mf{\Sigma}}))$ & 14.6645 & 14.8323 & 14.7237 & 14.7202 \\
    $M(poserr)$ & 0.9535 & 0.9535 & 0.9535 & 0.9535 \\
    $SD(poserr)$ & 0.0424 & 0.0424 & 0.0424 & 0.0424 \\
    $M(negerr)$ & 0.9559 & 0.9559 & 0.9559 & 0.9559 \\
    $SD(negerr)$ & 0.0421 & 0.0421 & 0.0421 & 0.0421 \\
    $M(nnerr)$ & 0.0007 & 0.0007 & 0.0007 & 0.0007 \\
    $SD(nnerr)$ & 0.0008 & 0.0008 & 0.0008 & 0.0008 \\
    $M(\psi_{ord})$ & 5.93  & 6.61  & 6.65  & 7.09 \\
    $M(\rho_{ord})$ & 2.73  & 3.56  & 3.39  & 3.32 \\
    \hline
%    $(c^{\psi}_{thr}$ & 5     &       &       &  \\
%    $(c^{\\rho}_{thr}$ & 10    &       &       &  \\
%    num_thr & 9     &       &       &  \\
    \end{tabular*}%
  \label{tab:t4}%
\end{table}%

\newpage

% Table generated by Excel2LaTeX from sheet '5 norm Latex'
\begin{table}[htbp]
\tiny
  \tabcolsep=0pt
  \caption{Simulation results: Setting 5 under the standard normal case. Estimates with $c^{\psi}_{thr}=5$ and $c^{\rho}_{thr}=5$, $n_{thr}=9$.}
   \begin{tabular*}{\columnwidth}{@{\extracolsep{\fill}}lcrcrrrrrrrrrrr@{}}
   \hline
          & \multicolumn{1}{l}{UNALCE-LD} & \multicolumn{1}{l}{ALCE-LD} & \multicolumn{1}{l}{UNALCE-F} & \multicolumn{1}{l}{ALCE-F} \\
   \hline
    $M(\widehat{r})$ & 6     & 6     & 6     & 6 \\
    $SD(\widehat{r})$ & 0     & 0     & 0     & 0 \\
    $M(\widehat{\theta})$ & 0.7848 & 0.7741 & 0.7907 & 0.7743 \\
    $SD(\widehat{\theta})$ & 0.0114 & 0.0108 & 0.0099 & 0.0128 \\
    $M(\widehat{\pi}_s)$ & 0.0020 & 0.0017 & 0.0029 & 0.0023 \\
    $SD(\widehat{\pi}_s)$ & 0.0017 & 0.0022 & 0.0021 & 0.0018 \\
    $M(\widehat{\rho}_{\widehat{\mf{S}}})$ & 0.0006 & 0.0005 & 0.0009 & 0.0007 \\
    $SD(\widehat{\rho}_{\widehat{\mf{S}}})$ & 0.0005 & 0.0004 & 0.0005 & 0.0005 \\
    $M(LL2)$ & 12.8676 & 12.854 & 12.8729 & 12.8525 \\
    $SD(LL2)$ & 2.5315 & 2.3597 & 2.575 & 2.4096 \\
    $M(SLM)$ & 0.7663 & 0.7635 & 0.7679 & 0.7597 \\
    $SD(SLM)$ & 0.3088 & 0.3085 & 0.3012 & 0.3120 \\
    $M(TL2)$ & 12.8594 & 12.8963 & 12.848 & 12.8913 \\
    $SD(TL2)$ & 2.4454 & 2.3027 & 2.4869 & 2.3557 \\
    $M(TLF)$ & 20.8945 & 21.0036 & 20.8753 & 20.9909 \\
    $SD(TLF)$ & 2.4450 & 2.3852 & 2.4612 & 2.4062 \\
    $M(\lambda(\widehat{\mf{L}}))$ & 9.9414 & 10.2549 & 9.8727 & 10.2164 \\
    $SD(\lambda(\widehat{\mf{L}}))$ & 2.9008 & 2.6856 & 2.9583 & 2.7332 \\
    $M(\lambda(\widehat{\mf{S}}))$ & 1.2268 & 1.2786 & 1.2007 & 1.2550 \\
    $SD(\lambda(\widehat{\mf{S}}))$ & 0.3069 & 0.3083 & 0.3028 & 0.2934 \\
    $M(\lambda(\widehat{\mf{\Sigma}}))$ & 9.8685 & 10.262 & 9.7756 & 10.2126 \\
    $SD(\lambda(\widehat{\mf{\Sigma}}))$ & 2.8087 & 2.6245 & 2.8593 & 2.6700 \\
    $M(poserr)$ & 0.9652 & 0.9652 & 0.9652 & 0.9652 \\
    $SD(poserr)$ & 0.0262 & 0.0262 & 0.0262 & 0.0262 \\
    $M(negerr)$ & 0.9611 & 0.9611 & 0.9611 & 0.9611 \\
    $SD(negerr)$ & 0.0310 & 0.0310 & 0.0310 & 0.0310 \\
    $M(nnerr)$ & 0.0003 & 0.0003 & 0.0003 & 0.0003 \\
    $SD(nnerr)$ & 0.0004 & 0.0004 & 0.0004 & 0.0004 \\
    $M(\psi_{ord})$ & 5.12  & 5.84  & 5.90   & 6.61 \\
    $M(\rho_{ord})$ & 4.61  & 5.02  & 4.70   & 5.17 \\
    \hline
%    $(c^{\psi}_{thr}$ & 5     &       &       &  \\
%    $(c^{\\rho}_{thr}$ & 5     &       &       &  \\
%    num_thr & 9     &       &       &  \\
    \end{tabular*}%
  \label{tab:normal5}%
\end{table}%

\newpage

% Table generated by Excel2LaTeX from sheet '5 t nu=5 Latex'
\begin{table}[htbp]
\tiny
  \tabcolsep=0pt
  \caption{Simulation results: Setting 5 under Student't case with $5$ dof.
  Estimates with $c^{\psi}_{thr}=3$ and $c^{\rho}_{thr}=5$, $n_{thr}=9$.}
   \begin{tabular*}{\columnwidth}{@{\extracolsep{\fill}}lcrcrrrrrrrrrrr@{}}
   \hline
          & \multicolumn{1}{l}{UNALCE-LD} & \multicolumn{1}{l}{ALCE-LD} & \multicolumn{1}{l}{UNALCE-F} & \multicolumn{1}{l}{ALCE-F} \\
   \hline
    $M(\widehat{r})$ & 6.06     & 6.10     & 6.16     & 6.18 \\
    $SD(\widehat{r})$ & 0.2387     & 0.2387     & 0.3685     & 0.3861 \\												
    $M(\widehat{\theta})$ & 0.7597 & 0.7443 & 0.7782 & 0.7520 \\														
    $SD(\widehat{\theta})$ & 0.0364 & 0.0391 & 0.0275 & 0.0325 \\
    $M(\widehat{\pi}_s)$ &  0.0013 & 0.0011 & 0.0016 & 0.0015 \\
    $SD(\widehat{\pi}_s)$ & 0.0014 & 0.0013 & 0.0015 & 0.0017 \\
    $M(\widehat{\rho}_{\widehat{\mf{S}}})$ & 0.0006 & 0.0005 & 0.0007 & 0.0006 \\															
    $SD(\widehat{\rho}_{\widehat{\mf{S}}})$ & 0.0005 & 0.0005 & 0.0005 & 0.0006 \\														
    $M(LL2)$ & 42.8189 & 41.1639 & 43.6323 & 41.629 \\
    $SD(LL2)$ & 21.3713 & 21.3396 & 21.3767 & 21.3577 \\														
    $M(SLM)$ & 3.1938 & 3.1368 & 3.0012 & 3.0203 \\
    $SD(SLM)$ & 1.5925 & 1.4788 & 1.4482 & 1.3249 \\
    $M(TL2)$ & 42.9352 & 41.2892 & 43.7066 & 41.7499 \\
    $SD(TL2)$ & 21.3700 & 21.3507 & 21.3871 & 21.3716 \\
    $M(TLF)$ & 57.3352 & 54.5298 & 58.6631 & 55.2326 \\
    $SD(TLF)$ & 21.0614	& 21.093 & 21.059 & 21.0582 \\
    $M(\lambda(\widehat{\mf{L}}))$ & 47.7139 & 42.9427 & 44.7736 & 40.2206 \\
    $SD(\lambda(\widehat{\mf{L}}))$ & 24.2546 & 25.4349 & 28.0912 & 27.7984 \\ 															
    $M(\lambda(\widehat{\mf{S}}))$ & 5.2042 & 5.2850 & 4.6958 & 5.0653 \\ 															
    $SD(\lambda(\widehat{\mf{S}}))$ & 2.0411 & 1.9345 & 1.8401 & 1.7646 \\
    $M(\lambda(\widehat{\mf{\Sigma}}))$ & 50.4726 & 47.2834 & 51.973 & 48.0925 \\
    $SD(\lambda(\widehat{\mf{\Sigma}}))$ & 21.4752 & 21.6242 & 21.3727 & 21.5275 \\
    $M(poserr)$ & 0.9778 & 0.9778 & 0.9778 & 0.9778 \\
    $SD(poserr)$ & 0.0233 & 0.0233 & 0.0233 & 0.0233 \\
    $M(negerr)$ & 0.9757 & 0.9757 & 0.9757 & 0.9757 \\
    $SD(negerr)$ & 0.0265 & 0.0265 & 0.0265 & 0.0265 \\
    $M(nnerr)$ & 0.0003 & 0.0003 & 0.0003 & 0.0003 \\
    $SD(nnerr)$ & 0.0004 & 0.0004 & 0.0004 & 0.0004 \\
    $M(\psi_{ord})$ & 5.58  & 5.98  & 6.04   & 6.49 \\
    $M(\rho_{ord})$ & 4.92  & 4.95  & 5.16   & 5.17 \\
    \hline
%    $(c^{\psi}_{thr}$ & 3     &       &       &  \\
%    $(c^{\rho}_{thr}$ & 5     &       &       &  \\
%    num_thr & 9     &       &       &  \\
    \end{tabular*}%
  \label{tab:t5}%
\end{table}%

\end{document}